\documentclass{amsbook}
\usepackage{amsmath,amsthm,latexsym,amscd,amssymb,graphicx}
\usepackage[all]{xy}

\newtheorem{prop}{Proposition}[section]
\newtheorem{lem}[prop]{Lemma}
\newtheorem{ddd}[prop]{Definition}
\newtheorem{theorem}[prop]{Theorem}
\newtheorem{cor}[prop]{Corollary}
\newtheorem*{theor}{Theorem}

\newcommand{\Gl}{\mathop{\mbox{\rm Gl}}} 
\newcommand{\ch}{\mathop{\mbox{\rm ch}}}
\newcommand{\End}{\mathop{\mbox{\rm End}}}

\newcommand{\vol}{\mathop{\mbox{\rm vol}}}
\newcommand{\ind}{\mathop{\mbox{\rm ind}}}

\newcommand{\spann}{{\rm span}\!}

\newcommand{\dom}{\mathop{\rm dom}}

\newcommand{\Coker}{\mathop{\rm Coker}}
\newcommand{\Pj}{{\mathcal P}}
\newcommand{\RM}{{\mathcal R}}
\newcommand{\tr}{\mathop{\rm tr}}
\newcommand{\trs}{{\rm tr}_s}
\newcommand{\id}{\mathop{\rm id}}
\newcommand{\trsi}{{\rm tr}_{\sigma}}

\newcommand{\Tr}{\mathop{\rm Tr}}
\newcommand{\Trsi}{{\rm Tr}_{\sigma}}
\newcommand{\Trs}{{\rm Tr}_s}
\newcommand{\Op}{\mathop{\rm Op}}
\newcommand{\diag}{\mathop{\rm diag}}

\newcommand{\N}{{\mathcal N}}
\newcommand{\U}{{\mathcal U}}

\newcommand{\A}{{\mathcal A}}
\newcommand{\B}{{\mathcal B}}
\newcommand{\Ai}{{\mathcal A}_{\infty}}
\newcommand{\ta}{\otimes {\mathcal A}}

\newcommand{\C}{C^{\infty}}
\newcommand{\Ol}[1]{\hat \Omega_{\le {#1}}\A_i}
\newcommand{\Kappa}{{\mathcal K}}
\newcommand{\Ok}{\hat \Omega_k}
\newcommand{\Oi}{\hat \Omega_*}
\newcommand{\ra}{\partial}

\newcommand{\aq}{\Leftrightarrow}
\newcommand{\ten}{\otimes}
\newcommand{\pten}{\otimes_{\pi}}
\newcommand{\pl}[1]{\varprojlim\limits_{#1}}

\newcommand{\ve}{\varepsilon}
\newcommand{\ov}{\overline}

\newcommand{\incl}{\hookrightarrow}
\newcommand{\cp}{\clubsuit}
\newcommand{\dirac}{\partial \!\!\!/_E}
\newcommand{\dirz}{\partial \!\!\!/_Z}
\newcommand{\dira}{\partial \!\!\!/}

\DeclareMathOperator{\re}{Re}
\DeclareMathOperator{\supp}{supp}
\DeclareMathOperator{\Ran}{Ran}
\DeclareMathOperator{\Ker}{Ker}
\DeclareMathOperator{\di}{d}
\DeclareMathOperator{\Mod}{mod}

\def\bbbr{{\rm I\!R}} 
\def\bbbn{{\rm I\!N}} 

\def\bbbc{{\rm I\!C}}

\def\bbbq{{\mathchoice {\setbox0=\hbox{$\displaystyle\rm Q$}\hbox{\raise
0.15\ht0\hbox to0pt{\kern0.4\wd0\vrule height0.8\ht0\hss}\box0}}
{\setbox0=\hbox{$\textstyle\rm Q$}\hbox{\raise
0.15\ht0\hbox to0pt{\kern0.4\wd0\vrule height0.8\ht0\hss}\box0}}
{\setbox0=\hbox{$\scriptstyle\rm Q$}\hbox{\raise
0.15\ht0\hbox to0pt{\kern0.4\wd0\vrule height0.7\ht0\hss}\box0}}
{\setbox0=\hbox{$\scriptscriptstyle\rm Q$}\hbox{\raise
0.15\ht0\hbox to0pt{\kern0.4\wd0\vrule height0.7\ht0\hss}\box0}}}}

\def\bbbz{{\mathchoice {\hbox{$\sf\textstyle Z\kern-0.4em Z$}}
{\hbox{$\sf\textstyle Z\kern-0.4em Z$}}
{\hbox{$\sf\scriptstyle Z\kern-0.3em Z$}}
{\hbox{$\sf\scriptscriptstyle Z\kern-0.2em Z$}}}}

\def\bbbc{{\mathchoice {\setbox0=\hbox{$\displaystyle\rm C$}\hbox{\hbox
to0pt{\kern0.4\wd0\vrule height0.9\ht0\hss}\box0}}
{\setbox0=\hbox{$\textstyle\rm C$}\hbox{\hbox
to0pt{\kern0.4\wd0\vrule height0.9\ht0\hss}\box0}}
{\setbox0=\hbox{$\scriptstyle\rm C$}\hbox{\hbox
to0pt{\kern0.4\wd0\vrule height0.9\ht0\hss}\box0}}
{\setbox0=\hbox{$\scriptscriptstyle\rm C$}\hbox{\hbox
to0pt{\kern0.4\wd0\vrule height0.9\ht0\hss}\box0}}}}

\numberwithin{section}{chapter}

\begin{document}
\frontmatter
\title{Noncommutative Maslov Index and Eta-Forms}

\author{Charlotte Wahl}

\address{Department of Mathematics, 460 McBryde, Virginia Tech,
 Blacksburg, VA 24061, USA}
 
\subjclass[2000]{19K56 (53D12, 58J28, 46L87)}
\keywords{Dirac operator, Maslov index, eta-invariant, higher index theory, $K$-theory}

\begin{abstract}
We define and prove a noncommutative generalization of a formula relating the Maslov index of a triple of Lagrangian subspaces of a symplectic vector space  to eta-invariants associated to a pair of Lagrangian subspaces. The noncommutative Maslov index, defined for modules over a $C^*$-algebra $\A$, is an element in $K_0(\A)$. The generalized formula calculates its Chern character in the de Rham homology of certain dense subalgebras of $\A$. The proof is a noncommutative Atiyah-Patodi-Singer index theorem for a particular Dirac operator twisted by an $\A$-vector bundle. We develop an analytic framework for this type of index problem. 
\end{abstract}

\maketitle

\setcounter{page}{4}
\tableofcontents

\mainmatter
\chapter*{Introduction}

The purpose of this paper is twofold. We establish a noncommutative generalization of a formula by Cappell, Lee and Miller \cite{clm} relating the Maslov index of a triple of Lagrangian subspaces of a symplectic vector space to certain $\eta$-invariants which can be associated to a pair of Lagrangian subspaces. The formula was generalized to families of Lagrangian subspaces by Bunke and Koch \cite{bk}. Our Maslov index will be defined for $C^*$-modules and will be an element in the $K$-theory of the $C^*$-algebra, whereas the $\eta$-invariants are replaced by noncommutative differential forms.  The proof of our formula is a noncommutative index theorem for a particular Dirac operator twisted by a $C^*$-vector bundle on a two-dimensional manifold with boundary and cylindric ends. The index theorem calculates the Chern character of index of the Dirac operator in the de Rham homology of certain dense subalgebras of the $C^*$-algebra. 

The second aim is to provide a precise analytic framework for the index theory of Dirac operators over a certain type of dense subalgebras of $C^*$-algebras based on heat kernel methods and the Quillen superconnection formalism. In particular these dense subalgebras are assumed to be projective limits of Banach algebras, so that we mostly deal with Banach space valued functions. Our results apply to higher index theory for invariant Dirac operators on covering spaces. Proofs in higher index theory based on the superconnection formalism have been given before \cite{lo3} \cite{lo1} \cite{lp1} \cite{lp2} \cite{wu}. Some analytical methods relevant for the general situation have been developed by Lott in \cite{lo2}. \\

The Maslov index $\tau(L_0,L_1,L_2)$ of a triple of Lagrangian subspaces $(L_0,L_1,L_2)$ of $\bbbr^{2n}$ endowed with the standard symplectic form $\omega$ is defined as the signature of the quadratic form $q$ on $L_0 \cap (L_1 + L_2)$ \cite{lv} given by
$$q(x_1+x_2):=\omega(x_2,x_1),~x_i \in L_i \ .$$ Its geometric significance comes from a gluing formula for signatures of manifolds with boundary \cite{wa}: Let $M^{4n}$ be an oriented manifold with boundary and let $N$ be a hypersurface in $M$ with boundary. Cutting along $N$ produces two (topological) manifolds $M_1, M_2$ with boundary. The images of $H^{2n-1}(N)$ and $H^{2n-1}(M_i),~i=1,2,$ in $H^{2n-1}(\ra N)$ are Lagrangian subspaces with respect to the symplectic form induced by the cup product. Up to sign the difference of the signatures
$\sigma(M)-\sigma(M_1)-\sigma(M_2)$ equals the Maslov index of these subspaces.

The $\eta$-invariant associated to a pair $(L_0,L_1)$ of Lagrangian subspaces  is defined as the $\eta$-invariant of the operator $D_I=I_0 \frac {d}{dx}$ on $L^2([0,1],\bbbr^{2n})$ with boundary conditions $f(0) \in L_0,~f(1) \in L_1$. Here $I_0$ is the skewsymmetric matrix representing the symplectic form $\omega$ with respect to the standard scalar product on $\bbbr^{2n}$. Then 
 $$\eta(L_0,L_1)=\frac{1}{\sqrt{\pi}}\int_0^{\infty}t^{-\frac 12} \Tr~ D_I e^{-tD_I^2} dt \ .$$ It is a regularization of the difference between the number of positive and negative eigenvalues of $D_I$. The $\eta$-invariant $\eta(L_0,L_1)$ occurs as a correction term in a gluing theorem for $\eta$-invariants \cite{b}. The formula 
$$\tau(L_0,L_1,L_2)=\eta(L_0,L_1)  + \eta(L_1,L_2) + \eta(L_2,L_0) \ , $$ which will be generalized to a noncommutative context, 
can be interpreted as reflecting the connection between the gluing problems for signatures and $\eta$-invariants via the Atiyah-Patodi-Singer index theorem.

The noncommutative Maslov index and the $\eta$-forms should play a similar role in gluing problems for higher signatures for manifolds with boundaries \cite{lp4}, \cite{llp}, \cite{llk}, higher $\eta$-invariants and $\rho$-invariants \cite{lo3}. In \S \ref{glue} we explain how our formula is related to a gluing formula for noncommutative $\eta$-forms for Dirac operators on the circle.

The proof by Cappell, Lee and Miller relies on the axiomatic
properties of the Maslov index, which are not expected to hold for a noncommutative generalization.
Bunke and Koch \cite{bk} formulated an index problem  whose solution, an Atiyah-Patodi-Singer type index theorem, yields the above formula. They interpreted the Maslov index as the index of a Dirac operator on a two-dimensional manifold with boundary and six cylindric ends isometric to $\bbbr^+ \times [0,1]$. The Lagrangian subspaces enter in the definition of the boundary conditions. The purpose of their approach was to prove a formula for families of Lagrangian subspaces.    

Our proof is a generalization of the approach of Bunke and Koch. Instead of a family of Dirac operators we consider a Dirac operator twisted by a $C^*$-vector bundle \cite{mf}. \\

Before stating the result we give a short introduction in the type of noncommutative index theorem considered in this paper: 
  
Family index theorems describe Fredholm operators  depending continuously on a parameter from
some compact space. The index of a family of operators  is an element in
$K$-theory of the base space. If the kernel and the cokernel are vector bundles, then the index is
 the difference of the classes of these bundles.

One may reformulate this situation by replacing the base space $B$ by the $C^*$-algebra
$C(B)$ and the family of operators by an operator on a $C(B)$-module. The index 
is then in the $C^*$-algebraic $K$-theory $K_0(C(B))$, which is naturally isomorphic to $K^0(B)$. 

Mi\v s\v cenko and Fomenko \cite{mf} formalized and generalized this point of view by defining Fredholm operators on
Hilbert $C^*$-modules for general $C^*$-algebras. The index of such a Fredholm operator is an element in the $K$-theory of the $C^*$-algebra. They also elaborated a theory of pseudodifferential operators over $C^*$-algebras. Important examples for geometric applications are Dirac operators associated to $C^*$-vector bundles.  

In order to formulate a noncommutative analogue of the family index theorem, which calculates the Chern character of the index of a family of Dirac operators in the 
cohomology of the base space, now assumed to be a manifold, analogues of differential forms, de Rham cohomology and the Chern character are needed. 
Karoubi \cite{kar} introduced a complex of differential forms $(\Oi\A,\di)$ associated to a
Fr\'echet algebra $\A$, its de Rham homology $H^{dR}_*(\A)$ and a Chern character
 $\ch:K_0(\A) \to
H^{dR}_*(\A)$, which is defined if in addition $\A$ is a local Banach algebra.  Unfortunately the de Rham homology of
a $C^*$-algebra does not behave well, in particular the de Rham homology of
a commutative unital $C^*$-algebra is in general not the
cohomology of the corresponding compact space. By considering the
de Rham homology of the
algebra $\C(B)$ instead of $C(B)$ one recovers the Chern character from differential geometry.
This can be interpreted as reflecting the fact that a differentiable structure on $B$ is needed for the definition of de Rham cohomology. 

In order to get a reasonable Chern character in the general situation one chooses a dense subalgebra $\Ai$ holomorphically closed under the functional calculus in $\A$.
Then $K_0(\Ai)$ is
canonically isomorphic to $K_0(\A)$, and the Chern character yields a
homomorphism $$\ch:K_0(\A) \to H^{dR}_*(\Ai) \ .$$
In this setting homological index theorems for Dirac operators associated to $C^*$-vector bundles can be formulated and -- at least formally -- a noncommutative version of Bismut's family index theorem makes sense \cite{bgv}. Note that formally the noncommutative situation is less complex than the family case since the Riemannian metric is fixed.

For the generalization of the heat kernel theory additional conditions have to be imposed on the algebra $\Ai$: 
It should be the limit of a projective system $\{\A_i\}_{i \in
  \bbbn_0}$  of involutive Banach algebras  with
$\A_0=\A$, such that there are dense embeddings $\A_{i+1} \incl \A_i$  and such that each $\A_i$ is closed with respect to the holomorphic functional
  calculus in $\A$. The motivating example is $\A_i:=C^i(B)$ for a closed manifold $B$. 
More generally, if $\delta$ is an involutive closed derivation on $\A$ with $\cap_i \dom \delta^i$ dense in $\A$,
then the projective system given by $\A_i:=\dom \delta^i$ with norm $$\|a\|_i:=\sum_{j =0}^i \|\delta^j(a)\|$$ fulfills the conditions. 
For a discrete finitely generated group $G$ the unbounded operator $D$ on $l^2(G)$ defined by $D1_g=l(g)1_g$, where $l$ is a word length function on $G$, induces a closable derivation $\delta(f)=[D,f]$ on $B(l^2(G))$. Furthermore $\bbbc G \subset \A_i:=\dom \ov{\delta}^i \cap C^*_r(G)$. The projective limit $\Ai \subset C^*_r(G)$ is closely related to the algebra employed by Connes and Moscovici in their proof of the Novikov conjecture for Gromov hyperbolic groups via higher index theory \cite{cm}. Using this setting our results can be applied to higher index theory. A proof of a higher Atiyah-Singer index theorem based on heat kernel methods was given by Lott \cite{lo1} and a higher Atiyah-Patodi-Singer index theorem was proved by Leichtnam and Piazza \cite{lp1} \cite{lp2}. Lott introduced higher $\eta$-invariants and $\rho$-invariants \cite{lo3}.  A motivation for higher index theory is the study of higher signatures and the Novikov conjecture. Index theory for Dirac operators over $C^*$-algebras in general has been applied to the study of manifolds with positive scalar curvature \cite{ps} \cite{ros}.

The proofs of the higher index theorems mentioned above rely on the comparison with the situation on the covering space where one can deal with operators on sections of complex vector bundles.
In the general setting the main difficulty lies in the fact that the calculus of regular operators on a Hilbert $C^*$-module and the calculus of pseudodifferential operators over a $C^*$-algebra are not sufficient for the study of the heat semigroup since we have to deal with vector bundles whose fibers are projective $\Ai$-modules. By proving that the de Rham homology behaves well under the projective limit we justify the fact that one can deal with Banach spaces  instead of Fr\'echet spaces. Then the theory of holomorphic semigroups can be used instead of the calculus of selfadjoint operators. We define appropriate operator spaces on $L^2$-spaces of vector valued functions, for example Hilbert-Schmidt operators and trace class operators. Duhamel's principle is used for the construction and the short time asymptotics of the heat kernels. The long time asymptotics is more intricate: here we have to get hold of the spectrum of the Dirac operator. We adapt a method developed by Lott, namely a restricted pseudodifferential operator calculus giving information about the resolvent set of the Dirac operator and the regularizing properties of the resolvents \cite{lo2}.  

Although the theory is developed only for a particular two-dimensional manifold and a trivial vector bundle, the relevant parts generalize to the Atiyah-Patodi-Singer index problem for Dirac operators over $C^*$-algebras, which will be considered elsewhere. \\

Now we can formulate the noncommutative version of the equation relating the Maslov index and the
$\eta$-invariants.
 
Let $\bbbc^{2n}$ be endowed with the skewhermitian form $\omega$ induced by the standard
symplectic form on $\bbbr^{2n}$ via the identification $\bbbc^{2n}=\bbbr^{2n} \ten_{\bbbr}\bbbc$.

For a $C^*$-algebra $\A$ a Lagrangian projection on $\A^{2n}$
is a selfadjoint projection $P \in M_{2n}(\A)$ fulfilling $PI_0=(1-I_0)P$. Two projections are called transverse if their sum is invertible.
 
As above one associates an $\A$-valued hermitian form $q$ to every triple $(\Pj_0,\Pj_1,\Pj_2)$ of pairwise
transverse Lagrangian projections. The class $\tau(\Pj_0,\Pj_1,\Pj_2) := [q] \in K_0(\A)$ is called Maslov index of $(\Pj_0,\Pj_1,\Pj_2)$.

An $\eta$-form $\eta(\Pj_0,\Pj_1) \in \Oi\Ai/\ov{[\Oi\Ai,\Oi\Ai]_s}$ can be associated to a pair of transverse Lagrangian projections $(\Pj_0,\Pj_1)$ with
$\Pj_i \in M_{2n}(\Ai),~i=0,1$.  Here  $[~,~]_s$ denotes the supercommutator. The $\eta$-form is defined
via a superconnection associated to the operator $D_I=I_0 \frac{d}{dx}$ on $L^2([0,1],\A^{2n})$ with boundary conditions $f(0)=\Pj_0f(0)$ and $f(1)=\Pj_1f(1)$. 

Then our main result is:

\begin{theor} For a triple $(\Pj_0,\Pj_1,\Pj_2)$ of pairwise
transverse Lagrangian projections with $\Pj_i \in M_{2n}(\Ai),~i=0,1,2,$
$$\ch\tau(\Pj_0,\Pj_1,\Pj_2)= [\eta(\Pj_0,\Pj_1) + \eta(\Pj_1,\Pj_2)
  +\eta(\Pj_2,\Pj_0)] \in H_*^{dR}(\Ai) \  .$$
\end{theor}

\vspace{0.5cm}

This paper is based on the author's PhD thesis.
I would like to thank  my supervisor
Ulrich Bunke for drawing my attention to the problem and for fruitful discussions, and Margit R\"osler for the introduction into the theory of
semigroups.

\section*{Summary}

The paper is organized in the following way:

As mentioned above, the Maslov index $\tau(\Pj_0,\Pj_1,\Pj_2)$ is the index of  a Dirac operator $D^+$ twisted by a $C^*$-vector bundle on
a two-dimensional spin manifold $M$ with six cylindric ends 
isometric to $[0,\infty)  \times [0,1]$. \\

In Chapter 1 the manifold $M$ is described and we explain in more detail the family index theorem of Bunke and
Koch \cite{bk}. Then the universal differential algebra $\Oi\Ai$, the de Rham homology $H^{dR}_*(\Ai)$ and the
  Chern character are introduced and investigated. It is shown that $H^{dR}_*(\Ai)$ is the projective limit of $H^{dR}_*(\A_i)$. Moreover Lagrangian
  projections and the Maslov index are defined.
  
In Chapter 2 we introduce the Dirac operator $D=D^+ \oplus D^-$ on $M$
whose  boundary conditions are defined by a triple of Lagrangian projections $(\Pj_0,\Pj_1,\Pj_2)$. Furthermore the operator $D_I$ is defined and its properties on the Hilbert
$\A$-module $L^2([0,1], \A^{2n})$ are studied. Then we show that $D^+$ is
Fredholm  between appropriate Hilbert $C^*$-modules and that its index
equals $\tau(\Pj_0,\Pj_1,\Pj_2)$. We define a compact perturbation $D(\rho)$ of $D$ such that $D(\rho)^+$ is surjective. Then the index of $D^+$ can be expressed in terms of the kernel of $D(\rho)$. 

Chapter 3 is devoted to heat semigroups and their integral kernels, in
particular those associated to $D_I$ and $D(\rho)$.

In Chapter 4 we introduce superconnections in order to define the
$\eta$-form. Now we can formulate the index theorem. The remainder of the
chapter is devoted to its proof. We introduce the rescaled superconnection $A(\rho)_t$ associated to $D(\rho)$ and study the family of operators $e^{-A(\rho)_t^2}$. We show that it is a
family of integral operators with smooth integral kernel and obtain estimates
for the integral kernel
for small $t$. Once the heat kernel theory is established, the proof of the index theorem itself is fairly standard. We follow the proof in \cite{bk}, which is modelled on Melrose's $b$-calculus \cite{me}, and compare the limit of a generalized supertrace of $e^{-A(\rho)_t^2}$ 
for $t \to \infty$, which is the Chern character of the index of $D^+$, with
its limit for $t \to 0$. Since the differences between the calculus for Dirac operators associated to complex vector bundles and the one for Dirac operators associated to a projective system of vector bundles as developed here are  subtle, the proof is given in detail.

In Chapter 5 the functional analytic framework is developed. The function spaces we deal with are introduced, as for example the
$L^2$-spaces of vector valued functions, and operators on them are studied. In particular we define and study appropriate notions of adjointable operators, Hilbert-Schmidt
 and trace class operators. Furthermore we recall the properties of  Fredholm
operators and regular operators on Hilbert $C^*$-modules, and collect some facts from holomorphic semigroup theory. 
The reader is advised to go through the definitions of this chapter first in order to get acquainted with the functional analytic setting.

\section*{Notation and conventions}
If not specified vector spaces and algebras are complex, manifolds are smooth.

We often deal with $\bbbz/2$-graded spaces. Then $[~,~]_s$ denotes
the supercommutator and $\trs$ the supertrace. In a graded context the tensor
products  are graded. For an ungraded vector space $V$, we denote by $V^+$ resp. $V^-$ the same
space endowed with a grading: all elements are homogeneous of positive resp. negative degree.

Tensor products denoted by $\otimes$ are completed. The way of completion is
indicated by 
a suffix in all but the two most common cases: In the case of
Hilbert $C^*$-modules $\ten$ means the Hilbert $C^*$-module tensor product,
and if one  of the spaces is a
nuclear locally convex space, then $\ten$ means $\otimes_{\pi}$ or $\otimes_{\ve}$.  The algebraic tensor product
is denoted by $\odot$.

By a  differentiable function on an open subset of
$[0,1]^n$ we understand a function that can be extended to a differentiable function on an open subset
of $\bbbr^n$. This induces the notion of a differentiable function on a manifold
with corners, in our case $M \times M$.

If $S,X$ are sets with $S \subset X$, then  the
characteristic function of $S$ is denoted by $1_S:X \to \{0,1\}$. If $X$ is a metric space, $y \in X$ and $S \subset X$ then
$d(y,S):=\inf\limits_{x \in S} d(x,y)$. For $S_1,S_2 \subset X$ we set $d(S_1,S_2):=\inf\limits_{x \in S_1} d(x,S_2)$.

If $E_i, ~i=1,2,$ is a vector bundle on a space $X_i,~i=1,2,$ and $p_i:X_1 \times X_2 \to X_i$ is the projection,  then $E_1 \boxtimes
E_2=p_1^*E_1 \ten p_2^*E_2$ on $X_1 \times X_2$.

We use the notions from \cite{bgv} in the context of spin geometry. 
In addition let a Dirac bundle be a selfadjoint Clifford module endowed with a
Clifford connection with respect to which the metric is parallel. Our sign conventions differ from those in \cite{bgv} since we deal with right modules over the algebra of differential forms.

The value of the constant $C$ used in estimates
may vary during a series of estimates without an explicit remark. 

\chapter{Preliminaries}

\section{The geometric situation}
\label{situat}

In this section  the two-dimensional spin manifold $M$ with boundary and cylindric ends, a Dirac bundle $E$ on it and the associated Dirac operator will be introduced. The open covering ${\mathcal U}(r,b)$ of $M$ constructed in this section will be used for cutting and pasting arguments later on. Furthermore we will fix a flat open set $F \subset M$ containing the boundary and the cylindric ends and trivializations of $TM|_F$  and of $E|_F$, which will be used in the definition of the boundary conditions for the Dirac operator in \S \ref{defiD}.\\

We begin by defining the manifold $M$.\\

For $k \in \bbbz/6$ let $Z_k$ be a copy of $\bbbr \times [0,1]$. With the euclidian metric and the standard
orientation $Z_k$ is an oriented Riemannian
manifold with boundary. Let
$(x^k_1,x^k_2)$ be the euclidian coordinates of $Z_k$.\\

For $r\ge -\frac 12$, $b \le \frac 13$ and $k \in \bbbz/6$ let 
$$F_k(r,b):=\{(x^k_1,x^k_2) \in~ ]r, \infty[ \times [0,1]
~\cup~ ]-1,r] \times \bigl([0, b[ \, \cup \, ]1-b,1]\bigr)\} \subset Z_k \ .$$ 
We define   
$$F(r,b):= \Bigl(\bigcup_{k \in \bbbz/6} F_k(r,b) \Bigr)/\sim$$
with $(x^k_1,x^k_2) \sim (-\frac 32 -x^{k-1}_1,1-x^{k-1}_2)$ for $(x^k_1,x^k_2) \in~ ]-1,-\frac 12[ \times
[0, b[$ and $k \in \bbbz/6$. \\

Then $F(r,b)$ inherits the structure of
an oriented Riemannian manifold from the sets $F_k(r,b)$.\\

The set $F(-\frac 12,\frac 13) \setminus \ov{F(-\frac 13,\frac 14)}$ is
diffeomorphic to the open ring
$B_1(0)^{\circ} \setminus B_{1/2}(0) \subset \bbbr^2$ via an oriented diffeomorphism $\phi$.
We define the manifold with boundary $$M:=F(-\tfrac 12,\tfrac 13) \cup_{\phi} B_1(0)^{\circ} \ .$$
For $r>-\frac 12$, $b \le \frac 13$ we identify 
$F(r,b)$ and $F_k(r,b)$ with the corresponding subsets in $M$.  The sets $F_k(r,b)$ are coordinate patches of $M$ with the coordinates $(x_1^k,x_2^k)$ from above.\\

Extend the orientation and metric from $F:=F(0,\frac 14)$  to the whole of
$M$ and endow $TM$ with the Levi-Civit\`a connection. In the following we identify $TM$ and $T^*M$.\\

For $r \ge 0$ and $b \le \frac 14$ we define an open covering ${\mathcal U}(r,b)=\{\U_k\}_{k \in J}$ of $M$ as follows.
Let $J$ be the union of $\bbbz/6$ with a one-element set $\{\cp\}$. For
$k \in \bbbz/6$ let $\U_k:=F_k(r,b)$ and let $\U_{\cp}:=M \setminus \ov{F(r+1,b/2)}$.\\

For $r\ge 0$ let $M_r:=M \setminus \ov{F(r,0)}$.\\

The connected components of $\ra M$ are labelled $\ra_k M,~ k \in \bbbz/6,$ in such a way
that $\ra_k M \cap F_k(r,b)\subset \{x_2^k=0\}$ and $\ra_{k+1}M \cap
F_k(r,b) \subset \{x_2^k=1\}$.\\

The manifold $M$ can be embedded diffeomorphically into $\bbbr^2$, even with a diffeomorphism that
is an isometry on the complement of $M_r$ for some $r>0$. The image of the embedding is
illustrated by figure 1. 

\begin{figure}
\includegraphics{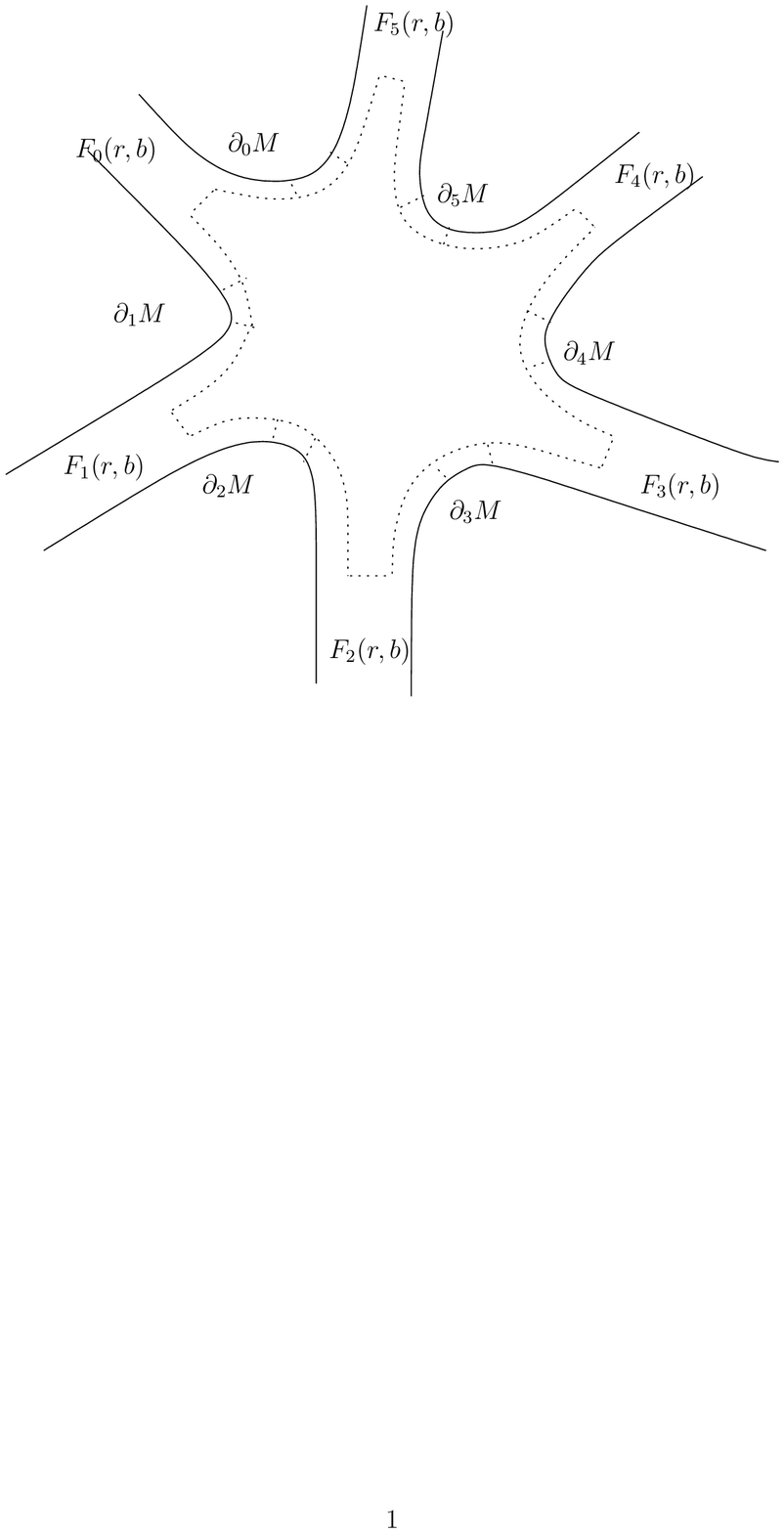}
\caption{The manifold $M$}

\end{figure}

Now we define the Dirac operator. Choose a spin structure on $M$ and fix $d \in \bbbn$. Let $S$ be the spinor
bundle endowed with the Levi-Civit\`a connection and a
parallel hermitian metric. Let $E$ be the
graded vector bundle $S \ten ((\bbbc^+)^d \oplus (\bbbc^-)^d)$. 
The hermitian metric on $S$
and the standard hermitian product on $\bbbc^d$  induce a hermitian metric
$\langle\cdot, \cdot \rangle $ on
$E$. 



Furthermore there are connections on $E$ and its dual $E^*$ induced by the Clifford connection
on $S$ and the de Rham differential. The Dirac operator
associated to the Dirac bundle $E$ is denoted by $\dirac$.\\

The oriented orthonormal frames $((-1)^kdx_1^k,(-1)^kdx_2^k)$ of $TM|_{F_k(0,
\frac 14)}$ patch together to an oriented orthonormal frame $(e_1,e_2)$ of
 $TM|_F$. The Clifford multiplication $c:TM \to \End E$ induces an even parallel endomorphism $$I=c(e_2)c(e_1):E|_F \to E|_F$$ defining a
skew-hermitian form  $$E|_F \times E|_F \to \bbbc,~(x,y)
\mapsto \langle x,Iy\rangle  \ .$$ 
Since the holonomy of $TM|_F$ is $4\pi \bbbz$ (measured with respect to any trivialization of
$TM$ on the whole of $M$), there are nonvanishing parallel sections of $S^+|_F$ and
$S^-|_F$. We choose a parallel unit section $s$ of $S^+|_F$ 
and  fix
once and for all the
 trivialization of $E|_F$  defined for $x \in F$ by
$$E_x^+=S_x^+ \ten (\bbbc^+)^d  \oplus S_x^- \ten (\bbbc^-)^d \to
(\bbbc^+)^{d} \oplus (\bbbc^+)^{d}\ ,$$
$$(s(x) \ten v) \oplus (ic(e_1)s(x) \ten w)  \mapsto (v,w) \ ,$$
and $$E_x^-=S_x^-\ten (\bbbc^+)^d  \oplus S_x^+ \ten (\bbbc^-)^d  \to
(\bbbc^-)^{d} \oplus (\bbbc^-)^{d} \ ,$$
$$(ic(e_1)s(x)\ten v) \oplus (s(x) \ten w)\mapsto (v,w) \ .$$
With respect to this trivialization the endomorphism $I|_{E^+}$  corresponds to
$$I_0:=\left(\begin{array}{cc} i &0 \\
0 & -i 
\end{array}\right) \in M_{2d}(\bbbc)$$
and $I|_{E^-}$ to $-I_0$.

\section{The family index theorem}
\label{commind}

In order to give a motivation for the definitions in the subsequent sections we sketch the corresponding family index theorem. The presentation follows \cite{bk}. \\
 
Let $\bbbc^{2d}$ be endowed with the standard hermitian product $\langle ~,~\rangle $ and the
skew-hermitian form $(x,y) \mapsto \langle x,I_0y\rangle $. 

Let $B$ be a compact space and let $(L_0,L_1,L_2)$ be a triple of pairwise transverse Lagrangian
subbundles of the trivial bundle $B \times \bbbc^{2d}$. For any $b
\in B$ the Lagrangian subspaces $L_i(b) \subset \bbbc^{2d}$
define parallel Lagrangian subbundles of $E^+|_F$ via the trivialization fixed
in the previous section. 

Let $D^+(b)$ be the Dirac operator associated to $E$ with
$$\dom D^+(b):=\{s \in \C_c(M,E^+)~|~ s(x) \in L_i(b) \text{ for } x \in \ra_i M \cup
\ra_{i+3}M,~i=0,1,2\} \ .$$
It turns out that for any $b \in B$ the kernel and cokernel of the closure of
$D^+(b)$ are finite dimensional and that the family
$\{D^+(b)\}_{b \in B}$ has a well-defined index in $K^0(B)$, which equals the generalized Maslov index of
$(L_0,L_1,L_2)$ defined as follows.

The triple $(L_0,L_1,L_2)$ induces a nondegenerate hermitian form $h$ on $L_0$: Let $v=v_1+v_2,~w=w_1+w_2 \in L_0(b)$ with $v_1,w_1
\in L_1(b);~ v_2,w_2 \in L_2(b)$, then 
$$h_b(v,w):=\langle v_2,I_0w_1\rangle  \ .$$
The generalized Maslov index is the element $[L_0^+]-[L_0^-] \in K^0(B)$ where
$L_0^+$ and $L_0^-$ are subbundles of $L_0$ with $L_0^+ \oplus L_0^-=L_0$ and
such that $h$ is positive on $L_0^+$ and negative on $L_0^-$.

If $B$ is a manifold and the
bundles are smooth, then by a generalization of the Atiyah-Patodi-Singer index theorem
the Chern character of the index bundle
can be expressed in terms of $\eta$-forms; the local term vanishes since $\ch(E/S)=0$. 

We outline the definition of the $\eta$-forms.

Let $\ra$ be the differentiation operator on $\C([0,1],\bbbc^{2d})$.

For $i \neq j$ and any $b \in B$ the operator $D_I(b):=I_0 \ra$ with domain 
$$\dom D_I(b):=\{s \in \C([0,1],\bbbc^{2d})~|~ s(0) \in L_i(b),~s(1) \in L_j(b) \}$$
is essentially selfadjoint and its closure has a bounded inverse on
$L^2([0,1],\bbbc^{2d})$.

There is a family of rescaled superconnections $A^I_t$
 associated to the family of operators $\sigma D_I$, where $\sigma$ is a formal parameter of degree $1$
 with $\sigma^2=1$. 

The $\eta$-form
$$\eta(L_i,L_j):= \frac{1}{ \sqrt \pi}\int_0^{\infty} t^{-\frac 12} \Trsi
D_Ie^{-(A^I_t)^2} dt \in \Omega^*(B) \ ,$$
with $\Trsi(a+\sigma b):=\Tr(a)$, is well-defined.

The statement of the index theorem is:
$$\ch(\ind D^+)= [\eta(L_0,L_1) + \eta(L_1,L_2) + \eta(L_2,L_0)] \in
H^*_{dR}(B) \ .$$

\vspace{0.3cm}

\section{The algebra of differential forms}
\label{univdg}

In this section we 
 study the noncommutative analogues of the algebra $\Omega^*(B)$, the Chern
character and the de Rham cohomology $H^*_{dR}(B)$. \\

\subsection{The universal graded differential algebra}

Let ${\mathcal B}$ be an involutive locally $m$-convex Fr\'echet algebra with
unit. In particular, the multiplication $\B \times \B \to \B,~(a,b) \mapsto ab$ is continuous, the group of invertible elements $\Gl(\B)$ in $\B$ is open and the map
$\Gl(\B) \to \Gl(\B),~a \mapsto a^{-1}$ is continuous \cite{ma}. In this section we recall the definition of the
topological universal graded differential algebra $\Oi\B$ and collect its main
properties  \cite{kar} \cite{cq}.

We write $\pten$ for the completed projective tensor product.

Let
$$\Ok {\mathcal B}:={\mathcal B} \pten (\pten^k({\mathcal B}/\bbbc)) \
$$ and $$\Oi {\mathcal B} :=\prod\limits_{k=0}^{\infty} \Ok {\mathcal B} \ .$$
With the following structures $\Oi{\mathcal B}$ is an involutive 
Fr\'echet locally $m$-convex graded differential algebra:

{\bf Product:} There is a graded continuous product on  $\Oi {\mathcal B}$ defined
for elementary tensors by
$$(b_0 \ten b_1 \ten \dots \ten b_k) (b_{k+1} \ten b_{k+2} \ten \dots
b_n):=\sum\limits_{j=0}^{k} (-1)^{k-j}   (b_0 \ten b_1 \ten \dots \ten b_jb_{j+1} \ten b_{j+2} \ten \dots
b_n) \ .$$

{\bf Differential:} There is a continuous differential $\di$ of degree one on the graded algebra $\Oi{\mathcal
B}$ defined on elementary tensors by 
$$\di (b_0 \ten b_1 \ten \dots \ten b_k):= 1 \ten b_0 \ten b_1 \ten \dots \ten
b_k \ .$$
It satisfies the graded Leibniz rule: For $\alpha \in \Ok\B$ and $\beta \in \Oi\B$  
$$\di (\alpha \beta)=(\di\alpha)\beta +(-1)^k \alpha(\di\beta) \ .$$
Furthermore in $\Ok\B$
 $$b_0 \ten b_1 \ten \dots \ten b_k=b_0\di b_1 \di b_2 \dots \di b_k \ .$$
If ${\mathcal B}$ is a Banach algebra, then $\di$ is a map of norm one.

{\bf Involution:} We extend the $*$-operation on $\B$ to a continuous involution
on $\Oi {\mathcal B}$
by  setting
$$(b_0 \ten b_1 \ten \dots \ten b_k)^*:=(1 \ten b_k^* \ten b_{k-1}^* \ten
\dots \ten b_1^*)b_0^* \ $$
or equivalently 
$$(b_0\di b_1 \di b_2 \dots \di b_k)^*=(\di b_k^* \di b_{k-1}^* \dots \di b_1^*) b_0^* \ .$$
For $\omega_1,\omega_2 \in \Oi {\mathcal B}$ 
$$(\omega_1 \omega_2)^*=\omega_2^*\omega_1^* \ ,$$
and for $\omega \in \Ok{\mathcal B}$ 
$$(\di \omega)^*=(-1)^k \di (\omega^*) \ .$$

Let $$\hat\Omega_{\le m}\B:=\Oi\B/\prod\limits_{k=m+1}^{\infty}
\Ok\B \ .$$ 
The above structures are well-defined on $\hat\Omega_{\le m}\B$ as well.

We identify $\hat\Omega_{\le m}\B$ as a graded vector space with the subspace $$\{\omega\in \Oi\B ~|~\omega^n=0 \text{ for } n>m\} \subset \Oi\B \ , $$ where $\omega^n$ denotes the homogeneous part of degree $n$ of $\omega \in \Oi\B$.

The Fr\'echet space $\ov{[\Oi\B,\Oi\B]_s}$, generated by the supercommutators in
$\Oi\B$, is preserved by $\di$ by  Leibniz rule. It follows that
$(\Oi\B/\ov{[\Oi\B,\Oi\B]_s},\di)$ is a complex.

\begin{ddd}  The {\sc de Rham
homology of $\B$} is 
$$H^{dR}_*(\B):=H_*(\Oi\B/\ov{[\Oi\B,\Oi\B]_s},\di) \ .$$
On the right hand side we take the topological homology, i.e. we
quotient out the closure of the range of $\di$ in order to obtain a Hausdorff space.
\end{ddd}

The map $\di$ induces maps $\di:(\Oi\B)^n \to (\Oi\B)^n$ and $\di:M_{n}(\Oi{\mathcal B}) \to
M_{n}(\Oi{\mathcal B}),~ A \mapsto \di(A)$ by applying $\di$ componentwise. Note the difference between $\di A= \di \circ A$ and
$\di(A)$. Sometimes we write $(\di A)$ for $\di(A)$.

For $A \in M_n(\Ok\B)$ 
$$\di A=(\di A) +(-1)^kA \di \ .$$ 
The trace $$\tr:M_n(\Oi{\mathcal B}) \to \Oi\B/\ov{[\Oi\B,\Oi\B]_s}$$ is defined by adding up the diagonal elements. It vanishes on supercommutators.\\

In the following definition and proposition a projection is not assumed to be selfadjoint. 

\begin{ddd}
\label{defch}
Let $P \in M_n(\B)$ be a projection. Then 
$$\ch(P):= \sum\limits_{k=0}^{\infty}(-1)^k  \frac{1}{k!} \tr P(\di (P))^{2k} \in \Oi\B/\ov{[\Oi\B,\Oi\B]_s}$$
 is the {\sc Chern character form of $P$}.
\end{ddd}

\begin{prop}
\label{homch}
\begin{enumerate}
\item The Chern character form is closed.
\item If $P:[0,1] \to M_n(\B)$ is a differentiable path of projections, then 
$\ch(P_1)-\ch(P_0)$ is exact.
\item If $\B$ is a local Banach algebra \cite{bl}, then the Chern character form
induces a homomorphism
$$\ch: K_0(\B) \to H^{dR}_*(\B),~ \ch([P]-[Q]):=\ch(P)-\ch(Q) \ ,$$
called the {\sc Chern character}.
\end{enumerate}
\end{prop}

\begin{proof}
First note that for a projection $P \in M_n(\B)$
$$0=\di((1-P)P)=(1-P)(\di P) - (\di P)P \ , $$
hence $(1-P)(\di P)=(\di P)P$ and $P(\di P)P=0$. Therefore $$P(\di P)^2=(\di P)^2P\ .$$
Analogous formulas hold for the derivative of a path of projections.

(1) 
follows from
\begin{eqnarray*}
\di \tr P(\di P)^{2k}&=&\tr(\di P)^{2k+1}\\
&=& \tr(1-P)(\di P)^{2k+1}(1-P) +\tr P(\di P)^{2k+1}P\\
&=& 0 \ .
\end{eqnarray*}

(2) 
We have that
\begin{eqnarray*}
(\tr P(\di P)^{2k})'
 &=& \tr P'(\di P)^{2k}+  \tr P((\di P)^{2k})' \\
&=& \sum\limits_{i=0}^{2k-1} \tr P(\di P)^i(\di P')(\di P)^{2k-i-1} \ .
\end{eqnarray*}
For $i$ even 
\begin{eqnarray*}
\lefteqn{\tr P(\di P)^i(\di P')(\di P)^{2k-i-1}}\\
&=& \tr (\di P)^i (\di (PP'))(\di P)^{2k-i-1}- \tr (\di P)^i (\di P)P'(\di P)^{2k-i-1}\\
&=&\tr (\di P)^i (\di (PP'))(\di P)^{2k-i-1}\\
&=&\di \tr P(\di P)^{i-1}(\di (PP'))(\di P)^{2k-i-1} \ .
\end{eqnarray*}
For $i$ odd the argument is similar.

(3) follows from (1) and (2).
\end{proof}

\subsection{Supercalculus}
\label{sutrace}
Let $\B$ be as before.

In this section all spaces and tensor products are  $\bbbz/2$-graded. If no grading is
specified we assume the grading to be trivial.

Let $V=V^+ \oplus V^-$ be a $\bbbz/2$-graded complex vector space with $\dim
V^+=m,~\dim V^-=n$ and consider $\Oi\B$ as a $\bbbz/2$-graded space with the grading induced by the form degree.

The space $V \ten \Oi\B$ is a free
$\bbbz/2$-graded right $\Oi\B$-module. It is furthermore a left supermodule of
the superalgebra $\End(V) \ten \Oi\B$. (Note that our setting differs from the corresponding one in \cite{bgv} since we deal with
right $\Oi\B$-modules. This leads to different signs.)

The supertrace $\trs:\End(V) \to \bbbc$ extends to a supertrace
$$\trs:\End(V) \ten \Oi\B  \to  \Oi\B/\ov{[\Oi\B,\Oi\B]_s} $$
$$\trs (T \ten \omega):= \trs(T) \omega \ .$$
By quotienting out the supercommutator we ensure that $\trs([T_1,T_2]_s)=0$
for $T_1, T_2 \in \End(V) \ten \Oi\B$.

The differential $\di$ acts on elements of $V \ten \Oi\B$ resp.  $\End(V) \ten \Oi\B$ by 
$$\di(A \ten \omega)=(-1)^{\deg A} A \ten \di \omega$$ for $A \in V^{\pm}$
resp. $A \in
\End^{\pm}(V)$ and $\omega \in \Oi\B$. 

Though the differential $\di$ is not a right $\Oi\B$-module map, the supercommutator $[\di,T]_s$ with $T \in \End(V) \ten \Oi\B$ is in  $\End(V) \ten \Oi\B$; namely 
$$[\di,T]_s= \di(T) \ .$$ 
Since $\trs(A \ten \omega)=0$ for $A\in
\End^-(V)$, 
we have that
$$\trs[\di,T]_s= \trs \di(T)=\di \trs T \ $$ for $T \in \End(V) \ten \Oi\B$. \\

See \S \ref{adop}, \S \ref{HS}, \S \ref{trclop} for notions of Hilbert-Schmidt operators, trace class operators and adjointable operators used in the following. 
 
If $M$ is a complete Riemannian manifold and $T$ is a trace class operator on $L^2(M,V \ten \hat\Omega_{\le \mu}\B)$, then
$$\Trs T:=\int_M \trs ~\ov R(T)(x)~dx \in
\hat\Omega_{\le \mu}\B/\ov{[\hat\Omega_{\le \mu}\B,\hat\Omega_{\le \mu}\B]_s} \ .$$
Note that $$\Trs[A,B]_s=0$$ if $A$ is an adjointable operator and $B$ a trace class operator or if $A$, $B$ are Hilbert-Schmidt operators.

\subsection{The algebras $\Ai$ and $\Oi\Ai$}
\label{einfai}
\label{einfOi}

Let $(\A_j,\iota_{j+1,j}:\A_{j+1} \to \A_j)_{j \in \bbbn_0}$ be a projective system of
involutive Banach algebras with unit satisfying the following conditions
(see \cite{lo2}, \S 2.1):

\begin{itemize} 
\item The algebra $\A:=\A_0$ is a $C^*$-algebra.
\item For any $j \in \bbbn_0$ the map $\iota_{j+1,j}:\A_{j+1} \to \A_j$ is injective.
\item For any $j \in \bbbn_0$ the map $\iota_j:\Ai:= \pl{i}\A_i \to \A_j$ has dense
range.
\item For any $j \in \bbbn_0$ the algebra $\A_j$ is stable with respect to the holomorphic functional calculus in
$\A$.  
\end{itemize}

The projective limit $\Ai$ is an involutive  locally
$m$-convex Fr\'echet algebra with unit.  

The motivating example is $\A_j=C^j(B)$ for a closed smooth manifold $B$.

\begin{prop}
\label{propai}

For  $j \in \bbbn_0$ and $n \in \bbbn$ we have:
\begin{enumerate}
\item The map $\iota_j:\Ai \to \A_j$ is injective. 

\item The algebras $M_n(\Ai)$ and $M_n(\A_j)$  are stable with respect to the holomorphic
functional calculus in $M_n(\A)$. 

\item The map $\iota_{0*}:K_0(\Ai) \to K_0(\A)$ is an isomorphism.
\end{enumerate}
\end{prop}

\begin{proof}
(1) follows immediately.  

(2)  follows from \cite{bo},
Prop. A.2.2. 

(3) follows from \cite{bo}, Th. A.2.1. 
\end{proof}

The
projective system $(\A_j,\iota_{j+1,j})_{j \in \bbbn_0}$ induces two projective
systems of involutive graded differential Fr\'echet algebras.

One of them is given by the maps
$$\iota_{j+1,j*}:\Oi\A_{j+1} \to \Oi\A_j \ .$$
Furthermore $(\hat\Omega_{\le m}\A_j)_{m,j \in \bbbn}$ is a projective
system of involutive Banach graded differential algebras. 

Their limits coincide:

The inclusion $\hat\Omega_{\le
m}\A_j  \to \Oi\A_j$
is  left inverse to the projection $\Oi\A_j \to \hat\Omega_{\le
m}\A_j$. The
induced maps between the projective limits are inverse to each other. It follows
that $$\pl{j}\Oi\A_j \cong \pl{j,m} \hat\Omega_{\le
m}\A_j \ .$$
Furthermore the inclusions $\iota_{j*}:\Oi\Ai \to \Oi\A_j$
induce a map 
$$\iota_*:\Oi\Ai \to \pl{j}\Oi\A_j \ .$$

\begin{prop}
\label{projlim}
There are the following canonical
isomorphisms of involutive Fr\'echet locally $m$-convex graded differential algebras:

\begin{enumerate}
\item $\Oi\Ai \cong \pl{j}\Oi\A_j  \cong \pl{j,m} \hat\Omega_{\le
m}\A_j \ .$
\end{enumerate}

There are the following canonical isomorphisms of graded Fr\'echet spaces:
\begin{itemize}
\item[(2)] $\Oi\Ai/\ov{[\Oi\Ai,\Oi\Ai]_s} \cong \pl{j,m} \hat\Omega_{\le
m}\A_j/\ov{[\hat\Omega_{\le m}\A_j,\hat\Omega_{\le m}\A_j]_s} \ ,$
\item[(3)] $H^{dR}_*(\Ai)\cong \pl{j} H_*^{dR}(\A_j) \cong \pl{j,m} H_*(\hat\Omega_{\le
m}\A_j/\ov{[\hat\Omega_{\le m}\A_j,\hat\Omega_{\le m}\A_j]_s},\di) \ .$
\end{itemize}
\end{prop}

\begin{proof}
(1) It is enough to prove that the right hand side and the left hand side are isomorphic as topological vector spaces. This
follows from the fact that projective limits and projective tensor products
commute (\cite{ko2}, 41.6).

(2) and (3) follow from the three technical lemmas below. 
\end{proof}  

The importance of this proposition for our purposes is the following: Since the analysis is easier on Banach spaces than on Fr\'echet spaces we will prove the
index theorem in $H_*(\hat\Omega_{\le
m}\A_j/\ov{[\hat\Omega_{\le m}\A_j,\hat\Omega_{\le m}\A_j]_s},\di)$, making sure that
the expressions in the index theorem behave
well under the projective limit. By the proposition this will prove the
index theorem in $H^{dR}_*(\Ai)$.

\begin{lem}
Let $(V_n,f_{n+1,n}:V_{n+1} \to V_n)_{n \in \bbbn}$ be a projective system of Banach spaces with projective limit $V_{\infty}$ . Let $(A_n,f_{n+1,n}|_{A_{n+1}})_{n \in \bbbn}$ be a projective subsystem of sets such that for any $n \in \bbbn$ the range of the induced map $f_n:A_{\infty}:=\pl{n}A_n \to A_n$ is dense in $A_n$. Then
$$\pl{n}\ov{A_n}=\ov{\pl{n}A_n} \ .$$
Hence a subset $S_1$ of $V_{\infty}$ is dense in $S_2 \subset V_{\infty}$ if and only if $f_nS_1 \subset V_n$ is dense in $f_nS_2 \subset V_n$ for any $n \in \bbbn$.
\end{lem}

\begin{proof} Without loss of generality we may assume that the maps $f_{n+1,n}$ have norm less than or equal to one.
Let $a \in \pl{n}\ov{A_n}$. Then for any $n \in \bbbn$ there is $b_n \in A_{\infty}$ with $|f_n(b_n)  - f_n(a)| \le \frac 1n$. The sequence $(b_n)_{n \in \bbbn}$ converges in $V_{\infty}$ to $a$. Hence $\ov{A_{\infty}}= \pl{n}\ov{A_n}$.
\end{proof}

\begin{lem}
\begin{enumerate}
\item Let $(V_n,f_{n+1,n})_{n \in \bbbn}, ~(W_n,g_{n+1,n})_{n \in \bbbn}$ be projective systems of Banach spaces.
Let $(d_n:V_n \to W_n)_{n \in \bbbn}$ be a morphism of projective systems with induced map $d_{\infty}: V_{\infty} \to W_{\infty}$. Then 
$$\Ker d_{\infty}= \pl{n} \Ker d_n \ .$$
\item In the situation of (1) assume furthermore that the induced maps $f_n:V_{\infty} \to V_n$ and $g_n:W_{\infty} \to W_n$ have dense range for all $n \in \bbbn$. Then
$$\ov{\Ran d_{\infty}}= \pl{n} \ov{\Ran d_n} \ .$$ 
\item Let $(A_n,f_{n+1,n})_{n \in \bbbn}$ be a projective system of Banach algebras such that for any $n \in \bbbn$ the range of $f_n:A_{\infty}:=\pl{n}A_n \to A_n$ is dense in $A_n$. Then 
$$\ov{[A_{\infty},A_{\infty}]}=\pl{n} \ov{[A_n,A_n]} \ .$$
\end{enumerate}
\end{lem}

\begin{proof}
(1) Let $x \in \Ker d_{\infty}$. Then $d_n f_n(x) = g_n d_{\infty}(x)=0$.

Conversely if $d_nf_nx=0$ for all $n \in \bbbn$, then by definition $d_{\infty}x=0$.

(2) follows from the previous lemma since $g_n\Ran d_{\infty}$ is dense in $\ov{\Ran d_n}$.

(3) follows from the previous lemma since  $f_n[A_{\infty},A_{\infty}]$ is dense in $\ov{[A_n,A_n]}$ for all $n \in \bbbn$.
\end{proof}

\begin{lem}
\label{projquot}
Let $(V_n,f_{n+1,n})_{n \in \bbbn}$ be a projective
  system of Banach spaces and let $(A_n,f_{n+1,n}|_{A_{n+1}})_{n \in \bbbn}$ be a projective
subsystem  such that $A_n$ is a closed subspace of $V_n$ for all $n \in \bbbn$.

Let $V_{\infty}:=\pl{n}V_n$ and $A_{\infty}:=\pl{n} A_n$ and assume furthermore
that the image of  $A_{n+1}$ is dense in $A_n$ for any $n \in \bbbn$.
 
Then, canonically,
$$V_{\infty}/A_{\infty} \cong \pl{n} V_n/A_n \ .$$
\end{lem}

\begin{proof}

For $n \in \bbbn$ let $f_n:V_{\infty} \to  V_n$ be the induced map. We prove
that the map
$$f_{\infty}:V_{\infty}/A_{\infty} \to  \pl{n} V_n/A_n$$ is an isomorphism:

For injectivity let $[v] \in V_{\infty}/A_{\infty}$ with $f_{\infty}[v]=0$. 
Then $f_n  v
\in A_n$ for all $n \in \bbbn$, so $v \in A_{\infty}$, thus $[v]=0$.

In order to prove surjectivity we assume that the norms of the maps
$f_{n+1,n}$ are less than
or equal to one, which can be obtained by rescaling the norms inductively.

Let $v \in  \pl{n} V_n/A_n$. Choose inductively $v_n \in V_n$ such that $[v_n] \in
V_n/A_n$ is the image of $v$ with respect to the map  
$\pl{n} V_n/A_n \to V_n/A_n$ and such that $$|f_{n+1,n}v_{n+1} -v_n | \le n^{-1} \ .$$ 
(Here we use that $f_{n+1,n}A_{n+1}$ is dense in $A_n$.)

For $n\ge j$ let $f_{n,j}:V_n \to V_j$ be the induced map.
For any $j \in \bbbn$ the sequence $(f_{n,j}(v_n))_{n \ge j}$ converges
  to an element $\tilde v_j \in V_j$. We have that $f_{j+1,j}\tilde v_{j +1}=
\tilde v_j $, hence there is
an element $\tilde v \in V_{\infty}$ with $f_j \tilde v = \tilde v_j$. Since $[\tilde v_j]=[v_j] \in V_j/A_j$, we have that $f_{\infty}[\tilde v]=v$.

\end{proof}

\section{Lagrangian projections}

Let $\A$ be a unital $C^*$-algebra.

In this section we define and study the analogues of Lagrangian subbundles and of the Maslov index bundle. 

Let $n \in \bbbn$.

\begin{ddd}
Two selfadjoint projections $P_1,P_2\in M_n(\A)$ are called {\sc transverse} if $$\Ran P_1 \oplus
\Ran P_2=\A^{n} \ .$$
\end{ddd}

We will often use the following transversality criterion: Two selfadjoint projections
$P_1,P_2$ are transverse if and only there exist $a,b \in \Gl(\A)$ such that
$aP_1+bP_2 \in M_n(\A)$ is invertible. This is equivalent to the invertibility of
$aP_1+bP_2$ for any $a,b \in \Gl(\A)$.

If $\A=C(B)$ for a compact space $B$, then the transversality of two projections is
equivalent to the transversality of
the corresponding subbundles.

\subsection{Definition and properties}
\label{lagr}
Let $\A^{2n}$ be endowed with the standard $\A$-valued scalar product and let 
$$I_0=\left(\begin{array}{cc}
i & 0\\
0 & -i 
\end{array}\right): \A^n \oplus \A^n \to \A^n \oplus \A^n \ .$$

\begin{ddd}
\label{lagrdef}
A {\sc Lagrangian projection on $\A^{2n}$} is a selfadjoint projection $P \in M_{2n}(\A)$ with
$$PI_0=I_0(1-P) \ .$$
\end{ddd}

If 
$\A=C(B)$ for some compact space $B$, then the Lagrangian projections are in one-to-one correspondence with the
subbundles of $B \times \bbbc^{2n}$ that are Lagrangian with respect to the
skew-hermitian form induced by $I_0$. 

Note furthermore that any Lagrangian projection $P$ on $\A^{2n}$ is also a Lagrangian
projection on the Hilbert $M_n(\A)$-module $M_n(\A)^2$.

Denote $$P_s:=\tfrac 12 \left( \begin{array}{cc}
1 & 1 \\
1 & 1
\end{array}\right) \ . $$

\begin{lem} 
\label{launit}
\begin{enumerate}
\item For every Lagrangian projection $P$ on $\A^{2n}$ there is a
unitary $p \in M_n(\A)$ such that
$$P=\tfrac 12 \left( \begin{array}{ccc}
1 & p^*\\
p & 1
\end{array}\right) \ .$$

\item For every Lagrangian projection $P$ on $\A^{2n}$ the unitary 
$$U=\left(\begin{array}{cc}
1 &0\\
0& p^* 
\end{array}\right)
\in M_{2n}(\A)$$ with $p$ as in (1) fulfills
$UI_0=I_0U$ and 
$UPU^*=P_s \ .$
\end{enumerate}
\end{lem}

\begin{proof}
(1) Since $P$ is selfadjoint, there are $a,b,c \in M_n(\A)$ with $a=a^*,c=c^*$
such that $P=\left( \begin{array}{ccc}
a & b\\
b^* & c
\end{array}\right)$. From $PI_0=I_0(1-P)$ it follows that
$$\left( \begin{array}{ccc}
ia & -ib\\
ib^* & -ic
\end{array}\right) = \left( \begin{array}{ccc}
i(1-a) & -ib\\
ib^* & -i(1-c)
\end{array}\right) \ ,$$
thus $a=c=\frac 12$. Furthermore from $P^2=P$ it follows that $2b$ is
unitary.

(2) is clear.
\end{proof}

Let $\Ai$ be as in the previous section.

\begin{lem}
\label{path}
Let $P\in M_{2n}(\Ai)$ be a Lagrangian projection of $\A^{2n}$ transverse to $P_s$.
Let $\tilde P \in M_{2n}(\bbbc)$ be a complex Lagrangian projection. Then for every $0<
\ve_1< \ve_2$ there is a
smooth path of unitaries $U: [0,\ve_2] \to M_{2n}(\Ai)$ such that
\begin{enumerate}
\item $U(0)PU(0)^*=\tilde P$, 
\item $U$ equals $1$ on a neighborhood of $\ve_2$, 
\item $U$ is constant on $[0,\ve_1]$,
\item $U$ is diagonal with respect to the decomposition $\A^{2n}=\A^n \oplus \A^n$.
\end{enumerate}   
\end{lem}

Note that (4) implies that $UI_0=I_0U$.

\begin{proof} 
It is enough to prove the assertion for $\tilde P=P_s$.

For $P$ let $p$ be as in the previous lemma.

Since $P$ and $P_s$ are transverse, $P-P_s$ is
invertible. It follows that $p-1$ and $p^*-1$ are invertible, so $\log p$ and $\log p^*$ are well-defined
if we choose the complement of $[0, \infty)$ in $\bbbc$ as a domain for the logarithm. 

Let
$\chi:[0,\ve_2] \to [0,1]$ be a smooth function with $\chi|_{[0,\ve_1]}=0$  and
$\chi(t)=1$ for $t \in [\frac{\ve_2-\ve_1}{2},\ve_2]$. The smooth path of
unitaries $$\gamma:[0, \ve_2] \to M_n(\A),~ \gamma(t):=\exp(2 \pi i\chi(t)
 + (1-\chi(t)) \log p^*)$$ connects $p^* $ with $1$. Moreover $\gamma(t)
\in M_n(\Ai)$ for all $t \in [0,\ve_2]$ since $M_n(\Ai)$ is stable under holomorphic function
calculus in $M_n(\A)$ by Prop. \ref{propai}.

Then $$U:=\left(\begin{array}{cc}
1 & 0 \\
0 & \gamma
\end{array} \right) $$
satisfies the conditions.
\end{proof}

\subsection{The Maslov index}
\label{maslov}

Let
$(P_0,P_1,P_2)$ be a triple of pairwise transverse Lagrangian
projections on $\A^{2n}$. For $x \in \A^{2n}$ write $x=x_1 + x_2$ with $x_i \in \Ran P_i,~
i=1,2$. 

The form
$$h: \Ran P_0 \times \Ran P_0 \to \A,~(x,y) \mapsto \langle x_2,I_0 y_1\rangle $$ is
hermitian and its radical vanishes \cite{wa}.  Since $x_1=P_1(P_1+P_2)^{-1}x$ and $x_2=P_2(P_1+P_2)^{-1}x$, the corresponding selfadjoint matrix is 
$$A:=P_0(P_1+P_2)^{-1}P_2I_0P_1(P_1+P_2)^{-1}P_0 \in M_{2n}(\A) \ .$$ 
Since 
\begin{eqnarray*}
A&=& P_0(P_1+P_2)^{-1}P_2I_0P_1(P_1+P_2)^{-1}P_0 \\
&=&P_0(P_1+P_2)^{-1}(P_1+ P_2)IP_1(P_1+P_2)^{-1}(P_0+P_2)\\
&=&(P_0+P_1)IP_1(P_1+P_2)^{-1}(P_0+P_1) \ ,
\end{eqnarray*} the range of $A$ is closed. Hence the hermitian form $h$
is non-singular and defines an element in $K_0(\A)$ \cite{ros}.

\begin{ddd} The {\sc Maslov index} $\tau(P_0,P_1,P_2) \in K_0(\A)$ of a triple of pairwise transverse Lagrangian
projections $(P_0,P_1,P_2)$ is the class of the hermitian form $h$ in $K_0(\A)$.
\end{ddd}

We can express the Maslov index in terms of $A$ as follows:
$$\tau(P_0,P_1,P_2)=[1_{\{x > 0\}}(A)]-[1_{\{x < 0\}}(A)] \in K_0(\A) \ .$$
 Note that the Maslov index is invariant under even permutations 
and changes sign under odd permutations. 

\begin{prop}
\label{mashomotop}
Let $P_i:[0,1] \to M_{2n}(\A),~i=0,1,2,$ be continuous paths of
Lagrangian projections such that $P_i(t) - P_j(t)$ is invertible for $i \neq
j$ and all $t \in [0,1]$. 

Then the Maslov index $\tau(P_0(t),P_1(t),P_2(t))$ does not depend on $t$.
\end{prop}

\begin{proof} The selfadjoint element $A(t) \in M_{2n}(\A)$ defined by
$(P_0(t),P_1(t),P_2(t))$ as above depends continuously on $t$ for all $t \in [0,1]$. It follows that the
projections $1_{\{x > 0\}}(A(t))$ and $1_{\{x < 0\}}(A(t))$ also depend
continuously on $t$, thus their $K$-theory classes are constant.
\end{proof}

Let $B$ be a compact space and let $(P_0,P_1,P_2)$ be a triple of pairwise transverse Lagrangian projections in $M_{2n}(C(B))$. Let $(L_0,L_1,L_2)$ be the corresponding
triple of  Lagrangian subbundles of $B \times \bbbc^{2n}$.
Then the Maslov index bundle $[L^+_0]-[L^-_0]$ defined in \S \ref{commind}  corresponds to
$\tau(P_0,P_1,P_2)$ under the canonical isomorphism $K^0(B) \cong K_0(C(B))$.\\

Now we study in some detail the Maslov index of a triple $(P_0,P_1,P_2)$ with $P_0=P_s$. The general case can be reduced
to this case by Lemma \ref{launit}.

The Cayley transform $a \mapsto \frac{a-i}{a+i}$, defined for selfadjoint $a \in M_n(\A)$, yields a bijective map $$a \mapsto P(a):=\tfrac 12 \left( \begin{array}{cc}
1 & \frac{a+i}{a-i} \\
\frac{a-i}{a+i} & 1
\end{array}\right) \ . $$
from the space of selfadjoint elements in $M_n(\A)$
to space of projections in $M_{2n}(\A)$ transverse to $P_s$.

\begin{lem}
Let  $a_1,a_2 \in M_n(\A)$ be selfadjoint. Then $P(a_1)$ and $P(a_2)$
are transverse projections if and only if $a_1-a_2$ is invertible.
\end{lem}

\begin{proof}
Let $U:= \frac{1}{\sqrt 2} \left(\begin{array}{cc}
1 & 1 \\
1 & -1
\end{array}\right) \ .$
Then  $$UP(a_j)U^*= (a_j^2+1)^{-1}\left(\begin{array}{cc}
a_j^2 & -ia_j \\
ia_j & 1
\end{array}\right) \ .$$

Now $P(a_1)$ and $P(a_2)$ are transverse if and only if
$$(a_1^2+1)UP(a_1)U^*-(a_2^2+1)UP(a_2)U^*$$ is invertible and this is the case if
and only if $a_1-a_2$ is invertible.
\end{proof}

\begin{lem}
\label{selfpath}
Let $(P_s,P_1,P_2)$ be a triple of pairwise transverse Lagrangian projections and let
$a_1,a_2 \in M_n(\A)$ be such that $P_i=P(a_i),~i=1,2$. Let $p^+:=1_{\{x>0\}}(a_1-a_2)$.
\begin{enumerate}
\item There are continuous paths $P_1,P_2:[0,2] \to M_{2n}(\A)$ of Lagrangian projections such that
$P_s,P_1(t),P_2(t)$ are pairwise transverse for all $t \in [0,2]$ and such that
$P_1(2)=P(2p^+-1)$ and $P_2(2)=P(1-2p^+)$.
\item $$\tau(P_s,P_1,P_2)=[p^+]-[1-p^+] \ .$$
\end{enumerate}
\end{lem}

\begin{proof}
(1) For $t \in [0, 1]$ define $a_1(t):=t(a_1-a_2)+(1-t)a_1$ and $a_2(t):=t(a_2-a_1)+(1-t)a_2$. Then 
$a_1(t)-a_2(t)=
(1+t)(a_1-a_2)$ is invertible, thus the projections $P(a_1(t))$ and $P(a_2(t))$ are transverse. Furthermore $a_1(0)=a_1$ and
$a_1(1)=a_1-a_2$, whereas $a_2(0)=a_2$ and $a_2(1)=a_2-a_1$. 

For $t \in [1,2]$ let $a_1(t)$ be a path of invertible selfadjoint elements with $a_1(1)=a_1-a_2$
and $a_1(2)=p^+-(1-p^+)$ and let $a_2(t)=-a_1(t)$.

Then the paths $P(a_1(t))$ and $P(a_2(t))$ satisfy the conditions.

(2) From (1) and Prop. \ref{mashomotop} it follows that
   $\tau(P_s,P_1,P_2)=\tau(P_s,P(2p^+-1),P(1-2p^+))$. 

The Cayley transform of $2p^+-1$ is $i(2p^+-1)$. By computing the
matrix $A$ we see that $1_{\{x>0\}}(A)=p^+$ and $1_{\{x<0\}}(A)=(1-p^+)$.
\end{proof}

\chapter{The Fredholm Operator and Its Index}

\section{The operator $D$ on $M$}
\label{defD}
 
\subsection{Definition of $D$}
\label{defiD}

Now we come back to the geometric situation described in \S\ref{situat}.

By taking the tensor product of the bundle  $E$ with the $C^*$-algebra $\A$ we obtain
an $\A$-vector bundle \cite{mf}. Furthermore we
consider the bundle  $E \ten
\Ol{\mu},~i, \mu \in \bbbn_0,$ of right $\Ol{\mu}$-modules. Keep in mind that $E$ can be trivialized on $M$
via a global orthonormal frame. Thus no theory of Banach space bundles is needed in this context.

The hermitian metric on $E$
extends to an $\A$-valued scalar product $\langle\cdot, \cdot\rangle $ on $E \ten \A$ and to
an $\Ol{\mu}$-valued non-degenerated product on $E \ten \Ol{\mu}$ (see \S\ref{adop} for
this notion). Furthermore parallel transport is defined on $E \ten \A$
resp. $E \ten \Ol{\mu}$.

By the trivialization of $E|_F$ fixed in \S \ref{situat} we identify 
$(E|_F \ten \A, I)_x,~x \in F,$ with $(\A^{4d} ,I_0 \oplus (-I_0))$ as a
Hilbert $\A$-module with  a skew-hermitian structure, and $(E\ten \Ol{\mu})_x,~x \in F,$ with $(\Ol{\mu})^{4d}$ as a right $\Ol{\mu}$-module with an $\Ol{\mu}$-valued non-degenerated product.\\

Recall that the Hilbert $\A$-module $L^2(M,E \ten \A)$ can be
 defined as the completion of $\C_c(M,E \ten \A)$ with respect to
 the norm induced by the $\A$-valued scalar product
$$\langle f,g\rangle :=\int_M \langle f(x),g(x)\rangle ~ dx \ .$$
(For fixing notation a short introduction to Hilbert
$\A$-modules can be found in \S\ref{HCmod}.) By Prop. \ref{orthsys} any orthonormal basis of the
Hilbert space $L^2(M,E)$ is an orthonormal basis of $L^2(M,E \ten \A)$ and
$L^2(M,E \ten \A)$ is isomorphic to
$l^2(\A)$. \\

We introduce the Schwartz space of sections of $E \ten \Ol{\mu}$:

For $k\in \bbbz/6$ define the Schwartz space $${\mathcal S}(Z_k,(\Ol{\mu})^{4d}):={\mathcal S}(\bbbr)\ten_{\pi}
\C([0,1],(\Ol{\mu})^{4d}) \ .$$

For  $r\ge 0,~0 \le b \le \frac 14$ choose a partition of unity $\{\phi_k\}_{k \in J}$ subordinate to the covering ${\mathcal
U}(r,b)$.

The embedding $F_k(r,b) \incl Z_k,~k \in \bbbz/6,$ and the trivialization of $E|_F$ induce a map 
$$\C(M,E \ten \Ol{\mu}) \to \C(Z_k,(\Ol{\mu})^{4d}),~f \mapsto \phi_k f \ .$$ 

As a vector space let ${\mathcal S}(M,E \ten \Ol{\mu})$   be the largest subspace of $\C(M,E \ten \Ol{\mu})$
such that for all $k \in \bbbz/6$ the maps
$${\mathcal S}(M,E \ten \Ol{\mu}) \to  {\mathcal S}(Z_k,(\Ol{\mu})^{4d}),~f \mapsto
\phi_kf $$ are well-defined. Endowed with the seminorms of $\C(M,E \ten \Ol{\mu})$ and  those of ${\mathcal S}(Z_k,(\Ol{\mu})^{4d})$ applied to $\phi_kf$ for $k \in \bbbz/6$ the space ${\mathcal S}(M,E \ten \Ol{\mu})$ is a Fr\'echet space.  
The topology does not depend on the choice of $r$ and $b$, neither on the choice of the
partition of unity.\\

To a triple $R=(P_0,P_1,P_2)$ of Lagrangian projections of
$\A^{2d}$ with $P_i \in M_{2d}(\Ai),~i=0,1,2,$ we associate boundary conditions on sections of $E \ten \Ol{\mu}$ in the following way: 

For $k \in \bbbn_0 \cup \{\infty\}$ let 
\begin{eqnarray*}
C^k_R(M,E \ten \Ol{\mu}) &:=& \{f \in C^k(M,E \ten \Ol{\mu})~|~(P_i \oplus P_i)
\dirac^l f(x)= \dirac^l f(x)\\
&& \qquad \text{ for } x \in
\ra_iM \cup \ra_{i+3}M,~i=0,1,2;~l \in \bbbn_0, l \le k\} 
\end{eqnarray*}
endowed with the subspace topology.

Analogously we define $C^k_{R0}(M,E \ten \Ol{\mu})$ and $C^k_{Rc}(M,E \ten \Ol{\mu})$.

Furthermore let ${\mathcal S}_R(M,E \ten \Ol{\mu})$ be the vector space ${\mathcal S}(M,E \ten \Ol{\mu})\cap \C_R(M,E
\ten \Ol{\mu})$  with the subspace topology of ${\mathcal S}(M,E \ten \Ol{\mu})$.\\

Now we introduce the operator $D$:

Fix a triple of pairwise transverse Lagrangian
projections $R=(\Pj_0, \Pj_1, \Pj_2)$ of $\A^{2d}$ with $\Pj_i \in
M_{2d}(\Ai),~i=0,1,2$. 
We define $D$ on $L^2(M,E \ten \A)$ as the closure of the Dirac operator $\dirac$ with domain
$\C_{Rc}(M, E \ten \A)$ and $D^+$ resp. $D^-$ as the restriction of $D$ to the sections of $E^+ \ten \A$ resp. $E^- \ten \A$.

Note that $D$ is symmetric.

\subsection{Comparison with $D_s$}
\label{stand}
Recall that $P_s=\frac 12\left(\begin{array}{cc} 1&1\\1&1 \end{array}\right) \in
M_{2d}(\bbbc) \ .$

Fix a triple of pairwise transverse Lagrangian projections
$(\Pj_0^s,\Pj_1^s,\Pj_2^s)$ with $\Pj_0^s=P_s$ and $\Pj_1^s,\Pj_2^s \in
M_{2d}(\bbbc)$. Let $D_s$ on
$L^2(M,E \ten \A)$ be the closure of the Dirac operator $\dirac$ with domain
$\C_{Rc}(M, E \ten \A)$ for the triple
 $R=(\Pj_0^s,\Pj_1^s,\Pj_2^s)$.

Let $W \in \C(M,\End^+ E \ten \Ai)$ be such that
\begin{itemize}
\item $WW^*=1$,
\item $W(x)(\Pj_i \oplus \Pj_i) W(x)^*=(\Pj^s_i \oplus \Pj^s_i)$ for all $x \in \ra_iM \cup
\ra_{i+3}M, ~i=0,1,2$,
\item  $W$ is parallel on $M \setminus F$ and on a neighborhood of $\ra M$,
\item for all $k \in \bbbz/6$ the restriction of $W$ to $F_k(0, \frac 14)$ depends only on the
coordinate $x_2^k$,
\item $W$ commutes with the Clifford multiplication.
\end{itemize}

Some of these properties are not needed in this section, but are important for the proof of the index theorem.

\begin{prop}
\label{DMreg}
\begin{enumerate}
\item We have that $WDW^*=D_s+ 
Wc(d W^*)$ with $c(d W^*)|_F:=c(e_2)\ra_{e_2}W^*$ and
$c(d W^*)|_{M \setminus F}:=0$. In particular, $Wc(d W^*) \in \C(M, \End^- E \ten \Ai)$. 

\item
The operator $D$ is regular and selfadjoint.
\end{enumerate}
\end{prop}

\begin{proof}
(1) For $R=(\Pj_0^s,\Pj_1^s,\Pj_2^s)$ and $f \in \C_{cR}(M,E \ten \A)$  we have that
   $(WDW^*f)|_{M\setminus F}=(D_sf)|_{M\setminus F}$ and
\begin{eqnarray*}
(WDW^*f)|_F &=& (D_s f)|_F+ W[c(e_2)\ra_{e_2},W^*]_s (f|_F) \\
&=&(D_sf)|_F + Wc(e_2)(\ra_{e_2}W^*)(f|_F) \ .
\end{eqnarray*}
(2) The restriction of $D_s$ to the Hilbert space $L^2(M,E)$ is selfadjoint \cite{bk}. Hence
$(1+D_s^2)$ has a bounded inverse on $L^2(M,E)$. It follows that
 the range of $(1+D^2_s)$ on $L^2(M,E\ten\A)$ is dense, thus $D_s$ is regular. By an analogous argument the operators $D_s \pm i$ have  dense range. From Lemma \ref{critself} it follows that $D_s$ is selfadjoint. 

By (1) the operator $D$ is a bounded perturbation of $D_s$, hence $D$
is selfadjoint and regular (see Prop. \ref{regpert}).
\end{proof}

The existence of a section $W$ fulfilling the properties above is proved in the following lemma and proposition:

\begin{lem}
\label{uniteq}

There is a parallel unitary section $U \in \C(M,\End E^+ \ten
\Ai)$ such that $U\Pj_0U^*=\Pj^s_0=P_s$.
\end{lem}

\begin{proof}
By the isomorphism $E^+ \ten \A \cong S^+ \ten (\A^+)^d  \oplus  S^- \ten (\A^-)^d$  any $A \in M_{2d}(\A)$
of the form $\left(\begin{array}{cc}a & 0 \\ 0 & b \end{array}\right)$ with $a,b \in M_d(\A)$ defines a parallel section of
$\End E^+ \ten \A$. 

By Lemma \ref{launit} there is a unitary $U \in M_{2d}(\Ai)$ of that form such that $U\Pj_0U^*=\Pj^s_0=P_s$.
\end{proof}

\begin{prop}
\label{constrW}
For any $0<b< \frac 14$ there is $W\in \C(M,\End^+ E \ten \Ai)$
with the properties above and such that $W$ is parallel on $\{ x \in M~|~d(x,\ra M)>b\}$.
\end{prop}

\begin{proof} 
By the previous lemma we may assume $\Pj_0=P_s$.

In the following we identify $\ra
M \times [0,b]$ with $\{x \in M~|~d(x, \ra M) \le b\}$. 

Since  for $i=1,2$ the projection $\Pj_i$ is transverse to $P_s$, there are,  by  Lemma \ref{path},
smooth  paths $W_i:[0,b] \to M_{2d}(\Ai)$ of unitaries with $[W_i,I_0]=0$ such that
$W_i(0)\Pj_iW_i^*(0)=\Pj_i^s$ and such that $W_i$ is equal to the identity in a neighborhood of
$b$ and constant on $[0,\frac b2]$. They induce maps $\tilde W_i:=\id \times W_i: (\ra_i M \cup
\ra_{i+3} M) \times [0, b] \to M_{2d}(\Ai)$. The map 
$$\cup_i (\tilde W_i \oplus \tilde W_i): \ra M \times
[0,b] \to M_{4d}(\Ai)$$ can be extended by $1$ to a
smooth section $W  \in \C(M,\End^+ E \ten \Ai)$. By construction $W$ has the
right properties.
\end{proof}

\section{The operator $D_I$ on $[0,1]$}

\subsection{Definition and comparison with $D_{I_s}$}
\label{DI}

Let $\ra$  be the differentiation operator on $[0,1]$.

For a pair $R=(P_0,P_1)$ of transverse Lagrangian projections with $P_j \in
M_{2d}(\Ai),\\ j=0,1,$ and $k \in \bbbn_0 \cup \{\infty\}$ we define the function space
\begin{eqnarray*}
\lefteqn{C^k_R([0,1],(\Ol{\mu})^{2d})}\\
&:=& \{f \in C^k([0,1],(\Ol{\mu})^{2d}) ~|~P_x(I_0\ra)^l f(x)=(I_0\ra)^lf(x) \\
&& \qquad \mbox{ for
}x=0,1; ~l \in \bbbn_0,~l \le k\} 
\end{eqnarray*}
and endow it with the subspace topology.\\

Given a pair $(P_0,P_1)$ of  transverse Lagrangian projections on   $\A^{2d}$ we write
$D_I$ for the closure of the operator $I_0 \ra$ on the Hilbert $\A$-module $L^2([0,1],\A^{2d})$ with domain
$\C_R([0,1],\A^{2d}).$  

We write $D_{I_s}$ if $(P_0,P_1)=(P_s,1-P_s)$.

\begin{prop}
\label{defWI}
Let $(P_0,P_1)$  be a pair of transverse Lagrangian projections  with $P_0,P_1 \in M_{2d}(\Ai)$. Then for $0<x_1<x_2<1$  there is $U
\in \C([0,1], M_{2d}(\Ai))$ such that
\begin{enumerate} 
\item $UU^*=1$, 
\item $UI_0=I_0U$,
\item $U(0)P_0U(0)^*=P_s$,
\item $U(1)P_1U(1)^*=1-P_s$,
\item $U$ is  constant on $[0,x_1]$ and on $[x_2,1]$. 
\end{enumerate}

\end{prop}

\begin{proof}
By Lemma \ref{launit} there is a unitary $U_0 \in M_{2d}(\Ai)$ with
$U_0I_0=I_0U_0$ and
$U_0P_0U_0^*=P_s$. Since $U_0P_1U_0^*$ is transverse to $P_s$, we can apply Lemma
\ref{path} to $P=U_0P_1U_0^*$ and $\tilde P=1-P_s$ in order to get a smooth path of unitaries
$U_1:[0,1] \to M_{2d}(\Ai)$ such that $U(t):=U_1(t)U_0$ has the right properties.
\end{proof}

\begin{prop}
\label{DIpert}
Let $(P_0,P_1)$ be a pair of transverse Lagrangian projections  of
$\A^{2d}$ with $P_0,P_1 \in M_{2d}(\Ai)$ and let $D_I$ be the associated operator.

\begin{enumerate}
\item Let $U$ be as in the previous proposition. Then  $$UD_IU^*=D_{I_s}+ UI_0
(\ra U^*)$$ with $UI_0
(\ra U^*) \in \C([0,1], M_{2d}(\Ai))$. 
 
\item
The operator $D_I$ is regular and selfadjoint.
\end{enumerate}
\end{prop}

\begin{proof}
(1) For $f \in \C_R([0,1],\Ai^{2d})$ with $R=(P_s,1-P_s)$ we have:
\begin{eqnarray*}
UD_IU^*f&=&U I_0\ra U^*f\\
&=&D_{I_s}f+ UI_0(\ra U^*)f \ .
\end{eqnarray*}

(2) follows as in Prop. \ref{DMreg}.
\end{proof}

\subsection{Generalized eigenspace decomposition}
\label{geneig}

We construct a decomposition of $L^2([0,1], \A^{2d})$
 into free finitely generated $\A$-modules preserved by $D_I$. For $d=1$
 these can be understood as analogues of eigenspaces.

Assume that the boundary conditions of $D_I$ are given by a pair $(P_0,P_1)$ with $P_0=P_s$. For the general case use Lemma \ref{launit}.

Let $p \in M_d(\A)$ be such that $P_1=\frac 12\left(\begin{array}{cc}1 & p^*\\ p &
1\end{array} \right)$.

The transversality of $P_0$ and $P_1$ implies that $1-p$ is invertible. Hence $\log p$ is defined for $\log:\bbbc\setminus [0,\infty) \to \bbbc$.

\begin{lem}  
\label{eigenv}
Assume that $d=1$.

Then for $k \in \bbbz$ the function
$$f_k(x)=\tfrac {1}{\sqrt 2}\left(\left(\begin{array}{c} 1\\ 0\end{array}\right)\exp[(-\tfrac 12 \log p + \pi k i)x] + \left(\begin{array}{c} 0\\ 1\end{array}\right)
\exp[(\tfrac 12 \log p - \pi k i)x] \right)$$ is in $\dom D_I$ and
$D_If_k=\lambda_kf_k$ with $\lambda_k:=-\frac 12 i \log p -\pi k$.

The system $\{f_k\}_{k\in \bbbz}$ is an orthonormal basis of the Hilbert $\A$-module
$L^2([0,1],\A^2)$.
\end{lem}

Note that $\lambda_kf_k=f_k \lambda_k$ and $\sigma(\lambda_k) \subset
]-\pi k,-\pi(k-1)[$.

\begin{proof} Clearly the system $\{f_k\}_{k\in \bbbz}$ is orthonormal in the sense of Def. \ref{deforthsys}.

In order to prove that these functions form a basis we first show that there
is an orthogonal projection onto the closure of their span. This implies that
the span is orthogonally complemented. The second step will be to see that the
complement is trivial. Then the claim follows from Prop. \ref{orthsys}.\\

The system $\{e_l^{\pm}(x):=v_{\pm}e^{2\pi i lx}\}_{l \in \bbbz,\pm}$ with
$v_+:=( 1, 0),~v_-:=( 0, 1)$ is an orthonormal basis of
$L^2([0,1],\A^2)$ by Prop. \ref{orthsys}.

We have that
\begin{eqnarray*}
\langle f_k,e^{\pm}_l\rangle  &=&\tfrac{1}{\sqrt 2} \int_0^1\exp((\mp \tfrac 12 \log p^* \mp \pi k i)x)\exp(2 \pi i lx)~dx \\
&=&\tfrac{1}{\sqrt 2}\left( \frac{ (-1)^k\exp(\mp \frac 12 \log p^*) - 1}{\mp
\frac 12 \log p^* + \pi i(2l \mp k)}\right) \ ,
\end{eqnarray*}
hence $\{\langle f_k,e^{\pm}_l\rangle \}_{k \in \bbbz} \in l^2(\A)$
for all $l \in \bbbz$ and for $\pm$. Thus
$$Pe_l^{\pm}:=\sum \limits_{k\in \bbbz}f_k \langle f_k,e_l^{\pm}\rangle  ~\in L^2([0,1],\A^2) \ .$$ 
The $\A$-linear extension of $P$ to the algebraic span of the functions $e_l^{\pm}$ has norm
one. It follows that its closure is an orthogonal projection $$P:L^2([0,1],\A^2) \to
\ov{\spann_{\A}\{f_k~|~k\in \bbbz\}} \ .$$ 
It remains to show that the
kernel of $P$ is trivial.

Let $g=(g_1,g_2) \in L^2([0,1],\A^2)$ with 
$\langle f_k,g\rangle =0$ for all $k \in \bbbz.$

Hence for all $k \in \bbbz$  
$$(*) \qquad \int_0^1\exp((- \tfrac 12 \log p^* - \pi k i)x)g_1(x) + \exp(( \tfrac
12 \log p^* + \pi k i)x)g_2(x)~dx=0 \ .$$
Since $\exp(- \frac 12 (\log p^*)~x)g_1(x)$  is in $L^2([0,1],\A)$, there is a unique $\{\lambda_l\}_{l \in \bbbz} \in l^2(\A)$ such that 
$$\sum\limits_{l \in \bbbz} \lambda_l e^{2\pi i lx}=\exp(- \tfrac 12 (\log p^*)~
x)g_1(x) \ .$$
Inserting this in $(*)$  and evaluating the integral for $k$ even leads to
$$\lambda_{k/2}+\int_0^1\exp(( \tfrac 12 \log p^* + \pi k i)x)g_2(x)dx=0 \ .$$
 It follows  that
$$\exp( \tfrac 12 (\log p^*)~ x)g_2(x)=\sum\limits_{l \in \bbbz} (-\lambda_l)
e^{-2\pi i lx} \ .$$  
Substituting again and evaluating $(*)$ for $k=2\nu+1$ with $\nu \in \bbbz$
we obtain
\begin{eqnarray*}
0&=&\int_0^1 \sum\limits_{l \in \bbbz} \lambda_l( e^{\pi i (2(l-\nu)-1)x} - e^{-\pi i (2(l-\nu)-1)x})dx\\ 
&=& -\sum\limits_{l \in \bbbz} \lambda_l \frac{4}{\pi i(2(l-\nu)-1)} \ .
\end{eqnarray*}
Note that for every
$l \in \bbbz$ the function 
$$a_l:\bbbz \to \bbbc,~\nu \mapsto \frac{2}{\pi
i(2(l-\nu)-1)}$$ is in $l^2(\bbbc)$. We claim that $\{a_l\}_{l\in \bbbz}$ is an
orthonormal basis of $l^2(\bbbc)$. Then $\{a_l\}_{l\in \bbbz}$
is an orthonormal basis of $l^2(\A)$ as well, thus $\lambda_l=0$ for
all $l \in \bbbz$ and hence $g_1=g_2=0$.\\

The  Fourier transform of $$a_0:\nu \mapsto \frac{-2}{\pi i(2\nu+1)}$$
 is $$h(x)=-ie^{-\pi ix}(1_{[0,1/2]}(x)- 1_{]1/2,1]}(x)) \in L^2(\bbbr/\bbbz) \ .$$
Since $a_l(\nu)= a_0(\nu-l)$, the Fourier transform of $a_l$ is
$h(x)e^{-2\pi i lx}
$. 

By $hh^*=1$ the system $\{h(x)e^{-2 \pi i lx}\}_{l \in \bbbz}$ is an orthonormal
basis of $L^2(\bbbr/\bbbz)$. This implies that $\{a_l\}_{l\in \bbbz}$ is an orthonormal basis of $l^2(\bbbc)$.
\end{proof}

For general $d$ there is a decomposition of $L^2([0,1]\A^{2d})$ into $\A$-modules of rank $d$:

\begin{prop}
\label{geneigenv}
For $k \in \bbbz$ let
$U_k \subset L^2([0,1],\A^{2d})$ be the right $\A$-module spanned by the column vectors of
$$\frac {1}{\sqrt 2}\left(\begin{array}{c} \exp[(-\frac 12 \log p + \pi k
i)x] \\
\exp[(\frac 12 \log p - \pi k i)x] \end{array} \right) \in \C([0,1],M_{2d \times d}(\Ai)) \
.$$
Each $U_k$ is a free right
$\A$-module of rank $d$ and 
$$L^2([0,1],\A^{2d})=\hat{\bigoplus_{k \in \bbbz}} U_k \ .$$
The sum is orthogonal.

Furthermore $U_k \subset \dom D_I$ and for $f \in U_k$ 
$$D_If=\left(\begin{array}{cc} \lambda_k &0 \\ 0& \lambda_k \end{array}\right)f
\in U_k$$
with $\lambda_k=-\frac 12 i \log p -\pi k$. 
\end{prop}

\begin{proof}
The projections $P_0,P_1\in M_{2d}(\A)$ are Lagrangian projections on
$M_d(\A)^2$ (see remark after Def. \ref{lagrdef}). 

Let $\tilde D_I$ be the closure of $I_0
\ra$ on $L^2([0,1],M_d(\A)^2)$ with domain $\C_R([0,1],M_{2d}(\A))$ with $R=(P_0,P_1)$. Then $\tilde D_I= \oplus^d D_I$ with
respect to the decomposition as right $M_d(\A)$-modules $$L^2([0,1],M_d(\A)^2)=L^2([0,1],\A^{2d})^d$$ induced by the decomposition of a matrix into its column vectors.

By the
previous proposition the  Hilbert $M_d(\A)$-module $L^2([0,1],M_d(\A)^2)$ has an orthonormal basis $\{f_k\}_{k
\in \bbbz}$ such that 
$$\tilde D_I f_k=\left(\begin{array}{cc} \lambda_k &0 \\ 0& \lambda_k
\end{array}\right)f_k \ .$$ 
For $k \in \bbbz$ let $P_k$ be the orthogonal projection onto the span of
$f_k$ in $L^2([0,1],M_d(\A)^2)$. It is  diagonal with respect to the decomposition
$L^2([0,1],M_d(\A)^2)=L^2([0,1],\A^{2d})^d$. 

Hence 
$$L^2([0,1],\A^{2d})=\hat{\bigoplus_{k \in
\bbbz}}P_k L^2([0,1],\A^{2d}) \ .$$
The assertion follows now since $P_k L^2([0,1],\A^{2d})=U_k$.
The module $U_k$ is free since $\exp[(-\frac 12 \log p + \pi k
i)x]$ is invertible in $M_d(\A)$ for all $x \in [0,1]$.
\end{proof}

\begin{cor} Let $\lambda \in \bbbc$.
The operator $D_I-\lambda$ has a bounded inverse on $L^2([0,1],\A^{2d})$ if and only if  $\exp(2 i\lambda) \notin \sigma(p).$

\end{cor}

\begin{proof}
By the previous
proposition 
$(D_I-\lambda)|_{U_k}$ is invertible if and only if $e^{2 i \lambda}
\notin \sigma(p)$. Furthermore for $\lambda$ with $e^{2 i \lambda}
\notin \sigma(p)$ the inverse of $(D_I-\lambda)|_{U_k}$ is uniformly bounded
in $k$. Hence the closure of $\oplus_k (D_I-\lambda)|_{U_k}$
has a bounded inverse by Cor. \ref{decinv}. In particular the closure of $\oplus_k D_I|_{U_k}$ is
selfadjoint. Since $D_I$ is a selfadjoint extension of $\oplus_k D_I|_{U_k}$, it follows that $D_I$ is the closure of $\oplus_k
D_I|_{U_k}$. Hence $D_I-\lambda$ has a bounded inverse if $\exp(2 i\lambda) \notin \sigma(p).$
\end{proof}
 
\section{The operator $D_Z$ on the cylinder}
\label{Diraccyl}
Let $X= \bbbr,\bbbr/\bbbz$. Endow $X \times [0,1]$
with the euclidean metric and a spin structure and let $(x_1,x_2)$ be the euclidian coordinates
of $X \times [0,1]$. Let $S=S^+ \oplus S^-$ be the spinor bundle on $X \times
[0,1]$ endowed with the Levi-Civit\`a connection and a parallel metric. Then $S \ten \bigl((\A^+)^d \oplus
(\A^-)^d \bigr)$ is Dirac bundle on $X \times [0,1]$ with the
$\A$-valued scalar product induced by the hermitian metric on $S$ and the standard $\A$-valued scalar product on
$\A^d$. Let $\dirz$ be the associated Dirac operator.

We choose a parallel unit section $s$ of $S^+$ and identify $S\ten \bigl((\A^+)^{d} \oplus (\A^-)^{d}\bigr)$ with  the trivial bundle $(X \times [0,1]) \times \bigl((\A^+)^{2d} \oplus
(\A^-)^{2d}\bigr)$
via the isomorphisms
$$S_x^+ \ten (\A^+)^d  \oplus S_x^- \ten (\A^-)^d \to
(\A^+)^{d} \oplus (\A^+)^{d}\ ,$$
$$(s(x) \ten v) \oplus (ic(dx_1)s(x) \ten w)  \mapsto (v,w) \ ,$$
and $$S_x^-\ten (\A^+)^d  \oplus S_x^+ \ten (\A^-)^d  \to
(\A^-)^{d} \oplus (\A^-)^{d} \ ,$$
$$(ic(dx_1)s(x)\ten v) \oplus (s(x) \ten w)\mapsto (v,w) \ .$$
for $x \in X \times [0,1]$.

Let $I:=c(dx_2)c(dx_1)=I_0 \oplus (-I_0).$ Then
$$\dirz=c(d x_1)(\ra_{x_1} + I \ra_{x_2}) \ .$$
Given a pair $(P_0,P_1)$ of transverse Lagrangian projections of $\A^{2d}$ with
$P_i \in M_{2d}(\Ai)$ we define $D_Z$ to be the closure of $\dirz$
 with domain $$\{f \in \C_c(X \times [0,1], \Ai^{4d})~|~ (P_i \oplus
P_i)f(x,i)=f(x,i) \mbox{ for all } x \in X, ~i=0,1\}\ $$ on the Hilbert
$\A$-module $L^2(X \times [0,1],\A^{4d})$.
 
We write  $H(D_Z)$ for the Hilbert $\A$-module whose underlying right $\A$-module is
$\dom D_Z$  and whose $\A$-valued scalar product is
$\langle f,g\rangle _{D_Z}=\langle f,g\rangle +\langle D_Zf,D_Zg\rangle $ (see \S \ref{regop}).

\begin{prop}
\label{inv}
\begin{enumerate}
\item
The operator
$D_Z$ is selfadjoint on $L^2(X \times [0,1],\A^{4d})$ and has a bounded inverse.
\item If $X=\bbbr/\bbbz$, then 
the inclusion
$\iota:H(D_Z) \to L^2(\bbbr/\bbbz \times [0,1],\A^{4d})$ is
compact.
\end{enumerate}
\end{prop}

\begin{proof}
The proof is analogous to the family case \cite{bk}.

(1) By Lemma \ref{launit} we may
assume $P_0=P_s$.
Recall from Prop. \ref{geneigenv} that the operator $D_I$ with boundary
conditions $(P_0,P_1)$ induces a decomposition $$L^2([0,1],(\A^+)^{2d})=\hat{\bigoplus_{l \in
\bbbz}} U_l$$ such that for each $l \in \bbbz$ there is $\lambda_l \in M_d(\A)$ with $D_If=\left(\begin{array}{cc}\lambda_l &0 \\ 0 & \lambda_l
\end{array}\right)f$ for $f \in U_l$. 

We define submodules
$U_{l,+}:=U_l$ and $U_{l,-}:=ic(dx_1)U_l$ of $L^2([0,1],\A^{4d})$. Then 
$$ \bigoplus_{l\in \bbbz} (U_{l,+} \oplus U_{l,-})
L^2(X)$$
is dense in $L^2(X \times [0,1],\A^{4d})$.

First consider the case $X=\bbbr$:

For $l \in \bbbz$ let  $\ra_{l,\pm}$ be the closure of the unbounded operator $(\ra_{x_1}  \pm
\lambda_l)$ on $U_{l,\pm}L^2(\bbbr)$ with domain $U_{l,\pm}{\mathcal S}(\bbbr)$  
and let 
$\ra_e$ be the closure of $c(d
x_1)\bigl(\bigoplus_{l\in \bbbz,\pm} \ra_{l,\pm}\bigr) \ .$

The operator $D_Z$ is an extension of $\ra_e$.

We claim that $\ra_e$ has a bounded inverse on $L^2(\bbbr \times
[0,1],\A^{4d})$.

The Fourier transform on $L^2(\bbbr)$ induces an automorphism on $U_{l,\pm} L^2(\bbbr)$. Conjugation by it transforms $\ra_{l,\pm}$
into multiplication by $$ix_1 \pm \left(\begin{array}{cc}\lambda_l &0 \\ 0 & \lambda_l
\end{array}\right) \ .$$ 
Since $\sigma(\lambda_l) \subset
\bbbr^*$, we see that the operator $\ra_{l,\pm}$ has a bounded inverse and the
norm of the inverse tends to zero for $l \to \pm \infty$. By Cor. \ref{decinv} it follows  that
the closure of $\oplus_{l,\pm} \ra_{l,\pm}$ has a bounded inverse. 

Hence the operator  $\ra_e$ has a bounded inverse as well. In particular it is
selfadjoint, thus  $D_Z=\ra_e$. For $X=\bbbr$ the assertion follows.

If $X=\bbbr/\bbbz$, then the spaces $U_{l,\pm}  L^2(\bbbr/\bbbz)$ decompose further into a
direct sum 
$$U_{l,\pm}  L^2(\bbbr/\bbbz)=\hat{\bigoplus_{k \in \bbbz}}~ V_{kl,\pm}$$
with  $V_{kl,\pm}:=e^{2\pi ikx_1}U_{l,\pm}$. Note that $V_{kl,\pm}$ is
isomorphic to $\A^d$ as a Hilbert $\A$-module. 

Let $\ra_{kl,\pm}\in B(V_{kl,\pm})$ be defined by
$$\ra_{kl,\pm}f:=(\ra_{x_1}  \pm
\lambda_l)f=\left(\begin{array}{cc}2\pi i k\pm
\lambda_l &0 \\ 0 &
2\pi i k\pm
\lambda_l \end{array} \right)f$$  
 and let $\ra_e$ be the closure of $c(d
x_1)\bigl(\bigoplus_{k,l,\pm} \ra_{kl,\pm}\bigr) \ .$

The operator $D_Z$ is an extension of $\ra_e$.

Since $|(2\pi i k\pm
\lambda_l)^{-1}|$ tends to zero for $k,l \to \pm \infty$, 
the closure of  $\oplus_{kl,\pm} \ra_{kl,\pm}$ has a compact inverse by Cor. \ref{decinv} and hence $\ra_e$ has a compact inverse as well.

Now (1) follows as above. 

(2) follows from the compactness of $D_Z^{-1}$ for $X=\bbbr/\bbbz$ since $\iota=D_Z^{-1}D_Z:H(D_Z) \to L^2(\bbbr/\bbbz \times [0,1],\A^{4d})$.
 \end{proof}

\section{The index of $D^+$}

Recall from \S \ref{defiD} that the boundary conditions of $D$ were defined by a triple $R=(\Pj_0,\Pj_1,\Pj_2)$ .
For an open precompact subset $U$ of $M$ we define $H^1_{R0}(U,E
\ta)$ to be the closure of
$\C_{Rc}(U, E \ta)$ in  $H(D)$ (see \S \ref{regop} for the definition of $H(D)$).

Note that $H(D)$ is isomorphic to $l^2(\A)$ as a Hilbert $\A$-module, since
$L^2(M, E \ta)$ is isomorphic to $l^2(\A)$ and since $L^2(M, E \ta)$ and
$H(D)$ are isomorphic by Lemma \ref{adj}.

\begin{lem}
The inclusion  $\iota:H^1_{R0}(M_r,E \ta)\to L^2(M,E \ta),
~r \ge 0,$ is compact.
\end{lem}

\begin{proof}
Let ${\mathcal
V}:=\{V_k\}_{k \in L}$ be an open covering of $M_r$ such that the index set $L$
is a finite subset of $\bbbn$  with $1 \in L$  and such that:
\begin{itemize}
\item $V_1=M \setminus \ov{F(r,\frac 16)}$; 
\item for $k>1$ there is
an isometry $V_k \cong  ]0,\frac 12[ \times [0,\frac 15[$. In particular $V_k$
is in the flat region and has exactly one boundary component.
\end{itemize}

First we prove that the maps $\iota_k:H^1_{R0}(V_k,E \ta) \to L^2(V_k,E \ta)$ are
compact for all $k \in L$. 

The compactness of $\iota_1$ follows from the
Sobolev embedding theorem (\cite{mf}, Lemma 3.3). 

For $k \neq 1$ let $i \in \{0,1,2\}$ be such that $\ra V_k \subset (\ra_i M \cup
\ra_{i+3}M)$ and set $P_k:=\Pj_i$. Let $D_k$ be the
 operator $D_Z$ on the bundle $(\bbbr/\bbbz
\times [0,1]) \times \A^{4d}$  with boundary conditions
given by the pair $(P_k,1-P_k)$.

The map $\iota_k$ is compact since it factors through the map $H(D_k) \to L^2(\bbbr/\bbbz \times
[0,1],\A^{4d})$, which is
compact by Prop. \ref{inv}.

Let $\{\phi_k\}_{k \in L}$ be a smooth partition of unity subordinate to the covering ${\mathcal
V}$
such that $\ra_{e_2}\phi_k(x)=0$ for all $x \in \ra M$ and $k \in L$ 
. Multiplication with $\phi_k$ is a bounded map from
$H^1_{R0}(M_r,E \ten \A)$ to $H^1_{R0}(V_k,E \ten \A)$, hence by
$\iota=\sum\limits_{k \in L} \iota_k
\phi_k$, the inclusion $\iota$ is compact.
\end{proof}

Let $H(D)^+$ be the subspace of $H(D)$ of homogeneous elements of degree zero.

\begin{prop}
\label{Dfred}
The operator $$D^+:H(D)^+ \to L^2(M,E^- \ten \A)$$ is a Fredholm operator.
\end{prop}

\begin{proof}
By constructing a
parametrix for $D^+$ we show that $D^+$ is Fredholm (see Prop. \ref{paramMF}). 

The construction of the parametrix is analogous to the construction in the family case \cite{bk}:

Choose a smooth partition of unity $\{\phi_k\}_{k \in J}$ subordinate to the
covering ${\mathcal U}(0, \frac 14)$ (defined in \S \ref{situat}) and a system of
smooth functions
$\{\gamma_k\}_{k \in J}$ on $M$ such that for all $k \in J$ 
\begin{itemize}
\item $\supp \gamma_k \subset \U_k$ and $\gamma_k \phi_k= \phi_k$, 
\item $\ra_{e_2} \gamma_k(x)=\ra_{e_2} \phi_k(x)=0$ for all $x \in \ra M$.
\end{itemize}

For $k \in \bbbz/6$ let $D_{Z_k}$ be the operator defined in \S\ref{Diraccyl} on
$L^2(Z_k,\A^{4d})$ with boundary conditions given by
 $(\Pj_{k \Mod 3},\Pj_{(k+1) \Mod 3})$.
 By Prop. \ref{inv} it is
invertible. Let
$$Q_k:=D_{Z_k}^{-1}:L^2(Z_k,\A^{4d}) \to H(D_{Z_k}) \ .$$ 

Since the symbol of
$D$ is elliptic  and since $\U_{\cp}$ is precompact, there is a
parametrix $Q_{\cp}:L^2(\U_{\cp},E\ta) \to H_{0}^1(\U_{\cp},E \ta)$ such that $\gamma_{\cp}(DQ_{\cp}-1)\phi_{\cp}$
resp. $\gamma_{\cp}(Q_{\cp}D-1)\phi_{\cp}$ is compact on $L^2(\U_{\cp},E\ta)$ resp. $H_{0}^1(\U_{\cp},E
\ta)$ \cite{mf}. Furthermore $Q_{\cp}$
 can be chosen to be an odd operator.

We claim that 
$$Q:=\sum_{k \in J} \gamma_k Q_k \phi_k:L^2(M,E \ta) \to H(D)$$
is a parametrix of $D$. Then it follows that $Q^-:L^2(M,E^-
\ta) \to H(D)^+$ is a parametrix of $D^+$. 

In the following calculations the
operators $D_{Z_k}$ and the restriction of $D$ to $\U_k$ are denoted by $D$ as well. Let  $\sim$
denote equality up to compact operators.

On $L^2(M,E \ten \A)$ 
\begin{eqnarray*}
DQ-1 &=& \sum_{k \in J} [D,\gamma_k]Q_k\phi_k + \sum_{k \in J} \gamma_k DQ_k \phi_k -1\\
&\sim & \sum_{k \in J} c(d \gamma_k)Q_k\phi_k \ .
\end{eqnarray*}

Since $c(d \gamma_k)Q_k\phi_k$ is bounded from $L^2(M,E \ta)$ to
$H^1_{R0}(M_r,E \ta)$ with $r>0$ big enough,  it
is a
compact operator on $L^2(M,E \ta)$ by the previous lemma. Hence $DQ-1$ is compact.

On $H(D)$  
\begin{eqnarray*}
QD-1&=& \sum_{k \in J} \gamma_kQ_k[D,\phi_k] + \sum_{k \in J} \gamma_k Q_kD \phi_k -1 \\
 &\sim& \sum_{k \in J} \gamma_k Q_k c(d \phi_k) \ .
\end{eqnarray*}

Here $c(d \phi_k):H(D) \to L^2(\U_k, E \ta)$ is a compact operator by the
previous lemma since $\supp c(d \phi_k)$ is compact.
Moreover $\gamma_k Q_k:L^2(\U_k,E \ta)
\to H(D)$ is bounded, hence $QD-1:H(D) \to H(D)$ is compact. 
\end{proof}

\begin{lem}
\label{homindD}
Let $P_0,P_1,P_2:[0,1] \to M_{2d}(\A)$ be continuous paths of Lagrangian
projections and assume that $P_0(t),P_1(t),P_2(t)$ are pairwise transverse for each $t
\in [0,1]$. Let $R(t)=(P_0(t),P_1(t),P_2(t))$. We define $D(t)$ to be  the closure of $\dirac$ on $L^2(M,E \ten \A)$ with domain
$\C_{R(t)c}(M, E \ten \A)$. 

Let $\ind D(t)^+$ be the index of the Fredholm operator $D(t)^+:H(D(t))^+ \to
L^2(M,E \ten \A)$.

Then 
$$\ind D(0)^+ =\ind D(1)^+ \in K_0(\A) \ .$$
\end{lem}

\begin{proof}
There is a continuous path of unitaries $[0,1] \to \C(M,\End E^+ \ten \A),~t \mapsto W_t$, such that $W_tD(t)^+W^*_t=D^+_s+W_tc(d W_t^*)$. 
The family $$D^+_s+W_tc(d W_t^*):H(D_s)^+ \to L^2(M,E^- \ten \A)$$ is a continuous path of Fredholm
operators, thus
$\ind W_0D(0)^+W^*_0=\ind W_1D(1)^+W^*_1$ by Prop. \ref{fredhomotop}. 
\end{proof}

\begin{prop}
\label{indeqmasl}

The index of $D^+:H(D)^+ \to L^2(M, E^- \ta)$ is 
$$\ind D^+=\tau(\Pj_0,\Pj_1,\Pj_2) \in K_0(\A) \ .$$
\end{prop}

\begin{proof}
The argument is analogous to the one in \cite{bk}.

By Lemma \ref{uniteq} we may assume $\Pj_0=P_s$.

Let $a_j:=i\frac{p_j+1}{p_j-1} \in M_{2d}(\A)$. Then  $\Pj_j=P(a_j),~j=1,2 ,$ in the notation of \S
\ref{maslov} .

Let $p^+:=1_{\{x>0\}}(a_1-a_2)$.

From Lemma \ref{selfpath} it follows that
$$\tau(\Pj_0,\Pj_1,\Pj_2)=[p^+]-[1-p^+] \in K_0(\A)$$
and that there are continuous paths $P_1,P_2:[0,2] \to M_{2d}(\A)$  of 
Lagrangian projections with $P_j(0)=\Pj_j,~j=1,2$, with $P_1(2)=P(2p^+-1)$ and
$P_2(2)=P(1-2p^+)$ and such that
$P_s,P_1(t),P_2(t)$ are pairwise transverse for all $t \in [0,2]$.  

For $t \in [0,2]$ let $D(t)$ be the Dirac operator on
$L^2(M,E \ten \A)$ whose boundary conditions are given by the triple
$(P_s,P_1(t),P_2(t))$. The previous lemma implies that $\ind D(0)^+=\ind D(2)^+$.

We show that the index of $D(2)^+$ equals $[p^+]-[1-p^+]$:

Let $Q_0=\frac 12\left(\begin{array}{cc} 1&1 \\
1& 1 \end{array}\right) \in M_2(\bbbc)$,
$Q_1=\frac 12 \left(\begin{array}{cc} 1&-i \\
i & 1 \end{array}\right)$ and $Q_2=\frac 12 \left(\begin{array}{cc} 1&i \\
-i & 1 \end{array}\right)$. 

Then $$P_1(2)=(Q_1 \ten p^+) \oplus (Q_2
\ten (1-p^+))$$ and $$P_2(2)=(Q_2 \ten p^+) \oplus (Q_1
\ten (1-p^+))$$ with respect to the decomposition 
$$E^+ \ten \A
= (S  \ten p^+\A^n) \oplus (S  \ten (1-p^+)\A^n) \ .$$
The Dirac operator respects the
decomposition. By \cite{bk} the Dirac
operator associated to the bundle $S\ten (\bbbc^+ \oplus \bbbc^-)$ has index 1
if  the boundary
conditions are given by the triple
$(Q_0,Q_1,Q_2)$, and index $-1$ if they are are given by
$(Q_0,Q_2,Q_1)$. It follows that $$\ind D(2)^+=[p^+]-[1-p^+] \ .$$
\end{proof}

\section{A perturbation with closed range}
\label{stoer}

Imitating the construction in \cite{bgv}, \S9.5, we define a perturbation of $D$ by a compact operator in order to obtain an operator with closed range. Then we can express the index of $D^+$ in terms of the kernel and cokernel of the
perturbed operator.

Choose an orthonormal basis $\{\psi_i\}_{i \in \bbbn} \subset L^2(M,E^-)$ such that
 $\psi_i \in \C_c(M,E^-)$ and $\supp \psi_i \subset M\setminus \ra M$ for all $i \in
 \bbbn$. By 
 Prop. \ref{orthsys} this is an
 orthonormal basis of $L^2(M,E^- \ta)$ as well.

 Since $D^+$ is a Fredholm operator,
 there is a projective $\A$-module $P \subset L^2(M,E^- \ta)$ and a closed
 $\A$-module  $Q \subset \Ran D^+$ such that $P \oplus Q=L^2(M,E^- \ta)$. By 
 Prop. \ref{dirsum} there is $\N\in \bbbn$ such that $L_{\N}:=\spann_{\A}\{\psi_i~|~i=1,\dots, \N\}$ fulfills $L_{\N} + P
 =L^2(M,E^- \ta)$. In particular it follows that
$$L_{\N} + \Ran D^+=L^2(M,E^- \ta) \ .$$

Let $M':=M \cup *$ be the disjoint union of $M$ and one isolated point. Let $E'^+$ be
the hermitian vector bundle $E^+ \cup (* \times
\bbbc^{\N})$ on $M'$, where we endow $\bbbc^{\N}$ with the standard hermitian product. Let $E'^-$ be the hermitian
bundle $E^- \cup (* \times \{0\})$ and let $E'=E'^+\oplus
E'^-$. Extend $D$ by zero to a selfadjoint odd operator $D'$ on  $L^2(M',E' \ten \A)$.

As $D$ is regular, $D'$ is regular as well.

Furthermore $D':H(D') \to L^2(M',E' \ten \A)$ is a Fredholm operator and  $$\ind D'^+ = \ind D^+ +
[\A^{\N}] \ .$$
We extend the standard basis of $\bbbc^{\N}$ by zero to a system $\{e_k\}_{k=1,
  \dots,\N}$ of sections of $E'^+$ and define a compact selfadjoint odd operator $K$ on $L^2(M',E' \ta)$ by
$$Kf:=\sum\limits_{k=1}^{\N}
e_k \langle \psi_k,f\rangle  + \psi_k\langle e_k,f\rangle  \ . $$
Set $D(\rho):=D'+\rho K$ for $\rho \in \bbbr$. Then 
$$\ind D(\rho)^+ = \ind D'^+ = \ind D^+ +
[\A^{\N}] \ .$$

Furthermore $D(\rho)$ is regular by Prop. \ref{regpert}.

By construction $D(\rho)^+$ is surjective for $\rho \neq 0$. Hence by
 Prop. \ref{kerproj} and Prop. \ref{fredclos}
 its kernel is a
projective submodule of $H(D')^+$ and the kernel of $D(\rho)^-$ is
trivial. 

Hence:

\begin{prop}
\label{propstoer}

For $\rho\in  \bbbr$ the operator $D(\rho)^+$ is a Fredholm operator with index independent of $\rho$. For
$\rho\neq 0$ 
$$\ind D^+=\ind D(\rho)^+ -[\A^{\N}]=[\Ker D(\rho)]-[\A^{\N}] \ .$$
\end{prop}

From now on we write $D,E,M$ for $D',E',M'$ and we extend the operator
$D_s$ by zero to the new manifold $M$. Furthermore we redefine the open covering ${\mathcal U}(r,b)$ from \S
\ref{situat} by including the isolated point in the set $\U_{\cp}$.

\chapter{Heat Semigroups and Kernels}

\section{Complex heat kernels}
\label{cplhk}

In this section we collect some properties of the heat
kernels associated to the operators $D_{I_s}$ and $D_s$, which were defined in \S \ref{DI} and \S \ref{stand}, and prove some technical lemmas for further reference. The results are proved by applying standard methods of the theory of partial differential equations to the particular geometric situation. The reader might skip this section at first reading.

\subsection{The heat kernel of $e^{-tD_{I_s}^2}$} 
\label{sgIcp}

Since $D_{I_s}$ is selfadjoint on $L^2([0,1],\bbbc^{2d})$, the operator $-D_{I_s}^2$
generates a semigroup $e^{-tD_{I_s}^2}$ on $L^2([0,1],\bbbc^{2d})$. In this section we determine the
corresponding family of integral
kernels by using the method of images (see \cite{ta}, Ch. 3, \S 7) and study its properties.

The space $L^2([0,1],\bbbc^{2d})$ decomposes into an orthogonal sum 
$$L^2([0,1],P_s\bbbc^{2d}) \oplus L^2([0,1],(1-P_s)\bbbc^{2d})$$
and the semigroup $e^{-tD_{I_s}^2}$ is diagonal with respect to this
decomposition.

We define an embedding
$$\tilde~:L^2([0,1],\bbbc^{2d}) \to L^2(\bbbr/4\bbbz,\bbbc^{2d})$$ by requiring
that $\tilde~$ is a right inverse
of the map $$L^2(\bbbr/4\bbbz,\bbbc^{2d}) \to L^2([0,1],\bbbc^{2d}),~f \mapsto f|_{[0,1]}$$
and that the image of $L^2([0,1],P_s\bbbc^{2d})$ resp. $L^2([0,1],(1-P_s)\bbbc^{2d})$ consists of functions that
are even resp. odd with respect to $y=0$ and $y=2$ and odd resp. even
with respect to $y=1$ and $y=3$.

Then $\tilde~$
maps $\C_R([0,1],\bbbc^{2d})$ with $R=(P_s,1-P_s)$ into $\C(\bbbr/4\bbbz,\bbbc^{2d})$ since 
$$\C_R([0,1],\bbbc^{2d})=\C_l([0,1],P_s\bbbc^{2d})\oplus \C_r([0,1],(1-P_s)\bbbc^{2d}) \ .$$
Here
\begin{eqnarray*}
\lefteqn{\C_l([0,1],P_s\bbbc^{2d})}\\
&:=&\{ f \in \C([0,1],P_s\bbbc^{2d})~|~ (i\ra)^{2k}f(1)=0,~ (i\ra)^{2k+1}f(0)=0~\forall k\in \bbbn_0  \} 
\end{eqnarray*}
and 
\begin{eqnarray*}
\lefteqn{\C_r([0,1],(1-P_s)\bbbc^{2d})}\\
&:=&\{ f \in \C([0,1],(1-P_s)\bbbc^{2d})~|~ (i\ra)^{2k}f(0)=0,~ (i\ra)^{2k+1}f(1)=0~\forall k\in \bbbn_0\} \ .
\end{eqnarray*}
The scalar heat kernel of
$e^{t\ra^2}$  on $\bbbr/4\bbbz$  is 
 $$H(t,x,y)=(4 \pi t)^{-\frac 12} \sum\limits_{k \in \bbbz}
e^{-\frac{(x-y+4k)^2}{4t}} \ .$$
For  $f \in L^2([0,1],(1-P_s)\bbbc^{2d})$ and $x \in [0,1]$ 
\begin{eqnarray*}
(e^{-tD_{I_s}^2}f)(x)&=&(e^{t\ra^2}\tilde f)(x)\\
&=&\int_0^1 H(t,x,y)f(y) dy +\int_1^2 H(t,x,y)f(2-y) dy \\
&&+ \int_2^3 H(t,x,y)(-f(y-2)) dy +\int_3^4 H(t,x,y)(-f(4-y)) dy \ .
\end{eqnarray*}
It follows that the action of $e^{-tD_{I_s}^2}$ on the space 
$L^2([0,1],(1-P_s)\bbbc^{2d})$ is  given by the scalar integral kernel
$$(x,y) \mapsto H(t,x,y)+H(t,x,2-y)-H(t,x,y+2)-H(t,x,4-y) \ .$$
Analogously we conclude that the action of $e^{-tD_{I_s}^2}$ restricted to
$L^2([0,1],P_s \bbbc^{2d})$ is  given by the integral kernel
$$(x,y) \mapsto H(t,x,y)-H(t,x,2-y)-H(t,x,y+2)+H(t,x,4-y) \ .$$
This yields the integral kernel $k_t$ of $e^{-tD_{I_s}^2}$. 

In the following we write $\C_R([0,1],M_{2d}(\bbbc))$ for the space of functions in
$\C([0,1],M_{2d}(\bbbc))$ with column vectors in $\C_R([0,1],\bbbc^{2d})$.

\begin{lem}
\label{smoothI}
The map 
$$(0,\infty) \to \C([0,1],\C_R([0,1],M_{2d}(\bbbc))),~t \mapsto \bigl(y \mapsto
k_t(\cdot,y)\bigr) \ ,$$
is smooth.

For $\phi,\psi \in \C([0,1])$ with $\supp \phi \cap \supp \psi = \emptyset$ the
map $t \mapsto \bigl(y\mapsto \phi(\cdot)k_t(\cdot,y)\psi(y)\bigr)$ can be extended by zero to a smooth map
from $[0,\infty)$ to $\C([0,1],\C_R([0,1],M_{2d}(\bbbc)))$.
\end{lem}

\begin{proof}
This follows from the corresponding properties of $H$.
\end{proof}

\begin{lem}
\label{hkIsest}
\label{hkIsest1}
Let $m,n \in \bbbn_0$. Then there is $C>0$ such that for
all $x,y \in [0,1]$ and all $t >0$ 
$$|\ra_x^m \ra_y^n k_t(x,y)| \le  C(1+t^{-\frac{m+n+1}{2}}) e^{-\frac{d(x,y)^2}{4t}} \ .$$
\end{lem}

\begin{proof}
The assertion follows from the explicit formula of $H$ above. When estimating the
derivatives we
take into account that for all $m \in \bbbn$ the function
$(x,y,t) \mapsto \frac{(x-y)^{2m}}{t^m}e^{-\frac{(x-y)^2}{4t}}$ can be
continuously extended by zero to $t=0$.
\end{proof}

\subsection{The heat kernel of $e^{-tD_s^2}$}

\label{est}

The operator $D_s$  is selfadjoint on the Hilbert space
$L^2(M,E)$. Hence $-D_s^2$ generates a semigroup on $L^2(M,E)$.
In this section we prove the existence of the integral kernel of $e^{-tD_s^2}$
and study its properties. 

In parallel  we study the semigroup $e^{-tD_Z^2}$  on
$L^2(Z,\bbbc^{4d})$ where $Z= \bbbr \times [0,1]$ and $D_Z$ is the
operator defined in \S \ref{Diraccyl}. Here we assume that the boundary conditions are given by a pair
$(P_0,P_1)$ with $P_0,P_1 \in M_{2d}(\bbbc)$.  
Then $D_Z$ is
selfadjoint on $L^2(Z,\bbbc^{4d})$. We will compare $e^{-tD_s^2}$ on the
cylindric ends with $e^{-tD_Z^2}$ with appropriate boundary conditions in order to get estimates for the integral kernel of $e^{-tD_s^2}$.

Since the proofs are standard, they are only sketched here.

Recall that a solution $u:
\bbbr \to \dom D_s$  of the initial-value problem $$\frac
{d}{dt}u(t)=iD_s u(t),~u(0)=f$$ with $f \in \C_{Rc}(M,E)$ is unique by an energy estimate. An analogous statement holds for a
solution of the corresponding problem for $D_Z$.

\begin{lem}
\label{finprop}
If $f \in \C_{Rc}(M, E)$, then $d(x, \supp
f)\le |t|$  for any $x \in \supp (e^{itD_s}f)$. 

An analogous result holds for $D_Z$ on $L^2(Z,\bbbc^{4d})$.
\end{lem}

This property is called ``finite propagation speed property''.

\begin{proof}
The proof relies on a cutting-and-pasting argument. 

For $j \in \bbbz/6$ let $V_j:=\{x \in M~|~ d(x, \ra_j M) < \frac 14\}$ and
let $W:= M \setminus \ra M$. These sets define an open covering of $M$.  We show that the finite propagation speed property holds on these sets for small times. From this we conclude that it holds on $M$ for all times.

By an oriented isometry we identify $V_j$ with
$\{0 \le x_2 < \frac 14\} \subset \bbbr^2$. Recall that $E|_{V_j}$ was identified with the trivial bundle with fiber $\bbbc^{4d}$ in \S \ref{situat}. 

The restriction of $\dirac$ to $V_j$ extends to
a translation invariant differential operator $\dira_{\bbbr^2}:\C_c(\bbbr^2,\bbbc^{4d}) \to
L^2(\bbbr^2,\bbbc^{4d})$.  Let $D_{\bbbr^2}$ be the closure of $\dira_{\bbbr^2}$.

We define an embedding
$$\tilde~: \C_{Rc}(V_j,\bbbc^{4d}) \incl \C_{c}(\bbbr^2,\bbbc^{4d})$$ intertwining the
operators $D_s$ and $D_{\bbbr^2}$ using the method of images similar to \S \ref {sgIcp}. 

Recall that the triple defining the boundary conditions of $D_s$ was denoted by $(\Pj^s_0,\Pj^s_1,\Pj^s_2)$. Let
\begin{eqnarray*}
\lefteqn{\C_l(V_j, \Pj^s_{j \Mod 3} \bbbc^{2d})}\\
&:=&\{ f \in \C_c(V_j, \Pj^s_{j \Mod 3} \bbbc^{2d})~|~
(\ra_{e_2}^{2k+1}f)(x)=0,~\forall k\in \bbbn_0,~\forall x \in \ra_j M\} 
\end{eqnarray*}
and
\begin{eqnarray*}
\lefteqn{\C_r(V_j,(1-\Pj^s_{j \Mod 3}) \bbbc^{2d})}\\
&:=&\{ f \in \C_{c}(V_j,(1-\Pj^s_{j \Mod 3})\bbbc^{2d})~|~
(\ra_{e_2}^{2k}f)(x)=0,~\forall k\in \bbbn_0 ,~\forall x \in \ra_j M \}  \ .
\end{eqnarray*}
Then  $\C_{Rc}(V_j,E^+)$ and $\C_{Rc}(V_j,E^-)$ decompose into a direct sum
$$\C_l(V_j,\Pj^s_{j \Mod 3}\bbbc^{2d}) \oplus
\C_r(V_j,(1-\Pj^s_{j \Mod 3})\bbbc^{2d})\ . $$

For  
$f \in \C_l(V_j,\Pj^s_{j \Mod 3}\bbbc^{2d})$ resp.  $f \in \C_r(V_j,(1-\Pj^s_{j \Mod
3})\bbbc^{2d})$ we define $\tilde f$ by
first extending $f$ by zero to the half plane $\{x_2\ge 0\}$ and then
reflecting 
 such that $\tilde f$ is even resp. odd with respect
to $\{x_2=0\}$. 

For $D_{\bbbr^2}$ the finite propagation speed property holds. Hence the
assertion of the lemma holds for all $f \in \C_{Rc}(V_j,E)$ with $\supp f\subset \{x \in
M~|~ d(x, \ra_j M) < \frac{3}{16}\}$ and for $|t|<\frac{1}{16}$.

For $f \in \C_c(W,E)$ with $\supp f \subset \{x \in
M~|~ d(x, \ra M) > \frac{1}{16}\}$ and for $|t|<\frac{1}{16}$  the assertion
holds by the standard theory of hyperbolic equations on open subsets of $\bbbr^2$.

Since every $f \in \C_c(M,E)$ can be written as $f=f_W+f_0+ \dots f_5$ with
$f_W \in\C_c(W,E)$ and $f_j \in \C_{Rc}(V_j,E)$, the
assertion holds for every $f \in \C_c(M,E)$ and for $|t|< \frac{1}{16}$, and by
the group property of $e^{itD_s}$ it follows for
all $t \in \bbbr$.

The proof for $D_Z$ is analogous.
\end{proof}

For $k \in \bbbn_0$ let $H^k(\bbbc,D_s)$ be the Hilbert space whose underlying vector space is $\dom
D_s^k$ and whose scalar product is given by $$\langle f,g\rangle _{H^k}:=\langle (1+D_s^2)^{\frac
k2}f , (1+D_s^2)^{\frac k2}g\rangle  \ .$$
Define $H^k(\bbbc,D_Z)$ analogously.

\begin{lem} 
\label{sob}
Let $k \in \bbbn,~k \ge 2$.
\begin{enumerate}
\item There is an embedding $H^k(\bbbc,D_Z) \to C^{k-2}(Z,\bbbc^{4d})$.
\item There is an embedding $H^k(\bbbc,D_s) \to C^{k-2}(M,E)$.
\end{enumerate}
\end{lem}

\begin{proof}
We sketch the proof of (2), the proof of (1) is analogous.

Let $D_{\bbbr^2}$ be as in the previous proof.

For fixed $r>0$ the constants in the  G\aa rding inequality for
the elliptic operator $(1+D_{\bbbr^2}^2)^k$ on balls $B_r(x),~x \in \bbbr^2,$
can be chosen
independent of $x$.

For $j \in \bbbz/6$ the embedding $\C_{Rc}(V_j,E) \incl
\C_{c}(\bbbr^2,\bbbc^{4d})$ 
defined in the previous proof intertwines the Dirac operators $D_s$ on $\C_{Rc}(V_j,E)$ and
$D_{\bbbr^2}$ on
$\C_{c}(\bbbr^2,\bbbc^{4d})$. Hence the  G\aa rding inequality for the operator $(1+D_s^2)^k$ on
balls $B_{1/8}(x) \subset M$
with $x \in \ra M$ holds with constants independent of $x$. 

For $r >0$ fixed and small enough we can also find global constants for the G\aa
rding inequality for $(1+D_s^2)^k$ on balls $B_r(x) \subset M$ with $B_r(x) \cap
\ra M =\emptyset$ since $M$ is of bounded geometry.

Then the assertion follows from the Sobolev embedding theorem.
\end{proof}

\begin{cor}
\label{smoothcpker}
The operators $e^{-tD^2_Z}$ on $L^2(Z,\bbbc^{4d})$ and $e^{-tD_s^2}$ on $L^2(M,E)$ are integral operators with smooth integral kernels.
\end{cor}

\begin{proof}
This follows from the previous two lemmas (see \cite{ro}, Lemma 5.6).
\end{proof}

\begin{lem}
\label{gaussest}
Let $f:[0,\infty) \times [0,\infty) \to \bbbr$ be a function and assume
that for every $\ve,\delta>0$ there is $C>0$ such that for all $r>\ve$ and $t > 0$   
$$f(r,t) \le Ce^{-\frac{(r-\ve/2)^2}{(4+\delta)t}} \ .$$
Then for all $\ve,\delta>0$ there is $C>0$ such that for all
$r>\ve$   and $t>0$
$$f(r,t) \le Ce^{-\frac{r^2}{(4+\delta)t}} \ .$$
\end{lem}

\begin{proof}
Choose  $0<a<1$  with $\frac{1-a}{4+\delta/2} >\frac{1}{4+\delta}$ and
let $m>\frac 2a$.

Then there is $C>0$ such that for all $r>\ve$ and $t>0$ 
$$f(r,t) \le Ce^{-\frac{(r-\ve/m)^2}{(4+\delta/2)t}} \ .$$
It follows that
\begin{eqnarray*}
f(r,t) &\le& Ce^{-\frac{(1-a)r^2}{(4+\delta/2)t}}e^{\frac{r}{(4+\delta/2)t}(-ar+\frac{2\ve}{m})}e^{-\frac{(\ve/m)^2}{(4+\delta/2)t}}\\
&\le& C e^{-\frac{r^2}{(4+\delta) t}} \ .
\end{eqnarray*}
In the last step we used the fact that $\frac{r}{(4+\delta/2)t}(-ar+\frac{2\ve}{m})
< 0$ for $r>\ve$.
\end{proof}

\begin{lem}
\label{estL2}
Let $N$ be closed manifold resp. let $N=M,Z$. If $N$ is a closed manifold, let $E_N$ be
a Dirac bundle on $N$ and let $D_N$ be the associated Dirac operator.
 If $N=M$ resp. $N=Z$, then let $E_N=E$ resp. $E_N=Z \times \bbbc^{4d}$ and let
$D_N=D_s$ resp. $D_N=D_Z$. Let $k_t$ be the integral kernel of $e^{-tD_N^2}$. 
 
For every $\ve,\delta>0$ there is $C>0$ such that for all $t>0,~r>\ve$ and $x \in
N$ 
$$ \int_{N \setminus B_r(x)} |k_t(x,y)|^2 ~dy \le C e^{-\frac{r^2}{(4+\delta)t}} \ .$$

Analogous estimates hold for the partial derivatives in $x$ and $y$ with respect
to unit vector fields on $N$.
\end{lem}

\begin{proof}
Let  $k \in 2\bbbn$ with $k>\frac{\dim N}{2}$.

Let $S(x,\ve):=\{u \in \C_c(B_{\ve}(x),E_N)~|~ \|(1+D_N^2)^{-\frac k2} u\| \le 1\}$.

Then by the Sobolev embedding theorem resp. by Lemma \ref{sob} there
is $C>0$ such that for all $x \in N$, $t>0$ and $r>\ve$ 
$$\int_{N \setminus B_r(x)} |k_t(x,y)|^2 ~dy \le C\sup\limits_{u \in S(x,\ve/2)}
\|e^{-tD_N^2}u\|_{N \setminus B_r(x)}^2 \ .$$ 
By a standard argument using the finite propagation speed property of $D_N$ (see the proof of \cite{cgt}, Prop. 1.1) it follows that
\begin{eqnarray*}
\int_{N \setminus B_r(x)} |k_t(x,y)|^2 ~dy &\le & C t^{-1/2}\int\limits_{r-\ve/2}^{\infty}|(1+(\frac{d}{ids})^2)^{k/2}e^{-s^2/4t}| ~ds\\
&=&C  \int\limits_{\frac{r-\ve/2}{\sqrt
t}}^{\infty}|(1+t^{-1}(\frac{d}{ids'})^2)^{k/2}e^{-s'^2/4}|
~ds' \ .
\end{eqnarray*}
Thus there is $l \in \bbbn$ such that 
\begin{eqnarray*}
\int_{N \setminus B_r(x)} |k_t(x,y)|^2 ~dy &\le& C (1+ t^{-l})\int\limits_{\frac{r-\ve/2}{\sqrt t}}^{\infty} (1+s'^l)e^{-s'^2/4} ~ds'\\
&\le& C   (1+ t^{-l}) e^{-\frac{(r-\ve/2)^2}{(4+\delta/2)t}}\\
&\le& C e^{-\frac{(r-\ve/2)^2}{(4+\delta) t}} \ .
\end{eqnarray*}
Then the assertion follows by applying the previous lemma to $$f(r,t):=\sup\limits_{x
\in N}\int_{N\setminus B_r(x)} |k_t(x,y)|^2 ~dy \ .$$ For the derivatives the argument is similar.
\end{proof}

\begin{lem}
\label{estpun}
Let $k_t$ be as in the previous lemma.

For any $\ve,\delta>0$ there is $C<\infty$ such that for all $x,y \in N$ with $d(x,y)>\ve$ and all $t>0$ 
$$|k_t(x,y)| \le C e^{-\frac{d(x,y)^2}{(4+\delta)t}} \ .$$
Analogous estimates hold for the partial derivatives in $x$ and $y$ with respect
to unit vector fields on $N$.
\end{lem}

\begin{proof}
Let $S(y,\ve)$ and $k \in 2\bbbn$ be as in the proof of the previous lemma.
By the Sobolev embedding theorem resp. Lemma \ref{sob} there is
$C>0$ such that for all $r>\ve$, $t>0$ and all $x,y \in N$ with $d(x,y)\ge r$
$$ |k_t(x,y)|\le C\sup\limits_{u \in S(y,\ve/4)}
\|(1+D_N^2)^{\frac k2}e^{-tD_N^2}u\|_{N \setminus B_{r-\ve/4}(y)}^2 \ .$$ 
As in the proof of the previous lemma this implies
$$|k_t(x,y)| \le C e^{-\frac{(r-\ve/2)^2}{(4+\delta) t}} \ .$$
Then the assertion follows from Lemma \ref{gaussest} with $f(r,t):=\sup\limits_{x,y \in N:~d(x,y)=r} |k_t(x,y)|$. For the derivatives the argument is similar.
\end{proof}

For the next lemma assume that $U \subset M$ is an open set for which one of the following properties holds:
\begin{enumerate}
\item  $U$ is precompact and $\ov{U}
\cap \ra M= \emptyset$, 
\item  there is $k \in \bbbz/6$ such that  $U \subset F_k(0,\frac 14)$.
\end{enumerate}

In the first case there is a closed manifold $N$ and a Dirac bundle $E_N$ on
$N$ such that there is a Dirac bundle isomorphism $E|_U \to
E_N$ whose base map is an isometry. We identify $E|_U$ with its image in
$E_N$. Then $D_s$ coincides with $D_N$ on $U$.

In the second case $U$ is a subset of $Z_k$ by  \S \ref{situat}. Let $D_{Z_k}$
be the operator $D_Z$ on
$Z_k \times \bbbc^{4d}$ with boundary conditions
given by $(\Pj_{k \Mod 3},\Pj_{(k+1) \Mod 3})$. Then $D_s$ coincides with $D_{Z_k}$
on $U$.

\begin{lem}
\label{estker}
Let $U$ be as in (1) resp. (2).
Let $k_t$ be the integral kernel of the heat semigroup of $D_s$ on $M$ and let $k'_t$
be the integral kernel of the heat semigroup of $D_N$ resp. $D_{Z_k}$.

For every $T>0$ and $\ve,\delta>0$ there is $C>0$ such that for all
$0<t<T,~ r>\ve$ and $x,y \in U$ with $B_r(x), B_r(y) \subset U$ 
$$|k_t(x,y)-k'_t(x,y)| \le Ce^{-\frac{r^2}{(4+\delta)t}} \ .$$
Analogous estimates hold for the partial derivatives with respect to unit vector
fields on $U$.
\end{lem}

\begin{proof}
The notation is as in the proof of Lemma \ref{estL2}.

The estimate follows from 
\begin{eqnarray*}
|k_t(x,y)-k'_t(x,y)| &=& \sup_{\phi \in S(x,\ve)}\sup_{\psi \in S(y,\ve/2)}|\langle \phi,e^{-tD_s^2}
 \psi\rangle -\langle \phi,e^{-tD_N^2} \psi\rangle |\\
&\le& Ct^{-\frac 12}\sup_{\phi \in S(x,\ve)}\sup_{\psi \in S(y,\ve/2)} |\langle \phi,
 \int_{\bbbr} e^{-s^2/4t}(e^{isD_s}-e^{isD_N})\psi\rangle |\\
&\le&C t^{-\frac
 12}\int\limits_{r-\ve/2}^{\infty}|(1+(\frac{d}{ids})^2)^ke^{-s^2/4t}| \ .
\end{eqnarray*}
Here we used that $(e^{isD_s}-e^{isD_N})\psi=0$ for $|s| \le r- \ve/2$ by the
finite propagation speed property (Lemma \ref{finprop}) and the uniqueness of solutions of
hyperbolic equations.
\end{proof}

\section{The heat semigroup on closed manifolds}
\label{cpman}

Let $\B$ be a Banach algebra with unit.

Let $N$ be a closed manifold of dimension $n$. Let $E_N$ be  a Dirac bundle
 on $N$ and let $D_N$ be  the associated Dirac operator. For simplicity (we will only need this case) assume
that $E_N$ is trivial as a vector bundle.

By Cor. \ref{L2ker} the associated heat kernel defines a family of bounded operators on
$L^2(N,E_N \otimes { \mathcal B})$. The operators are
smoothing, thus they restrict to a family of bounded operators on $C^m(N,E_N \ten
\B)$ for any $m \in \bbbn_0$. In order to show that the family extends to a
holomorphic semigroup we have to study its behavior for small times.

By Lemma \ref{closable} we can define $-D_N^2$ as a closed operator on $L^2(N,E_N\ten \B)$ by requiring that $\C(N,E_N\ten \B)$ is
a core of $-D_N^2$.  

For $t \to 0$  the heat kernel $k_t \in \C(N
\times N,E_N \boxtimes E_N)$ can be estimated as follows: 

\begin{lem}
Let $\ve>0$ be smaller than the injectivity radius of $N$ and let $\chi: [0,\infty) \to
[0,1]$ be a smooth monotonously decreasing function such that $\chi(r)=1$ for $r \le \ve/2$ and $\chi(r)=0$ for $r \ge \ve$. 

Let $A$ be a differential operator  of order $m$ on $\C(N,E_N)$.

Then there is $C >0$ such that for all $x,y \in N$ and for all $t>0$
$$|A_x k_t(x,y)| \le C+
Ct^{-(n+m)/2}e^{-d(x,y)^2/4t}\chi(d(x,y))\sum\limits_{i=0}^m d(x,y)^i
t^{-\frac i2} \ .$$                                                                                                                          
\end{lem}

\begin{proof}
This follows from \cite{bgv}, Prop. 2.46, and its proof.
\end{proof}

\begin{prop}
\label{cpmandiff}

Let $A$ be a differential operator  of order $m$ on $\C(N,E_N \ten \B)$. Then
there is $C>0$ such that the
action of the integral kernel $A_xk_t(x,y)$ on $L^2(N,E_N \ten \B)$ is bounded by
$C(1+t^{- m/2})$ for all $t>0$.
\end{prop}

\begin{proof}
Choose a finite open covering $\{U_{\nu}\}_{\nu \in I}$ of $N$ of normal
coordinate patches and assume that for every $x,y \in U_{\nu}$ the shortest
geodesic connecting $x$ and $y$ is in $U_{\nu}$.

Then there are $c_1,c_2>0$ such that for all $\nu \in I$
and all $x,y \in U_{\nu}$ 
$$c_1|x-y|_{\nu} \le  d(x,y) \le c_2|x-y|_{\nu} \ ,$$
where $| \cdot |_{\nu}$ denotes the euclidian distance on $U_{\nu}$ defined by the
coordinates.

Let $\{\phi_{\nu}\}_{\nu \in I}$ be a
partition of unity subordinate to the covering $\{U_{\nu}\}_{\nu \in I}$.  

Let $\ve>0$ be smaller than the
injectivity radius of $N$ and such that 
$\{x \in N~|~d(x,\supp \phi_{\nu}) \le \ve\} \subset U_{\nu}$ for every $\nu \in I$.
Let $\chi$ be as in the previous lemma. Then
$$\phi_{\nu}(x) \chi(d(x,y)) \le \phi_{\nu}(x)\chi(c_1|x-y|_{\nu})1_{U_{\nu}}(y) \ .$$ 
By the previous lemma there is $C>0$ such that for all $x,y \in N$ and $t>0$ the term
$|A_xk_t(x,y)|$ is bounded by
$$C+Ct^{-(n+m)/2} \sum_{\nu \in J}\phi_{\nu}(x)
\Bigl(e^{-c_1^2|x-y|_{\nu}^2/4t}\chi(c_1|x-y|_{\nu})\sum\limits_{i=0}^m
c_2^i|x-y|_{\nu}^i t^{-\frac i2}\Bigr)1_{U_{\nu}}(y) \ .$$
The $\nu$-th term of the outer sum is supported on $U_{\nu} \times U_{\nu}$. In the coordinates of
$U_{\nu}$ it is of
the form $\phi_{\nu}(x)f_t(x-y)1_{U_{\nu}}(y)$ with $f_t \in
L^1(\bbbr^n)$, and there is $C>0$ such that 
$$\|f_t\|_{L^1} \le Ct^{n/2}$$ for all $t>0$.

The assertion follows now from Prop. \ref{poskern} and Cor. \ref{conv}.
\end{proof}

\begin{prop}
\label{cpmanhsg}
\begin{enumerate}
\item
The family of integral kernels $k_t(x,y)$ defines a bounded strongly continuous semigroup on $L^2(N,E_N
\ten \B)$, which extends to a bounded holomorphic semigroup. Its generator is $-D_N^2$.
\item 
The family of integral kernels $k_t(x,y)$ defines a bounded strongly continuous
semigroup on $C^m(N,E_N \ten \B)$ for every $m \in \bbbn_0$.
\end{enumerate}
\end{prop}

\begin{proof}
(1) By the previous proposition the action of  the integral kernel
$k_t(x,y)$  on $L^2(N,E_N \ten \B)$ is uniformly
bounded for $t > 0$.  On $L^2(N,E_N)\odot \B$ it converges strongly to the
identity. Thus $k_t(x,y)$ induces  a bounded strongly continuous semigroup on $L^2(N,E_N
\ten \B)$. 

On $\C(N,E_N \otimes
{\mathcal B})$ the action
of the generator coincides  with the action of $-D_N^2$. Since $\C(N,E_N \otimes
{\mathcal B})$ is invariant under the semigroup and dense in $L^2(N,E_N
\ten \B)$, it is a core for the generator. Hence the generator is $-D_N^2$.

By the previous proposition there is $C>0$ such that on $L^2(N,E_N
\ten \B)$ for all $0<t<1$ 
$$\|D_N^2e^{-tD_N^2}\| <Ct^{-1} \ .$$  Since $\Ran e^{-tD_N^2} \subset \C(N,E_N \ten
{\mathcal B}) \subset \dom D_N^2$ for $t>0$, it follows, by Prop. \ref{critsg}, that $e^{-tD_N^2}$ extends
to 
a holomorphic semigroup.

The integral kernel of $D_N^2e^{-tD_N^2}$ is exponentially decaying in the
supremum norm for $t \to \infty$, hence $D_N^2e^{-tD_N^2}$ is exponentially
decaying as an operator on $L^2(N,E_N\ten
\B)$.  

By Prop. \ref{critsg} this shows that the holomorphic extension is bounded.

(2) follows from the fact that $k_t(x,y)$ defines a strongly continuous bounded
   semigroup on $C^m(N,E_N)$ by \cite{bgv}, Th. 2.30, and that $C^m(N,E_N
   \ten \B)\cong C^m(N,E_N)\ten_{\ve}\B$.
\end{proof}

It can be deduced from the asymptotic expansion of the heat kernel (\cite{bgv}, Th. 2.30) that $e^{-tD_N^2}$ extends even to a holomorphic semigroup on $C^m(N,E_N \ten {\mathcal B})$ -- we do not need this fact in the following.

\section{The heat semigroup on $[0,1]$}

\subsection{The semigroup $e^{-tD_I^2}$}
\label{Igen}

Let $(P_0,P_1)$ be a pair of transverse Lagrangian projections of $\A^{2d}$ with
$P_0,P_1 \in M_{2d}(\Ai)$. By Lemma \ref{closable} the operator $I_0 \ra$  with domain
$\C_R([0,1],(\Ol{\mu})^{2d})$ is closable on $L^2([0,1],(\Ol{\mu})^{2d})$.  In order to avoid indices its closure is denoted by $D_I$ (compare with the operator $D_I$ from \S \ref{DI}).

Let $U\in M_{2d}(\Ai)$ be a unitary with $UI_0=I_0U$ and $UP_0U^*=P_s$ and let  $p \in M_d(\Ai)$ be such 
that $$UP_1U^*=\tfrac 12 \left( \begin{array}{ccc}
1 & p^*\\
p & 1
\end{array}\right) \ .$$ 
The unitaries $U$ and $p$ exist by Lemma \ref{launit}.

\begin{prop}
\label{DIinv}
Let $\lambda \in \bbbc$ with $\exp(2 i\lambda) \notin \sigma(p)$. 
\begin{enumerate}
\item
The operator
$D_I-\lambda$ has a bounded inverse on $L^2([0,1],(\Ol{\mu})^{2d})$.
\item The inverse  $(D_I-\lambda)^{-1}$ maps
$C^l_R([0,1],(\Ol{\mu})^{2d})$ isomorphically to\\
$C^{l+1}_R([0,1],(\Ol{\mu})^{2d})$ for any $l \in \bbbn_0$.
\item The inverse  $(D_I-\lambda)^{-1}$ maps $L^2([0,1],(\Ol{\mu})^{2d})$
continuously to\\ $C([0,1],(\Ol{\mu})^{2d})$.
\end{enumerate} 
\end{prop}

\begin{proof}
The inverse of $D_I- \lambda$ is given by  
$$((D_I - \lambda)^{-1}f)(x)=\int_0^x I_0 e^{-I_0\lambda(x-y)}f(y)~dy + \int_0^1
e^{-I_0\lambda(x-y)} A(y) f(y) ~dy $$
with $$A(y)=U^*\frac{i}{e^{2i\lambda}-p} \left(\begin{array}{ccc}
p & e^{2i\lambda(1-y)}\\
pe^{2i\lambda y} & e^{2i \lambda}
\end{array}\right)U \ .$$
It is straightforward to check that this map fulfills (1), (2) and (3). 
\end{proof}

We define  $D_I$ as an unbounded operator on
$C^l_R([0,1],(\Ol{\mu})^{2d}),~l \in \bbbn_0,$ by setting 
$\dom D_I:=C^{l+1}_R([0,1],(\Ol{\mu})^{2d}).$
By the previous proposition $D_I$ is a closed operator on
$C^l_R([0,1],(\Ol{\mu})^{2d})$.

As before $D_{I_s}$ denotes $D_I$ with $R=(P_s,1-P_s)$.

We show that $-D_I^2$ generates a holomorphic semigroup on
$L^2([0,1],(\Ol{\mu})^{2d})$ and on $C^l_R([0,1],(\Ol{\mu})^{2d})$. This
will be done by first proving that $-D_{I_s}^2$ generates a holomorphic
semigroup and by
then applying Prop. \ref{sgpert}. 

Moreover  we will deduce a norm estimate of the semigroup $e^{-tD_I^2}$ for large $t$ from the knowledge of the resolvent set
of $-D_I^2$.

We will use the same method again when we study the heat semigroup on $M$.

\begin{lem}
\label{sgIcpex}
Assume that $R=(P_s,1-P_s)$.

The operator $-D_{I_s}^2$   is the generator of a
bounded holomorphic semigroup $e^{-tD_{I_s}^2}$ on $L^2([0,1],(\Ol{\mu})^{2d})$ and on
$C^l_R([0,1],(\Ol{\mu})^{2d})$ for any $l \in \bbbn_0$.
\end{lem}

\begin{proof}
Let $k_t$ be the integral kernel of $e^{-tD_{I_s}^2}$ (see \S \ref{sgIcp}) and let $S(t)$ be the
induced integral operator.

By Lemma \ref{hkIsest} and  Prop. \ref{poskern} the family $S(t)$
is uniformly bounded on $L^2([0,1],(\Ol{\mu})^{2d})$ for $t>0$ and the family
$D_{I_s}^2S(t)$ is bounded by $C(1+t^{-1})e^{-\omega t}$ for some $C,\omega>0$
and all $t>0$.

Since $S(t)$ converges strongly to the identity on
$L^2([0,1],\bbbc^{2d}) \odot \Ol{\mu}$ for $t \to 0$ and has the semigroup property, it is a strongly
continuous semigroup on $L^2([0,1],(\Ol{\mu})^{2d})$.
By Prop. \ref{critsg} it extends to a
bounded holomorphic
semigroup on $L^2([0,1],(\Ol{\mu})^{2d})$.

The kernel $H$ from \S \ref{sgIcp} defines a bounded holomorphic semigroup on\\ $C^k(\bbbr/4\bbbz, (\Ol{\mu})^{2d})$ as
well by \S \ref{fsp}. This
implies that $S(t)$ restricts to a bounded holomorphic semigroup on
$C_R^k([0,1],(\Ol{\mu})^{2d})$. In particular it follows that $\C_R([0,1],(\Ol{\mu})^{2d})$ is a core of
the generator of $S(t)$. Thus the generator is $-D_{I_s}^2$.
\end{proof}

\begin{prop}
\label{sgIgen}
The operator
$-D_I^2$  generates a holomorphic semigroup $e^{-tD_I^2}$ on
$L^2([0,1],(\Ol{\mu})^{2d})$ as well as on $C^k_R([0,1],(\Ol{\mu})^{2d})$ for
all $k \in \bbbn_0$, and there are $C,\omega>0$ such that for all $t\ge 0$
$$\|e^{-tD_I^2}\| \le Ce^{-\omega t}$$ on  $L^2([0,1],(\Ol{\mu})^{2d})$
and on $C^k_R([0,1],(\Ol{\mu})^{2d})$.
\end{prop}

\begin{proof}
The following arguments hold on $L^2([0,1],(\Ol{\mu})^{2d})$
as well as on $C^k_R([0,1],(\Ol{\mu})^{2d})$:

The operator $D_I-U^*D_{I_s}U$ is bounded by Prop. \ref{DIpert}. Since
$U^*D_{I_s}U$ has a bounded inverse by Prop. \ref{DIinv} and $-U^*D_{I_s}^2U$ generates a bounded holomorphic semigroup by
the previous lemma, we can apply Prop. \ref{sgpert}. It follows that
$-D_I^2$ generates a holomorphic
semigroup. 

By Prop. \ref{DIinv} there is $\omega >0$ such that the
spectrum of $D_I^2$ is a subset of $]\omega,
\infty[$. Hence by Prop. \ref{sec} there is $C>0$ with  $\|e^{-tD_I^2}\| \le Ce^{-\omega t} \ .$
\end{proof}

\subsection{The integral kernel} 
\label{hksemI}

In this section we prove the
existence of the heat kernel for $e^{-tD_I^2}$. 
By cutting and pasting we construct an approximation of the semigroup $e^{-tD_I^2}$ by integral
operators. Using Duhamel's principle (see Prop. \ref{duhprin}) we prove that the
error term is an integral operator as well. At the same time we obtain estimates
for the integral kernel of $e^{-tD_I^2}$.

The same method will be used later in order to study the heat kernel on $M$ (see \S \ref{hsDrho}) and the heat kernel associated to the superconnection (see \S \ref{supconMhk}).\\

Let $R=(P_0,P_1)$ be the boundary conditions of $D_I$.

Let $D_{I_0}$ resp. $D_{I_1}$ be  defined as $D_I$ with
boundary conditions given by $(P_0,1-P_0)$ resp. $(1-P_1,P_1)$.  

From Lemma \ref{launit} and \S \ref{sgIcp} it follows that
$e^{-tD_{I_k}^2},~k=0,1,$ is an integral operator for $t>0$. Its integral kernel
is denoted by $e^k_t(x,y)$.

Let $\phi_0:[0,1] \to [0,1]$ be a smooth function with $\supp \phi_0 \subset [0,
\frac 23[$ and $\supp (1-\phi_0) \subset]\frac 13,1]$ and let
$\phi_1:=(1-\phi_0)$.  Furthermore choose smooth functions
$\gamma_0,\gamma_1:[0,1] \to [0,1]$ with
\begin{itemize}
\item $\gamma_k|_{\supp \phi_k} = 1,~k=0,1 ,$
\item $\supp \gamma'_k \cap \supp \phi_k = \emptyset, ~k=0,1 ,$
\item $\supp \gamma_0 \subset [0,\frac 56]$ and  $\supp \gamma_1 \subset [\frac 16 ,1] .$
\end{itemize}

Write $E_t$ for the integral operator corresponding to the integral kernel
$$e_t(x,y):=\gamma_0(x) e^0_t(x,y) \phi_0(y) + \gamma_1(x) e^1_t(x,y) \phi_1(y) \
.$$
Set $E_0:=1$.

Then $E_t$ is  strongly continuous on
 $L^2([0,1],(\Ol{\mu})^{2d})$  as well as on\\
 $C_R^l([0,1],(\Ol{\mu})^{2d})$ at any $t\ge 0$. 

For  $f \in \C_R([0,1],(\Ol{\mu})^{2d})$ the map $[0, \infty) \to
L^2([0,1],(\Ol{\mu})^{2d}),~ t \mapsto E_tf$ is even differentiable. 
Hence
by Duhamel's principle (Prop. \ref{duhprin}) we have for $f \in
\C_R([0,1],(\Ol{\mu})^{2d})$  in
$L^2([0,1],(\Ol{\mu})^{2d})$:
$$(*) \qquad e^{-tD_I^2}f-E_tf= -\int_0^t e^{-sD_I^2}(\frac{d}{dt}+D_I^2)E_{t-s}f~ ds \ .$$

In the following the norm on
$M_{2d}(\A_i)$ is denoted by $|\cdot |$. 

We define $C^k_R([0,1],M_{2d}(\Ol{\mu})),~k \in \bbbn_0 \cup \{\infty\},$
as the space of functions in  $C^k([0,1],M_{2d}(\Ol{\mu}))$  whose column vectors
are in $C^k_R([0,1],(\Ol{\mu})^{2d})$. Then any bounded
operator on $C^k_R([0,1],(\Ol{\mu})^{2d})$ acts as a bounded operator on $C^k_R([0,1],M_{2d}(\Ol{\mu}))$ in
an obvious way.

\begin{prop}
\label{hkIgen}
For $t>0$ the operator $e^{-tD_I^2}$ is an integral operator. Let $k_t$ be its
integral kernel.

\begin{enumerate}
\item The map $$(0,\infty) \to \C([0,1],\C_R([0,1],M_{2d}(\Ai))),~t\mapsto \bigl(y
\mapsto k_t(\cdot,y)\bigr)$$
is well-defined and smooth.

\item  $k_t(x,y)=k_t(y,x)^*$.

\item For every $m,n \in \bbbn_0$ and every $\delta>0$ there is $C>0$ such that
$$|\ra_x^m\ra_y^n k_t(x,y)- \ra_x^m \ra_y^n e_t(x,y)| \le Ct \sum_{k=0,1}
e^{-\frac{d(y,\supp \gamma_k')^2}{(4+ \delta)t}} 1_{\supp \phi_k}(y)$$
for all $t>0$ and all $x,y \in [0,1]$.
\end{enumerate}
\end{prop}

\begin{proof} 
Let $f \in \C_R([0,1],(\Ol{\mu})^{2d})$. From $(*)$ it follows that
$$e^{-tD_I^2}f- E_t f = \sum_{k=0,1}  \int_0^t \int_0^1
e^{-sD_I^2}(\gamma_k'\ra+\ra\gamma_k')e^k_{t-s}(\cdot , y) \phi_k(y) f(y)~ dy ds
\ .$$
By Lemma \ref{smoothI} we can extend the map 
$$t \mapsto \bigl(y \mapsto (\gamma_k'\ra+\ra\gamma_k')e^k_t(\cdot
, y) \phi_k(y)\bigr)$$ by zero to a smooth map from
$[0,\infty)$ to $\C([0,1],\C_R([0,1],M_{2d}(\A_i))$. 

Since $e^{-sD_I^2}$ acts
as a uniformly bounded operator on $\C_R([0,1],M_{2d}(\A_i))$ by Lemma \ref{sgIcp}, it
follows that the operator on the right hand side is an integral operator with smooth integral
kernel. 

Hence $e^{-tD_I^2}$ is an integral operator with
 smooth integral kernel and (1) holds.

The selfadjointness of $e^{-tD_I^2}$ implies (2).

Since  $d(\supp \phi_k, \supp \gamma'_k)> \ve$ for some $\ve>0$, 
there is $C>0$, by Lemma \ref{hkIsest}, such that for all $x,y \in [0,1]$ and $t>0$ 
\begin{eqnarray*}
\lefteqn{|\ra_x^m\ra_y^n\bigl(k_t(x,y)- e_t(x,y)\bigr)|}\\
&\le& C\sum_{k=0,1} \int_0^t \|e^{-sD_I^2}(\gamma_k' \ra + \ra
\gamma_k')\ra_y^n\bigl(e^k_{t-s}(\cdot,y)\phi_k(y)\bigr) \|_{C^m}~ds~ 1_{\supp \phi_k}(y)\\
&\le& C \sum_{k=0,1} \int_0^t  e^{-\frac{d(y,\supp \gamma_k')^2}{(4+\delta)(t-s)}}ds~1_{\supp \phi_k}(y) \\
&\le& C\sum_{k=0,1} te^{-\frac{d(y,\supp  \gamma_k')^2}{(4+\delta)t}}~1_{\supp \phi_k}(y) \ .
\end{eqnarray*}

This shows statement (3).
\end{proof}

\begin{cor}
\label{gaussI}

Let $k_t(x,y)$ be the integral kernel of $e^{-tD_I^2}$.

For every $m,n \in \bbbn_0$ and $\delta,\ve>0$ we find $C>0$ such that for all
$x,y \in [0,1]$ with  $d(x,y)>\ve$ and $t>0$ 
$$|\ra_x^m\ra_y^n k_t(x,y)| \le C e^{-\frac{d(x,y)^2}{(4+\delta)t}} + Ct \sum_{k=0,1} e^{-\frac{d(y,\supp  \gamma_k')^2}{(4+\delta)t}}1_{\supp \phi_k}(y) \ .$$
\end{cor}

\begin{proof}
This follows from the previous proposition and Lemma \ref{hkIsest}.
\end{proof}

\begin{cor}
\label{hkIgenub}
Let $k_t(x,y)$ be the integral kernel of $e^{-tD_I^2}$.

Let  $\omega$ be as in Prop. \ref{sgIgen}.
For every $m,n \in \bbbn_0$ there is $C >0 $ such that for any $t>0$ and any $x,y \in [0,1]$ 
$$|\ra_x^m \ra_y^n k_t(x,y)| \le C(1+t^{-\frac{m+n+1}{2}})e^{-\omega t} \ .$$

\end{cor}

\begin{proof}
There is $C>0$ such that for all $x,y \in [0,1]$  and all $0<t<1$ 
$$ |\ra_x^m \ra_y^nk_t(x,y)-\ra_x^m \ra_y^ne_t(x,y)| \le Ct \sum_{k=0,1}
e^{-\frac{d(y,\supp \gamma_k')^2}{5t}}~1_{\supp \phi_k}(y)  \ ,$$
hence, by Lemma \ref{hkIsest},  
$$|\ra_x^m \ra_y^nk_t(x,y)| \le C(1+t^{-\frac{m+n+1}{2}}) \ .$$
For all $t>1$ and $y
\in [0,1]$ 
$$k_t(\cdot,y)=e^{-(t-1)D_I^2}k_1(\cdot,y) \ .$$
The assertion follows now since $(y \mapsto k_1(\cdot,y)) \in
C^n([0,1],C^m_R([0,1],M_{2d}(\A_i)))$ and since by Prop. \ref{sgIgen} the action of
$e^{-(t-1)D_I^2}$ on $C^m_R([0,1],M_{2d}(\A_i))$ is bounded
by $Ce^{-\omega t}$ for some $C,\omega>0$ and every $t>1$.
\end{proof}

The following facts will be needed for the definition of the $\eta$-form.

\begin{lem} Let $k_t$ be the integral kernel of $e^{-tD_{I_s}^2}$.
  
Then for all $x,y \in [0,1]$ and $t>0$ 
$$\tr(D_{I_s})_xk_t(x,y)=0.$$
\end{lem}

\begin{proof}
Let $S:=2P_s-1\in \Gl_{2d}(\bbbc)$. Then $S^2=1,~SI+IS=0$, $SP_s=P_s$ and $S(1- P_s)=-(1-P_s)$. 

This implies
 $SD_{I_s}e^{-tD_{I_s}^2} + D_{I_s}e^{-tD_{I_s}^2}S=0$.
Therefore  $$S (D_{I_s})_xk_t(x,y)+ (D_{I_s})_xk_t(x,y)S=0 \ ,$$
hence
$$\tr (D_{I_s})_xk_t(x,y)= \tr(-S(D_{I_s})_xk_t(x,y)S)= - \tr (D_{I_s})_xk_t(x,y) \ .$$
It follows that $\tr (D_{I_s})_xk_t(x,y)=0$.
\end{proof}

\begin{cor}
\label{tracet0}
Let $(D_Ik_t)$ be the integral kernel of $D_Ie^{-tD_I^2}$. We have, uniformly on $[0,1]$: 
$$\lim\limits_{t \to 0} \tr (D_Ik_t)(x,x)=0 \ .$$
\end{cor}

\begin{proof}
By the previous lemma $\tr(D_I)_x e_t(x,y)=0$ for all
$x,y \in [0,1]$. Then the assertion follows from the estimate in Prop. \ref{hkIgen}.
\end{proof}

\section{The heat semigroup on the cylinder}
\label{cyl}

Let $Z=\bbbr \times [0,1]$.

Let $R=(P_0,P_1)$  be a pair of pairwise transverse  Lagrangian projections of
$\A^{2d}$  with $P_0,P_1 \in M_{2d}(\Ai)$. 

In this section we study the action of $D_Z$ as an unbounded operator on\\
$L^2(Z,(\Ol{\mu})^{4d})$. If not specified the notation is as in \S \ref{Diraccyl}.\\

First we define the following  function spaces and operators:

For $k\in \bbbn_0 \cup \{\infty\}$ let
\begin{eqnarray*}
\lefteqn{C^k_R(Z,(\Ol{\mu})^{4d})}\\ &:=& \{f \in C^k(Z,(\Ol{\mu})^{4d})~ |~   (P_i \oplus P_i)(\dirz^lf)(x,i)=f(x,i)\\
&& \qquad \text{ for }x \in \bbbr;~i=0,1; ~l \in \bbbn_0, ~l \le k\} \ .
\end{eqnarray*}
Further suffixes, like $c$ or $0$ \dots, have their usual meaning.

We endow these spaces with the subspace topologies.

For a Fr\'echet space $V$ we define the Schwartz space
 $${\mathcal S}(Z,V):=
{\mathcal S}(\bbbr) \ten \C([0,1],V) \ .$$
Moreover let
${\mathcal S}_R(Z,(\Ol{\mu})^{4d})$ be ${\mathcal S}(Z,(\Ol{\mu})^{4d}) \cap
\C_R(Z,(\Ol{\mu})^{4d})$ as a vector space with the topology induced by ${\mathcal S}(Z,(\Ol{\mu})^{4d})$.\\

Let $D_Z$  as an unbounded operator on
$L^2(Z,(\Ol{\mu})^{4d})$ be the closure of  $\dirz$  with domain ${\mathcal S}_R(Z,(\Ol{\mu})^{4d})$. The existence of the closure follows from Lemma \ref{closable}.

Since at the moment it is not clear whether $D_Z^2$ is closed on
$L^2(Z,(\Ol{\mu})^{4d})$, we define $\Delta$  as the closure of
$\dirz^2=-\ra_{x_1}^2-\ra_{x_2}^2$ with domain ${\mathcal S}_R(Z,(\Ol{\mu})^{4d})$.\\

Let $\Delta_{\bbbr}$ be the closure of $-\ra_{x_1}^2$  with domain ${\mathcal S}_R(Z,(\Ol{\mu})^{4d})$.\\

Let $\tilde D_I$ be the closure of $I\ra_{x_2}$ as an unbounded operator
on $L^2(Z,(\Ol{\mu})^{4d})$  with domain ${\mathcal S}_R(Z,(\Ol{\mu})^{4d})$.

By Prop. \ref{DIinv} the operator $D_I$ has a bounded inverse on
$L^2([0,1],(\Ol{\mu})^{2d})$. By Lemma
\ref{identL2} the space $L^2(Z,(\Ol{\mu})^{4d})$
 can be identified with \\  $L^2(\bbbr,
L^2([0,1],(\Ol{\mu})^{4d}))$, hence $\tilde D_I$ has a bounded inverse on
$L^2(Z,(\Ol{\mu})^{4d})$. It follows that $-\tilde D_I^2$ is closed. 

By an analogous
argument $-\tilde D_I^2$
generates a bounded holomorphic semigroup on $L^2(Z,(\Ol{\mu})^{4d})$ with
integral kernel $k^I_t(x_2,y_2)\oplus k^I_t(x_2,y_2)$ for $t>0$.  Here $k^I_t$ is the integral kernel of $e^{-tD_I^2}$, which exists by
Prop. \ref{hkIgen}.\\

Furthermore we have a natural candidate for the integral kernel of a semigroup generated by
$-\Delta$, namely
$$k_t^Z(x,y):=\frac{1}{\sqrt{4 \pi
t}}~e^{-\frac{(x_1-y_1)^2}{4t}}
(k^I_t(x_2,y_2) \oplus k^I_t(x_2,y_2)) \ .$$\\ 
In the following let $\omega>0$ be as in Prop. \ref{sgIgen} such that
$$\|e^{-t\tilde D_I^2}\| \le Ce^{-\omega t}$$ 
 on $L^2(Z,(\Ol{\mu})^{4d})$ for all $t\ge 0$.

\begin{prop}
\label{sgZ}
\begin{enumerate}
\item
The integral kernel $k_t^Z(x,y)$
defines a  holomorphic semigroup
 on $L^2(Z,(\Ol{\mu})^{4d})$ with generator $-\Delta$.

\item
For every $m \in \bbbn_0$ the integral kernel $k_t^Z(x,y)$ defines a holomorphic semigroup on
$C^m_R(Z,(\Ol{\mu})^{4d})$, denoted by $e^{-t\Delta}$ as well.

\item
Let $A$ be a differential operator of order $m$ with coefficients in\\ $\C(Z,M_{4d}(\Ol{\mu}))$.
Then for the operator $Ae^{-t\Delta}$ on $L^2(Z,(\Ol{\mu})^{4d})$ as
well as for
$Ae^{-t\Delta}:C^n_R(Z,(\Ol{\mu})^{4d}) \to C^n(Z,(\Ol{\mu})^{4d}),~n \in
\bbbn_0,$ we have:

There is $C>0$ such that for all $t>0$   
$$\| Ae^{-t\Delta} \| \le C(1+t^{-m/2})e^{-\omega t} \ .$$

\end{enumerate}
\end{prop}

\begin{proof}
(1) By Prop. \ref{poskern} the
kernel $\frac{1}{\sqrt{4 \pi
t}}~e^{-\frac{(x_1-y_1)^2}{4t}}$ defines a uniformly bounded family of operators on
$L^2(Z,(\Ol{\mu})^{4d})$. Since it converges strongly to the identity on
$L^2 (Z) \odot (\Ol{\mu})^{4d}$ for $t \to 0$, it is a strongly continuous
semigroup. The space ${\mathcal S}_R(Z,(\Ol{\mu})^{4d})$ is invariant under the
action of the semigroup and the action of the generator on that space
is given by $\ra_{x_2}^2$. Hence the generator is
$-\Delta_{\bbbr}$. 

By checking the assumptions of Prop. \ref{critsg} we show that the semigroup
$e^{-t\Delta_{\bbbr}}$ extends to a holomorphic one:

The operator $(i\ra_{x_1})e^{-t\Delta_{\bbbr}}$ equals the convolution with
the function $$g(x_1):=\frac{-i}{\sqrt{ 4 \pi t}} \left(\frac{x_1}{2t}\right)
e^{-x_1^2/4t} \ .$$ 
Since  there is $C>0$ such that for
$0<t<1$ 
$$\|g\|_{L^1} \le C t^{-1/2} \ ,$$ it follows that
$$\|(i\ra_{x_1}) e^{-t\Delta_{\bbbr}}\| \le C t^{- 1/2} $$
on $L^2(Z,(\Ol{\mu})^{4d})$, thus for $0<t<1$
$$\|\Delta_{\bbbr} e^{-t\Delta_{\bbbr}}\| \le \|(i\ra_{x_1})
e^{-(t/2)\Delta_{\bbbr}}\|^2 \le C t^{-1} \ .$$
Hence $-\Delta_{\bbbr}$ generates a holomorphic semigroup on $L^2(Z,(\Ol{\mu})^{4d})$. 

Note that this estimate also holds on $C_R(Z,(\Ol{\mu})^{4d})$ showing that
$e^{-t\Delta_{\bbbr}}$ is a holomorphic semigroup on  $C_R(Z,(\Ol{\mu})^{4d})$
as well.

Since the semigroups $e^{-t\Delta_{\bbbr}}$ and $e^{-t\tilde D_I^2}$ commute with each other, their composition is a
holomorphic semigroup. The space
${\mathcal S}_R(Z,(\Ol{\mu})^{4d})$  is invariant under
the action of the semigroup and the generator acts on it as $\ra_{x_1}^2+
\ra_{x_2}^2$. Thus the generator is $-\Delta$. 

(2) Since for $f \in
C_R^n(Z,(\Ol{\mu})^{4d}), ~n \in \bbbn,$ we have that $(i\ra_{x_1})e^{-t\Delta}f=e^{-t\Delta}(i\ra_{x_1})f$ and $\tilde D_Ie^{-t\Delta}f= e^{-t\Delta} \tilde D_If$, 
 the assertion
can be reduced to the case $n=0$.
  
From Prop. \ref{sgIgen} it follows that the action of the integral kernel $k^I_t(x_2,y_2)
\oplus k^I_t(x_2,y_2)$  on
$C_R(Z,(\Ol{\mu})^{4d})$ extends to a holomorphic semigroup. Furthermore in (1) we showed that the integral kernel
$\frac{1}{\sqrt{4 \pi t}}e^{-\frac{(x_1-y_1)^2}{4t}}$ defines
a holomorphic semigroup on $C_R(Z,(\Ol{\mu})^{4d})$. Hence the kernel $k^Z_t(x,y)$ defines
a semigroup on $C_R(Z,(\Ol{\mu})^{4d})$ that extends to a holomorphic one.

(3) We can restrict 
   to the case $n=0$ by
the argument
in the proof of (2). 

In the following the operator norms can be understood with respect to
the action on $L^2(Z,(\Ol{\mu})^{4d})$ as well as with respect to the action on $C_R(Z,(\Ol{\mu})^{4d})$. 

The differential operator $A$ is a sum of operators $a_{hk}
\tilde D_I^h (i\ra_{x_1})^k$ with $a_{hk} \in \C(Z, M_{4d}(\Ol{\mu}))$ and $h+k \le m$.
We have
$$\tilde D_I^h (i\ra_{x_1})^k e^{-t\Delta}=\tilde D_I^h e^{-t \tilde D_I^2}(i\ra_{x_1})^k e^{-t\Delta_{\bbbr}} \ .$$
By Cor. \ref{semest} there is $C>0$ such that for $0<t$
$$\|\tilde D_I^he^{-t \tilde D_I^2}\| \le Ct^{- h/2}e^{-\omega t}  \ .$$
By the estimate in the proof of (1) 
 $$\|(i\ra_{x_1})^ke^{-t\Delta_{\bbbr}}
\| \le \|(i\ra_{x_1})e^{-(t/k)\Delta_{\bbbr}}\|^k \le C t^{- k/2}$$
 for $0<t<1$.

Now the assertion follows by taking into account that $e^{-t\Delta_{\bbbr}}$ is uniformly
 bounded.
 \end{proof}

\begin{cor}\label{DinvZ}
Let $\lambda \in \bbbc$ with $\re \lambda^2 < \omega$.
\begin{enumerate}
\item The operator $D_Z-\lambda$ is invertible on ${\mathcal S}_R(Z,(\Ol{\mu})^{4d})$.
\item
The operator $D_Z - \lambda$  is invertible on $L^2(Z, (\Ol{\mu})^{4d})$.
\item $\Delta=D_Z^2$.
\end{enumerate}
\end{cor}

\begin{proof}
(1) For any seminorm $p$ of ${\mathcal S}_R(Z,(\Ol{\mu})^{4d})$ there are $C>0$, $n \in \bbbn$ and a seminorm $q$
 such that 
$p(e^{-t\Delta_{\bbbr}}f) \le C(1+t^n)q(f)$ and $p(e^{-t\tilde D_I^2}f) \le
Ce^{-\omega t}q(f)$ for all $f \in {\mathcal S}_R(Z,(\Ol{\mu})^{4d})$.  Hence 
$e^{-t\Delta}$ restricts to a bounded operator on ${\mathcal S}_R(Z,(\Ol{\mu})^{4d})$, and
 the integral 
$$G(\lambda)=\int_0^{\infty}(D_Z+\lambda)e^{-t(\Delta-\lambda^2)} f dt $$ defines a bounded operator on ${\mathcal S}_R(Z,(\Ol{\mu})^{4d})$ inverting
$D_Z-\lambda$.

(2) The operator $G(\lambda)$ extends to a bounded operator on
   $L^2(Z,(\Ol{\mu})^{4d})$ since by the previous proposition there is $C>0$ such that on
   $L^2(Z,(\Ol{\mu})^{4d})$
$$\|(D_Z+\lambda)e^{-t(\Delta-\lambda^2)} \| \le C(1+t^{-\frac 12}) e^{-(\omega-\re
\lambda ^2) t}$$
for all $t>0$.

From (1) it follows that $G(\lambda)$ inverts $D_Z-\lambda$ on
$L^2(Z,(\Ol{\mu})^{4d})$.

(3) From (2) it follows that the
   operator $D_Z^2$ is closed, and from (1) that ${\mathcal S}_R(Z,(\Ol{\mu})^{4d})$ is a core for $D_Z^2$. Hence the closure of $D_Z^2$ equals $\Delta$. 
\end{proof}

\begin{prop}
\label{embed}

Let $\lambda \in \bbbc$ with $\re \lambda< \omega$.
\begin{enumerate}
\item  The operator $(D_Z^2-\lambda)^{-1}$ maps
$L^2(Z,(\Ol{\mu})^{4d})$ continuously to\\ $C(Z, (\Ol{\mu})^{4d})$.
\item Let $n \in \bbbn, n\ge 2$.
The operator $(D_Z^2-\lambda)^{-n}$ maps
$L^2(Z,(\Ol{\mu})^{4d})$ continuously to $C^{2n-3}(Z, (\Ol{\mu})^{4d})$.
\end{enumerate}
\end{prop}

\begin{proof}
(1)
For $\re \lambda < \omega$ we have on $L^2(Z,(\Ol{\mu})^{4d})$:
\begin{eqnarray*}
(D_Z^2-\lambda)^{-1} &=& \int_0^{\infty} e^{-t(D_Z^2-\lambda)} dt \\
&=& \int_0^{\infty} e^{\lambda t}e^{-t\Delta_{\bbbr}^2}e^{-t\tilde D_I^2} dt \ .
\end{eqnarray*}
By Prop. \ref{DIinv} the operator $D_I^{-1}:L^2([0,1],(\Ol{\mu})^{2d}) \to
C([0,1],(\Ol{\mu})^{2d})$ is bounded.
Thus the family of operators $$e^{-t \tilde D_I^2}=\tilde D_I^{-1} \tilde D_I
e^{-t\tilde D_I^2}:L^2(Z,(\Ol{\mu})^{4d}) \to L^2(\bbbr, C([0,1],(\Ol{\mu})^{4d}))$$ is
bounded by $C(1+t^{-\frac 12})e^{-\omega t}$ for all $t>0$. 

Furthermore the family 
$$e^{-t\Delta_{\bbbr}}:
L^2(\bbbr,C([0,1],(\Ol{\mu})^{4d})) \to C(\bbbr,C([0,1],(\Ol{\mu})^{4d}))$$ 
is bounded by $\sup\limits_{x_1 \in \bbbr}\|\frac{1}{\sqrt{4 \pi
t}}e^{-\frac{(x_1-y_1)^2}{4t}} \|_{L^2_{y_1}}$, hence by $Ct^{- 1/4}$ for some $C>0$.

Thus the integral converges as a bounded operator
from $L^2(Z,(\Ol{\mu})^{4d})$ to $C(Z, (\Ol{\mu})^{4d})$.

(2) We show that for any $k \in \bbbn_0$ with $k\le 2n-1$ the map
$$\tilde D_I^k (i\ra_{x_1})^{2n-1-k}(D_Z^2-\lambda)^{-n}:{\mathcal S}_R(Z,(\Ol{\mu})^{4d}) \to
L^2(Z,(\Ol{\mu})^{4d})$$ extends to a bounded operator on $L^2(Z,(\Ol{\mu})^{4d})$. Then the assertion follows from the first part.

We have for $k \in \bbbn_0$ with $k \le 2n+1$
\begin{eqnarray*}
\lefteqn{\tilde D_I^k (i\ra_{x_1})^{2n-k+1}(D_Z^2-\lambda)^{-n-1} }\\
&=& \int_0^{\infty} t^{n-1}\tilde D_I^k (i\ra_{x_1})^{2n-k+1} e^{-t(D_Z^2-\lambda)}~dt\\
&=&\int_0^{\infty} t^{n-1}\tilde D_I^k e^{-t(\tilde D_I^2-\lambda)}(i\ra_{x_1})^{2n-k+1} e^{-t\Delta_{\bbbr}} ~ dt \ .
\end{eqnarray*}
As bounded operators on $L^2(Z,(\Ol{\mu})^{4d})$ 
$$\|\tilde D_I^k e^{-t(D_I^2-\lambda)}\| \le Ct^{-k/2}e^{-(\omega -\lambda)t}$$
and $$\| (i\ra_{x_1})^{2n-k+1} e^{-t\Delta_{\bbbr}}\| \le C(1+t^{-(2n-k+1)/2}) \ ,$$   hence the integral converges.
\end{proof}

In the following $| \cdot |$ denotes the norm on $M_{4d}(\A_i)$. 

\begin{lem}
\label{gaussZ}
For every $\ve>0$ and $\alpha,\beta \in \bbbn_0^2$ there are $c,C>0$ such that
for all $x,y \in Z$ with $d(x,y)>\ve$ and all $t>0$ 
$$|\ra_x^{\alpha}\ra_y^{\beta} k^Z_t(x,y)| \le
C e^{-\frac{d(x,y)^2}{ct}} \ .$$
\end{lem}

\begin{proof}
For $m,n \in \bbbn_0$ there are $C,c>0$ such that 
$$|\ra_{x_2}^m\ra_{y_2}^n k^I_t(x_2,y_2)| \le Ce^{-\frac{(x_2-y_2)^2}{c t}}$$
for $|x_2-y_2| \ge \ve/2$ and $t>0$.
This follows from Cor. \ref{gaussI} for $t<1$ and from Cor. \ref{hkIgenub} for $t
\ge 1$. For $|x_2-y_2| \le \ve/2$ the left hand side is bounded by
$C(1+t^{-\frac{m+n+1}{2}})$ by Cor. \ref{hkIgenub}.
Similar estimates hold for $\frac{1}{\sqrt{4 \pi t}}~e^{-\frac{(x_1-y_1)^2}{4t}}$. A combination of these estimates implies the lemma.
\end{proof}

In the next lemma ${\mathcal S}_R(Z,M_{4d}(\A_i))$ is the subspace of  ${\mathcal
S}(Z,M_{4d}(\A_i))$ of functions whose columns are in ${\mathcal
S}_R(Z,\A_i^{4d})$. Then operators on ${\mathcal
S}_R(Z,\A_i^{4d})$ act on ${\mathcal S}_R(Z,M_{4d}(\A_i))$ columnwise. The space $\C_{Rc}(Z,M_{4d}(\A_i))$ is analogously defined.

\begin{lem}
\label{resker}
Let $\lambda \in \bbbc$ with $\re \lambda < \omega$.  Let $\xi_1,\xi_2 \in
\C(Z)$  with $\supp \xi_1 \cap \supp \xi_2 =\emptyset$ and assume that 
$\supp \xi_2$ is compact.

Then for any $n \in \bbbn$ the operator $\xi_1 (D_Z^2-\lambda)^{-n} \xi_2$ is an integral
operator. 

Let  $\kappa$  be its integral kernel. Then $(y \mapsto \kappa(\cdot,y)) \in
\C_{c}(Z,{\mathcal S}_R(Z,M_{4d}(\A_i)))$ and $(x \mapsto \kappa(x,\cdot)^*) \in {\mathcal S}(Z,\C_{Rc}(Z,M_{4d}(\A_i)))$.

In particular $\xi_1 (D^2-\lambda)^{-n} \xi_2$ maps $L^2(Z,(\Ol{\mu})^{4d})$
continuously to\\
${\mathcal S}_R(Z,(\Ol{\mu})^{4d})$.
\end{lem}

\begin{proof}
First we prove the claim for $n=1$.

Let $f \in \C_{Rc}(Z, (\Ol{\mu})^{4d})$.

In $L^2(Z,(\Ol{\mu})^{4d})$
\begin{eqnarray*}
\xi_1 (D_Z^2-\lambda)^{-1} \xi_2 f &=&\int _0^{\infty} \xi_1 e^{-t(D_Z^2 -
\lambda)} \xi_2  f ~dt  \\
&=&\int_0^{\infty} \int_Z\xi_1k^Z_t(\cdot,y)  e^{\lambda t}  \xi_2(y)  f(y) ~dy
dt \ .
\end{eqnarray*}
Let $\ve>0$ be such that $d(\supp \xi_1,\supp\xi_2) >\ve$.

By the previous lemma there are $c,C>0$ such that for all $x,y
\in Z$ and all $t>0$ 
\begin{eqnarray*}
|\xi_1(x)e^{\lambda t} k_t^Z(x,y) \xi_2(y)|
& \le & C1_{\{t>1\}}(t)|\xi_1(x)| e^{(\re \lambda -\omega) t} e^{-\frac{(x_1-y_1)^2}{4t}}
|\xi_2(y)|\\
&& + ~C1_{\{t \le 1\}}(t) |\xi_1(x)| e^{-\frac{d(x,y)^2}{ct}}
|\xi_2(y)| \ .
\end{eqnarray*}
Analogous estimates hold for the partial derivatives. Hence  we can
interchange the order of integration. It follows that $\xi_1 (D_Z^2-\lambda)^{-1}\xi_2$ is
an integral operator with integral kernel
$$\kappa(x,y):=\int_0^{\infty} \xi_1(x)e^{-\lambda t} k^Z_t(x,y) \xi_2(y) dt \ .$$   For $n=1$ the other
statements of the lemma also follow from the estimates.

For $n>1$ choose a smooth compactly supported function $\psi:Z \to [0,1]$ such that $\supp
\psi \cap \supp \xi_1= \emptyset$ and $\supp(1-\psi) \cap \supp
\xi_2=\emptyset$. Then
\begin{eqnarray*}
\xi_1 (D_Z^2-\lambda)^{-n} \xi_2 &=&
\xi_1(D_Z^2-\lambda)^{-1}\psi(D_Z^2-\lambda)^{-n+1} \xi_2 \\
&& +~
\xi_1(D_Z^2-\lambda)^{-1}(1-\psi)(D_Z^2-\lambda)^{-n+1} \xi_2 \ .
\end{eqnarray*}
By induction the lemma can be applied to
$\xi_1(D_Z^2-\lambda)^{-1}\psi$ and $(1- \psi)(D_Z^2-\lambda)^{-n+1} \xi_2 $. The statement of the lemma
follows for $\xi_1 (D_Z^2-\lambda)^{-n} \xi_2$ from this and the fact that by Cor. \ref{DinvZ} the operator $(D_Z^2-\lambda)^{-m}$ acts continuously on ${\mathcal S}_R(Z,M_{4d}(\A_i))$ for all $m \in \bbbn$. 
\end{proof}

\section{The heat semigroup on $M$}

\subsection{Definitions}
\label{hsDrho}
Recall the definition of the operator $D(\rho)^2$ on the Hilbert
$\A$-module $L^2(M,E \ten \A)$ in \S \ref{defiD} and \S \ref{stoer}. By Lemma \ref{closable} the operator $(\dirac + \rho K)^2$ with domain ${\mathcal S}_R(M, E \ten \Ol{\mu})$ is closable on $L^2(M,E \ten \Ol{\mu})$. Its closure will be denoted by $D(\rho)^2$ as well in order to simplify notation. 
 
So far the notation is misleading: It suggests that $D(\rho)^2$ is the
square of some unbounded operator on $L^2(M,E \ten \Ol{\mu})$. This is indeed the case as will become clear in
\S \ref{semMgen}.

Define $D_s^2$ as a closed operator on $L^2(M,E \ten \Ol{\mu})$ in an analogous way.\\

We will often make use of cutting and pasting arguments on $M$. We fix the setting:

Let $0<b_0 \le \frac 14$ be small enough and
$r_0>0$ large enough such that
$$\supp k_K \cap \bigl((F(r_0,b_0) \times M) \cup (M \times F(r_0,b_0))\bigr)= \emptyset
\ ,$$
with $F(r_0,b_0)$ as in \S \ref{situat}. Here $k_K$ denotes the integral kernel
of the operator $K$ from \S \ref{stoer}.

Let ${\mathcal U}(r_0,b_0)$ be the open covering defined in \S \ref{situat}
and modified in \S \ref{stoer}. 

Choose a smooth partition of unity $\{\phi_k\}_{k \in J}$
subordinate to ${\mathcal U}(r_0,b_0)$ and smooth functions $\{\gamma_k\}_{k \in J}$ on $M$ such
that for all $k \in J$  
\begin{itemize}
\item $\supp \gamma_k \subset \U_k$,
\item $\supp (1-\gamma_k) \cap \supp \phi_k= \emptyset$,
\item the derivatives 
$\ra_{e_2} (\phi_k|_F)$ and $\ra_{e_2} (\gamma_k|_F)$ vanish in a neighborhood of $\ra M$.
\end{itemize}
  
Let $E_N$ be a Dirac bundle on a compact spin manifold $N$ that is
trivial as a vector bundle and assume that there is a Dirac bundle isomorphism $E|_{\U_{\cp}} \to
E_N$, whose base map is an isometric embedding. Let $D_N$ be the associated
Dirac operator. We identify $\U_{\cp}$ with its
image in $N$ and $E|_{\U_{\cp}}$ with
its image in $E_N$. 

Since the support of $k_K$ is in $\U_{\cp} \times \U_{\cp}$, the restriction of
$D(\rho)=D+\rho K$ to $\U_{\cp}$ extends to an operator $D_N+\rho K$ on the sections
of $E_N$.

For $k \in \bbbz/6$ let $D_{Z_k}$ be the operator $D_Z$ on $L^2(Z, (\Ol{\mu})^{4d})$ from \S \ref{cyl} with boundary conditions
given by the pair $(\Pj_{k \Mod 3},\Pj_{(k+1) \Mod 3})$.

\subsection{The resolvents of $D(\rho)^2$}
\label{resD}
This section has three different aims:

Using a method of Lott (\cite{lo2}, \S 6.1.) we investigate the
resolvent set of $D(\rho)^2$ on $L^2(M,E \ten \Ol{\mu})$. 

Furthermore we prove a kind of Sobolev embedding theorem -- more precisely an
analogue of Lemma \ref{sob} for the operator $D(\rho)^2$ on $L^2(M,E \ten
\Ol{\mu})$.

Third we obtain more information about the kernel of $D(\rho)^2$ on $L^2(M,E \ten
\Ol{\mu})$, namely that there is a projection on it and that this projection is a
Hilbert-Schmidt operator with a smooth integral kernel.\\

Let $\omega>0$ be such that there is $C >0$ with  
$$\| e^{-tD_{Z_k}^2} \|
\le Ce^{-\omega t}$$ on $L^2(M,E \ten \Ol{\mu})$ for all $t\ge 0$ and all $k \in \bbbz/6$.

Let $\nu \in \bbbn$.
For $\lambda \in \bbbc$ with $\re \lambda <\omega$ we define a parametrix of
$(D(\rho)^2-\lambda)^{\nu}$:

By
Cor. \ref{DinvZ} we can set $Q_k(\lambda)=(D_{Z_k}^2-\lambda)^{-\nu}$ for $k \in \bbbz/6$. 

Let $Q_{\cp}(\lambda)$ be a local parametrix of $(D^2-\lambda)^{\nu}$ on $\U_{\cp}$ defined by the
symbol of $(D^2-\lambda)^{\nu}$ such that $\phi_{\cp}(Q_{\cp}(\lambda)(D^2-\lambda)^{\nu}-1)\gamma_{\cp}$ and $\phi_{\cp}((D^2-\lambda)^{\nu}Q_{\cp}(\lambda) -1)\gamma_{\cp}$ are integral operators with smooth integral kernels. 

The operator  $$Q(\lambda):=\sum_{k\in J} \phi_k Q_k(\lambda) \gamma_k $$ 
acts as a bounded operator on the spaces $L^2(M,E \ten \Ol{\mu})$ and ${\mathcal S}_R(M,E \ten \Ol{\mu})$ by \S \ref{pseudo}
and by  Cor. \ref{DinvZ}.  

\begin{lem}
For any $\rho \in \bbbr$ and $\lambda \in \bbbc$ with $\re \lambda< \omega$
the operator $Q(\lambda)(D(\rho)^2-\lambda)^{\nu}-1$ restricted to ${\mathcal S}_R(M,E \ten \Ol{\mu})$ is an integral operator $\Kappa$ with smooth integral kernel $\kappa \in L^2(M
\times M, (E \boxtimes 
E^*) \ten \A_i)$.

Furthermore $(x \mapsto \kappa(x, \cdot)^*) \in  {\mathcal S}(M, \C_{cR}(M, E\ten
\A_i)\ten E^*)$ and $(y \mapsto \kappa( \cdot, y)) \in \C_c(M, {\mathcal S}_R(M,E \ten
\A_i)\ten E^*)$. 

In particular $\Kappa$ extends to a bounded operator from $L^2(M, E \ten \Ol{\mu})$ to
${\mathcal S}_R(M, E \ten \Ol{\mu})$.
\end{lem}

\begin{proof} The difference
$(D(\rho)^2-\lambda)^{\nu}-(D^2-\lambda)^{\nu}$ on ${\mathcal S}_R(M,E \ten \Ol{\mu})$ is an integral operator with smooth
integral kernel whose support is contained in $\supp k_K$.
Hence we only need to investigate $Q(\lambda)(D^2-\lambda)^{\nu}-1$.

For any $k \in J$ choose a function $\xi_k \in \C_c(M)$
with values in $[0,1]$ and such that $\supp
\xi_k \subset \U_k,~ \xi_k|_{\supp d \gamma_k}=1$ and $\supp \phi_k \cap
\supp  \xi_k=\emptyset$. Furthermore assume that  $\ra_{e_2}(\xi_k|_F)$ vanishes in a
neighborhood of $\ra M$.

By induction we have that $[\gamma_k,(D^2-\lambda)^{\nu}]=\xi_k[\gamma_k,(D^2-\lambda)^{\nu}]$ since
\begin{eqnarray*}
\lefteqn{[\gamma_k,(D^2-\lambda)^{\nu}]}\\
&=&(D^2-\lambda)^{\nu-1}\bigl(c(d\gamma_k)D+D
c(d\gamma_k)\bigr) + [\gamma_k,(D^2-\lambda)^{\nu-1}](D^2-\lambda)\\
&=&\xi_k (D^2-\lambda)^{\nu-1}\bigl(c(d\gamma_k)D+D
c(d\gamma_k)\bigr)  + [\gamma_k,(D^2-\lambda)^{\nu-1}](D^2-\lambda) \ .
\end{eqnarray*}
In the following the operators $D_{Z_k},~k \in \bbbz/6,$ are denoted by $D$
as well. Furthermore $\sim$ means equality up to integral operators with smooth
compactly supported integral kernels.

Then on ${\mathcal S}_R(M,E \ten \Ol{\mu})$
\begin{eqnarray*} 
Q(\lambda)(D^2-\lambda)^{\nu}-1 &=&\sum_{k \in J} \phi_k Q_k(\lambda)[\gamma_k,(D^2-\lambda)^{\nu}] +\sum_{k \in J}\phi_k Q_k(\lambda)(D^2-\lambda)^{\nu}\gamma_k -1\\
&\sim& \sum_{k \in J} \phi_k Q_k(\lambda)\xi_k[\gamma_k,(D^2-\lambda)^{\nu}] \ .
\end{eqnarray*}
For all $k \in J$ the operator $\phi_k
Q_k(\lambda)\xi_k$ is an integral operator whose integral kernel has the
properties stated in the lemma. This holds for $k \in
\bbbz/6$ by Lemma \ref{resker}  and for $k=\cp$ by the
properties of pseudodifferential operators. Now the
assertion follows.
\end{proof}

\begin{prop}
\label{spec}
Let $\rho \in \bbbr$.
Let $\lambda \in \bbbc$ with $\re \lambda<\omega$ such that
$D(\rho)^2-\lambda$ has a bounded inverse on the Hilbert $\A$-module $L^2(M,E \ten \A)$.

Then
$D(\rho)^2-\lambda$ has a bounded inverse on
$L^2(M,E \ten \Ol{\mu})$.

The inverse $(D(\rho)^2-\lambda)^{-1}$
acts as a bounded operator on the space ${\mathcal S}_R(M, E \ten \Ol{\mu})$.
\end{prop}

\begin{proof}
Let $Q(\lambda)$ and $\Kappa$ be as in the previous lemma such that $Q(\lambda)(D(\rho)^2-\lambda)f=(1-\Kappa)f$ for $f \in \dom D(\rho)^2$.

We want to apply Prop. \ref{intker}. Since in general $1-\Kappa$ is not
invertible on $L^2(M, E\ten \A)$, we modify the parametrix:

Choose an integral kernel $s \in \C_c(M\times M,(E \boxtimes E^*) \ten \A_i)$
vanishing near $(\ra M \times M) \cup (M \times \ra M)$ such that in
$B(L^2(M,E \ten \A))$ 
$$\|\Kappa- S(D(\rho)^2-\lambda)\| \le \tfrac 12.$$
This choice 
is possible
since by assumption $(D(\rho)^2-\lambda)$ has a bounded inverse, hence also
$(D(\rho)^2-\ov{\lambda})$ has a bounded inverse on $L^2(M,E \ten \A)$.

It follows that $$(Q(\lambda)+S)(D(\rho)^2-\lambda)=1 -\bigl(\Kappa- S(D(\rho)^2-\lambda)\bigr)$$ has
a bounded inverse on $L^2(M,E \ten \A)$.

Prop. \ref{intker} implies that $1-\Kappa-S(D(\rho)^2-\lambda)$ is invertible on $L^2(M,E \ten
\Ol{\mu})$ as well. Thus
$$\bigl(1 -(\Kappa- S(D(\rho)^2-\lambda))\bigr)^{-1}(Q(\lambda)+S)$$
is a bounded operator on $L^2(M,E \ten \Ol{\mu})$, which is a right inverse for $(D(\rho)^2-\lambda)$. Hence $D(\rho)^2-\lambda$ is injective and bounded below. It remains to show that its range is dense. This follows from the fact that $D(\rho)^2-\lambda$ is invertible on $L^2(M, E \ten \A)$.

Since $Q(\lambda)$ acts continuously on ${\mathcal S}_R(M,E \ten \Ol{\mu})$ and $\Kappa$ maps $L^2(M, E \ten
\Ol{\mu})$ continuously to ${\mathcal S}_R(M,E \ten \Ol{\mu})$ by the previous
lemma, the operator 
$(D(\rho)^2-\lambda)^{-1}$ acts continuously on ${\mathcal S}_R(M,E \ten \Ol{\mu})$  by
\begin{eqnarray*}
(D(\rho)^2-\lambda)^{-1} &=& (1 -\Kappa)(D(\rho)^2-\lambda)^{-1} +
\Kappa(D(\rho)^2-\lambda)^{-1} \\
&=& Q(\lambda) + \Kappa(D(\rho)^2-\lambda)^{-1} \ .
\end{eqnarray*}
\end{proof}

\begin{prop}
\label{smooth}
Let $\rho \in \bbbr$. Let $\lambda \in \bbbc$ with $\re \lambda
 < \omega$ be
 such that $(D(\rho)^2-\lambda)$ has a bounded inverse on $L^2(M,E \ten
 \Ol{\mu})$.

Then for  $\nu \in \bbbn$, $\nu \ge 2$, the operator $(D(\rho)^2-\lambda)^{-\nu}$ maps
 $L^2(M,E \ten \Ol{\mu})$ continuously to $C^{2\nu-3}_R(M,E \ten \Ol{\mu})$. 
\end{prop}

\begin{proof}
Let $Q(\lambda)(D(\rho)^2-\lambda)^{\nu}=1-\Kappa$ as before, thus
$$(D(\rho)^2-\lambda)^{-\nu} =Q(\lambda)+
\Kappa(D(\rho)^2-\lambda)^{-\nu}\ . $$
By Prop. \ref{embed} and Lemma \ref{sobol} the operator $Q(\lambda)$
maps $L^2(M,E \ten \Ol{\mu})$ continuously to $C_R^{2\nu-3}(M,E \ten
\Ol{\mu})$. Furthermore $\Kappa$ is smoothing. 
\end{proof}

\begin{cor} 
\label{kerDschw}
The kernel of $D(\rho)^2$ on $L^2(M, E \ten \Ol{\mu})$ is a subspace of ${\mathcal S}_R(M,E \ten \Ol{\mu})$.
\end{cor}

\begin{proof}
Let $\lambda\neq 0$ be as in the previous proposition.

Then  $(D(\rho)^2-\lambda)^{-\nu}f=(-\lambda)^{-\nu}f$ for $f \in \Ker
D(\rho)^2$ and  every $\nu \in
\bbbn$. By the previous proposition it
 follows that the elements of $\Ker D(\rho)^2$ are
 smooth. 

For $k \in \bbbz/6$ and $f \in \Ker D(\rho)^2$ 
$$D_{Z_k}^2\phi_k
 f \in \C_{cR}(Z_k,(\Ol{\mu})^{4d}) \ .$$ From Cor. \ref{DinvZ} it follows that
$$\phi_k
 f=D_{Z_k}^{-2}(D_{Z_k}^2\phi_k
 f) \in {\mathcal S}_R(Z_k, (\Ol{\mu})^{4d}) \ .$$
Hence $f \in {\mathcal S}_R(M,E \ten \Ol{\mu})$.
\end{proof}

\begin{prop}
\label{Philsch}
Let $\rho \neq 0$.

Let $P$ be the projection onto the kernel of $D(\rho)$ on $L^2(M,E \ten \A)$.

Then  $P$ is a
finite Hilbert-Schmidt operator whose integral
kernel is of the form $\sum_{j=1}^m f_j(x)h_j(y)^*$ with $f_j,h_j \in \Ker D(\rho)\cap  {\mathcal S}_R(M,E \ten \Ai)$.

Furthermore on $L^2(M,E \ten \Ol{\mu})$ we have that $\Ker D(\rho)^2=PL^2(M,E \ten \Ol{\mu})$ and $\Ran D(\rho)^2=(1-P)L^2(M,E \ten \Ol{\mu})$. Hence there is a decomposition $$L^2(M,E \ten \Ol{\mu})=\Ker D(\rho)^2
\oplus \Ran D(\rho)^2$$
with respect to which
$$D(\rho)^2= 0 \oplus D(\rho)^2|_{\Ran D(\rho)^2} \ .$$
Moreover $D(\rho)^2|_{\Ran D(\rho)^2}$ is invertible.
\end{prop}

\begin{proof}
First consider the situation on $L^2(M,E \ten \A)$:
Since the range of $D(\rho)$ is closed, there is an orthogonal projection $P$ onto the
kernel of $D(\rho)$ by Prop. \ref{cloran}.
Furthermore $D(\rho)$ is selfadjoint, hence $\Ker D(\rho)=\Ker D(\rho)^2$. The range of
$D(\rho)^2$ is closed, thus zero is an isolated point in the spectrum of
$D(\rho)^2$ on $L^2(M,E \ten \A)$.

Hence, for $r$ small enough,
$$P= \frac{1}{2\pi i}\int_{|\lambda|=r} (D(\rho)^2-\lambda)^{-1} d\lambda \ .$$

From Prop. \ref{spec} it follows that zero is an isolated point in the spectrum of
$D(\rho)^2$ on $L^2(M,E \ten
\Ol{\mu})$ as well. Thus $P$ is well-defined as a bounded operator on
$L^2(M,E \ten \Ol{\mu})$.  By 
Prop. \ref{projker} it is a Hilbert-Schmidt operator whose integral kernel is as asserted.  

The remaining parts follow from the spectral theory for closed operators on
Banach spaces (\cite{da}, Th. 2.14).
\end{proof}

\begin{cor}
\label{DMinv}
Let $\re \lambda < \omega$.
\begin{enumerate}
\item Let $\rho \neq 0$ and let $P$ be the
orthogonal projection onto the kernel of $D(\rho)^2$.
If $D(\rho)^2+P-\lambda$ has a bounded inverse on $L^2(M, E \ten \A)$, then
$D(\rho)^2+P-\lambda$ has a bounded inverse on
$L^2(M,E \ten \Ol{\mu})$ and the inverse acts as a bounded operator on the space ${\mathcal S}_R(M, E \ten
\Ol{\mu})$ as well.

\item Let $P_0$ be the orthogonal projection onto $\Ker D_s^2$.
If $D_s^2 +P_0-\lambda$ has a bounded inverse on $L^2(M, E \ten \A)$, then
$D_s^2+P_0-\lambda$ has a bounded inverse on
$L^2(M,E \ten \Ol{\mu})$ and the inverse acts as a bounded operator on the space
${\mathcal S}_R(M, E \ten \Ol{\mu})$ as well.
\end{enumerate}

In particular there is $c>0$ such that $\{\re \lambda <c\}$ is in the resolvent
set of $D(\rho)^2+P$ resp. $D_s^2+P_0$.
\end{cor}

\begin{proof}
(1) From the previous proposition it follows that
$P(1-\lambda)^{-1}+(1-P)(D(\rho)^2-\lambda)^{-1}$ inverts
$D(\rho)^2+P-\lambda$ on $L^2(M,E \ten
\Ol{\mu})$. Since $P$ acts as a bounded operator on the space
${\mathcal S}_R(M, E \ten \Ol{\mu})$ by Cor. \ref{kerDschw} and $(D(\rho)^2-\lambda)^{-1}$
is bounded on ${\mathcal S}_R(M, E \ten \Ol{\mu})$  by Prop. \ref{spec}, the
operator $(D(\rho)^2+P-\lambda)^{-1}$ is bounded on ${\mathcal S}_R(M, E \ten
\Ol{\mu})$ as well. 

(2) follows analogously. 
\end{proof} 

\subsection{An approximation of the semigroup}
\label{appr}
By cutting and  pasting we construct a family of integral operators that 
behaves similar to a semigroup generated by $-D(\rho)^2$ for small times.
 
We work in the setting fixed in \S \ref{hsDrho}.

Let $e(\rho)^{\cp}_t(x,y)$ be the restriction of the integral kernel of
$e^{-t(D_N+\rho K)^2}$  to $\U_{\cp} \times \U_{\cp}$ 
and  for $k \in \bbbz/6$ let $e(\rho)^k_t(x,y)$ be the restriction
of the integral
kernel of $e^{-tD_{Z_k}^2}$ to $\U_k \times \U_k$. Extend these functions by
zero to $M \times M$. Clearly for $k \in \bbbz/6$
we have $e(\rho)^k_t(x,y)=e(0)^k_t(x,y)$.

We write $E(\rho)_t$ for the family of integral operators on $L^2(M,E \ten
\Ol{\mu})$   corresponding to the integral kernel $$e(\rho)_t(x,y):=\sum_{k \in J} \gamma_k(x) e(\rho)^k_t(x,y) \phi_k(y)$$  and set $E(\rho)_0:=1$.

Using the results in \S
\ref{cpman} and \S \ref{cyl} one deduces:

\begin{enumerate}
\item The family $E(\rho)_t$ is strongly continuous
in $t$ on $L^2(M,E \ten \Ol{\mu})$.
\item If $f \in \C_{cR}(M,E \ten
\Ol{\mu})$, then the map $[0,\infty) \to L^2(M,E \ten \Ol{\mu}),~ t \mapsto
  E(\rho)_t f$ is differentiable.
\item $\Ran E(\rho)_t \subset {\mathcal S}_R(M,E \ten \Ol{\mu})$  for $t>0$.
\item Let $A$ be a differential operator on $\C(M,E \ten \Ol{\mu})$ of order $m$ with bounded
coefficients. For $T>0$ there is $C>0$ such that for all $0<t<T$ on $L^2(M,E \ten \Ol{\mu})$
$$\| A E(\rho)_t \| \le C t^{-\frac m2} \ .$$
\item For any $m \in \bbbn_0$ and $T>0$ there is $C>0$ such that for
all $y \in \U_{\cp},~\rho \in [-1,1]$ and $0<t<T$ in $C^m(\U_{\cp},(E \ten E_y) \ten \A_i)$
$$\|e(\rho)^{\cp}_t(\cdot, y)-e(0)^{\cp}_t(\cdot,y)\|_{C^m} \le C t|\rho| \ .$$
An analogous estimate holds for the partial derivatives in $y$ with respect to unit
vector fields on $\U_{\cp}$.
\end{enumerate}
  
The last statement follows by Volterra development (Prop. \ref{boundpert}):
On $N$  
\begin{eqnarray*}
\lefteqn{e^{-tD_N(\rho)^2}-e^{-tD_N^2}}\\
&=&\rho t\sum\limits_{n=1}^{\infty}(-1)^n(\rho t)^{n-1} \int_{\Delta^n}e^{-u_0tD_N^2}
\bigl([D_N , K]_s+ \rho K^2\bigr) e^{-u_1tD_N^2} \dots \\
&& \qquad \dots e^{-u_ntD_N^2} ~du_0 \dots
du_n  \ ,
\end{eqnarray*}
and the sum is an integral operator whose integral kernel is uniformly bounded
in $0<t<T$ and $\rho \in [-1,1]$.

\subsection{The semigroup $e^{-tD_s^2}$}
\label{Mcp}

By Cor. \ref{smoothcpker}  the operator $e^{-tD_s^2}$ on the Hilbert
space $L^2(M,E)$ is
an integral operator with smooth integral kernel $k_t$ for $t>0$. In this section we show that $k_t$ defines a
strongly continuous semigroup on $L^2(M,E \ten \Ol{\mu})$ and that this
 semigroup extends to a holomorphic one.

For $D_s$ define $e^k_t(x,y)$ analogously to $e(\rho)^k_t(x,y)$
for $D(\rho)$ in the previous section and let 
$$e_t(x,y):=\sum_{k \in J} \gamma_k(x) e^k_t(x,y) \phi_k(y) \ .$$

The  corresponding family of
operators is denoted by $E_t$. We set $E_0:=1$. The properties of $E_t$ are
as described in the
previous section.

For $f \in \C_{cR}(M,E)$
$$e^{-tD_s^2}f-E_tf=-\int_0^te^{-sD_s^2}\left(\frac{d}{dt}+D_s^2\right)E_{t-s}f~ds
$$
by Duhamel's principle.

\begin{prop}
\label{complsemgr}
\begin{enumerate}
\item The heat kernel $k_t$ associated to $D_s$ defines a strongly continuous semigroup on $L^2(M,E \ten
\Ol{\mu})$ with generator $-D_s^2$, which extends to a bounded holomorphic semigroup.
\item Let $A$ be a differential operator of order $m \in \bbbn_0$  with bounded smooth
coefficients. Then for any $t>0$ the operator $Ae^{-tD_s^2}$ is bounded on
$L^2(M,E \ten \Ol{\mu})$ and for any
$T>0$ there is $C>0$ such that for $0<t < T$ 
$$\|Ae^{-tD_s^2}\| \le C t^{-\frac m2} \ .$$
\end{enumerate}
\end{prop}

\begin{proof}
First we show that for $T>0$ there is $C>0$ such that  for $0<t<T$ the difference $A_xk_t(x,y)-A_xe_t(x,y)$ is bounded by $Ct$ in $L^2(M \times
M,E \boxtimes E^*)$.

For $k \in J$ let $\chi_k \in \C_c(M)$ be a function with values in $[0,1]$, with compact
support in $\U_k$  and equal to one on a neighborhood of $\supp d \gamma_k$.

From Duhamel's principle it follows that
\begin{eqnarray*}
\lefteqn{A_xk_t(x,y)-A_xe_t(x,y)}\\
&=&-\sum_{k \in J}  \int_0^t \int_M  A_xk_s(x,r)[(\dirac)_r,c(d
\gamma_k(r))]_se^k_{t-s}(r,y) \phi_k(y) ~dr ds \ .
 \end{eqnarray*}
This can be re-written as
\begin{eqnarray*}
&&- \sum_{k \in J}  \int_0^t \int_M (1- \chi_k(x)) A_x k_s(x,r)  [(\dirac)_r ,c(d \gamma_k(r))]_se^k_{t-s}(r,y) \phi_k(y)~ dr ds
\\
&&- \sum_{k \in J}  \int_0^t \int_M \chi_k(x)A_x(k_s(x,r)-e^k_s(x,r))
[(\dirac)_r ,c(d \gamma_k(r))]_se^k_{t-s}(r,y) \phi_k(y) ~dr ds\\
&&- \sum_{k \in J}  \int_0^t \int_M \chi_k(x)A_x e^k_s(x,r) [(\dirac)_r ,c(d
\gamma_k(r))]_se^k_{t-s}(r,y) \phi_k(y) ~dr ds \ .
\end{eqnarray*}
Using Lemma \ref{estpun} and Lemma \ref{estker} we estimate the norms of the three terms in $E_x \ten E_y$:

Since $\supp d\gamma_k \cap \supp \phi_k=\emptyset$ and $\supp d\gamma_k \cap
\supp(1-\chi_k)=\emptyset$,  there is
$C>0$ such  that for $x,y \in M$ and $t > 0$ the norm of the first term is bounded by
$$C \sum_{k \in J} t(1- \chi_k(x)) e^{-\frac{d(x,\supp d
\gamma_k)^2}{(4+ \delta)t}} e^{-\frac{d(y,\supp d \gamma_k)^2}{(4+\delta)t}}
1_{\supp \phi_k}(y) \ ,$$
and such that for $x,y \in M$ and  $t>0$
 the norms of the second and third term are
bounded  by
$$Ct\chi_k(x)e^{-\frac{d(y,\supp d
\gamma_k)^2}{(4+\delta)t}} 1_{\supp \phi_k}(y) \ .$$
When estimating the third term we used that the action of
the integral kernel $e^k_s(x,r)\chi_k(r)$  is uniformly bounded for $k=\cp$  on $C^n(\U_{\cp},E \ten
E_y)$  by Prop. \ref{cpmandiff} and for $k
\in \bbbz/6$ on $C^n_R(\U_k,E
\ten E_y)$ by Prop. \ref{sgZ}.

Analogous estimates hold for the derivatives in $y$ with respect to unit vector
fields on $M$.

Hence $A_xk_t(x,y)-A_xe_t(x,y)$ is bounded   by
$Ct$ in $L^2(M \times
M,E \boxtimes E^*)$ for $0<t<T$
and some $C>0$. By Cor. \ref{L2ker} the corresponding family of operators on
$L^2(M,E \ten \Ol{\mu})$ is bounded by $Ct$ for $0<t<T$, hence $A_xk_t(x,y)$ defines a family of
bounded operators on $L^2(M,E \ten \Ol{\mu})$.

Write $S(t)$ for the integral operator induced by the integral kernel
$k_t(x,y)$.
 
By property (4) in \S \ref{appr} there is $C>0$ such that $\| A E_t \| \le C t^{-\frac m2}$  on $L^2(M,E \ten \Ol{\mu})$
for $0<t<T$ and some $C>0$, hence  
$$(*) \qquad \|AS(t)\| \le C t^{-\frac m2} \ . $$
The fact that $S(t)$ extends to a bounded holomorphic semigroup on $L^2(M,E \ten \Ol{\mu})$ is an almost
immediate consequence of $(*)$:

Since $E_t$ converges strongly to the identity on $L^2(M,E \ten \Ol{\mu})$ for
$t \to 0$, the operator $S(t)$
also does. Furthermore the kernels $k_t$ obey the semigroup law, hence $S(t)$ is a strongly continuous semigroup on $L^2(M,E \ten \Ol{\mu})$.

Note that the range of $S(t)-E_t$ is a subset of ${\mathcal S}_R(M,E \ten
 \Ol{\mu})$. Hence ${\mathcal S}_R(M,E \ten \Ol{\mu})$ is invariant under the action
 of $S(t)$. It
 follows that
$-D_s^2$ is the generator of $S(t)$.

From Prop. \ref{critsg} and the estimate $(*)$ applied to $A=D_s^2$ it follows   that the
 semigroup $S(t)=e^{-tD_s^2}$ extends to a
 holomorphic semigroup  on $L^2(M,E \ten \Ol{\mu})$.

In order to show that the holomorphic semigroup is bounded, let $P_0$ be the projection onto the kernel of $D_s^2$.

By
Cor. \ref{DMinv} there is $c>0$ such that $\{\re \lambda < c\}$ is in the
resolvent set of $D_s^2 + P_0$ on $L^2(M,E \ten \Ol{\mu})$. Hence by Prop. \ref{sec} the holomorphic
semigroup $e^{-t(D_s^2+P_0)}$ is
bounded. Thus $$e^{-tD_s^2}=e^{-t(D_s^2+P_0)}(1-P_0) + P_0$$ is bounded as well.  
\end{proof}

Recall that in \S \ref{hsDrho} we fixed the domain of $D_s^2$ as an unbounded
operator on $L^2(M,E \ten \Ol{\mu})$, but not the domain of
$D_s$. This is done now:

Let $D_s$ be the closure on $L^2(M,E \ten \Ol{\mu})$ of the Dirac operator
$\dirac$ with domain ${\mathcal S}_R(M, E \ten \Ol{\mu})$.

\begin{cor}
\label{Dsboundinv}

Let $P_0$ be the projection onto the kernel of $D_s^2$.

Let $\lambda \in \bbbc$ with $\re \lambda^2 < 0$.

Then the operators $D_s+P_0$ and $D_s - \lambda$  have a bounded inverse on $L^2(M,E \ten \Ol{\mu})$.
\end{cor}

\begin{proof}
By Cor.  \ref{DMinv} there is $c >0$ such that $\{\re \lambda \le c\}$ is in the
resolvent set of $D_s^2 + P_0$ on $L^2(M,E \ten \Ol{\mu})$. 

By the previous proposition and 
Prop. \ref{sec} it follows that there is $C>0$ such that for all $t>0$  on $L^2(M,E \ten \Ol{\mu})$
$$\|(D_s+P_0)e^{-t(D_s^2+P_0)}\| \le Ct^{-\frac 12} e^{-c t} \ .$$
Thus 
$$G:=\int_0^{\infty} (D_s+P_0)e^{-t(D_s^2+P_0)} dt$$ 
is a bounded operator on $L^2(M,E \ten \Ol{\mu})$. Furthermore
$G=(D_s+P_0)^{-1}$ on the Hilbert space $L^2(M,E)$. From Cor. \ref{DMinv} it follows that $G$ acts as a bounded operator on ${\mathcal S}_R(M, E)$. 
Hence $G$ inverts $D_s+P_0$ on ${\mathcal S}_R(M, E) \odot \Ol{\mu}$, thus 
$G$ inverts $D_s+P_0$ on $L^2(M,E \ten \Ol{\mu})$.

The proof of the fact that
$$ \int_0^{\infty}(D_s+\lambda)e^{-t(D_s^2-\lambda^2)}$$
inverts $(D_s-\lambda)$ for $\re \lambda^2<0$ is analogous.
\end{proof}

This answers a question from the beginning of \S \ref{hsDrho}:  the operator $D_s^2$ is indeed the square of $D_s$.

\subsection{The semigroup $e^{-tD(\rho)^2}$}
\label{semMgen}
This section is devoted to the study of the holomorphic semigroup generated by
$-D(\rho)^2$.  

First we prove its existence.

Define $D(\rho)$ on $L^2(M,E \ten \Ol{\mu})$ as the closure of the
operator $\dirac +
\rho K$ with domain ${\mathcal S}_R(M, E \ten \Ol{\mu})$.

Let $P_0$ be the orthogonal projection onto the kernel of $D_s^2$. Then
$D(\rho)$ is a bounded perturbation of $W^*(D_s+P_0)W$ with $W$ as in \S \ref{stand}. 
By the results of the previous section we can apply
 Prop. \ref{sgpert}  and conclude
 that $D(\rho)^2$, as defined in \S \ref{hsDrho}, is the square of $D(\rho)$ and
 furthermore that $-D(\rho)^2$ generates a holomorphic semigroup. This shows the first assertion of the following proposition.

For $\rho \neq 0$ let $P$ be the orthogonal projection onto the kernel of
$D(\rho)^2$.

\begin{prop}
\label{Mgensg}
\begin{enumerate}
\item Let $\rho \in \bbbr$.
The operator $-D(\rho)^2$ generates a holomorphic semigroup
$e^{-tD(\rho)^2}$ on $L^2(M,E \ten \Ol{\mu})$.  For $\rho \neq 0$ the semigroup is bounded holomorphic.
\item For $t \ge 0$ the operator $e^{-tD(\rho)^2}$ depends analytically on
$\rho$. For every $T>0$ there is $C>0$ such that 
$$\|e^{-tD(\rho)^2}\| \le C $$
for all $\rho \in [-1,1]$ and $0
\le t \le
T$.

\item For $\rho \neq 0$ there are $C,\omega>0$ such that for all $t\ge 0$ 
$$\|(1-P)e^{-tD(\rho)^2}\| \le Ce^{-\omega t} \ .$$
\end{enumerate}
\end{prop}

\begin{proof}
(1) Let $\rho \neq 0$. 
By Cor. \ref{DMinv} there is $c>0$ such that $\{\re \lambda
\le c\}$ is in the resolvent set of $D(\rho)^2 + P$ on $L^2(M,E \ten
\Ol{\mu})$, hence by Lemma \ref{sec} the holomorphic semigroup
$e^{-t(D(\rho)^2+P)}$ is
bounded by $Ce^{-c t}$ for some $C>0$, thus for $T>0$ there is $C>0$ such
that for $t>T$ 
$$\|D(\rho)^2e^{-tD(\rho)^2}\|=\|D(\rho)^2 e^{-t(D(\rho)^2+P)}\| \le Ce^{-c
t} \ .$$
Now Prop. \ref{critsg} implies that the semigroup $e^{-tD(\rho)^2}$
is bounded holomorphic. 

(2) follows from Prop. \ref{boundpert}.

(3) follows from $(1-P)e^{-tD(\rho)^2}=(1-P)e^{-t(D(\rho)^2+P)}$. 
\end{proof}

The following proposition shows that $e^{-tD(\rho)^2}$ is smoothing.

\begin{prop}
\label{Mgensmooth}

\begin{enumerate}
\item Let $n \in \bbbn_0$. For every $\rho \in \bbbr$ and every $t>0$ the operator $e^{-tD(\rho)^2}$ maps $L^2(M,E \ten \Ol{\mu})$
continuously to $C_R^n(M,E \ten \Ol{\mu})$. 
\item Let $n \in 2\bbbn,~ n \ge 4.$ For every $\rho\neq 0$ the family $e^{-tD(\rho)^2}:C_{cR}^{n}(M,E \ten \Ol{\mu})
\to C_R^{n-3}(M,E \ten \Ol{\mu})$ is uniformly bounded.
\item Let $n \in 2\bbbn,~n \ge 4$. For every $T>0$ the family $e^{-tD(\rho)^2}:C_{cR}^n(M,E \ten \Ol{\mu})
\to C_R^{n-3}(M,E \ten \Ol{\mu})$ is uniformly bounded in $0\le t<T$ and in $\rho
\in [-1,1]$.
\end{enumerate}
\end{prop}

\begin{proof}

We conclude (1) from
$$e^{-tD(\rho)^2}=(D(\rho)^2+1)^{-k}(D(\rho)^2+1)^k e^{-tD(\rho)^2} \ ,$$
taking into account that $(D(\rho)^2+1)^k e^{-tD(\rho)^2}$ is bounded on $L^2(M,E \ten
\Ol{\mu})$ for $t >0$ and that $(D(\rho)^2+1)^{-k}$ maps $L^2(M,E \ten \Ol{\mu})$ continuously
to $C_R^{2k-3}(M,E \ten \Ol{\mu})$ for $k \in \bbbn,~k \ge 2,$ by Prop. \ref{smooth}.

 (2) and (3) follow  from
$$e^{-tD(\rho)^2}=(D(\rho)^2+1)^{-k}
e^{-tD(\rho)^2}(D(\rho)^2+1)^k $$ by Prop. \ref{smooth}.
\end{proof}

\subsection{The integral kernel}

In order to  prove that the operator $e^{-tD(\rho)^2}$ is an integral
operator we use the same method as in \S \ref{hksemI}: Via Duhamel's principle we
compare $e^{-tD(\rho)^2}$ and the approximation $E(\rho)_t$, which was defined in \S \ref{appr}.
 
In the following $| \cdot|$ denotes the norm on the fibers of $(E \boxtimes E^*)\ten \A_i$.

\begin{prop}
\label{gaussM}

For every $\rho \in \bbbr$  and every $t>0$ the operator $e^{-tD(\rho)^2}$ is an integral operator with
smooth integral kernel. Let $k(\rho)_t(x,y)$ be its integral kernel. 

\begin{enumerate}
\item The map $(0,\infty) \to
\C(M \times M, (E \boxtimes E^*) \ten \Ai), ~t \mapsto k(\rho)_t$ is
smooth. 
\item $k(\rho)_t(x,y)=k(\rho)_t(y,x)^*$.
\item 
For every $T>0$ there are $c,C>0$ such that 
$$|k(\rho)_t(x,y)- e(\rho)_t(x,y)| \le Ct\bigl(|\rho| 1_{\U_{\cp}}(y) + \sum\limits_{k \in J}
e^{-\frac{d(y,\supp d\gamma_k)^2}{ct}}1_{\supp \phi_k}(y) \bigr)$$
for all $0<t<T$, $\rho \in [-1,1]$ and $x,y \in M$.
\item Let $\rho \neq 0$. There are $c,C>0$ such that 
$$|k(\rho)_t(x,y)- e(\rho)_t(x,y)| \le Ct  e^{-\frac{d(y,\U_{\cp})^2}{ct}}$$
for all $t>0$ and all $x,y \in M$.
\end{enumerate}

Statements analogous to (3) and (4) hold for the  partial derivatives in $x$ and
$y$ with respect to unit vector fields on $M$. 
\end{prop}

\begin{proof}
For (1) it is enough to prove an analogous statement for the operator
$e^{-tD(\rho)^2}- E(\rho)_t$.

Let $f \in \C_{cR}(M,E \ten \A_i)$. Then by Duhamel's principle
\begin{eqnarray*}
\lefteqn{e^{-tD(\rho)^2}f- E(\rho)_t f }\\
&=& -\sum_{k \in J}  \int_0^t \int_M
e^{-sD(\rho)^2}\bigl(\frac{d}{dt}+D(\rho)^2\bigr)\gamma_ke(\rho)^k_{t-s}(\cdot ,
y) \phi_k(y) f(y)~ dy ds \ .
\end{eqnarray*}
We write
\begin{eqnarray*}
\lefteqn{ \sum_{k \in J} \bigl(\frac{d}{d\tau}+D(\rho)^2\bigr)\gamma_ke(\rho)^k_{\tau}(\cdot , y) \phi_k(y)} \\
&=&- \sum\limits_{k \in J} [\dirac,c(d \gamma_k)]_se(0)^k_{\tau}(\cdot , y)
\phi_k(y) \\
&&-  
[\dirac,c(d \gamma_{\cp})]_s\bigl(e(\rho)^{\cp}_{\tau}(\cdot,y)-e(0)^{\cp}_{\tau}(\cdot,y)\bigr)
\phi_{\cp}(y)  \ .
\end{eqnarray*}
Note that the map 
$$(0,\infty) \to \C(M,\C_{cR}(M,E\ten \Ai) \ten E) \ ,$$
$$\tau \mapsto \bigl(y
\mapsto \sum_{k \in J} \bigl(\frac{d}{d\tau}+D(\rho)^2\bigr)e(\rho)^k_{\tau}(\cdot, y) \phi_k(y)\bigr)$$ is
smooth. We show that it extends smoothly by zero to $\tau=0$. From this and Prop. \ref{Mgensmooth} it follows that $e^{-tD(\rho)^2}- E(\rho)_t$ is an integral operator.
 
We estimate the terms on the right hand side of the previous equation.

From the estimates in Lemma \ref{estpun} and Lemma \ref{gaussZ} it follows that for any $m
\in \bbbn_0$ there are $C,c>0$ such that 
$$\|[\dirac,c(d \gamma_k)]_s e(0)^k_{\tau}(\cdot , y) \phi_k(y)\|_{C^m_{R}} \le C e^{-\frac{d(y,\supp
d\gamma_k)^2}{c\tau}}1_{\supp \phi_k}(y) $$ 
 in $C_R^m(M,E \ten
\A_i)\ten E^*_y$ for all $k \in J,~0 <\tau<T,~ \rho \in [-1,1]$ and $y \in M$.

Furthermore by \S \ref{appr}, property (5), 
 there is $C>0$ such that
$$\|[\dirac,c(d \gamma_{\cp})]_s\bigl(e(\rho)^{\cp}_{\tau}(\cdot,y)-e(0)^{\cp}_{\tau}(\cdot,y)\bigr)
\phi_{\cp}(y)\|_{C^m_R} \le C \tau |\rho| 1_{\U_{\cp}}(y) $$
 in $C_{R}^m(M,E\ten \A_i)\ten E_y$ for $0<\tau <T,~\rho \in
[-1,1]$ and $y \in M$. 

Analogous estimates hold for the derivatives in $y$ and also in $\tau$ since by the
heat equation the derivatives with respect to $\tau$ can be expressed in terms of the
derivatives with respect to $x$.
This shows (1).

The property
$k(\rho)_t(x,y)=k(\rho)_t(y,x)^*$ follows from the selfadjointness of $e^{-tD(\rho)^2}$. 

Statement (2) and (3) follow from the estimates using Prop. \ref{Mgensmooth}. For
the proof of (3) we also take into account that the kernel
$e(\rho)^{\cp}_t$ and its derivatives are uniformly bounded
in $t$ with $t>T$. 
\end{proof}

\begin{cor}
\label{boundsmt}
Let $\nu \in \bbbn_0$. For every $\rho \in \bbbr$ and $T>0$ there is $C>0$ such that
$$\|D(\rho)^{\nu}e^{-tD(\rho)^2}\| \le C t^{-\frac{\nu}{2}} $$
 for $0 <t<T$  on $L^2(M,E \ten \Ol{\mu})$.
\end{cor}

\begin{proof}
By Duhamel's principle, for $f \in \C_{Rc}(M,E\ten \Ol{\mu})$,
\begin{eqnarray*}
\lefteqn{D(\rho)^{\nu}e^{-tD(\rho)^2}f- D(\rho)^{\nu}E(\rho)_t f}\\
&=& -\sum_{k \in J}  \int_0^t \int_M
e^{-sD(\rho)^2}D(\rho)^{\nu}[D(\rho)^2,\gamma_k]_se(\rho)^k_{t-s}(\cdot , y) \phi_k(y) f(y)
~dy ds \ .
\end{eqnarray*}
There is $C>0$ such that this term is bounded in $L^2(M,E \ten \Ol{\mu})$ by $Ct$ for $0<t<T$. Furthermore by Prop. \ref{cpmandiff} and Prop. \ref{sgZ} there is $C>0$ such that on
$L^2(M,E \ten \Ol{\mu})$  $$\| D(\rho)^{\nu}E(\rho)_t\| \le C
t^{-\frac{\nu}{2}} \ .$$
The assertion follows.
\end{proof}

\begin{cor}
\label{semCm}

For every $\rho \neq 0$ and $m \in \bbbn$ the family of integral kernels
$k(\rho)_t$ defines a strongly continuous semigroup on $C^m_R(M,E
\ten \Ol{\mu})$ bounded by $C(1+t)^{\frac 32}$ for
 some $C >0$ and all $t>0$.

It is denoted by $e^{-tD(\rho)^2}$ as well.
\end{cor}

\begin{proof}
By the estimates in the proposition the integral kernel $k(\rho)_t- e(\rho)_t$
defines an operator on $C^m_R(M,E \ten
\Ol{\mu})$  bounded by $Ct^{3/2}$ for every $t>0$
and some $C >0$. 

We show that $E(\rho)_t$ is a strongly continuous uniformly bounded family of operators on $C^m_R(M,E \ten
\Ol{\mu})$.

For $k \in \bbbz/6$ the action of the integral
kernel $e(\rho)^k_t$ on $C^m_R(\U_k,E \ten \Ol{\mu})$ is strongly continuous and uniformly bounded in $t$ by Prop. \ref{sgZ}. 

By Prop. \ref{cpmanhsg} and  Prop. \ref{boundpert} the family $e^{-t(D_N+\rho K)^2}$
is a strongly continuous semigroup on $C^m(N,E_N \ten \Ol{\mu})$. It is bounded
since its integral kernel is uniformly bounded for $t>1$. Hence the action of $e(\rho)^{\cp}_t$ on
$C^m(\U_{\cp},E \ten \Ol{\mu})$ is strongly continuous and
 uniformly bounded.

Since $E(\rho)_t$  converges strongly to the
identity on $C^m_R(\U_k,E \ten \Ol{\mu})$ for $t \to 0$, so does $e^{-tD(\rho)^2}$. It clearly
satisfies the semigroup property.
\end{proof}

\begin{cor}
\label{gaussMdec}

Let $\rho \neq 0$ and $n \in \bbbn$.

For every $\ve>0$
we can find $C,c>0$ such that for 
$x,y \in M$ with $d(x,y)> \ve$ and $t>0$ 
$$|k(\rho)_t(x,y)| \le C  \bigl(e^{-\frac{d(x,y)^2}{ct}}+
te^{-\frac{d(y,\U_{\cp})^2}{ct}}\bigr) \ .$$
An analogous statement holds for the derivatives in $x$ and $y$ with respect to unit vector
fields on $M$.
\end{cor}

We refer to \S \ref{HS} for the notion of a Hilbert-Schmidt operator and the Hilbert-Schmidt
norm $\| \cdot \|_{HS}$ used in the next corollary.

\begin{cor}
\label{hkHS}
Let $\rho\neq 0$ and $\nu \in \bbbn_0$.

The operators
$1_{M_r}D(\rho)^{\nu}e^{-tD(\rho)^2}$ and $D(\rho)^{\nu}e^{-tD(\rho)^2}1_{M_r}$ with $r\ge 0$
are Hilbert-Schmidt operators.
\begin{enumerate}
\item For every  $T>0$ there is $C>0$ such that for any  $r\ge 0$ and $t>T$ 
$$\|1_{M_r}D(\rho)^{\nu}e^{-tD(\rho)^2}\|_{HS} \le C(1+r)$$
and $$\|D(\rho)^{\nu}e^{-tD(\rho)^2}1_{M_r}\|_{HS} \le C(1+r) \ .$$
\item  For every $\ve>0$ there is $C>0$ such that for every  $r,t>0$  
$$\|1_{M_r}D(\rho)^{\nu}e^{-tD(\rho)^2}(1-1_{M_{r+\ve}})\|_{HS} \le C(1+r)t^{1/2} $$
and $$\|(1-1_{M_{r+\ve}})D(\rho)^{\nu}e^{-tD(\rho)^2}1_{M_r}\|_{HS} \le C(1+r)t^{1/2}
\ . $$
\end{enumerate}
\end{cor}

\begin{proof}
(1) Since by Prop. \ref{Mgensg} the semigroup $e^{-tD(\rho)^2}$ is bounded, it follows
from Prop. \ref{HScomp} that there is $C>0$ such that for all $r>0$ and $t>T$
$$\|1_{M_r}D(\rho)^{\nu}e^{-tD(\rho)^2}\|_{HS}\le
C\|1_{M_r}D(\rho)^{\nu}e^{-TD(\rho)^2}\|_{HS} \ .$$
By the previous corollary there are $C,c>0$ such that 
$$|1_{M_r}(x)D(\rho)^{\nu}_xk(\rho)_T(x,y)| \le
C1_{M_{r}}(x)(e^{-cd(x,y)^2}+ e^{-cd(y,M_r)^2})$$ 
 for all $r>0$
  and  $x,y \in M$.
This yields the asserted estimate. The second estimate in (1) is proved analogously.

(2) follows from the previous corollary and (1).
\end{proof}

\chapter{Superconnections and the Index Theorem}

The notion of a superconnection 
 on a free finitely generated $\bbbz/2$-graded module, which we define now, generalizes the notion of a
 connection on a free module \cite{kar}.

In the family case  superconnections usually
act on infinite
dimensional bundles (\cite{bgv}, Ch. 9). In analogy the superconnections we  consider act
on modules with infinitely many generators. In that sense the following definition should
be merely seen as a motivation for the definitions of the superconnections to come.

\begin{ddd}
Let $\B$ be a locally $m$-convex Fr\'echet algebra and let $p,q \in \bbbn_0$.
 
Let $V:=(\bbbc^+)^p \oplus (\bbbc^-)^q$. Consider $V \ten \Oi\B$ as a
$\bbbz/2$-graded space.

A {\sc superconnection on $V \ten \B$} is an odd linear map
$$A:V \ten \Oi\B \to V \ten \Oi\B$$
satisfying Leibniz rule:

For $\alpha \in V^{\pm} \ten \Ok\B$ and $\beta \in \Oi\B$ 
$$A(\alpha \beta)= A( \alpha) \beta + (-1)^{\deg \alpha} \alpha \di\beta \ $$
where $\deg \alpha$ is the degree of $\alpha$ with respect to the
$\bbbz/2$-grading of $V \ten \Oi\B$.

The map $A^2$ is called the {\sc curvature of $A$}.  
\end{ddd}

As for a connection \cite{kar} the curvature of a superconnection is a right $\Oi\B$-module map.

\section{The superconnection $A^I_t$ associated to $D_I$}
\label{supconI}

\subsection{The family $e^{-(A^I_t)^2}$}

Let $C_1$ be the $\bbbz/2$-graded unital algebra generated by an odd element
$\sigma$ with $\sigma^2=1$. As a vector space $C_1$ is isomorphic to $\bbbc \oplus
\bbbc$ via the map $C_1\to \bbbc \oplus
\bbbc, a+b\sigma \mapsto (a,b)$. We endow $C_1$ with the scalar product induced by the standard hermitian scalar product on $\bbbc \oplus
\bbbc$.

We identify $L^2([0,1],C_1 \ten
(\Ol{\mu})^{2d})$ with $C_1 \ten L^2([0,1],
(\Ol{\mu})^{2d})$ and consider
$D_I$, as defined in \S \ref{Igen} with boundary conditions induced by a pair $(P_0,P_1)$,
as an unbounded operator on the $\bbbz/2$-graded $\A_i$-module $L^2([0,1],C_1 \ten \A_i^{2d})$.

Let $U \in \C([0,1],M_{2d}(\Ai))$ be as in Prop. \ref{defWI} with
$U(0)P_0U(0)^*=P_s$ and $U(1)P_1U(1)^*=1-P_s$. The map $U^*\di U$   
can
be seen as a flat superconnection on $L^2([0,1],C_1 \ten \A_i^{2d})$.
It preserves the space $C_1 \ten \C_R([0,1],(\Ol{\mu})^{2d})$.

We define
$$A_I:= U^*\di U + \sigma D_I$$ and for $t\ge 0$ 
$$A^I_t := U^*\di U + \sqrt t \sigma D_I \ .$$
Then $A_I$ is an odd map on $C_1 \ten \C_R([0,1], (\Ol{\mu})^{2d})$  fulfilling
Leibniz rule and is called a superconnection associated to $\sigma D_I$. The map $A^I_t$ is called the corresponding rescaled superconnection.

The curvature of $A_I$ is 
$$A_I^2=U^*\di^2U+U^*\di U\sigma D_I +\sigma D_I U^*\di U +D_I^2=D_I^2+\sigma [D_I,U^* \di U] \ .$$
From
\begin{eqnarray*}
 [D_I,U^* \di U] &=& -U^*[\di,UD_IU^*]U\\
&=&- U^*([\di,D_{I_s}]+[\di,UI_0(\ra U^*)])U\\
&=&- U^*\di(UI_0(\ra U^*))U =:R 
\end{eqnarray*}
it follows that
$A_I^2=D_I^2 + \sigma R$
with $R \in  \C([0,1],M_{2d}(\hat\Omega_1\Ai))$ vanishing near the boundary and fulfilling $R^*=-R$.
 
The curvature of the rescaled superconnection $A^I_t$ is 
$$(A^I_t)^2=tD_I^2 +\sqrt t \sigma R \ .$$
We see that the curvature and the rescaled curvature are right $\Ol{\mu}$-module maps.

Since $A_I^2$ is a bounded perturbation of $D_I^2$, it defines a holomorphic semigroup $e^{-tA_I^2}$ on $L^2([0,1],C_1 \ten(\Ol{\mu})^{2d})$.

In the following we restrict to $t\ge 0$:

By Volterra development 
\begin{eqnarray*}
e^{-tA_I^2}&=& \sum\limits_{n=0}^{\infty}(-1)^n t^n\int_{\Delta^n}e^{-u_0tD_I^2}\sigma
Re^{-u_1tD_I^2} \sigma R\dots e^{-u_ntD_I^2} ~du_0 \dots du_n\\
&=&\sum\limits_{n=0}^{\infty}\sigma^n (-1)^{\frac{(n+1)n}{2}}t^n I_n(t) 
\end{eqnarray*}
with
$$I_n(t):= \int_{\Delta^n}e^{-u_0tD_I^2}Re^{-u_1tD_I^2}R\dots e^{-u_ntD_I^2}
~du_0 \dots du_{n} \ .$$
Note that the series is finite on  $L^2([0,1],C_1 \ten (\Ol{\mu})^{2d})$.

It follows that
$$e^{-(A^I_t)^2}=\sum\limits_{n=0}^{\infty}\sigma^n (-1)^{\frac{(n+1)n}{2}}
t^{n/2}I_n(t) \ .$$
The operators $I_n(t)$ obey the following recursion relation for $n \ge 1$:
\begin{eqnarray*}
I_n(t)&=& \int_0^1 du_0 ~e^{-u_0tD_I^2}R \int_{(1-u_0)\Delta^{n-1}} e^{-u_1tD_I^2}R\dots R e^{-u_ntD_I^2} ~du_1 \dots du_{n}\\
&=& \int_0^1 du_0~ (1-u_0)^{n-1}e^{-u_0tD_I^2}R \int_{\Delta^{n-1}} e^{-(1-u_0)u_1tD_I^2}R\dots \\
&& \qquad \dots Re^{-(1-u_0)u_ntD_I^2} ~du_1 \dots du_{n}\\
&=& \int_0^1 du_0~ (1-u_0)^{n-1}e^{-u_0tD_I^2}R I_{n-1}((1-u_0)t) \ .
\end{eqnarray*}
Note that $e^{-(A^I_t)^2}$ is
selfadjoint on $L^2([0,1],C_1 \ten (\Ol{\mu})^{2d})$ in the sense of \S \ref{adop} and that $I_n(t)=(-1)^nI_n(t)^*$.
 
\subsection{The integral kernel of $e^{-(A^I_t)^2}$}

Since $e^{-tD_I^2}$ is a bounded semigroup on $C^m_R([0,1],(\Ol{\mu})^{2d})$
for every $m \in \bbbn_0$ by Prop. \ref{sgIgen}, the family $I_n(t):C_R^m([0,1],
\A_i^{2d})\to C_R^m([0,1], (\hat\Omega_n\A_i)^{2d})$ is uniformly bounded  in $t \ge 0$. 

In the following we write $| \cdot|$ for the norm on $M_{2d}(\Ol{\mu})$.

\begin{prop}
\label{supkerI}
For every $n \in \bbbn_0$ and $t>0$ the operator $I_n(t)$ is an integral
operator. Let 
 $p_t(x,y)^n$ be its integral kernel.

\begin{enumerate}
\item The map $$(0,\infty) \to \C([0,1], \C_R( [0,1], M_{2d}(\hat\Omega_n\Ai))),~ t \mapsto
\bigl(y \mapsto p_t(\cdot, y)^n\bigr) $$ is smooth.

\item  $p_t(x,y)^n=(-1)^n(p_t(y,x)^n)^*$.
 
\item For every $l,m,n \in \bbbn_0$ there are $C,\omega>0$ such that for
$t>0$ and $x,y \in [0,1]$ 
$$|\ra_x^l \ra_y^m p_t(x,y)^n | \le C (1+t^{-\frac{l+m+1}{2}})e^{-\omega t} \ .$$
\end{enumerate}
\end{prop}

\begin{proof}
 In degree $n=0$ the assertions hold by Prop. \ref{hkIgen} and
 Cor. \ref{hkIgenub}.

Let $k_t(x,y)$ be the integral kernel of $e^{-tD_I^2}$.

For  $f \in \C_R([0,1], (\Ol{\mu})^{2d})$  
$$(I_n(t)f)(x)=\int_0^1 \int_0^1 \int_0^1 (1-s)^{n-1}
k_{st}(x,r)R(r)p_{(1-s)t}(r,y)^{n-1}f(y)~dydr ds $$
by induction and by the recursion formula above.

We can
interchange the integration over $r$ and $y$.

For the proof of the existence of the integral kernel and of (1) it suffices to  show that the map 
$$(s,t) \mapsto \bigl(y \mapsto \int_0^1
(1-s)^{n-1}k_{st}(\cdot,r)R(r)p_{(1-s)t}(r,y)^{n-1}dr \bigr)$$ 
is a smooth map from $[0,1]
\times (0,\infty)$ to $\C([0,1] ,\C_R[0,1], M_{2d}(\hat\Omega_n\A_i)))$: 

For $s \ge \frac 12$ this follows from the fact that  by induction
the map 
$$(s,t) \mapsto \bigl(y \mapsto Rp_{(1-s)t}(\cdot,y)^{n-1}\bigr)$$ 
is a
smooth map from $[\frac 12,1]\times (0,\infty)$ to $\C\bigl([0,1],\C_R([0,1],
M_{2d}(\hat\Omega_n\A_i))\bigr)$. Furthermore the family $e^{-stD_I^2}$ is
uniformly bounded  on $C^l_R([0,1], M_{2d}(\hat\Omega_n\A_i))$ for any $l \in
\bbbn_0$ by Prop. \ref{sgIgen} and depends smoothly on $s,t$.

For $s\le \frac 12$ this holds since  
$$(s,t) \mapsto \bigl(x \mapsto R^* k_{st}(x,\cdot)^*\bigr)$$ is a
smooth map with values in $\C([0,1],\C_R([0,1], M_{2d}(\hat\Omega_1\A_i)))$ and since the action of
the family
$I_{n-1}((1-s)t)^*$ on $C^m_R([0,1], M_{2d}(\Ol{\mu})))$ depends smoothly on
$s,t$ and is uniformly bounded for any $m \in \bbbn_0$.

Assertion (2) holds since $I_n(t)^*=(-1)^nI_n(t)$.

The preceding arguments and the following facts imply the estimate in (3): 

By induction   there is $C>0$ such that the norm of $\bigl(y \mapsto
Rp_{(1-s)t}(\cdot,y)^{n-1}\bigr)$ in $C^m\bigl([0,1],C_R^l([0,1],
M_{2d}(\hat\Omega_n\A_i))\bigr)$ is bounded by $C(1+
t^{-\frac{l+m+1}{2}})e^{-\omega t}$ for $0 \le s\le \frac 12,~t>0$. Furthermore by
Cor. \ref{hkIgenub} the norm of $\bigl(x \mapsto R^*k_{st}(x,\cdot)^*\bigr)$ is bounded in $C^l([0,1],C_R^m([0,1],
M_{2d}(\hat\Omega_1\A_i))$ by $C(1+t^{-\frac{l+m+1}{2}})e^{-\omega t}$ for $s> \frac 12, ~t>0$.
\end{proof}
 
The proof of the previous proposition did not use the fact that
$D_I^2+\sigma R$ is the curvature of a superconnection. Hence an analogous argument shows
that $e^{-t(D_{I_s}^2+\sigma R)}$ is an integral operator whose integral kernel can be written as $$\sum_{n=0}^{\infty} (-1)^{\frac{(n+1)n}{2}}\sigma^nt^n p^s_t(x,y)^n$$ with $\bigl(y
\mapsto p^s_t(\cdot,y)^n \bigr) \in \C([0,1],\C_R([0,1],
M_{2d}(\hat\Omega_n\A_i))$ for $R=(P_s,1-P_s)$.

\begin{lem}
\label{supIsgauss}
For every $l,m,n \in \bbbn_0$ and $\ve,\delta >0$ there is $C>0$ such that for all $x,y \in [0,1]$
with $d(x,y) > \ve$ and all $t>0$  
$$|\ra_x^l \ra_y^m p^s_t(x,y)^n| \le Ce^{-\frac{d(x,y)^2}{(4+ \delta)t}} \ .$$
\end{lem}

\begin{proof} 
In degree $n=0$ the assertion holds by Lemma \ref{hkIsest}.

Let $\eta=\ve/4$. 

Let $\chi:\bbbr \to [0,1]$ be a smooth function with $\chi(x)=0$ for $x> \eta$
and $\chi(x)=1$ for $x \le \eta/2$.

Let $k_t$ be the integral kernel of $e^{-tD_{I_s}^2}$.

For $l,m \in \bbbn_0$ 
\begin{eqnarray*}
\lefteqn{\ra_x^l \ra_y^m p^s_t(x,y)^n}\\
&=&\int_0^1 \int_0^1 (1-s)^{n-1}\ra_x^lk_{st}(x,r)R(r)\chi(d(x,r))\ra_y^mp^s_{(1-s)t}(r,y)^{n-1}dr ds \\
&+& \int_0^1 \int_0^1
(1-s)^{n-1}\ra_x^lk_{st}(x,r)R(r)\bigl(1-\chi(d(x,r))\bigr) \dots \\
&& \qquad \dots \chi(d(r,y))\ra_y^m
p^s_{(1-s)t}(r,y)^{n-1}dr ds \\
&+& \int_0^1 \int_0^1 (1-s)^{n-1} \ra_x^l
k_{st}(x,r)R(r)\bigl(1-\chi(d(x,r))\bigr)\dots \\
&& \qquad \dots \bigl(1-\chi(d(r,y))\bigr) \ra_y^m
p^s_{(1-s)t}(r,y)^{n-1}dr ds \ .
\end{eqnarray*}

We begin by estimating the first term on the right hand side: By induction there is $C>0$ such that for $x,y \in [0,1]$
with $d(x,y) > \ve,~0<s<1$ and $t>0$  
$$\| R\chi(d(\cdot,x))\ra_y^mp^s_{(1-s)t}(\cdot,y)^{n-1}\|_{C_R^l}
\le Ce^{-\frac{(d(x,y)-\eta)^2}{(4+ \delta)t}} $$
in $C^l_R([0,1],M_{2d}(\hat\Omega_n\A_i))$.

Since the operator $e^{-stD_{I_s}^2}$ is
uniformly bounded on $C^l_R([0,1],M_{2d}(\hat\Omega_n\A_i))$, the first term is
bounded by $Ce^{-\frac{(d(x,y)-\eta)^2}{(4+ \delta)t}}$.

An analogous bound exists for the second term:  By
 Lemma \ref{hkIsest} there is $C>0$ such that for all $x,y \in [0,1]$ with
$d(x,y) > \ve$ and $0<s<1$ and $t>0$ in $C^m_R([0,1],M_{2d}(\hat\Omega_n\A_i))$
$$\|\bigl(\ra_x^l
k_{st}(x,\cdot)R\bigr)^*\bigl(1-\chi(d(\cdot,x))\bigr)\chi(d(\cdot,y))\|_{C_R^m}
\le C e^{-\frac{(d(x,y)-\eta)^2}{(4+ \delta)t}} \ .$$
Furthermore the integral kernel
$(y,r) \mapsto (\ra^m_y p^s_{(1-s)t}(r,y)^{n-1})^*$ induces a  uniformly bounded family of operators from
$C^m_R([0,1],M_{2d}(\Ol{\mu}))$ to $C_R([0,1],M_{2d}(\Ol{\mu}))$.

The third term is bounded by  
$C e^{-\frac{d(x,y)^2}{(4+\delta)t}}$ since by Lemma
\ref{hkIsest} 
$$|\ra^l_x k_{st}(x,r)R(r)\bigl(1-\chi(d(x,r))\bigr)| \le C e^{-\frac{d(x,r)^2}{(4+ \delta)st}}$$
 and by induction   
$$| \bigl(1-\chi(d(x,r))\bigr)
 \ra^m_y p^s_{(1-s)t}(r,y)^{n-1}| \le  C e^{-\frac{d(r,y)^2}{(4+ \delta)(1-s)t}}$$
for all $x,y,r \in [0,1]$, all $0<s<1$ and $t>0$.

Hence there is $C>0$ such that for all $x,y \in [0,1]$ and all $t>0$ 
$$|\ra_x^l \ra_y^m p^s_t(x,y)^n| \le C e^{-\frac{(d(x,y)-\ve/2)^2}{(4+ \delta)t}} \
.$$
The assertion follows now from Lemma \ref{gaussest}.

\end{proof}

As in \S \ref{hksemI} we apply Duhamel's principle in order to obtain an estimate for the
kernel $p_t(x,y)^n$: 

Recall the definitions of
$\phi_k,\gamma_k,~k=0,1,$ in \S \ref{hksemI}.

By Lemma \ref{launit} we can find unitaries $U_k \in M_{2d}(\Ai),~k=0,1,$ with $U_kI_0=I_0U_k$ and
$U_kP_kU_k^*=P_s$. 

Let $W_n^k(t)$ be the integral operator with integral kernel
 $$w^k_t(x,y)^n:=U_k^*p^s_t(x,y)^nU_k $$ 
 and denote by $W_n(t)$ the integral operator of
$$w_t(x,y)^n:=\gamma_0(x)w_t^0(x,y)^n \phi_0(y) + \gamma_1(x)w^1_t(x,y)^n
 \phi_1(y) \ .$$
 Set $W_0(0):=1$ and $W_n(0):=0$ for
$n \ge 1$. Then  $W_n(t)$ is a strongly continuous family of operators on
 $L^2([0,1],(\Ol{\mu})^{2d})$ for all $n \in \bbbn_0$. For $f \in
\C_R([0,1],(\Ol{\mu})^{2d})$ the map $[0,\infty) \to  L^2([0,1],(\Ol{\mu})^{2d}),~t \mapsto W_n(t)f \in L^2([0,1],(\Ol{\mu})^{2d})$
is even smooth.

Furthermore for $t>0$ the range of $W_n(t)$ is in $\C_R([0,1],(\Ol{\mu})^{2d})$.

Hence Duhamel's principle yields for $f \in \C_R([0,1], \A_i^{2d})$:
\begin{eqnarray*}
\lefteqn{\bigl(e^{-tA_I^2}-\sum\limits_{n=0}^{\infty}\sigma^n (-1)^{\frac{n(n+1)}{2}}
t^n W_n(t)\bigr)f}\\
&=&- \int_0^t e^{-sA_I^2}\bigl(\frac{d}{dt}+A_I^2\bigr)\sum\limits_{n=0}^{\infty} \sigma^n (-1)^{\frac{n(n+1)}{2}} (t-s)^n W_n(t-s) f ~ds \\
&=& \int_0^t e^{-sA_I^2}\sum_{k=0,1} \sum\limits_{n=0}^{\infty}\sigma^n (-1)^{\frac{n(n+1)}{2}} (t-s)^n [\gamma_k,D_I^2]_s W_n^k(t-s) \phi_k f ~ds \\
&=& \sum\limits_{n=0}^{\infty}\sigma^n  \sum\limits_{j=0}^n (-1)^{k(n,j)} \int_0^t s^{n-j}(t-s)^j
I_{n-j}(s) \dots \\
&& \qquad \dots \sum_{k=0,1}(\gamma_k'\ra + \ra \gamma_k') W_j^k(t-s) \phi_k f ~ds 
\end{eqnarray*}
with $k(n,j)=\frac{(n-j)(n-j+1)+j(j+1)}{2}$.

It follows that
\begin{eqnarray*}
\lefteqn{(*) \qquad \bigl(I_n(t)-W_n(t)\bigr)f} \\
&=& t^{-n} \sum\limits_{j=0}^n (-1)^{j(n-j)}\int_0^t s^{n-j}(t-s)^j I_{n-j}(s)\sum_{k=0,1}(\gamma_k'\ra + \ra \gamma_k')  W_j^k(t-s) \phi_k f ~ds \ .
\end{eqnarray*}

\begin{prop}
For every $l,m,n \in \bbbn_0$ and every $\delta>0$ there is $C>0$ such that 
$$|\ra_x^l\ra_y^m p_t(x,y)^n- \ra_x^l \ra_y^m w_t(x,y)^n| \le Ct \sum_{k=0,1} e^{-\frac{d(y,\supp \gamma'_k)^2}{(4+ \delta)t}}1_{\supp \phi_k}(y)$$
for all $t>0$ and all $x,y \in [0,1]$ .
\end{prop}

\begin{proof}
The proof is as for $e^{-tD_I^2}$ (see Prop. \ref{hkIgen}).

From the previous lemma it follows that for every $j\in \bbbn_0$ and
$k=0,1$ there is  $C>0$ such that
$$\|(\gamma_k'\ra + \ra \gamma_k')  \ra_y^m \bigl(w_{t-s}^k( \cdot, y )^j \phi_k(y)\bigr) \|_{C_R^l}\le
Ce^{-\frac{d(y,\supp \gamma'_k)^2}{(4+ \delta)t}}1_{\supp \phi_k}(y)$$
for all $y\in [0,1]$, all $t>0$ and $0<s<t$.
From this estimate, from the fact that $I_{n-j}(s)$ is a uniformly bounded
family of operators on $C_R^l([0,1],M_{2d}(\Ol{n}))$ and from the equation $(*)$ one deduces as in the proof of Prop. \ref{hkIgen} that
$$\|\ra_y^m p_t(\cdot,y)^n-\ra_y^m w_t(\cdot,y)^n \|_{C^l_R} \le Ce^{-\frac{d(y,\supp \gamma'_k)^2}{(4+ \delta)t}}1_{\supp \phi_k}(y)$$
for all $y\in [0,1]$, all $t>0$ and $0<s<t$.
\end{proof}

\begin{cor}
\label{gaussupconI}
For every $l,m,n \in \bbbn_0$ and $\ve>0$ there are $c,C>0$ such that
$$|\ra_x^l\ra_y^m p_t(x,y)^n| \le Ce^{-\frac{d(x,y)^2}{ct}} \ $$
for all $x,y \in [0,1]$ with $d(x,y)>\ve$ and all $t>0$.
\end{cor}

\begin{proof}
Note that the assertion is equivalent to the assertion that for every $l,m,n \in \bbbn_0$ and $\ve>0$ there are $c,C>0$ such that 
$$|\ra_x^l\ra_y^m p_t(x,y)^n| \le Ce^{-\frac{1}{ct}} \ $$
for all $x,y \in [0,1]$ with $d(x,y)>\ve$ and all $t>0$.
The estimate follows for $0<t<1$ from the previous proposition
and  Lemma \ref{supIsgauss} since $\supp  \gamma'_k \cap
\supp  \phi_k =\emptyset$. For $t>1$ it holds by the estimate in Prop. \ref{supkerI}.
\end{proof}

\subsection{The $\eta$-form}

In the following $(D_Ip_t)(x,y)^n$ denotes the integral kernel of $D_II_n(t)$.

\begin{lem}
For every $n \in \bbbn_0$ the integral 
$$\int_0^{\infty} t^{\frac {n-1}{2}} \int_0^1 \tr (D_Ip_t)(x,x)^n dx dt$$
is well-defined in $\Oi\Ai/\ov{[\Oi\Ai,\Oi\Ai]_s}$.
\end{lem}

\begin{proof}
The integral converges in $\Ol{\mu}/\ov{[\Ol{\mu},\Ol{\mu}]_s}$ for all $n,\mu,i
\in \bbbn_0$: For $n=0$ and
$t \to 0$
the convergence follows from Cor. \ref{tracet0}; for $n>0$ and $t
\to 0$ and for $n \in \bbbn_0$ and $t \to \infty$ it follows from Prop. \ref{supkerI}. For the convergence in $\Oi\Ai/\ov{[\Oi\Ai,\Oi\Ai]_s}$ apply Prop. \ref{projlim}.
\end{proof}

Let $\trsi(a+\sigma b):=\tr a$ for $a,b \in M_{2d}(\Oi\A_i)$ and let
$\Trsi$ be the corresponding trace on integral operators.

\begin{ddd}
\label{defeta}
The {\sc $\eta$-form of the superconnection $A_I$} is
$$\eta(A_I) := \frac{1}{\sqrt \pi} \int_0^{\infty} t^{-\frac 12} \Trsi D_Ie^{-(A^I_t)^2} dt \
 \in \Oi\Ai/\ov{[\Oi\Ai,\Oi\Ai]_s} \ .$$
\end{ddd}

Since
$$\eta(A_I)=\frac{1}{\sqrt \pi} \sum_{n=0}^{\infty} (-1)^n  \int_0^{\infty}
t^{n-\frac 12} \int_0^1 \tr (D_Ip_t)(x,x)^{2n} ~dx dt$$
the $\eta$-form is well-defined by the previous lemma.

For $\ve \in ]0,\frac 12]$ let $A_I(\ve):=U_{\ve}^* \di U_{\ve}+ \sigma D_I$ be
  a family of superconnections with $\supp
(U_{\ve}-1) \subset [0,\ve] \cup  [1- \ve,1] \ .$

\begin{prop}
\label{limeta}
The limit $\lim\limits_{\ve \to 0} \eta(A_I(\ve))$ exists and does not depend on the choice of $U_{\ve}$.
\end{prop}

\begin{proof}
Let $R_{\ve}:=[D_I,U_{\ve}^* \di U_{\ve}] \ .$
It suffices to prove that for $y,z \in [0,1]$ and $s,t>0$ the limit 
$$\lim\limits_{\ve \to 0}\int_0^1 k_s(y,x)R_{\ve}(x)k_t(x,z)~dx$$ exists and does not
depend on the choice of $U_{\ve}$. 

Let $f(x):=k_s(y,x)$ and $g(x):=k_t(x,z)$.

Then
\begin{eqnarray*}
\lefteqn{\int_0^1 f(x)R_{\ve}(x)g(x)~dx}\\
 &=& \int_0^1 f(x)[I_0 \ra_x, U_{\ve}^*(x) \di U_{\ve}(x)]g(x)~dx
\\
&=& \int_0^1 \ra_x\bigl(f(x)I_0 U_{\ve}^*(x)\di U_{\ve}(x) g(x)\bigl) ~dx\\
&&-\int_0^1
f'(x)I_0 U_{\ve}^*(x)\di U_{\ve}(x) g(x)~dx- \int_0^1 f(x)I_0 U_{\ve}^*(x)\di U_{\ve}(x) 
g'(x)~dx \\
&&
\\
&=& f(1)I_0 U_{\ve}^*(1) \di U_{\ve}(1) g(1)- f(0)I_0 U_{\ve}^*(0) \di
 U_{\ve}(0) g(0) \\
&&-\int_0^1
f'(x)I_0 U_{\ve}^*(x)\di U_{\ve}(x) g(x)~dx- \int_0^1 f(x)I_0 U_{\ve}^*(x)\di U_{\ve}(x) g'(x)~dx\\
&=&-\int_0^1
f'(x)I_0 U_{\ve}^*(x)\di U_{\ve}(x) g(x)~dx- \int_0^1 f(x) I_0 U_{\ve}^*(x)\di U_{\ve}(x) 
g'(x)~dx \ .
\end{eqnarray*}

Since for $\ve \to 0$ and $x \in (0,1)$ the term $U_{\ve}^*(x)\di U_{\ve}(x) g(x)$
resp. $U_{\ve}^*(x)\di U_{\ve}(x)
g'(x)$ converges to $\di g(x)$ resp. $\di g'(x)$, it follows that
$$\lim\limits_{\ve \to 0}\int_0^1 f(x)R_{\ve}(x)g(x)~dx=f(0)I_0 \di
g(0)-f(1)I_0 \di g(1) \ .$$ 
\end{proof}

\begin{ddd}
\label{defetatwo}
We call $$\eta(P_0,P_1):=\lim\limits_{\ve \to 0} \eta(A_I(\ve))$$ the 
{\sc $\eta$-form associated to the pair $(P_0,P_1)$}.
\end{ddd}

\section{The superconnection associated to $D_Z$}
\label{supcyl}

Let $(P_0,P_1)$ be a pair of transverse Lagrangian projections with $P_i\in
M_{2d}(\Ai),\\ ~i=1,2,$ and let  $D_Z$ be the operator on
$L^2(Z,(\Ol{\mu})^{4d})$ with boundary conditions $(P_0,P_1)$ as in \S \ref{cyl}.

Let $U \in \C([0,1],M_{2d}(\Ai))$ be as in Prop. \ref{defWI} with
$U(0)P_0U(0)^*=P_s$ and $U(1)P_1U(1)^*=1-P_s$. We consider $U$ as a function on
$Z$ depending only on the variable $x_2$ and set 
$$\tilde W:= U \oplus U\in \C(Z,M_{2d}(\Ai)) \ .$$
We call $A_Z:=\tilde W^*\di \tilde W
+ D_Z$ a superconnection associated to $D_Z$ and $A^Z_t:=\tilde W^*\di \tilde W
+ \sqrt t D_Z$ the corresponding rescaled superconnection.

The curvature is
\begin{eqnarray*}
A_Z^2 &=& D_Z^2+[\tilde W^*\di \tilde W,D_Z]_s\\
&=&D_Z^2 +\tilde W^* [\di,\tilde W c(d x_2)\ra_{x_2} \tilde W^*]_s\tilde W\\
&=& D_Z^2 - c(d x_1)\tilde W^*[\di,\tilde W I\ra_{x_2} \tilde W^*]_s\tilde W \\
&=&D_Z^2 - c(d x_1)\tilde W^*\di(\tilde W I\ra_{x_2} \tilde W^*)\tilde W \\
&=& D_Z^2 +c(d x_1)(R \oplus (-R))  
\end{eqnarray*}
with $R=-U^*\di(UI_0(\ra U^*))U$.

Let $\tilde R=c(d x_1)(R \oplus (-R)).$ Then $\tilde R \in
\C(Z,M_{4d}(\Ol{\mu}))$ and $\tilde R^*= \tilde R$.

We have a holomorphic semigroup $e^{-tA_Z^2}$ on $L^2(Z, (\Ol{\mu})^{4d})$.

From \S\ref{cyl} recall that 
$e^{-tD_Z^2}=e^{-t\tilde D_I^2}e^{-t\Delta_{\bbbr}}$.  

By Prop. \ref{supkerI} the operator
$$I_n(t)=\int_{\Delta^n}e^{-u_0tD_I^2}Re^{-u_1tD_I^2}R\dots e^{-u_ntD_I^2}
~du_0 \dots du_{n}$$ is an integral operator. Its integral kernel is denoted by $p^I_t(x,y)^n$.

Since $e^{-t\Delta_{\bbbr}}$ commutes with $\tilde R$ and $c(d x_1)$ commutes with
$e^{-tD_Z^2}$ and $R \oplus (-R)$, we obtain from Volterra development:
\begin{eqnarray*}
e^{-tA_Z^2}&=&\sum_{n=0}^{\infty}(-1)^nt^n\int_{\Delta^n}e^{-u_0tD_Z^2}\tilde  Re^{-u_1tD_Z^2} \tilde R\dots e^{-u_ntD_Z^2} ~du_0
\dots du_n\\
&=&\sum_{n=0}^{\infty}(-1)^nt^n e^{-t\Delta_{\bbbr}}\int_{\Delta^n}e^{-u_0t\tilde
D_I^2}\tilde  Re^{-u_1t \tilde D_I^2} \tilde R\dots e^{-u_nt \tilde D_I^2} ~du_0
\dots du_n \\
&=& \sum_{n=0}^{\infty}(-1)^nt^n c(dx_1)^n e^{-t\Delta_{\bbbr}} \bigl(I_n(t)
\oplus (-1)^nI_n(t)\bigr) \ .
\end{eqnarray*}
We define $$p_t^Z(x,y)^n:=
c(d x_1)^n \frac{1}{\sqrt{4 \pi
t}}~e^{-\frac{(x_1-y_1)^2}{4t}}
\bigl(p^I_t(x_2,y_2)^n \oplus (-1)^n p^I_t(x_2,y_2)^n\bigr) \ .$$
The integral kernel
of $e^{-tA_Z^2}$ is $$\sum\limits_{n=0}^{\infty}(-1)^nt^n p_t^Z(x,y)^n$$   and
the integral kernel of $e^{-(A_t^Z)^2}$  $$\sum\limits_{n=0}^{\infty}(-1)^nt^{\frac n2}
p_t^Z(x,y)^n \ .$$
Note that for all multi-indices $\alpha,\beta \in \bbbn_0^2$ 
$$\trs\ra^{\alpha}_x \ra^{\beta}_y p_t^Z(x,y)^n=0 \ .$$
Furthermore 
$$(D_Z)_xp_t^Z(x,y)^n=c(dx_1)(\ra_{x_1}+I
\ra_{x_2})p_t^Z(x,y)^n \ .$$
Since $\trs c(dx_1)\ra_{x_1}p_t^Z(x,y)^n$ vanishes for $x_1=y_1$ and $\trs c(dx_1)
I \ra_{x_2}p_t^Z(x,y)^n$ vanishes for all $x,y \in Z$, 
$$\trs (D_Z)_xp_t^Z(x,y)^n=0$$ if $x=y.$ 

Furthermore the integral kernel $p_t^Z(x,y)^n$ satisfies the following Gaussian estimate:

\begin{lem}
\label{gaussupZ}
Let $\alpha,\beta \in \bbbn_0^2$.  For every $\ve>0$ there are $c,C>0$ such that
for all $x,y \in Z$ with $d(x,y)>\ve$ and all $t>0$
$$|\ra_x^{\alpha}\ra_y^{\beta} p^Z_t(x,y)^n| \le
Ce^{-\frac{d(x,y)^2}{ct}} \ .$$
\end{lem}

\begin{proof}
The assertion follows from  Prop. \ref{supkerI} and Cor. \ref{gaussupconI}.
The arguments are as in the proof of Lemma \ref{gaussZ}.
\end{proof}

\section{The superconnection $A(\rho)_t$ associated to $D(\rho)$}

\subsection{The family $e^{-A(\rho)_t^2}$}

\label{supconM} 

In \S \ref{hsDrho} we fixed $r_0,b_0>0$ such that
$$\supp k_K \cap \bigl((F(r_0,b_0) \times M) \cup (M \times F(r_0,b_0))\bigr)= \emptyset \ .$$
Let $W \in \C(M, \End^+ E \ten \Ai)$ be as in \S \ref{stand} such that $W$ is parallel on $\{ x \in M~|~d(x,\ra M)>b_0\}$. Then $[W,K]_s=0$. 

We define a superconnection associated to $D(\rho)$ on the $\bbbz/2$-graded $\A_i$-module
$L^2(M,E \ten \A_i)$  by 
$$A(\rho):=W^*\di W+D(\rho) \ .$$
The corresponding rescaled superconnection is
$$A(\rho)_t:=W^*\di W+ \sqrt t D(\rho) \ .$$
The curvature of $A(\rho)$ is
\begin{eqnarray*}
A(\rho)^2 &=& W^*\di^2W+ [D(\rho),W^*\di W]_s +D(\rho)^2\\
&=&D(\rho)^2 +   W^* [\di,  Wc(d   W^*)]_s  W\\
&=:&D(\rho)^2+\RM \ .
\end{eqnarray*}
We used Prop. \ref{DMreg} and $[W,K]_s=0$.

Then $\RM \in \C(M, \End E \ten \Ol{1})$ with $\RM^*=\RM$. In the flat region 
\begin{eqnarray*}
\RM|_F&=& W^* [\di,  W c(e_2) \ra_{e_2}  W^*)]_s  W\\
&=& - c(e_1)  W^*\di(W I\ra_{e_2}  W^*)  W \ .
\end{eqnarray*}
Furthermore $\RM$ vanishes on $\{ x \in M~|~d(x,\ra M)>b_0\}$.

For every $k \in \bbbz/6$ the restriction of $\RM$ to $\U_k$ is of
the form $c(e_1)(R \oplus (-R))$ with $R \in \C(\U_k, M_{2d}(\Ol{1}))$ and $R$ is independent of the variable $x_2^k$.
 
The rescaled curvature  is
$$A(\rho)_t^2=tD(\rho)^2+\sqrt t \RM \ .$$
As expected, $A(\rho)^2$ and $A(\rho)_t^2$ are right $\Ol{\mu}$-module homomorphisms.

Since $A(\rho)^2$ is a bounded perturbation of $D(\rho) ^2$, it generates a holomorphic semigroup
$e^{-tA(\rho)^2}$ on $L^2(M,E \ten \Ol{\mu})$. 

In the following we assume that $t \ge 0$.

By  Volterra development 
\begin{eqnarray*}
e^{-tA(\rho)^2} &=& \sum\limits_{n=0}^{\infty}(-1)^nt^n\int_{\Delta^n}e^{-u_0tD(\rho)^2}
\RM e^{-u_1tD(\rho)^2} \RM\dots e^{-u_ntD(\rho)^2} ~du_0 \dots du_n\\
&=:&\sum\limits_{n=0}^{\infty}(-1)^nt^n I_n(\rho,t) \ .
\end{eqnarray*}
It follows that
$$e^{-A(\rho)_t^2}=\sum\limits_{n=0}^{\infty}(-1)^n t^{n/2}I_n(\rho,t) \ .$$

For $\rho \neq 0$ the family $t \mapsto I_n(\rho,t)$ is uniformly bounded on $L^2(M,E
\ten \Ol{\mu})$. By Cor. \ref{semCm} it acts as a strongly continuous family of
operators  on $C^m_R(M,E \ten
\Ol{\mu})$ and there are $C,l>0$ such that the action is bounded by $C(1+t)^l$.

\subsection{The integral kernel of $e^{-A(\rho)_t^2}$}
\label{supconMhk}

In this section we prove that $I_n(\rho,t)$ is an
integral operator for $t>0$. As usual we construct an
approximation of the family $I_n(\rho,t)$ by a family of integral operators and
compare it with $I_n(\rho,t)$ by Duhamel's principle. 

Let $\U(r_0,b_0)=\{\U_k \}_{k \in J}$, $\{\phi_k\}_{k \in
J}$ and $\{\gamma_k\}_{k \in J}$ be as in
\S \ref{hsDrho}. 

For $k \in \bbbz/6$ the function $W|_{\U_k}:\U_k \to M_{4d}(\Ai)$ does not depend on the coordinate
$x_1^k$. We extend it to a section
$\tilde W_k:Z_k \to M_{4d}(\Ai)$ independent of $x_1^k$ and define
the superconnection $A_{Z_k}:=\tilde W_k^* \di \tilde W_k + D_{Z_k}$, which coincides on $\U_k$ with the superconnection $A(\rho)$.

For $k \in \bbbz/6$ and $n \in \bbbn_0$ let
$w(\rho)^k_t(x,y)^n$  be the restriction of $p^{Z_k}_t(x,y)^n$ to $\U_k \times \U_k$.

Let furthermore 
$w(\rho)^{\cp}_t(x,y)^0$ be the restriction of the
integral kernel of $e^{-tD_N(\rho)^2}$ to $\U_{\cp} \times \U_{\cp}$ and set
$w(\rho)^{\cp}_t(x,y)^n=0$ for $n>0$. 

The reason for this is that $A(\rho)^2$
 equals $D(\rho)^2$ on $\U_{\cp}$.

We extend $w(\rho)^k_t(x,y)^n$ by zero to
$M \times M$ and set $$w(\rho)_t(x,y)^n:=\sum_{k \in J} \gamma_k(x)
w(\rho)^k_t(x,y)^n \phi_k(y) \ .$$
Write $W_n(\rho,t)$ for the corresponding integral operator. It is a bounded
operator on $L^2(M,E \ten \Ol{\mu})$ and on $C^m_R(M,E \ten \Ol{\mu})$.

Set
$W_0(\rho,0)=1$ and $W_n(\rho,0)=0$ for $n>0$.

Then for $f \in L^2(M,E \ten \Ol{\mu})$ the family $W_n(\rho,t)f \in L^2(M,E
\ten \Ol{\mu})$ depends continuously on $t$ for all $t \in [0,\infty)$,  and for $f \in
\C_{Rc}(M,E \ten \A_i)$ even smoothly.

For $f \in \C_{Rc}(M,E \ten \A_i)$ Duhamel's principle yields: 
\begin{eqnarray*}
\lefteqn{\bigl(e^{-tA(\rho)^2}-\sum\limits_{n=0}^{\infty}(-1)^n t^n W_n(\rho,t)\bigr)f } \\
&=&-\int_0^t e^{-sA(\rho)^2}(\frac{d}{dt}+A(\rho)^2)\sum_{k \in J} \sum\limits_{n=0}^{\infty} (-1)^n(t-s)^n  \gamma_k W_n^k(\rho,t-s) \phi_k f ~ds 
\end{eqnarray*}

\begin{eqnarray*}
&=& \int_0^t e^{-sA(\rho)^2}\sum_{k \in J} [\gamma_k,D(\rho)^2]_s\sum\limits_{n=0}^{\infty} (-1)^n(t-s)^n  W_n^k(\rho,t-s) \phi_k f ~ds \\
&=& \sum\limits_{n=0}^{\infty}(-1)^n\int_0^t \sum\limits_{m=0}^n s^{n-m}(t-s)^m I_{n-m}(\rho,s) \dots \\
&& \qquad \dots \sum_{k \in J} [\gamma_k,D(\rho)^2]_s W_m^k(\rho,t-s) \phi_k f ~ds \ . \end{eqnarray*}
Hence
\begin{eqnarray*}
\lefteqn{\bigl(I_n(\rho,t)-W_n(\rho,t)\bigr)f} \\
&=&- t^{-n}\int_0^t \sum\limits_{m=0}^n s^{n-m}(t-s)^m I_{n-m}(\rho,s)\sum_{k \in J}[\gamma_k,D(\rho)^2]_s  W_m^k(\rho,t-s) \phi_k f ~ds \ .
\end{eqnarray*}
In the following proposition $|\cdot |$ denotes the norm on the fibers of $(E \boxtimes E^*)
\ten \Ol{\mu}$.

 \begin{prop}
\label{hksupcon}
The operator $I_n(\rho,t)$ is an
integral operator for $t>0$. Let $p(\rho)_t(x,y)^n$ be its integral kernel. 

\begin{enumerate}
\item The map $(0,\infty) \to
\C(M \times M, (E \boxtimes E^*) \ten \hat\Omega_n \Ai), t \mapsto p(\rho)_t^n$ is
smooth. 

\item $p(\rho)_t(x,y)^n=(p(\rho)_t(y,x)^n)^*$. 

\item
For every $T>0$ there are $C,c>0$ such that 
$$|p(\rho)_t(x,y)^n-w(\rho)_t(x,y)^n| \le Ct~\bigl(|\rho| 1_{\U_{\cp}}(y) + \sum_{k \in J} 
e^{-\frac{d(y,\supp d\gamma_k)^2}{ct}}1_{\supp\phi_k}(y) \bigr)$$
for all $0<t<T$, $\rho \in [-1,1]$ and all $x,y \in M$.
\item Let $\rho \neq 0$. Then there are $C,c>0$ and $j \in \bbbn$ such that 
$$|p(\rho)_t(x,y)^n- w(\rho)_t(x,y)^n| \le Ct(1+t)^j 
e^{-\frac{d(y, \U_{\cp})^2}{ct}}$$
for all $t>0$ and all $x,y \in M$.
\end{enumerate}
Statements analogous to (3), (4) hold for the partial derivatives of
$p(\rho)_t(x,y)^n$ in $x$ of $y$ with respect to unit vector fields on $M$.
\end{prop}

\begin{proof}
The proof follows the proof of Prop. \ref{gaussM}. 

In order to show the
existence of the integral kernel and (1) we need only investigate
$I_n(\rho,t)-W_n(\rho,t)$.

For $f \in \C_{Rc}(M,E \ten \A_i)$
\begin{eqnarray*}
\lefteqn{\bigl(I_n(\rho,t)-W_n(\rho,t)\bigr) f}\\
 &=& t^{-n}\int_0^t\int_M  \sum\limits_{m=0}^n s^{n-m}(t-s)^m
 I_{n-m}(\rho,s)\sum_{k \in J} [c(d \gamma_k),D]_s  w(\rho)_{t-s}^k( \cdot, y
 )^m \phi_k(y) f(y)~dy ds  \ .
\end{eqnarray*}
For $\rho \in \bbbr$ and $t>0$ the family $I_{n-m}(\rho,s)$ is uniformly bounded on $C_R^{\nu}(M,E \ten
\Ol{\mu})$ in $s<t$.

The function
$$\tau \mapsto \bigl(y \mapsto \sum_{k \in J}  [c(d \gamma_k),D]_s w(\rho)_{\tau}^k( \cdot, y )^m
\phi_k(y)\bigr)$$
is smooth from $(0,\infty)$ to  $C^l(M,C_{Rc}^{\nu}(M,E \ten
\Ol{m+1})\ten E^*)$ for any $l,\nu \in \bbbn_0$. 

If $k \in \bbbz/6$, then  by Lemma \ref{gaussupZ}  there
are $c,C>0$ such that for all $y \in M$ and $0 < \tau$ 
$$ \|[c(d \gamma_k),D]_s w_{\tau}^k( \cdot, y )^m \phi_k(y)\|_{C^{\nu}} 
\le C e^{-\frac{d(y,\supp d\gamma_k)^2}{c\tau}}1_{\supp\phi_k}(y) \ .$$
Furthermore $w(\rho)_{\tau}^{\cp}(x,y)^m=0$ for $m>0$. The kernel
$w(\rho)^{\cp}_{\tau}(x,y)^0$ is equal to the kernel $e(\rho)^{\cp}_{\tau}(x,y)$ in the
proof of Prop. \ref{gaussM} and was estimated there. It follows that for $T,\delta>0$
there is $C>0$ such that
$$\|[c(d \gamma_k),D]_s  w(\rho)_{\tau}^{\cp}( \cdot, y )^0 \phi_{\cp}(y)
\|_{C^{\nu}} \le C\bigl(e^{-\frac{d(y,\supp d\gamma_{\cp})^2}{(4+ \delta)\tau}}
+ \tau |\rho|\bigr)
1_{\supp\phi_{\cp}}(y)$$
 for $y \in M$, $\rho \in [-1,1]$ and $0 < \tau<T$.

Analogous estimates hold for the derivatives in $y$ and also in $\tau$ by the heat equation.
This shows that the function above extends smoothly to $\tau=0$. 
The existence and (1) and (3) follow now by the usual arguments.

Assertion (2) follows from $I_n(\rho,t)^*=I_n(\rho,t)$.

(4) follows from the estimates by taking into account that for every $\rho \neq 0$ there is $C>0$ and $j \in \bbbn$
such that the norm of $I_n(\rho,t)$ on $C_R^l(M,E \ten \Ol{\mu})$ is bounded by $C(1+t)^j$ for all $t>0$.
\end{proof}

We deduce the following estimates for later use:

\begin{cor}
\label{hksupconMgauss}

Let $\rho \neq 0$. For every $\ve >0$ and $m,n \in \bbbn_0$ there are $c,C>0$ and $j \in \bbbn$ such that for $t >0$ and $x,y \in M$ with $d(x,y) > \ve$

$$|D(\rho)_x^mp(\rho)_t(x,y)^n| \le C(1+t)^j\bigl(e^{-\frac{d(x,y)^2}{ct}} +
e^{-\frac{d(y,\U_{\cp})^2}{ct}} \bigr).$$
\end{cor}

\begin{proof}
This follows from the proposition and Lemma \ref{gaussupZ}.
\end{proof}

\begin{cor}
\label{compcyl}

Let $k \in \bbbz/6$.
 
\begin{enumerate}
\item For every $T>0$ and $m,n \in \bbbn_0$ there are $c,C>0$ such that for all $x,y
\in \U_k$, for $0 <t<T$ and $\rho \in [-1,1]$
$$|D(\rho)_x^mp(\rho)_t(x,y)^n-(D_{Z_k})_x^m p^{Z_k}_t(x,y)^n|\le
Ce^{-\frac{d(y,\U_{\cp})}{ct}}  \ .$$
\item For every $\rho \neq 0$ and $n \in \bbbn_0$ there are $c,C>0$ and $j \in \bbbn$ such that for all $x,y
\in \U_k$ and $t>0$ 
$$|D(\rho)_x^mp(\rho)_t(x,y)^n-(D_{Z_k})_x^m p^{Z_k}_t(x,y)^n| \le
C(1+t)^j e^{-\frac{d(y,\U_{\cp})}{ct}}  \ .$$
\end{enumerate}
\end{cor}

\section{The index theorem and its proof}

\subsection{The generalized supertrace}
\label{suptr}

In the following  $\trs$ denotes the supertrace on the fibers of $(E \ten E^*) \ten \Ol{\mu}$ and $\Trs$ the corresponding supertrace for trace class operators (see \S \ref{sutrace}). 

Let $\chi:\bbbr \to [0,1]$ be a smooth function with $\chi(x)=1$ for $x \le 0$ and
$\chi(x)=0$ for $x \ge 1$. Let  $\phi_r:M \to
[0,1],~ \phi_r(x):=\chi(d(M_r,x))$ for $r>0$.

\begin{ddd} 
Let $K$ be a bounded operator on $L^2(M,E \ten \Ol{\mu})$ such that
$\phi_r K \phi_r$ is a trace class operator for all $r\ge 0$. Then we define the generalized supertrace
$$\Trs  K: =\lim\limits_{r \to \infty} \Trs \bigl(\phi_r K\phi_r) $$
if the limit exists.
\end{ddd}

For trace class operators the generalized supertrace coincides with the usual one. Note that in general it does not vanish on supercommutators.

\begin{prop}
\label{trwelldef}
For
$\rho \in \bbbr$ and $t>0$ the generalized supertraces
$$\Trs e^{-A(\rho)_t^2}$$ and $$\Trs D(\rho)e^{-A(\rho)_t^2}$$ exist.
\end{prop}

\begin{proof} We show the assertion for $e^{-A(\rho)_t^2}$, the proof for $D(\rho)e^{-A(\rho)_t^2}$ is
analogous.

 The operator 
$\phi_r e^{-A(\rho)_t^2}\phi_r$ is trace class for $t>0$ since 
$$\phi_r e^{-A(\rho)_t^2}\phi_r=(\phi_r
e^{-{A(\rho)_t^2}/2})(e^{-{A(\rho)_t^2}/2}\phi_r)$$ is a
product of Hilbert-Schmidt operators. The
operators $\phi_r I_n(\rho,t)\phi_r$ are also trace class for $t>0$. 

We show that $\trs p(\rho)_t(x,x)^n$ is in
$L^1(M,\Ol{\mu}/\ov{[\Ol{\mu},\Ol{\mu}]_s})$.

By Cor. \ref{compcyl} there are $c,C>0$ such that 
$$|p(\rho)_t(x,x)^n-p^{Z_k}_t(x,x)^n| \le Ce^{-cd(x,\U_{\cp})^2}$$ for
all $x \in \U_k$. 

Since $\trs p_t^{Z_k}(x,x)^n=0$ by \S \ref{supcyl}, 
$$|\trs p(\rho)_t(x,x)^n| \le Ce^{-cd(x,\U_{\cp})^2} \ .$$
Now the assertion follows.
\end{proof}

\subsection{The limit of $\Trs e^{-A(\rho)_t^2}$ for $t \to \infty$}

This section is devoted to the proof of

\begin{theorem}
\label{liminfty}
Let $\rho \neq 0$.

Let $P_0$ be the projection onto the kernel of $D(\rho)$. Then for $T>0$ there
is $C>0$ such that for all $t>T$ 
 $$|\Trs e^{-A(\rho)_t^2}- \sum_{n=0}^{\infty}(-1)^n\frac{1}{n!}\Trs (P_0W\di W^*P_0)^{2n}| \le Ct^{-\frac 12}$$
and 
$$|\Trs D(\rho)e^{-A(\rho)_t^2}| \le C t^{-1}$$
in $\Ol{\mu}/\ov{[\Ol{\mu},\Ol{\mu}]_s}$. 

\end{theorem}

Note that $(P_0W\di W^*P_0)^{2n}=W (W^*P_0W)(\di (W^*P_0 W))^{2n}W^*$ by Lemma \ref{chHS}, hence  $(P_0W\di W^*P_0)^{2n}$ is a  trace
class operator on $L^2(M,E \ten \Ol{\mu})$. 

The proof is subdivided into some lemmata. 

Throughout the section $\rho \neq 0$ is fixed.

In the following $| \cdot|$ denotes the norm on
$\Ol{\mu}/\ov{[\Ol{\mu},\Ol{\mu}]_s}$ resp. the fiberwise norm of $(E \boxtimes
E^*) \ten \Ol{\mu}$ (depending on the context), and $\| \cdot \|$ denotes the
operator norm of $B(L^2(M,E \ten \Ol{\mu}))$. Furthermore we make use of the Hilbert-Schmidt norm $\|~ \|_{HS}$ and the trace class norm $\|~\|_{TR}$, which are defined in \S \ref{HS} and \S \ref{trclop}.

\begin{lem} Let $\nu=0,1$. For any $T>0$ there are $\ve,C>0$ such that for
all $t>T$ 
$$| \Trs D(\rho)^{\nu} e^{-A(\rho)_t^2}- \Trs \phi_t^2 D(\rho)^{\nu}
e^{-A(\rho)_t^2} | \le Ce^{-\ve t} \ .$$
\end{lem}

\begin{proof} We prove the case $\nu=0$, the case $\nu=1$ can be proved analogously.

By Cor. \ref{compcyl} there are $c,C,r>0$ such that  
$$|p(\rho)_t(x,x)^n-p^{Z_k}_t(x,x)^n| \le C(1+t)^je^{-\frac{d(x,M_r)^2}{ct}}$$ for $t>0$
and $x \in \U_k$ with $k \in \bbbz/6$.

Hence for $T>0$ there are $c,C>0$ such that for all $x \in M$ and $t>T$  
$$|\trs p(\rho)_t(x,x)^n| \le Ce^{-\frac{d(x,M_r)^2}{ct}}$$  and thus there are
$C,\ve>0$ such that for all $t>r$ 
\begin{eqnarray*}
|\Trs  e^{-A(\rho)_t^2}- \Trs \phi_t^2  e^{-A(\rho)_t^2} | 
&=&| \Trs (1-\phi_t^2) e^{-A(\rho)_t^2}| \\
&\le&  C \int_{t}^{\infty} e^{-\frac{(r'-r)^2}{ct}}dr'\\
&\le&C  e^{-\ve t} \ .
\end{eqnarray*}
\end{proof}

Let $P_1:=1-P_0$.

\begin{lem}
Let $\nu=0,1$ and $k \in \bbbn_0$. Then the integral
$$ \int_{\Delta^k} 
\phi_t D(\rho)^{\nu}P_1e^{-u_0tD(\rho)^2}\RM  P_1e^{-u_1tD(\rho)^2} \RM
\dots   \RM  P_1e^{-u_ktD(\rho)^2} \phi_t ~du_0 \dots du_k$$
converges in the space of trace class operators and
for $T>0$ there are $C,\ve>0$ such that for
all $t>T$ the norm
$$\|\int_{\Delta^k} 
\bigl(\phi_t D(\rho)^{\nu}P_1e^{-u_0tD(\rho)^2}\RM  P_1e^{-u_1tD(\rho)^2} \RM
\dots   \RM  P_1e^{-u_ktD(\rho)^2} \phi_t\bigr) ~du_0 \dots du_k \|_{TR}$$ is bounded by $C e^{-\ve t}$.
\end{lem}

\begin{proof}

Note that for any $(u_0,\dots u_k) \in \Delta_k$ we can find $i\in \{0,1, \dots, k\}$
such that $u_i \ge \frac{1}{k+1}$.

We begin by showing that for any  $T>0$ there are $C,\ve>0$ such that
for $\frac{1}{k+1} \le u_i \le 1 $, for $0 < u_0, \dots ,u_{i-1}\le 1$ and for $t>T$ the family
$$\phi_t D(\rho)^{\nu} P_1e^{-u_0tD(\rho)^2}\RM  P_1e^{-u_1tD(\rho)^2} \RM  \dots  \RM  P_1e^{-\frac{u_i}{2}tD(\rho)^2}$$
is a family of Hilbert-Schmidt operators with Hilbert-Schmidt norm bounded by
 $C u_0^{-\frac{\nu}{2}}e^{-\ve t}$. If not specified the estimates in the
following hold for $\frac{1}{k+1}\le u_i \le 1$, for $0 < u_0, \dots ,u_{i-1}\le 1$ and
for $t>T$.

\newpage

We have that  
\begin{eqnarray*}
\lefteqn{\phi_t  D(\rho)^{\nu} P_1 e^{-u_0tD(\rho)^2}\RM  P_1e^{-u_1tD(\rho)^2} \RM  \dots  \RM  P_1e^{-\frac{u_i}{2}tD(\rho)^2}}\\
&=&\phi_t D(\rho)^{\nu} P_1e^{-u_0tD(\rho)^2}\RM  P_1e^{-u_1tD(\rho)^2} \RM  \dots  \RM  P_1e^{-\frac{u_i}{2}tD(\rho)^2} \phi_{(t+6)} \\
&+& \phi_t D(\rho)^{\nu}P_1 e^{-u_0tD(\rho)^2}\RM  
\dots P_1 e^{-u_{i-1}tD(\rho)^2}\RM \phi_{(t+6)}P_1 e^{-\frac{u_i}{2}tD(\rho)^2}  (1-\phi_{t+6})   \\
\lefteqn{ \!\!\!\!\! + \phi_t D(\rho)^{\nu}P_1e^{-u_0tD(\rho)^2}\RM  
 \dots P_1 e^{-u_{i-1}tD(\rho)^2}\RM (1-\phi_{(t+6)}) P_1
e^{-\frac{u_i}{2}tD(\rho)^2} (1-\phi_{(t+6)}) \ .}
\end{eqnarray*}
Consider the first term on the right hand side. 
By Cor. \ref{hkHS} the Hilbert-Schmidt norm of $e^{-\frac{u_i}{4}tD(\rho)^2} \phi_{(t+6)}$  is  bounded by $Ct$.

Furthermore by Prop. \ref{Mgensg} and Cor. \ref{boundsmt} there are $\ve,C>0$ such that the operator norm of 
$$\phi_t D(\rho)^{\nu} P_1e^{-u_0tD(\rho)^2}\RM  P_1e^{-u_1tD(\rho)^2} \RM  \dots  \RM  P_1e^{-\frac{u_i}{4}tD(\rho)^2}$$
is bounded by $Cu_0^{-\frac{\nu}{2}}e^{-\ve t}$.

Hence (see Prop. \ref{HScomp}) the first term is a family of Hilbert-Schmidt operators with Hilbert-Schmidt norm bounded by $Cu_0^{-\frac{\nu}{2}} e^{-\ve t}$.

In the second term the factor
$$\phi_{(t+6)}
P_1 e^{-\frac{u_i}{2}tD(\rho)^2} (1-\phi_{(t+6)})=(\phi_{(t+6)}
e^{-\frac{u_i}{4}tD(\rho)^2})P_1 e^{-\frac{u_i}{4}tD(\rho)^2}(1-\phi_{(t+6)})$$
is a Hilbert-Schmidt
operator bounded by $Ce^{-\ve t}$ for some $C,\ve>0$. 
Hence the second term is bounded in the Hilbert-Schmidt norm by
$Cu_0^{-\frac{\nu}{2}}e^{-\ve t}$.

The estimate of the third term requires more effort.
We prove by induction on $j \in \bbbn$ that there is $C>0$ such that
$$\phi_t D(\rho)^{\nu} P_1e^{-u_0tD(\rho)^2}\RM  P_1e^{-u_1tD(\rho)^2} \RM 
\dots P_1e^{-u_jtD(\rho)^2}  (1-\phi_{(t+6)})$$
is a Hilbert-Schmidt operator with
Hilbert-Schmidt norm bounded by $Cu_0^{-\frac{\nu}{2}}(1+t)$ for $t>T$ and $0 <
u_0, \dots, u_j \le 1$.
Then it follows that the third term is uniformly bounded by $ C
e^{-\ve t}$ for some $C,\ve>0$  since $P_1
e^{-\frac{u_i}{2}tD(\rho)^2}$ is exponentially
decaying for $t \to \infty$ by  Prop. \ref{Mgensg}.

For $j=0$ the assertion follows from Cor. \ref{hkHS} by
\begin{eqnarray*}
\lefteqn{\phi_t D(\rho)^{\nu}P_1e^{-u_0tD(\rho)^2}(1-\phi_{(t+6)})}\\
&=&\phi_t
D(\rho)^{\nu}e^{-u_0tD(\rho)^2}(1-\phi_{(t+6)})-\phi_t
D(\rho)^{\nu}P_0(1-\phi_{(t+6)}) \ .
\end{eqnarray*}
Now assume the assertion is true for $j-1$. We have that
\begin{eqnarray*}
\lefteqn{\phi_t P_1D(\rho)^{\nu}e^{-u_0tD(\rho)^2}\RM 
P_1  \dots  \RM  P_1e^{-u_jtD(\rho)^2}
(1-\phi_{(t+6)})} \\
&=& \phi_t P_1D(\rho)^{\nu}e^{-u_0tD(\rho)^2}P_1\RM  P_1 \dots \\
&&\qquad \dots P_1 e^{-u_{j-1}tD(\rho)^2}P_1 \RM \phi_{(t+3)}   P_1e^{-u_j tD(\rho)^2}
(1-\phi_{(t+6)}) \\
&+& \phi_t P_1D(\rho)^{\nu}e^{-u_0tD(\rho)^2}P_1 \dots \\
&& \qquad \dots P_1 e^{-u_{j-1}tD(\rho)^2}P_1 (1-\phi_{(t+3)}) \RM P_1
e^{-u_j tD(\rho)^2} (1-\phi_{(t+6)}) \ .
\end{eqnarray*}
Both terms on the right hand side are bounded in the Hilbert-Schmidt norm by
 $Cu_0^{-\nu/2}(1+t)$ for all $t>T$ and $0 <
u_0, \dots, u_j \le 1$: the first term
since $\phi_{(t+3)}  P_1e^{-u_j tD(\rho)^2}
(1-\phi_{(t+6)})$ is bounded by $C(1+t)$ by Cor. \ref{hkHS}, the second term by induction.  

This shows the claim from the beginning of the proof.

An analogous argument yields that for any $T>0$ there are $C,\ve>0$ such that
 the family
$$P_1e^{-\frac{u_i}{2}tD(\rho)^2}\RM P_1e^{-u_{i+1}tD(\rho)^2} \RM  \dots   P_1e^{-u_ktD(\rho)^2}\phi_t$$
is a family of Hilbert-Schmidt operators with Hilbert-Schmidt norm  bounded by
 $C e^{-\ve t}$ for all $u_i>\frac{1}{k+1}$, ~$0 \le u_{i+1}, \dots ,u_k \le 1$ and  $t>T$.

It follows that the integral 
$$\int_{\Delta^k} \phi_t  P_1D(\rho)^{\nu}e^{-u_0tD(\rho)^2}\RM
P_1e^{-u_1tD(\rho)^2} \RM  \dots  \RM  P_1e^{-u_ktD(\rho)^2}\phi_t ~du_0
\dots du_k $$
converges in the trace class norm and is  bounded by $C e^{-\ve t}$  for some $C,\ve>0$ and all $t>T$.
\end{proof}

We have that
$$ \phi_t D(\rho)^{\nu}e^{-A(\rho)^2_t}\phi_t =\sum\limits_{k=0}^{\infty}(-1)^k t^{k/2}
~ \phi_t D(\rho)^{\nu}I_k(\rho,t)\phi_t \ .$$
The decomposition $e^{-u_itD(\rho)^2} = P_0  + P_1e^{-u_itD(\rho)^2}$ induces a
decomposition of $D(\rho)^{\nu}I_k(\rho,t)$
into a sum of $2^{k+1}$ terms. 
For $0\le j \le k+1$ let $P_{jk}^{\nu}(t)$ be the sum of those terms with exactly $j$ factors of the form 
$P_1e^{-u_itD(\rho)^2} \ .$ 

Thus $$D(\rho)^{\nu}e^{-A(\rho)^2_t}=\sum\limits_{k=0}^{\infty}(-1)^k t^{\frac k2}~
\sum\limits_{j=0}^{k+1}  P_{jk}^{\nu}(t) \ .$$
Note that Prop. \ref{Philsch} implies that $P_0$ is a trace class operator, hence the operators $P_{jk}^{\nu}(t)$ with $j \neq k+1$ are trace class for $t \neq 0$.

From $P_0\RM P_0=P_0[W^*\di W, D(\rho)]_sP_0=0$ it follows that $P_{jk}^{\nu}(t)=0$ for $j <
\frac k2$. 

Furthermore for $k$ even 
\begin{eqnarray*}
P^{\nu}_{\frac k2 k}(t) &=& \int_{\Delta^k} D(\rho)^{\nu}P_0\RM P_1e^{-u_1tD(\rho)^2} \RM P_0 \RM P_1 e^{-u_3tD(\rho)^2}
\RM P_0\dots \\
&& \quad \dots P_0\RM P_1 e^{-u_{k-1}tD(\rho)^2} \RM P_0  ~ du_0 \dots du_k \ .
\end{eqnarray*}
Since $D(\rho)P_0=0$, it follows that $P^1_{\frac k2 k}(t)=0$. 

Moreover by the previous lemma 
$$\|\phi_t P_{(k+1)k}^{\nu}(t)\phi_t \|_{TR} \le C e^{-\ve t} $$
for $t$ large.

Now we study the behavior of the remaining cases for large $t$.

\begin{lem}  Let $\nu=0,1$ and $j,k \in \bbbn_0$ with $\frac k2 \le j\le k$. For every $T>0$ there is $C>0$ such that  for $t>T$
 $$\|P_{jk}^{\nu}(t)\|_{TR} \le Ct^{-j} \ .$$
If $k$ is odd, then for every $n \in \bbbn$ and $T>0$ there is $C>0$ such that for
$t>T$ 
$$|\Trs \phi_t^2 P_{\frac{k+1}{2}{k}}^{\nu}(t)| \le Ct^{-n} \ .$$
\end{lem}

\begin{proof}
For $j \le k$ the operator $P_{jk}^{\nu}(t)$ is a sum of terms of the form 
$$\int_{\Delta^k} (A(u_0, \dots u_i,t)P_0)(P_0 B(u_{i+1}, \dots u_k,t))~du_0
\dots du_k \ ,$$
where $A$ and $B$ are continuous families of bounded operators on $L^2(M,E \ten
\Ol{\mu})$ for $u_0 \neq 0$.

Since $P_0$ is a Hilbert-Schmidt operator by Prop. \ref{Philsch}, Prop. \ref{HScomp} implies that
$$\| \int_{\Delta^k} A(u_0, \dots u_i,t)P_0\bigr)\bigl(P_0 B(u_{i+1}, \dots
u_k,t)~du_0 \dots du_k \|_{TR} $$
$$ \le C \|P_0\|^2_{HS}\int_{\Delta^k} \|A(u_0, \dots u_i,t)\| \| B(u_{i+1}, \dots
u_k,t)\| ~du_0 \dots du_k \ .$$
Let $\omega>0$ be such that there is $C>0$ with $\|P_1e^{-tD(\rho)^2}\| \le
Ce^{-\omega t}$ for all $t>0$. Then  
\begin{eqnarray*}
\lefteqn{ \|P^{\nu}_{jk}(t)\|_{TR}}\\
 &\le& C \int_0^1 du_0~ (u_0t)^{-\nu/2} e^{-\omega u_0 t}
 \int_{(1-u_0)\Delta^{k-1}} \exp(-\sum\limits_{i=1}^{j-1} \omega u_i t)  ~du_1 \dots du_k\\
&=&C\int_0^1 du_0~(u_0t)^{-\nu/2} e^{-\omega u_0 t} \int_0^{1-u_0} \!\!\! e^{-\omega s t} \vol(s\Delta^{j-2})\vol\bigl((1-u_0-s)\Delta^{k-j}\bigr)~ ds\\
&=& C \int_0^1 du_0~(u_0t)^{-\nu/2} e^{-\omega u_0 t} \int_0^{1-u_0}
 e^{-\omega s t}~ \frac{s^{j-2}(1-u_0-s)^{k-j}}{(j-2)!(k-j)!}ds\\
&=& C t^{-j} \int_0^t dy~ y^{-\nu/2}e^{-\omega y} \int_0^{t-y}  e^{-\omega x}~ \frac{x^{j-2}(1-y/t-x/t)^{k-j}}{(j-2)!(k-j)!}~dx\\
&\le& Ct^{-j}\Bigl( \int_0^{\infty}  y^{-\nu/2}e^{-\omega
 y}~dy \Bigr) \Bigl( \int_0^{\infty}  e^{-\omega x}~
 \frac{x^{j-2}}{(j-2)!(k-j)!}~dx \Bigr)\\
&\le &C t^{-j} \ .
\end{eqnarray*}
This shows the first statement.

For $k$ odd 
\begin{eqnarray*}
\lefteqn{\Trs \phi_t^2 P^{\nu}_{\frac{k+1}{2}k}(t)}\\
&=&\int_{\Delta^k} \Trs \phi_t^2 D(\rho)^{\nu}P_0\RM P_1e^{-u_1tD(\rho)^2}
\RM P_0  \dots P_0\RM P_1 e^{-u_ktD(\rho)^2}   ~du_0 \dots du_k \\
&+& \int_{\Delta^k} \Trs D(\rho)^{\nu} P_1e^{-u_0tD(\rho)^2} \RM P_0 \RM P_1
 \dots  e^{-u_{k-1}tD(\rho)^2} \RM P_0 \phi_t^2 ~du_0 \dots du_k
\ .
\end{eqnarray*}
By Prop. \ref{Philsch} we can estimate $\|(1-\phi_t^2)P_0\|_{HS}$ and $\|P_0(1-\phi_t^2)\|_{HS}$ by $Ct^{-n}$.
The second estimate follows then from the
cyclicity of the supertrace since $P_0P_1=P_1P_0=0$.
\end{proof}

From the estimates so far obtained the second assertion of
the theorem follows.

Furthermore the previous lemmas imply that for $t>T$ there is $C>0$
such that
$$|\Trs
\phi_t^2 e^{-A(\rho)_t^2}-\sum\limits_{n=0}^{\infty}  t^n
~\Trs \phi_t^2 P^0_{n(2n)}(t)|\le Ct^{-\frac 12} \ .$$ 
Hence it remains to study the behavior of $t^n P_{n(2n)}(t)$ for $t \to \infty$.

Recall that
\begin{eqnarray*}
P^0_{n(2n)}(t)&=&\int_{\Delta^{2n}}  P_0\RM P_1e^{-u_0tD(\rho)^2} \RM P_0 \RM P_1 e^{-u_1 tD(\rho)^2}
\RM P_0 \dots \\
&& \quad \dots P_0\RM P_1 e^{-u_{n-1}tD(\rho)^2}\RM P_0 ~du_0 \dots du_{2n} \ .
\end{eqnarray*}
By
Prop. \ref{Philsch} for any $j \in \bbbn$ there is $C>0$ such
that $$\|P_0-\phi_t^2 P_0\|_{HS} \le Ct^{-j} \ ,$$
hence $$|\Trs \phi_t^2 P^0_{n(2n)}(t)  -  \Trs P^0_{n(2n)}(t)) | \le Ct^{-j} \ .$$
In the next lemma we show that for $T>0$ there is $C>0$ such that
$$\|t^n P^0_{n(2n)}(t) -(-1)^n \frac{1}{n!} (P_0WdW^* P_0)^{2n}\| \le Ct^{-1}$$ in
$B(L^2(M,E \ten \Ol{\mu}))$ for $t >T$. This implies that there is $C>0$ such that for $t>T$  
\begin{eqnarray*}
\lefteqn{|t^n \Trs  P^0_{n(2n)}(t) -(-1)^n \frac{1}{n!} \Trs (P_0WdW^* P_0)^{2n}| \qquad \qquad}\\
&\le&  C\|P_0\|_{HS}^2\|t^n P^0_{n(2n)}(t)-(-1)^n \frac{1}{n!} (P_0WdW^*
P_0)^{2n}\|  \le C t^{-1}\ .
\end{eqnarray*}
Both estimates combined yield the first assertion of the theorem.

\begin{lem}
Let $k,n \in \bbbn_0$ with $n \le k$. For $t \to \infty$ the term 
$$t^n \int_{\Delta^k}P_0\RM P_1 e^{-u_0tD(\rho)^2}P_1\RM  P_0\RM P_1
e^{-u_1tD(\rho)^2} \RM P_0\dots $$
$$\dots P_0\RM P_1 e^{-u_{n-1}tD(\rho)^2} \RM P_0~du_0
\dots du_k$$
converges in $B(L^2(M,E \ten \Ol{\mu}))$  to
$$(-1)^n \frac{1}{(k-n)!}~ (P_0 W \di W^* P_0)^{2n} \ $$
with $O(t^{-1})$.
\end{lem}

\begin{proof} 
For $i,j \in \{0,1\}$ write  ${}_id_j:=P_i W \di W^* P_j$.
 
By $\RM =[W^*\di W,D(\rho)]_s$ 
\begin{eqnarray*}
P_1\RM P_0 &=&P_1D(\rho) W^*\di W P_0 \ ,\\
P_0\RM P_1 &=& P_0 W^*\di W D(\rho) P_1 \ ,
\end{eqnarray*} 
thus $$P_0\RM P_1 e^{-tD(\rho)^2} \RM P_0={}_0d_1 D(\rho)^2e^{-tD(\rho)^2} {}_1d_0 \ .$$ 
This term is uniformly bounded for $t \to 0$. Hence the integral $$\int_0^t {}_0d_1 D(\rho)^2e^{-sD(\rho)^2} {}_1d_0 ~ds$$ converges and equals
${}_0d_1 (e^{-tD(\rho)^2}-1)~ {}_1d_0.$

For $n \in \bbbn$ and $k \ge n$ 
$$t^n \int_{\Delta^k}P_0\RM P_1 e^{-u_0tD(\rho)^2}\RM  P_0\RM P_1
e^{-u_1tD(\rho)^2} \RM P_0\dots $$
$$\dots P_0\RM P_1 e^{-u_{n-1}tD(\rho)^2} \RM P_0~du_0
\dots du_k$$
$$=t^n \int_{\Delta^k} {}_0d_1D(\rho)^2e^{-u_0tD(\rho)^2} {}_1d_0{}_0d_1
D(\rho)^2e^{-u_1tD(\rho)^2} {}_1d_0\dots $$
$$\dots {}_0d_1
D(\rho)^2e^{-u_{n-1}tD(\rho)^2}{}_1d_0~   ~du_0 \dots du_k \ .$$
Set $$D_n:=\{(u_0,u_1, \dots u_{n-1}) \in \bbbr^n~|~ \sum\limits_{i=0}^{n-1}u_i \le 1;~0\le u_i \le
1,~i=0, \dots, n-1\} \ .$$
By integration on $u_n, \dots, u_k$ the previous term equals
$$(*) \qquad t^n\int_{D_n}
 \dots {}_0d_1
D(\rho)^2e^{-u_{n-1}tD(\rho)^2}{}_1d_0~ \vol((1-\sum\limits_{i=0}^{n-1}
u_i)\Delta^{k-n}) ~du_0 \dots du_{n-1} \ .$$
We claim that $(*)$ converges to $\frac{1}{(k-n)!}
({}_0d_1{}_1d_0)^{n}$ with $O(t^{-1})$.
Then it follows from $(W^*\di W)^2=0$  that ${}_0d_1{}_1d_0=-{}_0d_0^2$. This shows the
assertion of the lemma.

For $n=1$ and $k=1$ the term $(*)$ equals $1$. 

For $n=1$ and $k >1$ the claim follows
by induction since by partial integration  
\begin{eqnarray*}
\lefteqn{\int_0^1  {}_0d_1
D(\rho)^2e^{-u_0tD(\rho)^2} {}_1d_0  \vol((1-
u_0)\Delta^{k-1}) ~du_0} \\
&=&{}_0d_1
e^{-tD(\rho)^2} {}_1d_0+\int_0^1  {}_0d_1
D(\rho)^2e^{-u_0tD(\rho)^2} {}_1d_0  \vol((1-
u_0)\Delta^{k-2}) ~du_0
\end{eqnarray*}
because of
$$\partial_{u_0} \vol((1-u_0)\Delta^{k-1})=-\vol((1-u_0)\Delta^{k-2}) \ .$$
Furthermore ${}_0d_1
e^{-tD(\rho)^2} {}_1d_0$ decays exponentially for $t \to \infty$.

For general $k$ and $n$  the term $(*)$ equals, by partial integration on $u_{n-1}$, 
$$t^n\int_{\Delta^{n-1}}  ~  \dots {}_0d_1 \Bigl[-t^{-1} e^{-xtD(\rho)^2} \vol((1-\sum\limits_{i=0}^{n-2} u_i -x)\Delta^{k-n}) \Bigr]^{u_{n-1}}_0 {}_1d_0~ ~du_0 \dots du_{n-1}$$

$$+ \quad t^{n-1} \int_{D_n}   \dots {}_0d_1 e^{-u_{n-1}tD(\rho)^2} {}_1d_0~ \partial_{u_{n-1}} \vol((1-\sum\limits_{i=0}^{n-1} u_i)\Delta^{k-n}) ~du_0 \dots du_{n-1} \ .$$
Note that the first integral vanishes for $x=u_{n-1}$. 

We obtain
$$t^{n-1} \int_{D_{n-1}}  \dots {}_0d_1 D(\rho)^2e^{-u_{n-2} tD(\rho)^2} {}_1d_0{}_0d_1{}_1d_0  \vol((1-\sum\limits_{i=0}^{n-2} u_i)\Delta^{k-n}) ~du_0 \dots du_{n-2}
$$

$$- t^{n-1} \int_{D_{n-1}} \dots
{}_0d_1 e^{-u_{n-1}tD(\rho)^2} {}_1d_0  \vol((1-\sum\limits_{i=0}^{n-1}
u_i)\Delta^{k-n-1}) ~du_0 \dots du_{n-1} \ .$$
There are $C,\omega>0$ such that the last term is bounded by
$$Ct^{n-1}\int_0^1e^{-s \omega t}
\vol(s\Delta^{n-1})\vol((1-s)\Delta^{k-n-1})~ds \ ,$$ hence it vanishes  with     
$O(t^{-1})$ for $t \to \infty$ (by a calculation as the proof of the previous lemma).

By induction the first term converges with $O(t^{-1})$ to
$$\frac{1}{(k-1-(n-1))!} ({}_0d_1{}_1d_0)^{n}
=\frac{1}{(k-n)!}({}_0d_1{}_1d_0)^{n} \ . $$

\end{proof}

 \subsection{The limit of $\Trs e^{-A_t^2}$ for $t \to 0$}

Recall the definition of $\N$ from \S \ref{stoer}.

\begin{theorem}
\label{lim0}

\begin{enumerate}
\item
$\quad \lim\limits_{t \to 0} \Trs e^{-A_t^2} =\N \ .$
\item $\quad \lim\limits_{t \to 0} \Trs De^{-A_t^2}=0 \ .$
\end{enumerate}
\end{theorem}

\begin{proof}
(1) By Prop. \ref{hksupcon} 
$$\lim\limits_{t \to 0} \Trs e^{-A_t^2}= \sum\limits_{n=0}^{\infty}(-1)^n
\lim\limits_{t \to 0}~ t^{\frac
n2}~ \sum\limits_{k \in J}\int_{\U_k} \gamma_k(x)~\trs
w(0)^k_t(x,x)^n \phi_k(x) ~dx \ .$$ 
Now $\trs w(\rho)^k_t(x,x)^n=0$ for $k \in \bbbz/6$ and for all $n
\in \bbbn_0$ and $t > 0$ by \S \ref{supcyl}.

Furthermore $w(0)^{\cp}_t(x,y)^n=0$ for $n>0$. 

Recall that $\U_{\cp}$
contains the isolated point $*$. Since
$w(0)^{\cp}_t(x,y)^0$ is the integral kernel of $e^{-tD_N^2}$,
$$\lim\limits_{t \to 0} \trs w(0)^{\cp}_t(x,x)^0= \N1_{*}(x) $$ 
 in
$C(\U_{\cp})$  by the local index
theorem (\cite{bgv}, Th. 4.2) and by
 $\ch(E/S)=\ch((\bbbc^+)^d)+ \ch((\bbbc^-)^d)=0$.

(2) Prop. \ref{hksupcon} implies that
$$\lim\limits_{t \to 0} \Trs De^{-A_t^2}=\sum\limits_{n=0}^{\infty} (-1)^n \lim\limits_{t \to 0}~ t^{\frac
n2}~ \Trs D W_n(0,t) \ .$$
Since $DW_0(0,t)$ is an odd trace class operator, its supertrace vanishes. Furthermore $\Trs  D W_n(0,t)=0$  by \S \ref{supcyl}.
\end{proof}

\subsection{$\dfrac{d}{dt}\Trs e^{-A(\rho)_t^2}$ and $\dfrac{d}{d\rho}\Trs e^{-A(\rho)_t^2}$}
\label{deriv}

\begin{lem} 

\begin{enumerate}
\item
$$\frac{d}{dt}\Trs e^{-A(\rho)^2_t}=-\Trs
\frac{dA(\rho)_t^2}{dt}e^{-A(\rho)^2_t} \ .$$
\item
$$\frac{d}{d\rho}\Trs e^{-A(\rho)^2_t}=-\di\Trs
\frac{dA(\rho)_t}{d\rho}e^{-A(\rho)_t^2} \ .$$
\end{enumerate}
 \end{lem}

\begin{proof}
(1) First we calculate $\frac{d}{dt}e^{-A(\rho)^2_t}$.

Consider the holomorphic semigroup $e^{-t'(D(\rho)^2+z\RM )}$ depending on
the parameter $z$. By the semigroup law  
$$\frac{d}{dt'}~ e^{-t'(D(\rho)^2+z\RM )}= -(D(\rho)^2+z\RM )e^{-t'(D(\rho)^2+z\RM )}$$
and by Duhamel's formula (Prop. \ref{boundpert})
$$\frac{d}{dz}~e^{-t'(D(\rho)^2+z\RM )}=-\int_0^{t'} e^{-(t'-s)(D(\rho)^2+z\RM )}\RM e^{-s(D(\rho)^2+z\RM )}
ds \ .$$
It follows that
\begin{eqnarray*}
\frac{d}{dt}~e^{-A(\rho)_t^2}
&=&\frac{d}{dt}~e^{-t(D(\rho)^2+ t^{-1/2} \RM )}\\
&=& \frac{d}{dt'}~e^{-t'(D(\rho)^2+t^{-1/2}\RM )}|_{(t'=t)}
-\frac 12 t^{-3/2} \frac{d}{dz}~e^{-t(D(\rho)^2+z\RM )}|_{(z=t^{-1/2})}\\ 
&=&-(D(\rho)^2+t^{-1/2}\RM )e^{-t(D(\rho)^2  +t^{-1/2}\RM )}\\
&&+~\frac 12
t^{-3/2}\int_0^{t} e^{-(t-s)(D(\rho)^2+t^{-
1/2}\RM )}\RM e^{-s(D(\rho)^2+t^{-1/2}\RM )} ds\\
&=&-t^{-1} A(\rho)_t^2e^{-A(\rho)_t^2} + \frac 12~ t^{-1/2} \int_0^1
e^{-(1-s)A(\rho)_t^2}\RM  e^{-sA(\rho)_t^2}ds \ .
\end{eqnarray*}
Using this equation we prove that   
$$(*_1) \qquad |\frac{d}{dt}~\trs p(\rho)_t(x,x)^n| \le Ce^{-cd(x,M_r)^2}$$ for all $x$ uniformly in $t$ for $t \in [t_1,t_2]$ with $t_1,t_2>0$. This yields that
$$\frac{d}{dt}~\Trs e^{-A(\rho)^2_t}=\Trs \frac{d}{dt}e^{-A(\rho)^2_t} \ .$$
By Cor. \ref{compcyl} and the fact that the pointwise
supertrace of the integral kernel $(A^{Z_k}_t)^2e^{-(A^{Z_k}_t)^2}$ vanishes for
$k \in \bbbz/6$ the pointwise supertrace of the integral kernel of
$A(\rho)_t^2e^{-A(\rho)_t^2}$ can be estimated by  $Ce^{-cd(x,M_r)^2}$ uniformly in $t$ on compact subsets of $(0,\infty)$.

The integral kernel of
$\int_0^1
e^{-(1-s)A(\rho)_t^2}\RM  e^{-sA(\rho)_t^2}ds$ is the sum over $m,n \in \bbbn_0$ of
the terms
$$(-1)^{m+n}t^{\frac{m+n}{2}} \int_0^1(1-s)^{m/2}s^{n/2}  \int_M p(\rho)_{(1-s)t}(x,y)^m \RM (y) p(\rho)_{st}(y,x)^n~dy
ds \ .$$ 
We decompose 
\begin{eqnarray*}
\lefteqn{p(\rho)_{(1-s)t}(x,y)^m\RM (y) p(\rho)_{st}(y,x)^n}\\
&=&p(\rho)_{(1-s)t}(x,y)^m \RM (y) \bigl(p(\rho)_{st}(y,x)^n- w(\rho)_{st}(y,x)^n\bigr) \\ 
&& + 
\bigl(p(\rho)_{(1-s)t}(x,y)^m- w(\rho)_{(1-s)t}(x,y)^m\bigr)\RM (y) w(\rho)_{st}(y,x)^n \\
&&+ w(\rho)_{(1-s)t}(x,y)^m \RM (y) w(\rho)_{st}(y,x)^n \ .
\end{eqnarray*}
By Prop. \ref{hksupcon} and the fact that  the operator $I_m(\rho,(1-s)t)$ is uniformly bounded on
$C_R(M,E \ten \Ol{\mu})$ for $0 \le (1-s)t \le t_2$ there are $C,c,r>0$ such that  for all $x \in
M$ and $0\le s \le 1$ and $t_1\le t\le t_2$
$$|\int_M  p(\rho)_{(1-s)t}(x,y)^m \RM (y) \bigl(p(\rho)_{st}(y,x)^n- w(\rho)_{st}(y,x)^n\bigr) dy| \le  Ce^{-cd(x,M_r)^2} \ .$$
An analogous estimate holds  for the second term.

By \S \ref{supcyl} and since $\RM|_{\U_{\cp}}=0$
 $$\trs w(\rho)_{(1-s)t}(x,y)^m \RM(y)
w(\rho)_{st}(y,x)^n =0 \ .$$
We conclude that there are $r,c,C>0$ such that for $x \in M_r$, $0\le s \le 1$ and $t_1 \le
t\le t_2$ 
$$(*_2) \qquad |\trs \int_M p(\rho)_{(1-s)t}(x,y)^m \RM(y) p(\rho)_{st}(y,x)^n ~dy|  \le  Ce^{-cd(x,M_r)^2}\ .$$
Now $(*_1)$ follows.

The next step is to show that
$$\Trs \int_0^1
e^{-(1-s)A(\rho)_t^2}\RM e^{-sA(\rho)_t^2}ds =\Trs \RM e^{-A(\rho)_t^2}$$
or equivalently 
$$\int_M \int_0^1(1-s)^{m/2}s^{n/2}  \int_M
\trs p(\rho)_{(1-s)t}(x,y)^m \RM(y) p(\rho)_{st}(y,x)^n~dy
dsdx$$
$$= \int_0^1 (1-s)^{m/2}s^{n/2} \int_M \int_M
\trs p(\rho)_{(1-s)t}(x,y)^m \RM(y) p(\rho)_{st}(y,x)^n~dx
dyds \ .$$
We can interchange the integration over $s$ and $x$  by
the estimate $(*_2)$.

Fix $s$ and $t$. Consider once more the decomposition of $$p(\rho)_{(1-s)t}(x,y)^m \RM(y) p(\rho)_{st}(y,x)^n$$ from above. From Cor. \ref{hksupconMgauss} it follows that for $\ve>0$ there are $r,c,C>0$
such that for all $x,y \in M$ 
\begin{eqnarray*}
\lefteqn{|p(\rho)_{(1-s)t}(x,y)^m \RM(y) \bigl(p(\rho)_{st}(y,x)^n- w(\rho)_{st}(y,x)^n\bigr)|} \\
&\le& C\bigl(e^{-cd(x,y)^2} + 
e^{-cd(x,M_r)^2} \bigr)
e^{-cd(y,M_r)^2} \ .
\end{eqnarray*}
The second term can be estimated in an analogous manner and the supertrace of
the third term vanishes as we saw.

Hence we can interchange
$dx$ and $dy$.

This shows (1) by
\begin{eqnarray*}
\Trs \frac{d}{dt}~e^{-A(\rho)_t^2}&=&\Trs(-t^{-1}
A(\rho)_t^2 + \frac 12 t^{- 1/2}\RM) e^{-A(\rho)_t^2}\\
&=&-\Trs
\frac{dA(\rho)_t^2}{dt}~e^{-A(\rho)^2_t} \ .
\end{eqnarray*}

(2) As above, since $A(\rho)_t^2=t(D^2+ \rho[D,K]_s + \rho^2K^2) +\sqrt t \RM$, Duhamel's formula (Prop. \ref{boundpert}) and the chain rule 
 imply that 
$$\frac {d}{d\rho} ~e^{-A(\rho)_t^2}=-\int_0^1
e^{-(1-s)A(\rho)_t^2}\frac{dA(\rho)_t^2}{d\rho} ~e^{-sA(\rho)_t^2} ~ds \ .$$
Note that $\frac{dA(\rho)_t^2}{d\rho}$ is a trace class operator, hence we conclude immediately that
$$\frac {d}{d\rho}~ \Trs e^{-A(\rho)_t^2}=\Trs \frac {d}{d\rho}~
e^{-A(\rho)_t^2}=-\Trs \frac{dA(\rho)_t^2}{d\rho} ~e^{-A(\rho)_t^2}\
.$$
Since $A(\rho)_t$ is a bounded perturbation of $\sqrt t WD_sW^*$, we have, by 
Cor. \ref{Dsboundinv} and Prop. \ref{sgpert}, that $A(\rho)_t$ and
$e^{-A(\rho)_t^2}$ commute.

It follows that
\begin{eqnarray*}
 \frac{dA(\rho)_t^2}{d\rho}
e^{-A(\rho)_t^2} &=& [A(\rho)_t, \frac{dA(\rho)_t}{d\rho}e^{-A(\rho)_t^2}]_s\\
&=& [W^*\di W, \frac{dA(\rho)_t}{d\rho}e^{-A(\rho)_t^2}]_s+\sqrt t [D(\rho),
\frac{dA(\rho)_t}{d\rho}e^{-A(\rho)_t^2}]_s \ .
\end{eqnarray*}
The first term of the last line
is an integral operator with integral kernel
$W^*d(k)W$ where $k$ is the integral kernel of
$W\frac{dA(\rho)_t}{d\rho}e^{-A(\rho)_t^2}W^*$. Hence the supertrace of
the first term equals $$\di\Trs  \frac{dA(\rho)_t}{d\rho}e^{-A(\rho)^2_t} \
.$$
Now consider the second term. Let $P$ be the orthogonal
projection onto the range of $K$. It is an even Hilbert-Schmidt operator with smooth complex integral kernel. Since
$\frac{dA(\rho)_t}{d\rho}=\sqrt t K$,
$$\Trs [D(\rho),
\frac{dA(\rho)_t}{d\rho}e^{-A(\rho)_t^2}]_s=\Trs \sqrt t [D(\rho)P,Ke^{-A(\rho)_t^2}]_s
\ .$$
Since $D(\rho)P$ and $Ke^{-A(\rho)_t^2}$ are Hilbert-Schmidt operators, the
supertrace vanishes.
\end{proof}

Let
$D_{I_k}$  be the operator $D_I$ from \S \ref{Igen} with boundary conditions
given by the pair $(\Pj_{k \Mod 3},\Pj_{k+1 \Mod 3})$.

In \S \ref{supconMhk} we defined $A^{Z_k}_t=\tilde W_k^* \di \tilde W_k +\sqrt
t D_{Z_k}$ such that $A_t^{Z_k}$ coincides with $A(\rho)_t$ on $\U_k$. There is $U_k \in \C([0,1],M_{2d}(\Ai))$ such
that  $\tilde W_k(x_1,x_2) =U_k(x_2) \oplus U_k(x_2)$. Let $A^{I_k}_t=U_k^* \di
U_k +\sqrt t \sigma D_{I_k}$.

\begin{lem}
$$\frac{d}{dt} \Trs e^{-A(\rho)^2_t}=\frac{1}{\sqrt {4\pi t}} \sum\limits_{k \in \bbbz/6}  \Trsi D_{I_k} e^{-(A^{I_k}_t)^2} 
-\frac{1}{2 \sqrt t} \di \Trs D(\rho)e^{-A(\rho)_t^2} \ .$$
\end{lem}

\begin{proof}
Since $\frac{dA(\rho)^2_t}{dt}=[\frac{dA(\rho)_t}{dt},A(\rho)_t]_s ,$ it follows
from the previous lemma that
$$\frac{d}{dt} \Trs e^{-A(\rho)_t^2}=-\Trs
[\frac{dA(\rho)_t}{dt},A(\rho)_t]_se^{-A(\rho)_t^2} \ .$$
Furthermore 
\begin{eqnarray*}
\lefteqn{\!\!\!\!\!\!\!\!\!\!\!\!\!\!\!\!\!\!\!\!\!\!\!\! -[\frac{dA(\rho)_t}{dt},A(\rho)_t]_se^{-A(\rho)_t^2}}\\
&=& -[A(\rho)_t,\frac{dA(\rho)_t}{dt}e^{-A(\rho)_t^2}]_s \\
&=&-\frac {1}{2 \sqrt t}[W^*\di W,D(\rho)e^{-A(\rho)_t^2}]_s- \frac
12[D(\rho),D(\rho)e^{-A(\rho)_t^2}]_s \ .
\end{eqnarray*}
The supertrace of the first supercommutator in the last line equals
$$-\frac {1}{2 \sqrt t}\di \Trs D(\rho)e^{-A(\rho)_t^2} \ .$$
Now consider the second term.

We have that
\begin{eqnarray*}
\lefteqn{-\frac 12\Trs [D(\rho),D(\rho)e^{-A(\rho)_t^2}]_s}\\
&=& -\frac 12 \lim\limits_{r \to \infty}\left(\Trs  \phi_r
D(\rho)^2e^{-A(\rho)_t^2} \phi_r + \Trs \phi_r
D(\rho)e^{-A(\rho)_t^2}D(\rho) \phi_r \right) \ .
\end{eqnarray*}
Since  for $\nu \in \bbbn_0$ the operators $\phi_r
D(\rho)^{\nu}e^{-A(\rho)_t^2/2}$ and $e^{-A(\rho)_t^2/2}D(\rho)^{\nu}
\phi_r$ are  Hilbert-Schmidt operators,  it follows that 
\begin{eqnarray*}
\lefteqn{- \frac
12\Trs [D(\rho),D(\rho)e^{-A(\rho)_t^2}]_s}\\
&=& \!\!\!\!\! - \frac
12\lim\limits_{r \to \infty}( \Trs e^{-A(\rho)_t^2/2}\phi_r^2
D(\rho)^2e^{-A(\rho)_t^2/2}- \Trs e^{-A(\rho)_t^2/2} D(\rho) \phi_r^2
D(\rho)e^{-A(\rho)_t^2/2}) \\
&=& \frac
12\lim\limits_{r \to \infty} \Trs e^{-A(\rho)_t^2/2} c(d\phi_r^2) D(\rho)e^{-A(\rho)_t^2/2}\\
&=& \frac
12\lim\limits_{r \to \infty} \Trs c(d\phi_r^2) D(\rho)e^{-A(\rho)_t^2} \ .
\end{eqnarray*}
For $r>0$ we define the function
$\chi_r: Z \to \bbbr,~\chi_r(x):=\chi^2(x_1-r)$ where $\chi$ is the function from the beginning of \S \ref{suptr}.  Cor. \ref{compcyl} implies that
$$\frac 12 \lim_{r \to \infty}\Trs c(d\phi_r^2) D(\rho)e^{-A(\rho)_t^2}=
\frac 12 \lim_{r \to \infty}\sum_{k\in \bbbz/6} \Trs c(d\chi_r) D_{Z_k}e^{-(A^{Z_k}_t)^2}\
\ .$$ 
Recall from \S \ref{supcyl} that the integral
kernel of $e^{-(A^{Z_k}_t)^2}$ is $\sum\limits_{n=0}^{\infty}(-1)^nt^{\frac n2}
p_t^{Z_k}(x,y)^n$ with $$p_t^{Z_k}(x,y)^n= c(dx_1)^n\frac{1}{\sqrt{4 \pi
t}}~e^{-\frac{(x_1-y_1)^2}{4t}}
\bigl(p^{I_k}_t(x_2,y_2)^n\oplus (-1)^np^{I_k}_t(x_2,y_2)^n\bigr) \ .$$ 
The integral kernel of $c(d\chi_r) D_{Z_k}e^{-(A^{Z_k}_t)^2}$ is 
$$\sum_{n=0}^{\infty}(-1)^{n+1} t^{n/2}  \chi'(x_1-r)
(\ra_{x_1}+I\ra_{x_2})p_t^{Z_k}(x,y)^n \ .$$
It follows that
\begin{eqnarray*}
\lefteqn{\Trs c(d\chi_r) D_{Z_k}e^{-(A^{Z_k}_t)^2}}\\
&=&\frac{1}{\sqrt{4 \pi t}} \sum_{n=0}^{\infty}(-1)^nt^{n/2}   \int_0^1 \trs ~c(dx_1)^n
\bigl((D_{I_k}p^{I_k}_t)(x_2,x_2)^n  \oplus \\
&& \qquad \oplus 
(-1)^{n+1}(D_{I_k}p^{I_k}_t)(x_2,x_2)^n)~ dx_2  \ .
\end{eqnarray*}
Comparison with Def. \ref{defeta} and the subsequent remark yields that
$$\trs c(dx_1)^n
\bigl((D_{I_k}p^{I_k}_t)(x_2,x_2)^n \oplus 
(-1)^{n+1}(D_{I_k}p^{I_k}_t)(x_2,x_2)^n\bigr)$$
$$=2i^n\trsi \sigma^n
(D_{I_k}p^{I_k}_t)(x_2,x_2)^n \ .$$ 
It follows that
\begin{eqnarray*}
\Trs c(d\chi_r) D_{Z_k}e^{-(A^{Z_k}_t)^2}
&=&\frac{2}{\sqrt{4 \pi t}} \sum_{n=0}^{\infty}(-1)^n t^n   \int_0^1 \tr~ 
(D_{I_k}p^{I_k}_t)(x_2,x_2)^{2n} dx_2  \\
&=& \frac{1}{\sqrt{\pi t}}  \Trsi D_{I_k}e^{-(A^{I_k}_t)^2} \ ,
\end{eqnarray*}    
hence
$$- \frac
12\Trs [D(\rho),D(\rho)e^{-A(\rho)_t^2}]_s= \sum_{k \in \bbbz/6}
\frac{1}{\sqrt{4\pi t}}  \Trsi D_{I_k}e^{-(A^{I_k}_t)^2} \ .$$
\end{proof}

\subsection{The index theorem}
\label{indtheorem}
In this section the notation is as in Prop. \ref{propstoer}.

Let $A_{I_k}$ be as in the previous section.

\begin{theorem}
\label{indth}  
$$\ch(\ind D^+)=\frac 12[\sum_{k \in \bbbz/6} \eta(A_{I_k})] \in H_{ev}^{dR}(\Ai) \ .$$

\end{theorem}

Here we understand $\ind D^+$ as a class in $K_0(\Ai)$ via the isomorphism
$K_0(\A) \cong K_0(\Ai)$ induced by the
injection $\Ai \incl \A$.

\begin{proof}
Let $\rho \neq 0$.

By Prop. \ref{propstoer}
$$\ind D^+= [\Ker D(\rho)^2]-[\A^{\N}]$$
in $K_0(\A)$.
Let $P_0$ be the projection onto the kernel of $D(\rho)^2$.

From Prop. \ref{Philsch} and Prop. \ref{projker} it follows that 
$$\ind D^+=
[\Ran_{\infty} P_0]-[\Ai^{\N}]= [\Ran_{\infty} WP_0W^*]-[\Ai^{\N}]$$
in $K_0(\Ai)$,
hence in $H_*^{dR}(\Ai)$ 
$$\ch(\ind D^+)=\ch[\Ran_{\infty} WP_0W^*] - \N  .$$
By Prop. \ref{projker} in $H_*^{dR}(\Ai)$ 
\begin{eqnarray*}
\ch[\Ran_{\infty} WP_0W^*] &=& \sum\limits_{n=0}^{\infty}(-1)^n
\frac 1{n!} \Trs (WP_0W^*\di WP_0W^*)^{2n} \\
&=& \sum\limits_{n=0}^{\infty}(-1)^n
\frac 1{n!} \Trs (P_0W^*\di WP_0)^{2n} \ . 
\end{eqnarray*}
Let $T>0$. By Th. \ref{liminfty} and Th. \ref{lim0} in $\Oi\Ai/\ov{[\Oi\Ai,\Oi\Ai]_s}$
\begin{eqnarray*}
 \sum\limits_{n=0}^{\infty}(-1)^n
\frac {1}{n!} \Trs (P_0W^*\di WP_0)^{2n}-\N 
 &=& \lim\limits_{t \to \infty} \Trs
e^{-A(\rho)_t^2}-\lim\limits_{t \to 0} \Trs e^{-A_t^2}\\
&=& \int_T^{\infty}\frac {d}{dt}  \Trs e^{-A(\rho)_t^2} dt \\
&&
+ \int_0^{\rho}\frac {d}{d\rho'}  \Trs e^{-A(\rho')_T^2} d\rho' \\
&& + \int_0^T \frac {d}{dt} \Trs e^{-A_t^2}dt \ .
\end{eqnarray*}
By the results of \S \ref{deriv}, the estimate of Th.
\ref{liminfty} and Th. \ref{lim0} the integrals converge and we have:
\begin{eqnarray*}
\sum\limits_{n=0}^{\infty}(-1)^n
\frac {1}{n!} \Trs (P_0W^*\di WP_0)^{2n}-\N  &=&\sum\limits_{k \in \bbbz/6} \int_0^{\infty} \frac{1}{\sqrt {4\pi
t}} \Trsi D_{I_k} e^{-(A^{I_k}_t)^2} dt \\ 
&&- \di\int_T^{\infty} \frac{1}{2 \sqrt t}  \Trs D(\rho)e^{-A(\rho)_t^2}dt\\
&&-\di \int_0^{\rho} \Trs
\frac{dA(\rho')_T}{d\rho'}e^{-A(\rho')_T^2}d \rho' \\
&& - \di \int_0^T \frac{1}{2 \sqrt t} \Trs D(0)e^{-A_t^2}dt \ .\\
\end{eqnarray*}
The assertion follows. 
\end{proof}

\begin{cor}
\label{mainres}
In $H_*^{dR}(\Ai)$
$$\ch \tau(\Pj_0,\Pj_1,\Pj_2)=[\eta(\Pj_0,\Pj_1)+ \eta(\Pj_1,\Pj_2) + \eta(\Pj_2,\Pj_0)] \ .$$
\end{cor}

\begin{proof}
For $\rho>0$ we define a family $A(\rho,b)$ of superconnections associated to
$D(\rho)$.
For $0<b<b_0$ let $W(b) \in \C(M, \End^+ E \ten \Ai)$ be as in \S \ref{supconM} such that $W(b)$ is parallel on $\{x \in M~|~d(x,\ra M) \ge b\}$.  As in
\S \ref{supconM} we set $$A(\rho,b):=W(b)^*\di W(b)+D(\rho) \ .$$
For every $k \in \bbbz/6$ this induces a family of superconnections
$A_{I_k}(b)$. By Prop. \ref{limeta} and the subsequent definition
$$\lim\limits_{b \to 0} \eta(A_{I_k}(b))=\eta(\Pj_{k \Mod 3},\Pj_{k+1 \Mod 3}) \ .$$
The assertion follows now by the previous theorem and Prop. \ref{indeqmasl}.
\end{proof} 

\section{A gluing formula for $\eta$-forms on $S^1$}
\label{glue}

In the following we sketch a reinterpretation of Cor. \ref{mainres} as a gluing formula for $\eta$-forms on $S^1$ (which we identify with $\bbbr/\bbbz$ as a Riemannian manifold).

Given $u \in U(\Ai^n)$ we define a projective system of vector bundles on $S^1$:
$${\mathcal L}_i(u)=([0,1] \times \A_i^n)/\! \sim$$ with $(0,v) \sim (1,uv)$. The standard $\A$-valued scalar product on $\A^n$ induces a hermitian $\A$-valued metric on ${\mathcal L}(u):={\mathcal L}_0(u)$. We identify a smooth section of ${\mathcal L}_i(u)$ with a smooth function $f:\bbbr \to \A_i^n$ satisfying $f(x+1)=uf(x)$. The trivial connection $f \mapsto f' dx$ on $\bbbr \times \A^n$ induces a connection on ${\mathcal L}(u)$ preserving the metric. The associated Dirac operator is denoted by $\dira_{{\mathcal L}(u)}$ (as are its closures on the spaces $L^2(S^1,{\mathcal L}_i(u))$ in the following). 

Assume now $1 \notin \sigma(u)$. Then there is a path $w:[0,1] \to U(\Ai^n)$ with $w(0)=1,~w(1)=u$, and the map $f \mapsto w^*f$ defines a trivialization of ${\mathcal L}_i(u)$. Hence $\dira_{{\mathcal L}(u)}$ can be identified with the operator $i\ra + w^*w'$ on the trivial bundle. One can easily deduce that $-\dira_{{\mathcal L}(u)}^2$ generates a holomorphic semigroup on $L^2(S^1,{\mathcal L}_i(u))$ with integral kernel in $\C(S^1 \times S^1,M_n(\Ai))$. Since $\dira_{{\mathcal L}(u)}^2$ is invertible on $L^2(S^1,{\mathcal L}_i(u))$, the integral kernel vanishes exponentially for $t \to \infty$. Locally there is a trivialization of ${\mathcal L}(u)$ with respect to which $\dira_{{\mathcal L}(u)}$ equals $i\ra$. By Duhamel's principle it follows that $\Tr\dira_{{\mathcal L}(u)}e^{-t(\dira_{{\mathcal L}(u)})^2}$ converges to zero in $\Ai/\ov{[\Ai,\Ai]}$ for $t \to 0$.

We define a superconnection $A:=w\di w^* + \dira_{{\mathcal L}(u)}$ and the associated rescaled superconnection $A_t:=w\di w^* + \sqrt t \dira_{{\mathcal L}(u)}$. Then $e^{-A_t^2}$ is a well-defined integral operator with smooth integral kernel. Furthermore $\Tr \dira_{{\mathcal L}(u)}e^{-A_t^2}$ has a limit for $t \to 0$, and
the $\eta$-form 
$$\eta(A_w):= \frac{1}{\sqrt{\pi}} \int_0^{\infty}t^{-1/2}\Tr \dira_{{\mathcal L}(u)}e^{-A_t^2} ~dt \in \Oi\Ai/\ov{[\Oi\Ai,\Oi\Ai]_s}$$ is well-defined. As in Def. \ref{defetatwo} we can eliminate the dependence of the trivialization $w$ by defining $$\eta(\dira_{{\mathcal L}(u)})=\lim\limits_{\ve \to 0} \eta(A_{w_{\ve}})$$ where $w_{\ve}$ is a family of trivializations with $w_{\ve}' \subset [0,\ve] \cup [1-\ve,1]$. The existence of the limit follows from the proof of the next proposition.

In the following let $$P(u):=\frac 12 \left(\begin{array}{cc} 1&u^* \\u & 1\end{array}\right) \ .$$ 

\begin{prop} Let $u_0, u_1$ be unitaries such that $u_0-u_1$ is invertible. Then up to exact forms 
$$\eta(\dira_{{\mathcal L}(u_0^*u_1)})= \eta(P(u_0),P(u_1)) \ .$$
\end{prop}

\begin{proof}
For a unitary $U \in M_n(\Ai)$ 
$$\eta(U^*P(u_0)U,U^*P(u_1)U)=\eta(P(u_0),P(u_1)) \ .$$
Since for $U:=\left(\begin{array}{cc} 1 & 0 \\ 0 & u_0 \end{array}\right)$ we have that $U^*P(u_0)U=P(1)$ and 
$U^*P(u_1)U =P(u_0^*u_1)$, we may assume $u_0=1$. Let $u:=u_1$.

Let $V:\bbbr \to U(\Ai^n)$ be a smooth path with $V(x)=u$ for $x \le 0$ and $V(x)=-1$ for $x \ge 1$ and with $1 \notin \sigma(V)$; we define the bundle $$L_i(V) :=(\bbbr \times [0,1]) \times \A_i^n/\! \sim$$ where $(x_1,0,v) \sim (x_1,1,V(x_1)v)$. The smooth sections of $L_i(V)$ can be identified with smooth functions $f:\bbbr \times \bbbr \to \A^n$ fulfilling $f(x_1,x_2+1)=V(x_1)f(x_1,x_2)$.  Let $W:\bbbr \times [0,1] \to U(\Ai^n)$ be smooth with $W(x_1,0)=1,~W(x_1,1)=V(x_1)$ and such that $W(x_1,x_2)$ is independent of $x_2$ in a small neighborhood of $0$ resp. $1$ and independent of $x_1$ for $x_1<0$ resp. $x_1>1$. Then $f \mapsto W^*f$ is a bundle isomorphism between $L_i(V)$ and the trivial bundle $(\bbbr \times S^1) \times \A_i^n$. Since $dx_1 \ra_1 + W^*dx_2 \ra_2W= dx_1\ra_1 + dx_2 \ra_2 + dx_2 W^* (\ra_2 W)$ is a connection preserving the metric on the trivial bundle, the operator $W dx_1\ra_1W^* + dx_2\ra_2$ is a connection preserving the metric on $L_i(V)$. It agrees with the trivial connection for $|x_1|$ large. 

Let $\dira_{L_i(V)}$ be the Dirac operator associated to the connection $W dx_1\ra_1W^* + dx_2\ra_2$ on the bundle $L_i(V)$. The index of $\dira^+_{L_i(V)}$ vanishes, which can be seen as follows. For $t \in [0,1]$ let $W_t:\bbbr \times [0,1] \to U(\Ai^n)$ be a homotopy with $W_0=1$ and $W_1=W$ and independent of $x_1$ on the complement of a compact set.  Then $dx_1 \ra_1 + W_t^*dx_2 \ra_2W_t$ is a homotopy of connections on the trivial bundle interpolating between the trivial connection and the connection $dx_1 \ra_1 + W^*dx_2 \ra_2W$. The index of the Dirac operator associated to the trivial connection vanishes, hence the index of the Dirac operator associated to $dx_1 \ra_1 + W^*dx_2 \ra_2W$ vanishes as well. 

The local term in the index theorem is determined by the superconnection $$B_1:=\di + W\dira_{L_i(V)}W^*=\di + c(dx_1)\ra_1 + c(dx_2) \ra_2 + c(dx_2)  W(\ra_2 W^*)$$ on the trivial bundle $S \ten \A^n$ on $\bbbr \times S^1$. The contribution from the cylindric ends is given by $\frac 12(\eta(\dira_{{\mathcal L}(-1)})-\eta(A))$.

Compare this with the index theorem for the Dirac operator $\dirz$ on the sections of $(\bbbr \times [0,1]) \times ((\A^+)^{2n} \oplus (\A^-)^{2n})$
with boundary conditions defined by $P(1) \oplus P(1)$ and $P(V(x_1)) \oplus P(V(x_1))$. Let $$\tilde W:=\diag(1,W,1,W) \in M_{4n}(\Ai) \ .$$ 
Then $$\tilde W^* \dirz\tilde W= \dirz+ \tilde W^*c(dx_1) (\ra_1 \tilde W) + \tilde W^*c(dx_2)(\ra_2 \tilde W) \ ,$$  and the boundary conditions transform to $P(1) \oplus P(1)$ and $P(-1) \oplus P(-1)$.

An argument analogous to the one above shows that the index of $\tilde W^* \dirz^+ \tilde W$ vanishes. Since $\tilde W^*c(dx_1)(\ra_1 \tilde W)$ is compactly supported, the index of $\tilde W^* \dirz \tilde W$ equals the index of $\dirz + \tilde W^*c(dx_2)(\ra_2 \tilde W)$, which is the Dirac operator associated to the connection  $dx_1\ra_1 + dx_2 \ra_2 + dx_2 \tilde W^* (\ra_2 \tilde W)$. Hence we can use the superconnection $B_2:=\di + c(dx_1)\ra_1 + c(dx_2) \ra_2 + c(dx_2) \tilde W^* (\ra_2 \tilde W)$ for the index theorem, and it follows that the local term here equals the local term in the index theorem for $\dira_{L(V)}$. Since the proposition is true for $u=-1$ by \cite{clm}, Prop. 6.3, the assertion follows. 
\end{proof}

\begin{cor} If $u_0,u_1,u_2  \in M_n(\Ai)$ are unitaries such that $u_i-u_j$ is invertible for $i \neq j$, then
$$\ch \tau(P(u_0),P(u_1),P(u_2)) = [\eta(\dira_{{\mathcal L}(u_0^*u_1)})+ \eta(\dira_{{\mathcal L}(u_1^*u_2)}) + \eta(\dira_{{\mathcal L}(u_2^*u_0)}]  \in H_*^{dR}(\Ai) \ .$$
\end{cor}
 
Assume that $P_i-P_j$ and $P_i-(1-P_j)$ are invertible for $i \neq j$ and define 
\begin{eqnarray*}
\tau_I(P_0,P_1,P_2) &:=& \tau(P_0,P_1,P_2)+ \tau(P_0,1-P_1,P_1)\\
&& + ~\tau(P_1,1-P_2,P_2) + \tau(P_2,1-P_0,P_0) \in K_0(\A)\ .
\end{eqnarray*}
Then from the corollary, the previous formula and the fact that $1-P(u)=P(-u)$ it follows after some straightforward calculations:

\begin{cor} In $H_*^{dR}(\Ai)$
$$\ch \tau_I(P(u_0),P(u_1),P(u_2))= [\eta(\dira_{{\mathcal L}(-u_0^*u_1)})+ \eta(\dira_{{\mathcal L}(-u_1^*u_2)}) + \eta(\dira_{{\mathcal L}(-u_2^*u_0)}]  \ .$$
\end{cor}

In the following we explain how this formula is related to a gluing formula for $\eta$-invariants (\cite{b}, \S 1.7 and Cor. 1.20).

Consider the Clifford bundle $\A^n \times [0,\pi]$ on the manifold with boundary $[0,\pi]$ with Clifford multiplication $i$ and Dirac operator $i\ra$. Its restriction to the boundary is $\A^n \times (\{0\} \cup \{\pi\})$ and the Clifford multiplication by the inward pointing normal vector is given by the bundle morphism $I_0=i \cup (-i)$. Via the standard $\A$-valued scalar product it induces a skew-hermitian form on the bundle $\A^n \times (\{0\} \cup \{\pi\})$. We identify the sections of this bundle with $\A^{2n}$. The image of the kernel of the Dirac operator with respect to restriction to the boundary is the Lagrangian submodule $L=\{(x,x)~|~ x \in \A^n\}$.

Now let
$$\phi_j:=(1 \cup (u_j)): \A^n \times (\{0\} \cup \{\pi\}) \to  \A^n \times (* \cup *) \ .$$ Then $\phi_j(L)$ is the range of the projection $P(u_j)$. 
   
The bundle ${\mathcal L}(-u_j^*u_{j+1})$ on $S^1$ can be obtained by gluing
$$(\A^n \times [0,\pi]) \cup_{\phi_j^{-1}\circ I_0 \circ \phi_{j+1}} (\A^n \times [0,\pi]) \ ;$$
and  by gluing the operator $i\ra$ on one copy of $[0,\pi]$ with $-i \ra$ on the second copy one obtains $\dira_{{\mathcal L}(-u_0^*u_1)}$.

From the proposition it follows that up to exact forms $$\eta(\dira_{{\mathcal L}(-u_j^*u_{j+1})})=\eta(P(u_j),1-P(u_{j+1})) \ .$$
Since the projection $\phi_i(L)$ in $\A^{2n}$ is given by the projection $P(u_j)$, we have  in $H_*^{dR}(\Ai)$:
$$\ch \tau_I(\phi_0(L),\phi_1(L),\phi_2(L)) = [\eta(\dira_{{\mathcal L}(-u_0^*u_1)})+ \eta(\dira_{{\mathcal L}(-u_1^*u_2)}) + \eta(\dira_{{\mathcal L}(-u_2^*u_0)}] \ .$$
For $\A=\bbbc$ this is a particular case of the well-known general gluing formula.

\chapter{Definitions and Techniques}

\section{Hilbert $C^*$-modules}

\subsection{Bounded operators}
\label{HCmod}

Let $\A$ be a unital $C^*$-algebra with norm $| ~\cdot~|$.

In order to fix notation we recall some facts about Hilbert $\A$-modules. References are \cite{la} and \cite{wo}.

\begin{ddd}
A {\sc pre-Hilbert $\A$-module} is a right $\A$-module $H$ with an
$\A$-valued scalar product  $\langle ~,~\rangle :H \times H \to \A$; i.e.
\begin{enumerate}
\item $\langle ~,~\rangle $ is $\A$-linear in the second variable,
\item $\langle x,y\rangle =\langle y,x\rangle ^*$ for all $x,y \in H$,
\item $\langle x,x\rangle  \ge 0$ for all $x \in H$,
\item if $\langle x,x\rangle =0$, then $x=0$.
\end{enumerate}
If $H$ is complete with respect to the norm $\|v\|:=|\langle v,v\rangle |^{1/2}$, then $H$ is called
a {\sc Hilbert $\A$-module}.
\end{ddd}

The completion of a pre-Hilbert $\A$-module is a Hilbert $\A$-module.

The right $\A$-module $\A^n$ is a Hilbert $\A$-module when endowed with the standard $\A$-valued scalar product 
$$\langle a,b \rangle := \sum_{i=1}^n a_i^* b_i \ .$$
The right $\A$-module $\{(a_n)_{n \in \bbbn} \subset \A~|~\sum\limits_{n
=1}^{\infty}a_n^*a_n \mbox{ converges}\}$ endowed with the $\A$-valued scalar
product $$\langle (a_n)_{n \in \bbbn},(b_n)_{n \in \bbbn}\rangle :=\sum_{n=1}^{\infty}
a_n^*b_n$$ is a Hilbert $\A$-module and is denoted by $l^2(\A)$.  Sometimes we use $\bbbz$ as index set.

Let $M$ be a measure space and let $\langle ~,~\rangle $ be the standard $\A$-valued scalar product on $\A^n$. Then the Hilbert $\A$-module
$L^2(M,\A^n)$ is defined in the following way: By 
$$\langle f,g\rangle _{L^2}=\int_M \langle f(x),g(x)\rangle dx$$ an $\A$-valued scalar product is defined on
the quotient of the 
space of simple functions on $M$ with values in $\A^n$ by the space of simple functions
 vanishing almost everywhere. Hence the quotient is a pre-Hilbert
$\A$-module. Its completion is the Hilbert $\A$-module  $L^2(M,\A^n)$.

Let $H$ be a Hilbert $\A$-module.

A submodule $U \subset H$  is called complemented if 
$$U^{\perp}=\{x \in
H~|~\langle x,u\rangle =0~ \forall u \in U\}$$ satisfies $U \oplus
U^{\perp}=H \ .$ 
Any projective submodule in $H$ is complemented.  

Let $H_1,H_2$ be Hilbert $\A$-modules.  The elements of
\begin{eqnarray*}
B(H_1,H_2) &=& \{ T:H_1 \to H_2~|~ T \mbox{ is continuous and }\exists T^*:H_2 \to H_1 \mbox{ with }\\
&& \langle Tv,w\rangle _{H_2}=\langle v,T^*w\rangle _{H_1}~ \forall v \in H_1,~w \in H_2\} 
\end{eqnarray*}
are called bounded operators from $H_1$ to $H_2$. They form a Banach space with
respect to the operator norm. With the composition as a product $B(H_1):=B(H_1,H_1)$ is a $C^*$-algebra. Note
that the existence of an adjoint must be required.

A continuous $\A$-module map $K:H_1 \to H_2$ is called compact if it can be approximated by a linear combination of operators of the form $x \mapsto z\langle y,x\rangle$ with $y\in H_1,z \in H_2,$ in the operator norm topology. 

Every compact operator is
adjointable.

A projection in $B(H)$ is compact if and only if it is a projection onto a
projective submodule of $H$.

If the range of $T \in B(H_1,H_2)$ is complemented, we call its complement the cokernel
$\Coker T$. Clearly a necessary condition for its existence is that the range of $T$ is
closed. The following proposition shows that it is sufficient:

\begin{prop}
\label{bouncloran}
Suppose that $T \in B(H_1,H_2)$ has closed range. Then

\begin{enumerate}
\item $\Ker T$ is a complemented submodule of $H_1$,
\item $\Ran T$ is a complemented submodule of $H_2$,
\item $T^*:H_2 \to H_1$ also has closed range.
\end{enumerate}
\end{prop}

\begin{proof}
\cite{la}, Th. 3.2.
\end{proof}

\subsection{Fredholm operators}
\label{fredop}

Let $H_1,H_2$ be Hilbert $\A$-modules isomorphic to $l^2(\A)$. 

\begin{ddd}
\label{deffredMF}

An operator $F \in B(H_1,H_2)$ is {\sc Fredholm} if 
there are decompositions $H_1=M_1 \oplus N_1$ and $H_2=M_2 \oplus N_2$ with the
following properties:
\begin{enumerate}
\item $N_1,N_2$ are projective $\A$-modules and $M_1,M_2$ are closed $\A$-modules.
\item The operator $F$ is diagonal: $F=F_M \oplus F_N$ with $F_M:M_1 \to M_2$ and $F_N:N_1 \to N_2$.
\item The component $F_M:M_1 \to M_2$ is  an isomorphism. 
\item The projection onto $N_i$ along $M_i$ is adjointable for $i=1,2$.
\end{enumerate}

The {\sc index of $F$} is defined as
$$\ind F:=[N_1]-[N_2] \in K_0(\A) \ .$$
\end{ddd}

\begin{prop}
\label{paramMF}
A selfadjoint operator $F \in B(H_1,H_2)$ is Fredholm if and only if there
exists an $\A$-linear continuous, not necessarily adjointable, map $G:H_2 \to H_1$ such that $FG-1$
and $GF-1$ are compact.
\end{prop}

\begin{proof} 
Analogous to \cite{mf}, Theorem 2.4.
\end{proof}

From the proposition it follows that
  if $F$ is Fredholm, then for any compact operator $K$ the operator $F+K$ is Fredholm. 

\begin{prop}
If $F:H_1 \to H_2$ is a Fredholm operator and $K:H_1 \to H_2$ is a compact operator, then 
$$\ind F=\ind(F+K) \ .$$
\end{prop}

\begin{proof}
\cite{mf}, Lemma 2.3.
\end{proof}

Another important property of Fredholm operators is the following:

\begin{prop}
\label{fredclos}
If $F\in B(H_1,H_2)$ is Fredholm and $\Ran F$ is closed, then $\Ker F$ and
$\Coker F$ are projective modules. 

Hence
$$\ind F=[\Ker F]-[\Coker F] \ .$$
\end{prop}

\begin{proof}
Let $P_{\Ker F} \in B(H_1)$ resp. $P_{\Coker F} \in B(H_2)$ be the orthogonal projection onto the
kernel resp. cokernel of $F$. They exist by Prop. \ref{bouncloran}. We have to
prove that they are compact.

There is $E \in B(H_2,H_1)$ such that $EF=1-P_{\Ker F}$ and $FE=1-P_{\Coker F}$.

Let $G$ be a parametrix of $F$. It follows that
$$P_{\Ker F}=1-EF=(1-GF)(1-EF)$$
and $$P_{\Coker F}=1-FE=(1-FE)(1-FG)$$
are compact operators.
\end{proof}

\begin{prop}
\label{fredhomotop}
If $F:[0,1] \to B(H_1,H_2)$ is a continuous path of Fredholm operators, then the
map $[0,1] \to K_0(\A),~t \mapsto \ind F(t)$ is constant.
\end{prop}

\begin{proof} see \cite{wo}, Prop. 17.3.4.
\end{proof}

\subsection{Regular operators}
\label{regop}

In this section some basic facts about unbounded operators on Hilbert $\A$-modules are collected. Most of them and more can be found in \cite{la}.

Let $H$ be a Hilbert $\A$-module with $\A$-valued scalar product $\langle ~,~\rangle $. Let
$D:\dom D \to H$ be a  densely defined operator on $H$.

\begin{lem}
If the adjoint $D^*$ of $D$ is densely defined, then $D$ is closable.
\end{lem}

\begin{proof}
Let $(f_n)_{n \in \bbbn}$ be a sequence in $\dom 
D$ such that $(f_n,Df_n)$   converges to $(0,f)$ in $H \oplus H$ for $n \to
\infty$. Then for every $g \in \dom D^*$ 
$$\langle f,g\rangle =\lim\limits_{n \to \infty}\langle Df_n,g\rangle =\lim\limits_{n \to \infty}\langle f_n,D^*g\rangle =0 \ .$$
Since $\dom D^*$ is dense in $H$, it follows that $f=0$.
\end{proof}

If $D$ is closed, then $$\langle f,g\rangle _D:=\langle f,g\rangle  + \langle Df,Dg\rangle $$ is an $\A$-valued scalar
product on $\dom D$, with respect to which $\dom D$ is a Hilbert $\A$-module, denoted by $H(D)$.

\begin{lem}
\label{ker}
Assume that  $D$ is closed.
\begin{enumerate}
\item
Suppose that $D$ has a densely defined adjoint $D^*$. Then $\Ker D^* =(\Ran D)^{\perp}$. 
\item
$\Ker D$ is complemented in $H(D)$ if and only if $\Ker D$ is complemented in $H$.
\end{enumerate} 
\end{lem}

\begin{proof}
(1) Since for $f \in \Ker D^*$ and $h \in \dom D$   
$$\langle f,Dh\rangle =\langle D^*f,h\rangle =0 \ ,$$
the $\A$-module $\Ker D^*$ is a submodule of $(\Ran D)^{\perp}$.

For $g \in (\Ran D)^{\perp}$ the $\A$-linear functional
$$\dom D \to \A,~f \mapsto \langle g, Df \rangle $$ vanishes. Thus
$g \in \dom D^*$ and $D^*g=0$.

(2) For $g \in \Ker D$ and $f \in \dom D$ the conditions $\langle f,g\rangle _{D}=0$ and
$\langle f,g\rangle =0$ are equivalent. Hence if $\Ker D$ is complemented in $H(D)$, then $\Ker D$ is complemented in
$H$. 

Conversely, if $\Ker D$ is complemented in $H$, we can decompose $g \in
H(D)$ in a sum $g=g_1+g_2$ with $g_1 \in \Ker D$ and $g_2 \in (\Ker
D)^{\perp}$. From $\Ker D \subset H(D)$ it follows that $g_2=g-g_1 \in H(D)$, hence
$g_2 \in (\Ker D)^{\perp_{H(D)}}$.
\end{proof}

Recall that $D$ is called regular if it is closed with densely defined adjoint $D^*$ and if $1+D^*D$ has dense
range, or equivalently if it is closed with densely defined adjoint and if its
graph is complemented in $H \times H$. 

If $D$ is regular, then $1+D^*D$ has a bounded inverse.

In the following the adjoint of an operator $A \in B(H(D),H)$ is denoted by
$A^T \in B(H,H(D))$ in order to distinguish it from the adjoint $A^*$
of $A$ as an unbounded operator on $H$.

\begin{lem}
\label{adj}
Assume that $D$ is closed.
\begin{enumerate}
\item
The operator $D$ is regular if and only if the inclusion  $\iota:H(D) \to H$ is in $B(H(D),H)$ and $(1+D^*D)$ is selfadjoint. Then $\iota^T=(1+D^*D)^{-1} \in B(H,H(D))$ and 
$(1+D^*D)^{-\frac 12}:H \to H(D)$ is
an isometry.

\item Assume that $D$ is regular and selfadjoint. Then $D \in B(H(D),H)$ and 
$D^T=D(1+D^2)^{-1} \ .$
\end{enumerate}
\end{lem}

\begin{proof}
 If $D$ is regular, then for $v \in H(D)$  and $w\in H$
\begin{eqnarray*}
\langle \iota v,w\rangle  &=&\langle v,w\rangle  \\
&=&\langle v,(1+D^*D)(1+D^*D)^{-1}w\rangle \\
&=&\langle v,(1+D^*D)^{-1}w\rangle  +\langle Dv,D(1+D^*D)^{-1}w\rangle \\
&=& \langle v,(1+D^*D)^{-1}w\rangle _D \ .
\end{eqnarray*}
This shows $\iota^T=(1+D^*D)^{-1}$.

Now the converse direction:

Let $v \in H$.
Then for any $w \in \dom(1+D^*D)$ 
$$\langle v,w\rangle =\langle \iota^Tv,w\rangle _D=\langle \iota^Tv,(1+D^*D)w\rangle  \ .$$
Since $(1+D^*D)$ is selfadjoint, it follows that $\iota^Tv \in \dom (1+D^*D)$ and  $(1+D^*D)\iota^Tv=v$. 

This shows that $(1+D^*D)$ is surjective and that $\iota^T$ is a right inverse
of $(1+D^*D)$. Since $(1+D^*D)$ is
bounded below, it is injective as well. It follows that
$(1+D^*D)$ is invertible and $\iota^T$ is its inverse.

The remaining parts are  immediate.

\end{proof}

\begin{prop}
\label{regpert}
Let $D_0$ be a regular selfadjoint operator and assume $D=D_0 +V$ with $V \in B(H)$.

\begin{enumerate} 
\item Then $D$ is closed.
\item The identity map induces a
continuous isomorphism from $H(D_0)$ to $H(D)$.
 \item $D \in B(H(D_0),H)$.
\item Suppose that $V$ is selfadjoint. Then $D$ is regular.
\end{enumerate}
\end{prop}

\begin{proof}

(1) From $\dom D^{**}=\dom D_0^{**}=\dom D_0=\dom D$ it follows that
$D=D^{**}$. Thus $D$ is closed. 

Assertion (2) follows from the fact that there is $C>0$ such that for all
$f \in H(D_0)$ 
\begin{eqnarray*}
\|f\|^2_D &\le& \|f\|^2 + \|D_0f\|^2+ \|\langle Vf,D_0f\rangle \|+ \|\langle D_0f,Vf\rangle \|+ \|Vf\|^2\\
&\le& C(\|f\|^2+\|D_0f\|^2) +2\|Vf\| \|D_0f\|\\
&\le& C\|f\|^2_{D_0} \ .
\end{eqnarray*}
We applied Cauchy-Schwarz inequality.

(3) By (2) the operator $D:H(D_0) \to H$ is continuous. By the previous lemma the adjoint of  $D=D_0 + V \iota :H(D_0) \to  H$ is 
$$D^T=D_0(1+D_0^2)^{-1} + (1+D_0^2)^{-1} V^*: H \to  H(D_0) \ .$$ 
(4) By (3) the operator  $D+i:H(D_0) \to H$ is an adjointable bounded
operator. By \cite{la}, Lemma 9.7, the range of $D+i$ is closed. Thus it is
complemented by Prop. \ref{bouncloran}. From Lemma \ref{ker} it follows  
that the
cokernel of $D+i$ agrees with the kernel of $D-i$. By \cite{la}, Lemma 9.7, the operator $D-i$
is injective. It follows that $\Coker(D+i)=\{0\}$. By \cite{la}, Lemma 9.8, this shows that $D$ is regular. 
\end{proof}

It can be shown that the condition of selfadjointness for $D_0$ and $V$ in statement (4) of the previous proposition can be dropped.

\begin{prop}
\label{cloran}
Let  $D$ be a regular and selfadjoint operator on $H$ with  closed range. 
\begin{enumerate}
\item The cokernel of $D$ exists and $\Ker D =\Coker D$. In particular $\Ker D$ is  complemented.
\item The $\A$-module $\dom D \cap \Ran D$ is dense in $\Ran D$ and  $\dom D =
\Ker D \oplus (\dom D \cap \Ran D)$, thus $D= 0 \oplus D|_{\Ran D}$ and
$D|_{\Ran D}$ is invertible.
\end{enumerate}
\end{prop}

\begin{proof}
(1) Since, by  Lemma \ref{adj}, $D \in B(H(D),H)$ and since the range of $D$ is
closed,  Prop. \ref{bouncloran} implies that the range is complemented. Its complement is $\Coker D$. Since $D$ is
selfadjoint,  $\Ker D= \Coker D$ by Lemma \ref{ker}.

(2) Let $P:H \to \Ran D$ be the orthogonal projection. From 
$(1-P)(\dom D) \subset \Ker D \subset \dom D$ 
we conclude that $P(\dom D) \subset \dom D$.
The assertion follows because $P(\dom D)$ is dense in  $P(H)=\Ran D.$
\end{proof}

We will need the following $\bbbz/2$-version of the previous proposition:

If $H=H^+ \oplus H^-$ is $\bbbz/2$-graded, then we call a closed operator $D$ on $H$ even
resp. odd if $\dom D$ decomposes in $(\dom D)^+ \oplus 
(\dom D)^-$ and if the action of $D$
is even resp. odd. 
 
\begin{prop}
\label{kerproj}
Let $H$ be a $\bbbz/2$-graded Hilbert $\A$-module and let $D$ be an odd regular selfadjoint
operator on $H$.

Suppose that $D^+:(\dom D)^+ \to H^-$ is surjective.

Then the range of $D$ is closed and complemented. Furthermore $\Ker
 D^+=\Ker D=\Coker D= \Coker D^-$ and this module is complemented.

\end{prop}

\begin{proof}
Since $D^+$ is surjective, $D^-$ is injective and so $\Ker D^+=\Ker D$.

Let $P_+$ be the orthogonal projection onto $H^+$. Since $D$ is odd, $DP_+=P_+D$. 

By Lemma \ref{adj} the operator $DP_+:H(D) \to H$ is adjointable with adjoint
$P_+D(D^2+1)^{-1}$.  From $P^+D|_{H(D)^{\pm}}=D^{\pm}$ it follows that $D^-(D^2+1)^{-1}:H^- \to H(D)^+$ is the
adjoint of $D^+:H(D)^+ \to H^-$.

Since $D^+$ is surjective, $\Ker D^+$ is complemented in $H(D)^+$ and the range of
  the adjoint $D^-(D^2+1)^{-1}:H^- \to H(D)^+$ is closed. Furthermore $$D^-(D^2+1)^{-1}=(D^2+1)^{-\frac 12}D^-(D^2+1)^{-\frac
  12} \ ,$$ and 
 $(D^2+1)^{-1/2}:H^{\pm} \to H(D)^{\pm}$ is an isomorphism by Lemma \ref{adj},  hence $\Ran D^-$ is closed, too. 
\end{proof}

\begin{prop}
\label{specD}
Let $D$ be a regular selfadjoint operator on $H$.
\begin{enumerate}
\item
For all $\lambda \in
\bbbc \setminus \bbbr$  the operator $D-\lambda$ is invertible.
\item
Assume that the range of $D$ is closed and let $P$ be the projection onto the
kernel of $D$. Then there is $c>0$ such that the spectrum of $(D+P)$ is
contained in
$\bbbr \setminus ]-c,c[$ and the spectrum of $D$ is contained in\\ $(\bbbr
\setminus ]-c,c[) \cup \{0\}$.
\end{enumerate}
\end{prop}

\begin{proof}
This follows from the functional calculus for regular operators (\cite{la},
Th. 10.9) and from the decomposition in Prop. \ref{cloran}.
\end{proof}

The following criteria for selfadjointness are useful:

\begin{lem}
\label{critself}

Let $D$ be a symmetric regular operator such that the range of $D+i$ and
of $D-i$ is dense in $H$. Then $D$ is selfadjoint.
\end{lem}

\begin{proof}
The operators  $D+i$ and of $D-i$ have closed range (see \cite{la}, Lemma 9.7). It follows that they
have a bounded inverse on $H$. Then they are adjoint to each other, thus $D$ is
selfadjoint.
\end{proof}

\begin{lem}
\label{invreg}
Assume that $D$ is symmetric and has an inverse $D^{-1} \in B(H)$. Then $D$ is regular.
\end{lem}

\begin{proof}
Since $D$ is symmetric, the
adjoint is densely defined. From $D^{-1} \in B(H)$ it follows that the graph of $D^{-1}$ is complemented,
hence the graph of $D$ is complemented as well.  Hence $D$ is regular.
\end{proof}

\subsection{Decompositions of Hilbert $C^*$-modules}

Let $H$ be a Hilbert $\A$-module with $\A$-valued scalar product $\langle ~,~\rangle $.
Let $J=\{1, \dots, m\} \subset \bbbn$ resp. $J=\bbbn$. If $J=\bbbn$, then set $m =\infty$.

\begin{ddd}
\label{deforthsys}
A system $\{f_k\}_{k \in J} \subset H$ is called
{\sc orthonormal} if for all $k,l \in J$ 
$$\langle f_k,f_l\rangle =\delta_{kl} \ .$$
It is called an {\sc orthonormal basis of $H$} if for all $f \in H$ there is $(a_n)_{n \in J} \subset \A$ such that $f=\sum\limits_{n =1}^m
 f_n a_n$.
\end{ddd}

Since $a_n=\langle f_n,f\rangle $, the coefficients are uniquely defined by the system.

\begin{prop}
\label{orthsys}
Let $\{f_k\}_{k \in
J}$ be an orthonormal system in $H$ whose span is dense in $H$. Then it is
an orthonormal basis of $H$ and the map $$f \mapsto (\langle f_n,f\rangle )_{n \in J}$$ is an isomorphism from $H$ to $\A^m$ if $m<\infty$ and to $l^2(\A)$ else.
\end{prop}

\begin{proof}
Let $P_n$ be the orthogonal projection onto the span of the first $n$ vectors of the system $\{f_k\}_{k \in
J}$. On $\spann_{\A}\{f_k~|~k \in
J\}$ the projection $P_n$ converges
strongly to the identity for $n \to \infty$. Since
$\|P_n\|=1$ for all $n \in \bbbn$, it follows that $P_n$ converges strongly to 
the identity on $H$. 
\end{proof}

\begin{lem}
\label{sumop}

Let  $\{U_i\}_{i \in \bbbn}$ be a family of pairwise orthogonal closed subspaces
of $H$ such that $\oplus_{i \in \bbbn} U_i$ is dense in $H$. Let $\{T_i\}_{i \in
\bbbn}$ be a family  of operators with 
$T_i \in B(U_i)$ and assume that there is $c\in \bbbr$ such that
$\|T_i\|\le c$  for all $i \in
\bbbn$.

Then the closure $T$ of the operator $\oplus_{i \in \bbbn}T_i$ is in $B(H)$ and 
$\|T\| \le c.$
\end{lem}

\begin{proof}
The spectral radius of an operator $A\in B(H)$ is denoted by $r(A)$.

Write $T(n)$ for the restriction of $T$ on $\oplus_{i =1}^n U_i$. Then for all $n \in \bbbn$
$$\|T(n)\|^2=r(T(n)^*T(n))=\max\limits_{1\le i\le n}r(T_i^*T_i)=
\max\limits_{1\le i\le n} \|T_i^*T_i\| \le c^2\ .$$
For $v \in \dom T$ there is $n \in \bbbn$ such that $v \in \oplus_{i =1}^n U_i$. Then $Tv=T(n)v$ and thus 
$$\|Tv\|=\|T(n) v\| \le c\|v\| \ .$$
It follows that the closure of $\oplus_{i \in \bbbn}T_i$ is a
bounded operator on $H$. The adjoint is given by the closure of $\oplus_{i \in \bbbn}T_i^*$.
\end{proof}

\begin{cor}
\label{decinv}
Let  $\{U_i\}_{i \in \bbbn}$ be a family of pairwise orthogonal closed subspaces
of $H$ such that $\oplus_{i \in \bbbn} U_i$ is dense in $H$. Let $\{T_i\}_{i \in
\bbbn}$ be a family of operators such that $T^{-1}_i \in B(U_i)$ and
assume that there is $c\in \bbbr$ such that $\|T_i^{-1}\|\le c$ for all $i \in \bbbn$.
Then the closure $T$ of the map $\oplus_{i \in \bbbn}T_i$ is invertible with inverse in $B(H)$.
\end{cor}

\begin{proof}
The operator $\oplus_{i \in \bbbn}T_i^{-1}$ is inverse to $\oplus_{i \in \bbbn}T_i$. It fulfills the conditions of the
previous lemma, hence its closure is a bounded operator on $H$. It is the inverse of
the closure of $T$.
\end{proof} 

\begin{prop}
\label{dirsum}
Let $\{e_i\}_{i \in \bbbn}$ be the standard basis of $l^2(\A)$. Let $M$ be a closed
and $N$ a projective submodule of $l^2(\A)$ such that $l^2(\A)=M \oplus N$. Let
$P$ be the projection onto $N$ along $M$ and let $P_n$ be the orthogonal projection
onto $L_n:=\spann_{\A}\{e_i~|~i =1, \dots, n\}.$ Assume that $P$ is adjointable.

For $n \in \bbbn$ such that $$\|P(1-P_n) \|\le \frac 12$$ it holds:

\begin{itemize}
\item[(i)]
The $\A$-module $N':=P_n(N)$ is projective and the maps
$$P_n:N \to N' \text{ and } P:N' \to N $$ 
are isomorphisms.

\item[(ii)] $l^2(\A)=M \oplus N' \ .$
\end{itemize}

\end{prop}

Note that for $n \in \bbbn$ large enough $\|P(1-P_n) \|\le \frac 12$
 since $P$ is a compact operator.

\begin{proof}
From $\|P(1-P_n) \|\le \frac 12$ it follows that
$$\|1_N-(PP_n)|_N\|\le \frac 12 \ .$$
Hence $(PP_n)|_N:N \to N$ is invertible.

The module $N':=P_n(N)$ is closed and finitely generated, thus projective. Furthermore, since $(PP_n)|_N$ is an isomorphism, the maps
$$P_n:N \to N' \text{ and } P:N' \to N$$
are isomorphisms as well.

It remains to show $N' \oplus M =l^2(\A)$.

The intersection $N' \cap M$ is trivial: If $x \in N' \cap M$, then
$Px=0$, hence, as $P:N' \to N$ is an isomorphism, $x=0$.  

Let now $x \in l^2(\A)$. Since $PP_n:N \to N$ is invertible, there is $y \in N$ such that
$Px=PP_ny$.

Then
\begin{eqnarray*}
x&=& (1-P)x+Px\\
&=& (1-P)x+PP_ny\\
&=&(1-P)x + P_ny -(1-P)P_ny\\
&=& (1-P)(x-P_ny)+P_ny \ .
\end{eqnarray*}
Since $(1-P)(x-P_ny) \in M$ and $P_ny \in N'$, it follows that $x \in M \oplus N'$.
\end{proof}

\section{Operators on spaces of vector valued functions}

\subsection{Vector valued functions and tensor products}
\label{fsp}
Let $V$ be a Fr\'echet space.

In the following we are interested in spaces $\Gamma(M,V)$ where $\Gamma=C,\C,C^k,L^p$ and $M$ is endowed with the appropriate structure.
It is desirable to have an isomorphism $\Gamma(M, \bbbc) \ten V \cong \Gamma(M,V)$ extending the inclusion $\Gamma(M,\bbbc) \odot V \incl \Gamma(M,V)$  because this ensures that every bounded
operator on $\Gamma(M,\bbbc)$ extends to a bounded operator on $\Gamma(M,V)$ (at least if the tensor product is an $\ve$- or $\pi$-tensor product). Examples where such an isomorphism exists are listed below whereas the example $\Gamma=L^2$ where in general such an isomorphism cannot be found is studied in detail in the following sections.
 
Proofs can be found in \cite{tr}. If not specified
the functions are assumed to be complex valued. 

\begin{itemize}
\item
Let $M$ be a compact topological space. Then  $$C(M)
\ten_{\ve}V \cong C(M,V) \ .$$ 
For compact spaces $M,N$ 
$$C(M \times N) \cong C(M) \ten_{\ve} C(N) \ .$$

\item
Let $U \subset \bbbr^n$ be  open and precompact. For all $m \in \bbbn_0$  $$C^m_0(U)
\ten_{\ve}V \cong C^m_0(U,V) \ .$$ 

\item
Let $M$ be a closed smooth manifold. For all $m \in \bbbn_0$  $$C^m(M)
\ten_{\ve}V \cong C^m(M,V) \ .$$ 

\item Let $U \subset \bbbr^n$ be  open and precompact. Then $\C_0(U)$ is nuclear, in particular 
$$\C_0(U) \ten_{\pi} V \cong  \C_0(U)
\ten_{\ve}V \cong  \C_0(U,V) \ .$$
\item
Let $M$ be a closed smooth manifold. Then $\C(M)$ is nuclear, in particular 
$$\C(M) \ten_{\pi} V \cong \C(M)
\ten_{\ve}V \cong  \C(M,V) \ .$$
For closed smooth manifolds $M,N$ 
$$\C(M \times N) \cong \C(M) \ten \C(N) \ . $$

\item The space of Schwartz functions ${\mathcal S}(\bbbr)$ is nuclear, in
particular  $${\mathcal S}(\bbbr) \ten_{\pi}V \cong {\mathcal S}(\bbbr)
\ten_{\ve}V \cong {\mathcal S}(\bbbr,V) \ .$$ 
\end{itemize}
The isomorphisms are given by the unique extension of the inclusion of the algebraic tensor product.

\subsection{$L^2$-spaces and integral operators}
\label{Lp}

Let $E$ be a Banach space with norm $| \cdot |$. Let $\End E$ be the Banach
algebra of bounded operators on $E$. We denote the operator norm on $\End E$ by
$| \cdot |$ as well.  

\begin{ddd}
Let $M$ be a measure space and $p \in \bbbn$. 

The Banach space
$L^p(M,E)$ is defined to be the completion of the quotient of the space of simple $E$-valued
functions on $M$ by the subspace of functions vanishing almost everywhere with
respect to the norm
$$\|f\|_{L^p}:= \left(\int_M |f(x)|^p dx\right)^{\frac 1p} \ .$$
\end{ddd}

In order to avoid confusion we make the following convention: If $E=\A^n$ for a
$C^*$-algebra $\A$, then $L^2(M,E)$ denotes the Hilbert $\A$-module defined in
\S\ref{HCmod} and not the space just defined. In general these spaces do not coincide.

\begin{lem}
\label{identL2}
Let $M_1, M_2$ be $\sigma$-finite measure spaces.
Then the map
$$L^2(M_1 \times M_2,E) \to L^2(M_1, L^2(M_2,E)), ~f \mapsto (x \mapsto f(x,
\cdot)) $$
is an isometric isomorphism.
\end{lem}
            
\begin{proof}
The lemma follows from Fubini.
\end{proof}

\begin{prop}
\label{poskern}
Let $M$ be a measure space.

Let $k:M \times M \to \End E$ be a measurable function such
that  the integral
kernel $|k(x,y)|$ defines a bounded operator $|K|$ on $L^2(M)$ with norm $\| |K|\|$. Then $k$ defines
a bounded operator on $L^2(M,E)$ with norm less than or equal to $\| |K| \|$.
\end{prop}

\begin{proof}
For a simple function $f:M \to E$ 
$$\|\int_M k(\cdot ,y)f(y)~dy\|_{L^2} \le \|\int_M |k(\cdot ,y)| |f(y)| ~dy \|_{L^2}\le \| |K| \|~ \|f\|_{L^2}.$$
\end{proof}

\begin{cor}
\label{L2ker}
Let $M$ be a measure space.

There is a norm-decreasing map
$$L^2(M \times M, \End E) \to B(L^2(M,E),L^2(M,E))$$
$$k \mapsto \Bigl(f\mapsto Kf:=\int_M k(\cdot,y)f(y)~dy \Bigr) \ .$$
\end{cor}

\begin{cor}
\label{conv}

The convolution induces a continuous map
$$L^1(\bbbr^n, \End E) \to B(L^2(\bbbr^n,E)),~ f \mapsto (g \mapsto f*g) \ .$$
\end{cor}

\begin{proof}
The convolution with $f \in L^1(\bbbr^n, \End E)$ is an integral operator with
integral kernel $f(x-y)$. Since  $|f| \in L^1(\bbbr^n)$ for $f \in L^1(\bbbr^n, \End E)$, the convolution with $|f|$ is bounded
on $L^2(\bbbr^n)$. Then the assertion follows from the previous proposition. 
\end{proof}

\begin{lem}
\label{trans}
For any $f \in L^p(\bbbr^n,E)$ the translation map $$\tau f:\bbbr^n \to
L^p(\bbbr^n,E),~ y \mapsto \tau_y f \ ,$$ defined by $$\tau_yf(x):=f(x-y) \ ,$$
is continuous.
\end{lem}

\begin{proof}
The proof is analogous to the case $E=\bbbc$, see \cite{co}, Ch. VII, Prop. 9.2.
\end{proof}

\subsection{Adjointable operators on Banach spaces}
\label{adop}

Let $\B$ be an involutive Banach algebra with unit. In this section all operators are assumed to be right $\B$-module maps. Let $E$ be a Banach right $\B$-module with norm $|\cdot |$.

\begin{ddd} 
\label{nondegprod} 
A {\sc $\B$-valued non-degenerated product on $E$} is a map
$\langle ~,~\rangle :E \times E \to \B$ that is $\bbbc$-linear in the second variable and has the
following properties:
\begin{enumerate}
\item $\langle v,wb\rangle =\langle v,w\rangle b$ for all $v,w \in E,~b \in \B$,
\item $\langle v,w\rangle =\langle w,v\rangle ^*$ for all $v,w \in E$,
\item if $\langle v,w\rangle =0$ for all $w \in E$, then $v=0$,
\item there is $C>0$ such that $|\langle v,w\rangle | \le C|v|_E~|w|_E$ for all $v,w \in E$.
\end{enumerate}
\end{ddd}

Let $E$ be endowed with a $\B$-valued non-degenerated product $\langle ~,~\rangle $.

\begin{ddd} 
A bounded operator
$T:E \to E$ is called {\sc adjointable} if there is a map $T^*:E \to E$ satisfying 
$$\langle v,Tw\rangle =\langle T^*v,w\rangle $$ for all $v,w \in E$. 
\end{ddd}

\begin{lem}
Let $T:E \to E$ be adjointable. 
\begin{enumerate}
\item The adjoint $T^*$ is unique.
\item $T^*$ is a right $\B$-module map.
\item $T^*$ is bounded.
\item  $T^{**}=T$.
\item  $(ST)^*=T^*S^*$.
\end{enumerate}
\end{lem}

\begin{proof} 
(1) Let $T^*$ be an adjoint of $T$ and $\Gamma(T^*)$ its graph. Then
$$\Gamma(T^*)\subset G:=\{(x,y) \in E \times E~|~\langle y,w\rangle +\langle -x,Tw\rangle =0~ \forall w \in
E\} \ .$$ 
Let  $v \in E$. There is a unique $v_1 \in E$ with $(v,v_1) \in G$ since from
$$\langle v_1,w\rangle =\langle v,Tw\rangle =\langle v_2,w\rangle  ~\forall w \in E$$  
it follows $\langle v_1-v_2,w\rangle =0$ for all $w \in E$ and therefore $v_1=v_2$.
This shows $\Gamma(T^*)=G$.

(2) If $(x,y),(v,w) \in \Gamma(T^*)$ and $b \in \B$, then
$(xb+v,yb+w) \in \Gamma(T^*)$ by the proof of (1). 

(3)  Since $\Gamma(T^*)$ is closed, the operator $T^*$ is bounded.

(4) From (1) it follows that $\Gamma(T^{**})=\Gamma(T)$.

(5) $\langle (ST)^*v,w\rangle =\langle v,STw\rangle =\langle S^*v,Tw\rangle =\langle T^*S^*v,w\rangle $.
\end{proof}

\begin{lem}
\label{closable}
Let $T$ be a densely defined operator on $E$. If there exists a densely defined operator $S$, called {\sc formal adjoint} of $T$, such that
$$\langle f,Tg\rangle =\langle Sf,g\rangle $$ for all $f \in \dom S,~g \in \dom T$, then $T$ is closable.
\end{lem}

\begin{proof}
The set
$$\{(x,y)\in E \oplus E ~|~\langle f,y\rangle  - \langle Sf,x\rangle =0 ~ \forall f \in \dom S \}$$ is the graph of a closed extension of $T$.
\end{proof}

The standard $\B$-valued non-degenerated product on $\B^n$ is given by
$$\langle v,w\rangle :=\sum_{i=1}^n v_i^* w_i \ .$$
Since the endomorphism set of $\B^n$ can be identified with $M_n(\B)$, all elements of $\End(\B)$ are adjointable and taking the adjoint is a
bounded linear map.

If $M$ is a measure space, the standard $\B$-valued non-degenerated product on $L^2(M,\B^n)$ is defined by
$$ \langle f,g\rangle _{L^2}:=\int_M \langle f(x),g(x)\rangle  dx \ . $$
We check condition (3) of Def. \ref{nondegprod}: 

Let $\B'$ be the topological dual of $\B$ endowed with the weak topology.

We use that every $\lambda \in \B'$
induces a map  $\lambda:\B^n \to \bbbc^n$ by componentwise application.

If $f \in L^2(M,\B^n)$ with $\langle f,g\rangle _{L^2}=0$ for all $g
\in L^2(M,\B^n)$, then  in particular for all $g \in L^2(M,\bbbc^n)$ and $\lambda \in \B'$
$$\int_M \lambda(f(x)^*)g(x) dx=0 \ ,$$
hence $\lambda(f(x)^*)$ vanishes almost everywhere. Since $\B'$ is
separable, it follows $f=0$ in $L^2(M,\B^n)$.

\subsection{Hilbert-Schmidt operators}
\label{HS}

Let the notation be as in the previous section. Assume that $M$ is a $\sigma$-finite measure
space. 

\begin{lem}
Let  $k \in L^2(M \times
M,M_n(\B))$ and let $K$ be the corresponding integral operator on $L^2(M,\B^n)$. Then $k$ is uniquely defined by $K$.
\end{lem}

\begin{proof}
It is enough to show that $k$ vanishes in $L^2(M \times
M,M_n(\B))$ if $K=0$.

Applying $\lambda \in \B'$
componentwise yields maps $\lambda:M_n(\B) \to M_n(\bbbc)$ and $\lambda:\B^n \to
\bbbc^n$. 

For $f \in L^2(M,\bbbc^n)$ we have that almost everywhere
$$\lambda(\int_M k(x,y)f(y)~dy)=\int_M \lambda(k(x,y))f(y)~dy =0 \ .$$ 
It follows that $\lambda(k(x,y))=0$ in  $L^2(M \times
M,M_n(\bbbc))$. Since $\B'$ is separable, the set of all
$(x,y) \in M\times M$ such that there is $\lambda \in \B'$ with
$\lambda(k(x,y)) \neq 0$ has measure zero. On the complement of this set $k$ vanishes.
\end{proof}

\begin{ddd} 
A {\sc Hilbert-Schmidt operator on $L^2(M,\B^n)$} is an integral operator with
integral kernel in $L^2(M \times M,M_n(\B))$. Let $A$ be a Hilbert-Schmidt
operator on $L^2(M,\B^n)$ with integral kernel $k_A \in L^2(M \times
M,M_n(\B))$. We define
$$\|A\|_{HS}:=\|k_A\| \ ,$$
where the norm on the right hand side is taken in $L^2(M \times M,M_n(\B))$.

The normed space of Hilbert-Schmidt operators  on $L^2(M,\B^n)$ is denoted by $HS(L^2(M,\B^n))$.
\end{ddd} 

Note that $HS(L^2(M,\B^n))$ is a Banach algebra and that the inclusion $$HS(L^2(M,\B^n)) \to
 B(L^2(M,\B^n))$$ is bounded. Prop. \ref{HScomp} below shows that $HS(L^2(M,\B^n))$ is
 a left $B(L^2(M,\B^n))$-module.

All operators in $HS(L^2(M,\B^n))$ are adjointable. The integral kernels of $A \in HS(L^2(M,\B^n))$
and $A^*$ are related by
$k_{A^*}(x,y)=k_A(y,x)^*$. It follows that
taking the adjoint is a bounded map on $HS(L^2(M,\B^n))$.

\begin{prop} 
\label{HScomp}

Let $A \in B(L^2(M,\B^n)),~ K \in HS(L^2(M,\B^n))$.
\begin{enumerate}
\item Then $AK \in HS(L^2(M,\B^n))$. Furthermore there is $C>0$, independent of $A$ and $K$, such that $$\|AK\|_{HS} \le C \|A\|
\|K\|_{HS} \ .$$  
\item If $A$ is adjointable, then $KA \in HS(L^2(M,\B^n))$.
Furthermore  there is $C>0$, independent of $A$ and $K$,
with $$\| KA\|_{HS} \le C\|A^*\|
\|K\|_{HS} \ .$$
\end{enumerate}
\end{prop}

\begin{proof}
(1) There is an isomorphism  $$L^2(M \times M, M_n(\B)) \cong L^2(M\times M, \B^n)^n $$
that is equivariant with respect to the left $M_n(\B)$-action  on both spaces. 
Furthermore the map
$$L^2(M \times M,\B^n) \to L^2(M,L^2(M,\B^n)),~ k \mapsto (y \mapsto
k(\cdot,y))$$ is an isomorphism by Lemma \ref{identL2}. The operator $A$ induces a bounded map on
$L^2(M,L^2(M,\B^n))$, namely  $k \mapsto (y \mapsto
A k(\cdot,y))$, clearly its norm is less than or equal to the norm of $A$ on $L^2(M,\B^n)$.

(2) The map $K \mapsto KA$ is a composition of the following maps on $HS(L^2(M,\B^n))$:
$$K \stackrel{*}{\mapsto} K^* \stackrel{A^*}{\mapsto} A^*K^* =(KA)^*
\stackrel{*}{\mapsto} KA \ .$$
By (1)  and the fact that taking the adjoint is bounded on $HS(L^2(M,\B^n))$
these maps are bounded.
\end{proof}

\begin{ddd} 
\begin{enumerate}
\item Let $\langle ~,~\rangle :\B^n \times \B^n \to \B$ be the standard $\B$-valued non-degenerated product.
For $e \in \B^n$ define the map
$$e^*:\B^n \to \B,~v \mapsto
\langle e,v\rangle  \ .$$
\item
An integral operator $A$ on $L^2(M,\B^n)$ is called {\sc finite} if there is $k \in
\bbbn$ and there are
functions $f_j,h_j \in L^2(M,\B^n),~j=1 \dots k,$ such that
$$k_A(x,y)=\sum_{j=1}^k f_j(x)h_j(y)^* \ .$$
\end{enumerate}
\end{ddd}

\subsection{Trace class operators}
\label{trclop}

Let $M$ be a $\sigma$-finite Borel space. Assume that there is a uniformly bounded sequence $(K_m)_{m \in \bbbn} \subset  B(L^2(M,\B^n))$ converging strongly to the identity such that each $K_m$ is an  integral operator with  continuous compactly supported integral kernel. 

This condition is fulfilled for example if $M$ is a complete Riemannian manifold: Using Prop. \ref{poskern} and the fact that the heat kernel $k_t(x,y)$ of the scalar Laplacian $\Delta$ is positive, one can deduce that $k_t(x,y)$ defines a strongly continuous semigroup $e^{-t\Delta}$ on  $L^2(M,\B^n)$. Then $K_m:=\phi_m e^{-\frac 1m \Delta}\phi_m$ is an appropriate sequence, where $(\phi_m)_{m \in \bbbn}$ is a sequence in $C_c(M)$ that converges to the identity on compact subsets of $M$.

The composition of operators  induces a continuous map
$$\mu:HS(L^2(M,\B^n)) \ten_{\pi} HS(L^2(M,\B^n)) \to B(L^2(M,\B^n)) \ .$$

\begin{ddd}
A {\sc trace class operator on $L^2(M,\B^n)$} is an element in the range of $\mu$. The space of trace class operators  is denoted by $\Tr(L^2(M,\B^n))$. We identify it with $$\bigl(HS(L^2(M,\B^n)) \ten_{\pi} HS(L^2(M,\B^n))\bigr)/ \Ker \mu$$ and endow it with the quotient norm.
\end{ddd}  

It follows from Prop. \ref{HScomp} that $\Tr(L^2(M,\B^n))$ is a left $B(L^2(M,\B^n))$-module and that there is a well-defined action of adjointable operators from the right.

\begin{prop}
The map $$R:HS(L^2(M,\B^n)) \ten_{\pi} HS(L^2(M,\B^n)) \to L^1(M,M_n(\B)) \ ,$$ 
$$(A,B) \to (x \mapsto \int_M k_A(x,y) k_B(y,x)~dy)$$ descends to a continuous map 
$$\ov R:\Tr(L^2(M,\B^n)) \to L^1(M,M_n(\B)) \ .$$
\end{prop}

\begin{proof}
Note that $R$ is continuous, hence all we have to show is that $\ov{R}$ is well-defined. Since the sequence $K_m$ is uniformly bounded, the map $F \mapsto K_mFK_m$ on $HS(L^2(M,\B^n)) \ten_{\pi} HS(L^2(M,\B^n))$ converges strongly to the identity. The assertion follows if we can show that $\mu(F)=0$ implies $R(K_mFK_m)=0$. 

The map $F \mapsto K_mFK_m$ is continuous as a map from $HS(L^2(M,\B^n)) \ten_{\pi} HS(L^2(M,\B^n))$ to $C_c(M,L^2(M,M_n(\B))) \ten_{\pi} C_c(M,L^2(M,M_n(\B)))$. The composition with the continuous map $C_c(M,L^2(M,M_n(\B))) \ten_{\pi} C_c(M,L^2(M,M_n(\B))) \to C_c(M \times M,M_n(\B))$ induced by the product $$L^2(M,M_n(\B)) \times L^2(M,M_n(\B)) \to M_n(\B),~(f,g) \mapsto \int_M f(x)g(x)~dx$$ agrees with the map that sends $F$ to the continuous integral kernel $k_{\mu(K_mFK_m)}$ of $\mu(K_mFK_m)$. We compose further with the restriction to the diagonal and get a continuous map 
$$HS(L^2(M,\B^n)) \ten_{\pi} HS(L^2(M,\B^n)) \to C_c(M,M_n(\B)) \subset L^1(M,M_n(\B)) \ ,$$ which agrees with the map $F \mapsto R(K_mFK_m)$. This shows that the map $F \mapsto R(K_mFK_m)$ factors through the map $$HS(L^2(M,\B^n)) \ten_{\pi} HS(L^2(M,\B^n)) \to C_c(M \times M,M_n(\B)),~ F \mapsto k_{\mu(K_mFK_m)} \ .$$ Now $\mu(F)=0$ implies $\mu(K_mFK_m)=0$, hence $k_{\mu(K_mFK_m)}=0$ by the uniqueness of the integral kernel in $C_c(M \times M,M_n(\B))$. 
\end{proof}

Let $$\tr:M_n(\B) \to \B/\ov{[\B,\B]}$$ be the trace defined, as usual, by adding up the diagonal elements and
let $$\Tr(T):=\int_M \tr ~\ov R(T)(x)~dx$$ for $T \in \Tr(L^2(M,M_n(\B)))$.
Then
for Hilbert-Schmidt operators $A,B$ 
$$\Tr(AB)= \Tr(BA) \ .$$

From the fact that any trace class operator can be approximated by finite sums of products of Hilbert-Schmidt operators it follows that the equation holds also for $A \in \Tr(L^2(M,\B^n))$ and $B$ any adjointable operator.

\subsection{Pseudodifferential operators}

\label{pseudo}

Let $E$ be a Banach space. 

Let $U$ be an open precompact subset of $\bbbr^n$.
Recall the notion of a symbol of order $m$ on $U$:

\begin{ddd} 

A function $a \in \C(U \times \bbbr^n,M_l(\bbbc))$ is called a  {\sc symbol of order
$m\in \bbbr$} if it is compactly supported in the first variable and if for all
multi-indices $\alpha,\beta \in \bbbn_0^n$ the expressions
 $$\sup\limits_{x\in U,~\xi \in \bbbr^n} (1+|\xi|)^{-m+|\beta|}|\ra_x^{\alpha}\ra_{\xi}^{\beta}a(x,\xi)| $$
are finite. 

These are seminorms on the space $S^m(U,M_l(\bbbc))$ of symbols of order $m$ on $U$,
which turn $S^m(U,M_l(\bbbc))$ into a Fr\'echet space.
\end{ddd}

In order to simplify formula involving Fourier transform it is convenient to rescale the
Lebesgue measure on $\bbbr^n$ by setting
$d'x:=(2\pi)^{-\frac n2}dx$.

We consider $L^2(U,E^l)$ as a subspace of $L^2(\bbbr^n,E^l)$ in the following.

The Fourier transform is bounded from $L^1(\bbbr^n,E)$ to $C_0(\bbbr^n,E)$.

A symbol $a \in S^m(U,M_l(\bbbc))$ defines a continuous operator
$$\Op(a):\C_c(U,E^l) \to \C_c(U,E^l),~ (\Op(a)f)(x) = \int_{\bbbr^n}e^{ix \xi}
a(x, \xi) \hat f(\xi)~ d'\xi $$
called a pseudodifferential operator of order $m$.

Due to the fact that the Fourier transform is in general not continuous on $L^2(\bbbr^n,E^l)$, the continuity properties of pseudodifferential operators acting on vector valued functions  are in general weaker that those one gets for $E= \bbbc$. They are still weaker if we allow for symbols with values in $\End E$.

\begin{lem}
\label{sobol}
\begin{enumerate}
\item
For $m <  -\frac{n}{2}$ and $a \in S^m(U,M_l(\bbbc))$ the operator $\Op(a)$
extends to a bounded operator on $L^2(U,E^l)$ and the map
$$\Op:S^m(U,M_l(\bbbc)) \to B(L^2(U,E^l))$$
is continuous.

\item Let $m < -\frac{n}{2}$ and $\nu,k \in \bbbn_0$ with $k<-\frac{n}{2}-m$.  Then for $a \in S^m(U,M_l(\bbbc))$ the operators 
$$\Op(a):C^{\nu}_0(U,E^l) \to C^{\nu+k}_0(U,E^l)$$
and 
$$\Op(a):L^2(U,E^l) \to C^k_0(U,E^l)$$
are continuous.

\end{enumerate}
\end{lem}

\begin{proof}
(1) Let $m < -\frac{n}{2}$.

The Fourier transform in $\xi$ induces a bounded map
$$S^m(U,M_l(\bbbc)) \to \C_c(U,L^2(\bbbr^n,M_l(\bbbc))),~a \mapsto (x \mapsto
\hat a(x, \cdot)) \ .$$
If $a \in S^m(U,M_l(\bbbc))$, we extend $\Op(a)$ to $L^2(U,E^l)$ by
$$(\Op(a)f)(x) := \int_{\bbbr^n}\hat a(x,z) f(-x-z) d'z \ .  $$
The map
$$S^m(U,M_l(\bbbc)) \to B(L^2(U,E^l)),~a \mapsto \Op(a)$$ is well-defined and continuous since
 $$\bbbr^n \to L^2(\bbbr^n,E^l),~x \mapsto (z \mapsto f(-x-z))$$ is continuous by Lemma \ref{trans}, hence $\Op(a)f \in C_0(U,E^l)$  and
$$\|\Op(a)f\|_{C_0} \le \sup_{x \in U} \| \hat a(x, \cdot )\|_{L^2} \|f\|_{L^2} \ .$$
(2)  First let $m< -\frac n2$ and $k=0$.

If $f \in C^{\nu}_0(U,E^l)$, then
$x \mapsto (z \mapsto f(-x-z))$ is in $C_0^{\nu}(\bbbr^n,L^2(\bbbr^n,E^l))$.

It follows as above that $\Op(a)f \in C_0^{\nu}(U,E^l)$ and
$$\|\Op(a)f\|_{C^{\nu}} \le C\sup_{|\alpha|\le \nu}\sup_{x \in U} \|\ra_x^{\alpha} \hat a(x,\cdot)\|_{L^2} \|f\|_{C^{\nu}} \ .$$
Now assume that the assertion holds for $k-1$ and all $a
\in S^m(U,M_l(\bbbc))$ with $m < -\frac n2 -k +1$. 

We prove the assertion for $k$ and $a \in S^m(U,M_l(\bbbc))$ with $m < -\frac n2 -k$:

If $\alpha \in \bbbn_0^n$ with $|\alpha|=1$, then the map $$f \mapsto \ra^{\alpha}(\Op(a)f)$$ is a
pseudodifferential operator of degree $m+1$. By induction
it  is a
bounded operator
from $C^{\nu}(U,E^l)$ to $C^{\nu+k-1}(U,E^l)$. It follows that $\Op(a)$ is continuous from $C^{\nu}(U,E^l)$ to
$C^{\nu+k}(U,E^l)$.

An analogous induction argument shows that $\Op(a)$ is continuous from $L^2(U,E^l)$ to
$C_0^k(U,E^l)$. For $k=0$ this was proved in (1).
\end{proof}

\section{Projective systems and  function spaces}

The projective systems $(\A_i)_{i \in \bbbn_0}$ and $(\Ol{\mu})_{i,\mu \in \bbbn_0}$
from \S \ref{einfai} induce projective systems of spaces
$\bigl(L^2(M,\A_i^l)\bigr)_{i \in \bbbn_0}$ and
$\bigl(L^2(M,(\Ol{\mu})^l)\bigr)_{i,\mu \in \bbbn_0}$. 
 
Recall the convention fixed in \S \ref{Lp}:
The space $L^2(M,\A^l)$ is the Hilbert $\A$-module defined
in \S \ref{HCmod}. For $\mu \in \bbbn_0$ and $i\in \bbbn$  the space
$L^2(M,(\Ol{\mu})^l)$ was defined in \S \ref{Lp}.

In the following we investigate the behavior of some particular classes of
operators on $L^2(M,(\Ol{\mu})^l)$ under the projective limit.

\subsection{Integral operators}
\label{intop}

Hilbert-Schmidt operators on $L^2(M,\A_i^l)$ have the
property that they extend to bounded operators on $L^2(M,(\hat\Omega_{\le
\mu}\A_j)^l)$ for all $j \in \bbbn_0$ with $j \le i$ and all $\mu \in
\bbbn_0$. We investigate how the spectrum
depends on $j,\mu$.

We extend a method developed by Lott (\cite{lo2},
\S 6.1.) for  closed manifolds to certain non-compact
manifolds with boundary (in particular for the manifold defined in \S \ref{situat}). 

Let $[0,1]^n$ be endowed with a measure of the form $hdx$ where $h$ is a
positive continuous function on $[0,1]^n$ and $dx$ is
the Lebesgue measure.

In the proof of the following proposition we use that there exists a Schauder basis of
$C([0,1]^n)$ that is orthonormal in 
$L^2([0,1]^n)$ (here and in the following we understand $L^2([0,1]^n)$ with respect to $hdx$). For
$h=1$ a Franklin system \cite{se}
 yields such a basis $\{f_n\}_{n \in \bbbn}$, then for general $h$ the
system $\{h^{-\frac 12}f_n\}_{n \in \bbbn}$ is one. 

The proposition still holds true if $[0,1]^n$ is replaced by a compact Borel space for which such a basis exists.

\begin{prop}
\label{intker}
\begin{enumerate}
\item Let $[0,1]^n$ be endowed with a measure $hdx$ as above.
Let $k \in C([0,1]^n\times [0,1]^n, M_l(\A_i))$ and let
$K$ be the corresponding integral operator. 
Assume that 
$1-K$ is invertible in $B(L^2([0,1]^n,\A^l))$. Then $1-K:L^2([0,1]^n,(\Ol{\mu})^l) \to L^2([0,1]^n,(\Ol{\mu})^l)$ is invertible.

\item 
Let $M$ be a Riemannian manifold of dimension $n$, possibly with
boundary. Suppose that  there is a covering $\{K_m\}_{ m \in \bbbn}$ of $M$ with $K_m$ compact, $K_m \subset K_{m+1}$ and such that $K_m$ is diffeomorphic to
$[0,1]^n$ for every $m \in \bbbn$. Let $k \in L^2(M
\times M, M_l(\A_i))\cap C(M
\times M, M_l(\A_i))$ and assume furthermore that $x \mapsto k(x, \cdot )$ and $y \mapsto
k(\cdot , y)$ are in $C(M,L^2(M,M_l(\A_i)))$. 

If $1-K$ is invertible in $B(L^2(M,\A^l))$, then
$1-K$ is invertible in $B(L^2(M,(\Ol{\mu})^l))$. 
\end{enumerate} 
\end{prop}

\begin{proof}
(1) Choose a basis of $C([0,1]^n)$ which is
 orthonormal with respect to $h dx$ and 
let $P_N$ denote the projection onto the first $N$ basis vectors. 
The integral kernel of $P_N$ is in $L^2([0,1]^n \times [0,1]^n)$,
hence $P_N$ acts continuously on $L^2([0,1]^n,( \Ol{\mu})^l)$.

We decompose $L^2([0,1]^n,( \Ol{\mu})^l)$ into the direct sum 
$$P_N L^2([0,1]^n,( \Ol{\mu})^l) \oplus (1-P_N)L^2([0,1]^n,( \Ol{\mu})^l)$$ and write
$$1-K=\left(\begin{array}{cc} a & b\\
c& d
\end{array} \right)$$
with respect to the decomposition.

If we find $N$ such that $d$ is invertible on $(1-P_N)  L^2([0,1]^n,(\Ol{\mu})^l)$ and  prove that then $a-bd^{-1}c$ is invertible on $P_N  L^2([0,1]^n,( \Ol{\mu})^l)$,
we can conclude that  $(1-K)$ is invertible by the equality 
$$\left(\begin{array}{cc} a & b\\
c& d
\end{array} \right)
=\left(\begin{array}{cc} 1 & bd^{-1}\\
0 & 1
\end{array} \right)
\left(\begin{array}{cc} a-bd^{-1}c & 0\\
0& d
\end{array} \right)
\left(\begin{array}{cc} 1 & 0\\
d^{-1}c & 1
\end{array} \right) \ .$$
First we show that $$d=(1-P_N)(1-K)(1-P_N)$$ is invertible for $N$ big
enough. By Prop. \ref{HScomp}  the operator
$(1-P_N)K(1-P_N)$ is a Hilbert-Schmidt operator. Its integral
kernel is continuous. 
 
For $N \to \infty$ the projections $P_N$ converge strongly to the identity on
$C([0,1]^n)$. Since $$C([0,1]^n\times [0,1]^n,
M_l(\A_i))\cong C([0,1]^n) \ten_{\ve}C([0,1]^n) \ten_{\ve} M_l(\A_i) \ ,$$
 there is $N$ such that the norm of the integral kernel of 
$(1-P_N)K(1-P_N)$ is smaller than $\frac 12$ in $C([0,1]^n \times [0,1]^n
,M_l(\A_i))$.

For that $N$ the series $$(1-P_N)+\sum\limits_{\nu=1}^{\infty}((1-P_N)K(1-P_N))^{\nu}$$
converges as a bounded operator on $(1-P_N)L^2([0,1]^n, ( \Ol{\mu})^l)$ and inverts
$d$. Hence $a-bd^{-1}c$ is well-defined.

Via the basis we identify $a-bd^{-1}c$  with an
element of $M_{Nl}(\A_i)$. Since $1-K$ is invertible on $L^2(M,\A^l)$, the matrix $a-bd^{-1}c$ is
invertible in $M_{Nl}(\A)$. By Prop. \ref{propai} it follows that $a-bd^{-1}c$ is
invertible in $M_{Nl}(\A_i)$ as well.

(2) Let $m \in \bbbn$ be such that in $L^2(M \times M, M_l(\A_i))$
$$\|(1-1_{K_m}(x))k(x,y) (1-1_{K_m}(y))\| \le \frac 12 \ .$$
Write 
$$(1-K)= \left(\begin{array}{cc} a & b\\
c& d
\end{array} \right)$$
with respect to the decomposition $$L^2(M,  ( \Ol{\mu})^l)=1_{K_m}L^2(M,  ( \Ol{\mu})^l) \oplus (1-1_{K_m})L^2(M,  ( \Ol{\mu})^l) \ .$$
By the choice of $m \in \bbbn$ the entry $d=1-(1-1_{K_m})K(1-1_{K_m})$ is invertible on
$(1-1_{K_m})L^2(M, (\Ol{\mu})^l)$.

We prove that $a-bd^{-1}c$ is invertible on $L^2(K_m,(\Ol{\mu})^l)$ and then the assertion follows as in the proof of (1).

On $L^2(K_m,(\Ol{\mu})^l)$  
$$a-bd^{-1}c=1_{K_m}-(1_{K_m}K1_{K_m}+bd^{-1}c) \ .$$
The operator $1_{K_m}K1_{K_m}+bd^{-1}c$ is an integral
operator on $L^2(K_m,(\Ol{\mu})^l)$ with continuous integral kernel: The integral kernel of $b$ is $1_{K_m}(x)k(x,y) (1-1_{K_m}(y))$ and $x \mapsto
1_{K_m}(x)k(x, \cdot) (1-1_{K_m})$ is in  $C(K_m,L^2(M, M_l(\A_i)))$. For the
integral kernel $(1-1_{K_m}(x))k(x,y) 1_{K_m}(y)$ of $c$ we have $y \mapsto (1-1_{K_m})k(\cdot,y) 1_{K_m}(y) \in
C(K_m,L^2(M, M_l(\A_i)))$. It follows that $bd^{-1}c$ is an integral operator with
continuous kernel on $K_m \times K_m$. Clearly the integral kernel of $1_{K_m}K1_{K_m}$ is continuous
as well.

Since $a-bd^{-1}c$ is invertible on $L^2(K_m,\A^l)$ and  since the measure on
$K_m$ pulled back by an orientation preserving diffeomorphism $[0,1]^n \to K_m$ is of the form $h
dx$, we conclude by (1) that $a-bd^{-1}c$ is invertible on $L^2(K_m, ( \Ol{\mu})^l)$ as well. 
\end{proof}

\begin{cor}
\label{specint}
Let $k$ be an integral kernel as in part (2) of the proposition and let $K$ be the corresponding
integral operator. Then for $\lambda \in \bbbc^*$ the operator $K-\lambda$ is invertible on $L^2(M, ( \Ol{\mu})^l)$ if $K-\lambda$ is invertible  on $L^2(M,\A^l)$.
\end{cor}

\begin{proof}
For $\lambda \in\bbbc \setminus \{0\}$ the integral kernel $k/\lambda$ fulfills
the conditions of the lemma. Thus if $\lambda-K=\lambda(1-K/\lambda)$ is
invertible on $L^2(M,\A^l)$, then
$\lambda-K$ is invertible on $L^2(M, ( \Ol{\mu})^l)$ as well.
\end{proof}

\subsection{The Chern character}
\label{kern}
\begin{ddd}
Let $V$ be a $\bbbz/2$-graded finite dimensional vector space.
Let $K$ be an integral operator on $L^2(M,V \ten \Ol{\mu})$ with integral kernel $k:M \times M \to
\End(V) \ten \Ol{\mu}$. Then we define $\di (K)$ to be the integral operator on $L^2(M,V
\ten \Ol{\mu})$
with integral kernel $\di(k(x,y))$. (The action of
$\di$ on $\End(V) \ten \Ol{\mu}$ was described in \S\ref{sutrace}.)
\end{ddd}

Note that if $K$ is of degree $n$ with respect to the $\bbbz/2$-grading on
$L^2(M,V \ten \Ol{\mu})$, then $$\di (Kf)=\di (K)f+
(-1)^nK\di f\ .$$

\begin{lem}
\label{chHS}
Let $M$ be a $\sigma$-finite measure space.

Let  $P$ be a Hilbert-Schmidt operator on $L^2(M,(\Ol{\mu})^l)$ with integral
kernel in $L^2(M \times M, M_l(\A_i))$ and  assume that $P^2=P$. Then
$$P\di P \di P=P(\di (P))^2 \ .$$
\end{lem}

\begin{proof}
As for matrices (see the beginning of the proof of Prop. \ref{homch}) we have that $(\di P)=P(\di P)+(\di P)P$ and
 $P(\di P)P=0$.

It follows that
\begin{eqnarray*}
P\di P\di P &=& P(\di P)\di P=P(\di P)^2P \\ 
&=&P(\di P)(\di P)-P(\di P)P(\di P) \\
&=&P(\di (P))^2 .
\end{eqnarray*}
\end{proof}

\begin{lem}
Let $P:[0,1] \to HS(L^2(M,(\Ol{\mu})^l))$ be a differentiable path of Hilbert-Schmidt operators with integral kernels in $L^2(M \times
M,M_l(\A_i))$, and  assume that $P(t)^2=P(t)$ for every $t \in [0,1]$.

Then for $k \in \bbbn_0$ in $\Ol{\mu}/\ov{[\Ol{\mu},\Ol{\mu}]}$
$$\Tr P(1)(\di P(1))^{2k}- \Tr P(0)
(\di P(0))^{2k}$$ is exact.
\end{lem}

\begin{proof}
As for matrices (see Prop. \ref{homch}).
\end{proof}

If $T$ is a trace class operator on $L^2(M,(\Ol{\mu})^l)$ restricting to a trace class operator
  on $L^2(M,(\hat\Omega_{\le
\nu}\A_j)^l)$ for all $j,\nu \in \bbbn$ with $j \ge i$ and $\nu \ge \mu$,
then by taking the projective limit we can consider $\Tr(T)$ as in element in $\Oi\Ai/\ov{[\Oi\Ai,\Oi\Ai]_s}.$ 
In the next lemma we show that under certain conditions the formula for the Chern character can be generalized to projections which are Hilbert-Schmidt operators.

\begin{prop}
\label{projker}

Let $M$ be a Riemannian manifold of dimension $d$, possibly with boundary. Suppose that there is a covering $\{K_m\}_{ m \in \bbbn}$ of $M$ with $K_m$ compact, $K_m \subset K_{m+1}$ and such that $K_m$ is
diffeomorphic to $[0,1]^d$ for all $m \in \bbbn$. 

Let $P \in B(L^2(M, \A^l))$ be a projection onto a projective submodule of
 $L^2(M, \A^l)$. Assume further that for any $i \in \bbbn$ it restricts to a bounded
 projection on $L^2(M, \A_i^l)$ and that $P(L^2(M,\A_i^l)) \subset
 C(M,\A_i^l)$. Let $$\Ran_{\infty} P:=\bigcap\limits_{i \in \bbbn}
 P(L^2(M,\A_i^l))   \ .$$ 

\begin{enumerate}
\item The projection $P$ is a Hilbert-Schmidt
operator with integral kernel of the form $\sum_{j=1}^m f_j(x)h_j(y)^*$ with $f_j,h_j \in \Ran_{\infty} P$.
\item The intersection $\Ran_{\infty} P$ is a projective
$\Ai$-module. The classes $[\Ran P] \in K_0(\A)$ and $[\Ran_{\infty} P] \in K_0(\Ai)$
correspond to each other under the canonical isomorphism $K_0(\A) \cong
K_0(\Ai)$. 
\item In $H^{dR}_*(\Ai)$ 
$$\ch[\Ran_{\infty} P]=\sum\limits_{n=0}^{\infty}(-1)^n 
\frac 1{n!}  \Tr (P\di P)^{2n}  \ .$$

\end{enumerate}
\end{prop}

\begin{proof}

(1)
Let $\{e_n\}_{n \in \bbbn} \subset \C_0(M,\bbbc^l)$ be an orthonormal basis of $L^2(M,\bbbc^l)$.

The orthogonal projection $P_n$ onto the span of the first $n$ basis vectors is
a Hilbert-Schmidt operator with integral kernel
in $\C_0(M \times M ,M_l(\bbbc)).$

In particular $P_n \in B(L^2(M,\A_i^l))$ for any $i \in \bbbn$.

First we consider the situation on $L^2(M,\A^l)$:

Since $P$ is compact, there is $n \in \bbbn$ such that on $L^2(M,\A^l)$ $$\|P(P_n-1)\| \le \frac 12 \ .$$

Then by Prop. \ref{dirsum} the map 
$PP_nP:\Ran P \to \Ran P$ is an isomorphism.

It follows that
$\Ker PP_nP=(\Ran P)^{\perp}=\Ker P$ and
therefore 

$$P =1- P_{\Ker P P_nP} \ .$$
Here $P_{\Ker P P_nP}$ denotes the orthogonal projection onto
 $\Ker P P_nP$.

We can find $r>0$ such
that 
$B_r(0)\setminus \{0\}$ is in the resolvent set of $PP_nP$.

Then 
\begin{eqnarray*}
P &=& 1- P_{\Ker P P_nP} \\
&=& 1- \frac{1}{ 2 \pi i }\int_{|\lambda|=r} (\lambda-PP_nP)^{-1} d \lambda \\
&=& \frac{1}{2 \pi i} \int_{|\lambda|=r} \left(\lambda^{-1}-  (\lambda-PP_nP)^{-1} \right) d \lambda\\
&=& -\frac{1}{2 \pi i}  PP_nP \int_{|\lambda|=r}\lambda^{-1}(\lambda -PP_nP)^{-1}
d \lambda \ . 
\end{eqnarray*}
The integral
kernel of $PP_nP$ fulfills the conditions of Prop. \ref{intker}. From Cor. \ref{specint}  we conclude that the spectrum of $PP_nP$ in
$B(L^2(M,\A_i^l))$ is independent of $i \in \bbbn_0$, thus 
$$R:=\int_{|\lambda|=r}\lambda^{-1}(\lambda -PP_nP)^{-1}
d \lambda$$ is a bounded operator on
$L^2(M,\A_i^l)$ for all $i \in \bbbn_0$.

This and the equation $PR=PRP$ show that $P$ is an integral operator with  integral kernel
$$k_P(x,y)=- \frac{1}{2 \pi i}  \sum\limits_{j=1}^n Pe_j(x)(P R^* Pe_j(y))^* \ .$$
The integral kernel is in $L^2(M \times M, M_l(\A_i))$ for all $i \in
\bbbn_0$ and is of the form we asserted.

(2) Let $P_n$ as above with $\|(P_n-1)P\| \le \frac 12$ in $B(L^2(M,\A^l))$. 

The operator $(1+(P_n-1)P)$
is invertible in $B(L^2(M, \A^l))$, hence by Prop. \ref{intker} it is invertible in $B(L^2(M, \A_i^l))$ for any $i \in \bbbn_0$. From $P_nP=(1+(P_n-1)P)P$
it follows that $$P_nP(L^2(M, \A_i^l)) \cong \Ran P(L^2(M, \A_i^l))$$ for all $i \in \bbbn_0$.
Furthermore $Q:=(1+(P_n-1)P)P(1+(P_n-1)P)^{-1} \in B(L^2(M, \A_i^l))$ is a projection onto $P_nP(L^2(M, \A_i^l))$. 
 Identify $P_n(L^2(M, \A_i^l))$ with $\A_i^n$ via the basis and let $Q'
\in M_n(\Ai)$ be the restriction of $Q$ to $P_n(L^2(M, \A_i^l))$. From $\Ran_{\infty} P \cong Q'(\Ai^n)$ it follows that $[\Ran P]=[Q']$ in $K_0(\A)$ and $[\Ran_{\infty}
P]=[Q']$ in $K_0(\Ai)$. This shows the assertion.

(3) Let $P_n, Q$ and $Q'$ be as in the proof of (2).  
In $H_*^{dR}(\Ai)$ 
\begin{eqnarray*}\ch[\Ran_{\infty}
P]&=&\ch(Q')\\
&=&\sum\limits_{k=0}^{\infty}
(-1)^k\frac{1}{k!} \tr (Q'\di Q')^{2k}\\
&=&\sum\limits_{k=0}^{\infty}
(-1)^k\frac{1}{k!} \Tr  ((P_nQP_n)\di (P_nQP_n))^{2k}\\
&=&\sum\limits_{k=0}^{\infty}
(-1)^k\frac{1}{k!} \Tr (P_nQ)(\di P_nQ)^{2k}\\
&=&\sum\limits_{k=0}^{\infty}
(-1)^k\frac{1}{k!} \Tr(Q (\di Q)^{2k}) \ .
\end{eqnarray*}
Since $$H:[0,1] \to B(L^2(M,(\Ol{\mu})^{l})) \ ,$$
$$H(t)= (1+t(P_n-1)P)P(1+t(P_n-1)P)^{-1}$$
is a differentiable path of finite projections with $H(0)=P$ and
$H(1)=Q$, the
difference $\Tr (Q\di Q)^{2k}-\Tr(P\di P)^{2k}$ is  exact in $\Oi\Ai/\ov{[\Oi\Ai,\Oi\Ai]_s}$  by the previous
lemma.

This shows the assertion.
\end{proof}

\section{Holomorphic semigroups}
\label{semigr}

\subsection{General facts}

Let $X$ be a Banach space.

In order to fix notation we collect some well-known facts about holomorphic semigroups on $X$. A general reference is \cite{rr}.\\

For $\delta \in ]0, \pi]$ let $\Sigma_{\delta}:=\{ \lambda \in \bbbc^*~|~ |\arg
\lambda|<\delta\}$.

Recall that a family $S:\Sigma_{\delta} \cup \{0\} \to B(X)$ is called a holomorphic semigroup if it is a semigroup, i.e. $S(0)=1$ and $S(t+s)=S(t)S(s)$ for all $s,t \in \Sigma_{\delta} \cup \{0\}$, moreover if it is  strongly continuous on $\Sigma_{\delta'} \cup \{0\}$ for all $\delta'<\delta$ and holomorphic on $\Sigma_{\delta}$.\\

Given a holomorphic semigroup $S(t)$, the unbounded operator $Z$ with $Zx:= \frac{d}{dt}(S(t)x)|_{t=0}$ whenever defined is called the generator of $Z$. It is densely defined and closed. Any subset of its domain that is dense in $X$ and invariant under the action of the semigroup is a core for $Z$.\\
 
Recall that a densely defined operator $Z$ on $X$ generates a holomorphic semigroup $e^{tZ}$ if and only
if there is $\omega \in \bbbr$ such that  $Z+\omega$ is
$\delta$-sectorial for some $\delta \in ]0,\pi/2]$, i.e. if and only if
$\Sigma_{\delta + \pi/2}$ is a subset of the resolvent set $\rho(Z+\omega)$ 
and for any $\ve$ with $0<\ve<\delta$
there is $C > 0$ such that for all $\lambda \in \Sigma_{\ve+\pi/2}$ 
$$\|(Z+\omega-\lambda)^{-1}\| \le \frac{C}{|\lambda|} \ .$$ 
Then there is $C>0$ such that for all $t>0$
$$\|e^{tZ} \| \le Ce^{-\omega t} \ .$$

\begin{lem}
\label{secl}
Let $Z$ be a densely defined operator on $X$ and let $\delta>0$ be such that

\begin{itemize}
\item[(i)]
$\Sigma_{\delta+\pi/2}\cup \{0\} \subset \rho(Z) \ ,$
\item[(ii)] for every $\alpha<\delta$ there are $c \in \bbbr$ and $r>0$ such that $$\|(Z+c-\lambda)^{-1}\| \le
\frac{C}{|\lambda|}$$ for $\lambda \in \Sigma_{\alpha +\pi/2}$ with $|\lambda|> r$ and $\lambda-c \in \rho(Z)$.
\end{itemize}

Then $Z$ is $\delta$-sectorial.
\end{lem}

\begin{proof} Let $0<\ve<\delta$.

Let $\alpha \in (\ve,\delta)$ and $R>r$ big enough such that any $\lambda \in \Sigma_{\ve +\pi/2}$ with $|\lambda|>R$ satisfies $\lambda-c \in  \Sigma_{\alpha +\pi/2}$. 
By assumption we find $C>0$ such that for all $\lambda \in \Sigma_{\ve +\pi/2}$ with $|\lambda| > R$
\begin{eqnarray*}
\|(Z -\lambda)^{-1}\| &\le& \frac{C}{|\lambda +c|}\\
&\le&\left(\frac{C}{|\lambda|}\right) \left(\frac{|\lambda+c| + |c|}{|\lambda+c|} \right) \\
&\le& \left(\frac{C}{|\lambda|}\right) \left(1+
\frac{|c|}{(R-|c|)}\right)\\
&\le & \frac{C}{|\lambda|} \ .
\end{eqnarray*}
Since $\ov{\Sigma_{\ve +\pi/2}} \cap \{|z| \le R\}$ is a compact subset of the resolvent set of $Z$, this implies the assertion.

\end{proof}

We have the following relation between the spectrum of a
$\delta$-sectorial operator $Z$ and the behavior of the holomorphic semigroup
$e^{tZ}$ for $t \to \infty$.

\begin{prop}
\label{sec}

If $Z$ is a $\delta$-sectorial operator and there is $\omega>0$ such that 
$$\{ \re \lambda > -\omega \} \subset \rho(Z) \ ,$$ then 
 for any $\omega'<\omega$ there is $C>0$ such that for all $t \ge 0$  
$$\|e^{tZ} \| \le Ce^{-\omega' t} \ .$$ 
\end{prop}

\begin{proof}
For any $\omega'<\omega$ there is $0<\delta' \le \delta$ such that
$\Sigma_{\delta'+ \pi/2} \cup \{0\} \subset \rho(Z+\omega')$. Now it follows from  the
previous lemma that $Z+\omega'$ is $\delta'$-sectorial.
\end{proof}

\begin{prop}
\label{critsg}
If  $e^{tZ}$ is a strongly continuous semigroup with generator $Z$ such that
$\Ran e^{tZ} \subset \dom Z$ for all $t \in (0, \infty)$ and if there are $T>0$
and $C>0$ such that for $t \le T$
$$\|Ze^{tZ}\| \le Ct^{-1} \ ,$$ 
then $e^{tZ}$ extends to a holomorphic semigroup.

If the estimate holds for all $t >0$, then the extension is bounded holomorphic.
\end{prop}

\begin{proof}

The assertion follows immediately from \cite{da}, Th. 2.39.
\end{proof}

Part (1) of the next proposition is known under the name Volterra development and
the formula in part (2) is called Duhamel's formula:

\begin{prop}
\label{boundpert}
Let $Z$ be the generator of a strongly continuous semigroup and let $M, \omega
\in \bbbr$ be such that $\|e^{tZ}\| \le
Me^{\omega t}$ for all $ t \ge 0$.
\begin{enumerate}
\item Let $R \in
B(X)$. 
Then $Z+R$ is the generator of a strongly continuous semigroup and for all $t
\ge 0$ 
$$e^{t(Z+R)}=\sum\limits_{n=0}^{\infty}t^n\int_{\Delta^n}e^{u_0tZ}
Re^{u_1tZ}  R\dots e^{u_ntZ} ~du_0 \dots du_n$$
with $\Delta^n= \{u_0+ \dots +u_n=1; u_i \ge 0,~i=0, \dots, n\}$. 
Furthermore 
$$\|e^{t(Z+R)}\| \le Me^{(\omega+M\|R\|)t} \ .$$
\item Let $R_1, \dots, R_n \in B(X)$.
For $t \ge 0$ the map 
$$\bbbc^n \to B(X),~(z_1, \dots,z_n)  \mapsto e^{t(Z+z_1R_1+ \dots + z_n R_n)}$$
is analytical and for $i=0, \dots n$ 
$$\frac{d}{dz_i}e^{t(Z+z_1R_1+ \dots + z_n R_n)}=\int_0^t e^{(t-s)(Z+z_1R_1+ \dots + z_n R_n)}R_i e^{s(Z+z_1R_1+ \dots + z_n R_n)}~ds \ .$$
\end{enumerate}
\end{prop}

\begin{proof}
(1) follows from \cite{da}, Th. 3.1 and the proof of it.

(2) The analyticity follows from (1). 

For the formula it is enough to consider $n=1$. Let $R:=R_1$.

For $z_0 \in \bbbc$ by (1)
$$e^{(Z+zR)t}=\sum\limits_{n=0}^{\infty}(z-z_0)^n t^n\int_{\Delta^n}e^{u_0t(Z+z_0R)}
Re^{u_1t(Z+z_0R)}  R\dots e^{u_nt(Z+z_0R)} ~du_0 \dots du_n \ . $$
This implies that
$$\frac{d}{dz}e^{(Z+zR)t}|_{z_0}= t\int_0^1 e^{u_0t(Z+z_0R)}
Re^{(1-u_0)t(Z+z_0R)}du_0 \ .$$
\end{proof}

The following proposition is known in the literature as Duhamel's principle:

\begin{prop}
\label{duhprin}
Let $Z$ be the generator of a strongly continuous semigroup on $X$. Let $u \in
C^1([0, \infty ),X)$ such that $u(t) \in \dom Z$ and $\frac{d}{dt}u(t)-Zu(t) \in \dom Z$ for all $t \in
[0,\infty)$. Then 
$$e^{tZ}u(0)-u(t)=-\int_0^t e^{sZ}\bigl(\frac{d}{dt}-Z\bigr)u(t-s)ds \ .$$
\end{prop}

\begin{proof} see \cite{ta}, Appendix A, (9.37) and (9.38). \end{proof}

\subsection{Square roots of generators and perturbations}

Assume that $D$ is a densely defined closed operator on a Banach space $X$ with
bounded inverse and such that
$-D^2$ is $\delta$-sectorial. 

There are well-defined fractional powers
$(D^2)^{\alpha}$ for $\alpha \in \bbbr$ (\cite{rr}, \S 11.4.2). These are densely defined closed
operators that coincide for $\alpha \in \bbbz$ with the usual powers and satisfy $$(D^2)^{\alpha + \beta}f=(D^2)^{\alpha}  (D^2)^{\beta}f$$ for all $\alpha,\beta \in \bbbr$ and $f \in \dom (D^2)^{\gamma}$ with $\gamma= \max\{\alpha, \beta, \alpha + \beta\}.$
For $\alpha\le 0$ the operator $(D^2)^{\alpha}$ is bounded and depends  of $\alpha$ in a strongly
continuous way. 

In particular it follows  that for  $\alpha \ge 0$ the operator
$(D^2)^{-\alpha}$ is a bounded inverse of $(D^2)^{\alpha}.$

Define $|D|:=(D^2)^{\frac 12}.$ 

By \cite{ka}, Th. 2, the operator  $-|D|$ is $(\delta+\frac{\pi/2-\delta}{2})$-sectorial and can be
expressed in terms of the resolvents of $D^2$. 

Note that for every $n \in
\bbbn$ the domain of $|D|^n$ is a core of $|D|$ and $\dom D^n$ is a core of
$D$. 

\begin{lem}
\begin{enumerate}
\item
Let $P$ be a closed densely
defined operator on $X$ and assume that $P^2$ is densely
defined and  $P^2|_M=P|_M$ for some dense subset $M$ of $\dom P^2$. 

Then $P$ is a bounded projection.

\item
Let $I$ be a closed densely
defined operator on $X$ and assume that $I^2$ is densely defined and that
there is a dense subset $M$ of $\dom I^2$ with $I^2|_M=1|_M$. 

Then $I$ is a bounded involution.
\end{enumerate}
\end{lem}

\begin{proof}
(1) If $f \in M$, then $f= (1-P)f + Pf$ and $(1-P)Pf=P(1-P)f=0$. We conclude that $M\subset \Ker P + \Ker
(1-P) \subset \dom P$.  If $P$ is closed, then $\Ker P + \Ker
(1-P)$ is closed, hence $\Ker P + \Ker
(1-P)=X$, thus $\dom P=X$.

(2) The operator $P:= \frac 12(1+I)$ is a closed projection on $X$ in the sense of
(1). Thus $P$ is bounded. It follows that $I$ is bounded
as well.
\end{proof}

\begin{prop}
\label{betrag}
The closure of the operator $|D|^{-1}D: \dom D \to X$ is a bounded involution $I$
on $X$ and
\begin{enumerate}
\item $\dom D=\dom |D|$ and $I(\dom D)\subset \dom D$.
\item $|D|=ID=DI$ and $D=I|D|=|D|I$.
\end{enumerate}
\end{prop}

\begin{proof}  The operator $D^{-1}$ commutes with the resolvents of $D^2$. It
follows that $|D|^{-1}D^{-1}=D^{-1}|D|^{-1}$. Hence  $|D|^{-1}(\dom D) \subset
\dom D$ because of $\dom D=D^{-1} X$ , so $\dom D^2 \subset \dom (|D|^{-1}D)^2$. 

If $f
\in \dom D$, then $|D|^{-1} Df=D|D|^{-1}f$. For $f \in \dom D^2$ it
follows that $(|D|^{-1}D)^2f=f$.

Let $(f_n)_{n \in \bbbn}$ be a sequence in $\dom D$ converging to zero. If $|D|^{-1}Df_n =D|D|^{-1}f_n$  converges, the limit is
zero since $D$ is closed. Hence $|D|^{-1}D$ is closable. By the
previous lemma it extends to a bounded involution.

Since for $f \in \dom D^2$ we have that $DIf=|D|f$ and since $DI$ is closed, the composition $IDI:\dom |D| \to X$ is well-defined and closed. It coincides with $D$ on $\dom D^2$, hence it is
a closed extension of $D$. It follows that $\dom D \subset \dom |D|$. The
inclusion $\dom |D|\subset \dom  D$ is shown analogously. 

The equations are clear on $\dom D^2$, which is a core for $D$ and $|D|$.
\end{proof}

The operator $P:=\frac 12(1+I)$ with $I$ as in the previous proposition is a bounded
projection on $X$.
By the proposition  $P\dom D \subset \dom D$ and $P$ commutes with $D$
and $|D|$.

From $ID=|D|$ and $I|D|=D$ it follows that $PD=-P|D|$. Thus with respect to the decomposition $X=PX
\oplus (1-P)X$ 
$$D=\left(\begin{array}{cc}
PDP & 0\\
0 & -(1-P)D(1-P)
\end{array}\right)
$$ and 
$$|D|=\left(\begin{array}{cc}
PDP & 0 \\
0 & (1-P)D(1-P)
\end{array}\right) \ .
$$
By taking
into account that $-|D|$ is $(\delta+\frac{\pi/2-\delta}{2})$-sectorial it follows for
the resolvent set of $D$:

\begin{prop}
\label{compspec}
For $\lambda \in \bbbc$ 
$$\{\lambda,- \lambda\}  \subset \rho(D) \aq
\{\lambda,-\lambda\} \subset \rho(|D|)\ .$$ 
Thus if $\lambda \in \bbbc$ with
$-\lambda^2 \in \Sigma_{\pi/2+\delta}$,
then $\lambda \in
\rho(D)$.

Furthermore for every $\delta'<\delta$ there is $C>0$ such that for all $\lambda$ with
 $-\lambda^2 \in \Sigma_{\delta'+\pi/2}$ 
$$\|(D-\lambda)^{-1}\| \le \frac{C}{|\lambda|} \ .$$
\end{prop}

\begin{cor} 
\label{semest} 
Let $\omega \ge 0$ be such that there is $C>0$ with $\|e^{-tD^2}\| \le C
e^{-\omega t}$ for all $t\ge 0$.
\begin{enumerate}
\item
For every $\alpha \in \bbbr$ and $\omega'<\omega$ there is $C>0$ such that for
all $t >0$ 
$$\| |D|^{\alpha} e^{-tD^2}\| \le Ct^{-\alpha/2}e^{-\omega't} \ .$$
\item For every $n \in \bbbn$ and $\omega'< \omega$ there is $C>0$ such that for all
$t>0$  
$$\|D^ne^{-tD^2}\| \le C t^{-n/2}e^{-\omega' t} \ .$$
\end{enumerate}
\end{cor}

\begin{proof}
The first assertion is \cite{rr}, Lemma 11.36, and the second one follows from
the first one by $D=I|D|$ and $DI=ID$.
\end{proof}

\begin{prop}
\label{sgpert}
Let $A$ be a bounded operator and let $\delta'<\delta$.
\begin{enumerate}
\item There is $R>0$ such that $D+A-\lambda$ has a bounded inverse if $|\lambda|>R$ and $-\lambda^2 \in
  \Sigma_{\delta'+\pi/2}$. 
\item There is $\omega >0$ such that $-(D+A)^2+\omega$ is
  $\delta'$-sectorial. 
\item $D+A$ commutes with $e^{-t(D+A)^2}$. 
\end{enumerate}
\end{prop}

\begin{proof}
By Prop. \ref{compspec} there is $M>0$ such that for all $\lambda$ with
 $-\lambda^2 \in \Sigma_{\delta'+\pi/2}$
$$\|(D-\lambda)^{-1}\| \le \frac{M}{|\lambda|}.$$
Hence the Neumann series
$$(D+A-\lambda)^{-1}=(D-\lambda)^{-1}\sum_{n=0}^{\infty}(A(D-\lambda)^{-1})^n$$
converges for $|\lambda|>M \|A\|$ and $-\lambda^2 \in \Sigma_{\delta' + \pi/2}$.

This shows (1).

If $|\lambda|>2M\|A\|$ and $-\lambda^2 \in \Sigma_{\delta'+\pi/2}$, then
\begin{eqnarray*}
\|(D+A-\lambda)^{-1}\|&=&\|\sum_{n=0}^{\infty}(D-\lambda)^{-1}\left(A(D-\lambda)^{-1}\right)^n\|\\
 &\le& \sum_{n=0}^{\infty}\|A\|^n\|(D-\lambda)^{-1}\|^{n+1}\\
&\le& \sum_{n=0}^{\infty}\frac{M^{n+1}\|A\|^n}{|\lambda|^{n+1}}\\
&=&\frac{M}{|\lambda|}(1-\frac{M\|A\|}{|\lambda|})^{-1}\\
&\le&\frac{2M}{|\lambda|} \ .
\end{eqnarray*}
Let $\mu \in \{|\mu|>4M^2\|A\|^2 \}\cap \Sigma_{\delta'+\pi/2}$. If
$\lambda \in \bbbc$ with $-\lambda^2=\mu$, then $\lambda \in \rho(D+A)$, hence
 the resolvent
$$(-(D+A)^2-\mu)^{-1}=-(D+A-\lambda)^{-1}(D+A+\lambda)^{-1}$$ exists and is
bounded by
$$ \|(-(D+A)^2-\mu)^{-1}\| \le \frac{4M^2}{|\mu|} \ .$$
There is $\omega >4M^2\|A\|^2$ such that
$$\Sigma_{\delta'+\pi/2} \cup \{0\} \subset \bigl(\{|\mu|>4M^2\|A\|^2 \}\cap
\Sigma_{\delta'+\pi/2}\bigr) -\omega$$
and thus
$$\Sigma_{\delta'+\pi/2}\cup \{0\}  \subset \rho\bigl(-(D+A)^2 + \omega\bigr) \ .$$
Assertion (2)  follows now from Lemma \ref{secl}.

(3) follows from the fact that $e^{-t(D+A)^2}$ can be expressed in terms of the
   resolvents of $(D+A)^2$, which commute with $D+A$. 
\end{proof}

\backmatter

\end{document}